\documentclass{amsart}

\usepackage{setspace}

\usepackage{amssymb}   

\usepackage{amsmath}

\usepackage{mathrsfs}

\usepackage{stmaryrd}

\usepackage{amsthm}
\usepackage{newlfont}
\usepackage{amscd}
\usepackage{mathtools}
\usepackage{bm}
\usepackage{tikz-cd}
\usepackage{graphicx}

\usepackage{hyperref}
\usepackage{framed}
\usepackage{comment}

\usepackage{enumitem}

\usepackage[all]{xy}

\usepackage[normalem]{ulem}

\usepackage{textcomp}

\usepackage{amsbsy}

\usepackage[T1]{fontenc}

\newtheorem{thm}{Theorem}[section]

\newtheorem{thm-defn}[thm]{Theorem/Definition}
\newtheorem{lem}[thm]{Lemma}
\newtheorem{prop}[thm]{Proposition}
\newtheorem{cor}[thm]{Corollary}

\theoremstyle{definition}
\newtheorem{defn}[thm]{Definition}

\newtheorem{eg}[thm]{Example}

\newtheorem{construction}[thm]{Construction}
\newtheorem{set-up}[thm]{Set-up}
\newtheorem{notation}[thm]{Notation}

\theoremstyle{remark}
\newtheorem{rem}[thm]{Remark}

\numberwithin{equation}{section}

\usepackage{relsize}
\usepackage[bbgreekl]{mathbbol}
\usepackage{amsfonts}

\DeclareSymbolFontAlphabet{\mathbb}{AMSb}
\DeclareSymbolFontAlphabet{\mathbbl}{bbold}
\newcommand{\Prism}{{\mathlarger{\mathbbl{\Delta}}}}

\DeclareMathOperator{\Gal}{Gal}

\DeclareMathOperator{\Spec}{Spec}
\DeclareMathOperator{\Spa}{Spa}
\DeclareMathOperator{\Spf}{Spf}

\DeclareMathOperator{\Frac}{Frac}
\DeclareMathOperator{\Ker}{Ker}
\DeclareMathOperator{\Fil}{Fil}
\DeclareMathOperator{\Hom}{Hom}
\DeclareMathOperator{\End}{End}
\DeclareMathOperator{\Mor}{Mor}

\newcommand{\CR}{\mathrm{CR}}
\newcommand{\QCoh}{\mathrm{QCoh}}

\newcommand{\rank}{\mathrm{rank}}
\newcommand{\Vect}{\mathrm{Vect}}
\newcommand{\Loc}{\mathrm{Loc}}
\newcommand{\perf}{\mathrm{perf}}
\newcommand{\str}{\mathrm{str}}

\newcommand{\pr}{\mathrm{pr}}

\newcommand{\an}{\mathrm{an}}
\newcommand{\aff}{\mathrm{aff}}

\newcommand{\cris}{\mathrm{cris}}
\newcommand{\CRIS}{\mathrm{CRIS}}
\newcommand{\st}{\mathrm{st}}

\newcommand{\dR}{\mathrm{dR}}
\newcommand{\et}{\mathrm{\acute{e}t}}
\newcommand{\proet}{\mathrm{pro\acute{e}t}}
\newcommand{\Perfd}{\mathrm{Perfd}}

\newcommand{\ra}{\rightarrow}
\newcommand{\N}{\mathbf{N}}
\newcommand{\F}{\mathbf{F}}
\newcommand{\Z}{\mathbf{Z}}
\newcommand{\Q}{\mathbf{Q}}

\newcommand{\A}{\mathbf{A}}
\newcommand{\B}{\mathbf{B}}
\newcommand{\OA}{\mathbf{O}\mspace{-2mu}\mathbf{A}}
\newcommand{\OB}{\mathbf{O}\mspace{-2mu}\mathbf{B}}
\newcommand{\Ainf}{\mathbf{A}_{\mathrm{inf}}}

\newcommand{\calOB}{\mathcal{O}\mathbb{B}}

\renewcommand{\log}{\mathrm{log}}

\newcommand{\fkn}{\mathfrak{n}}

\newcommand{\fkM}{{\mathfrak M}}

\newcommand{\fkS}{{\mathfrak S}}

\newcommand{\calC}{\mathcal{C}}

\newcommand{\calE}{\mathcal{E}}
\newcommand{\calF}{\mathcal{F}}
\newcommand{\calG}{\mathcal{G}}

\newcommand{\calI}{\mathcal{I}}
\newcommand{\calJ}{\mathcal{J}}

\newcommand{\calM}{\mathcal{M}}

\newcommand{\calO}{\mathcal{O}}

\newcommand{\calX}{\mathcal{X}}

\newcommand{\bA}{\mathbb{A}}
\newcommand{\bB}{\mathbb{B}}

\newcommand{\bL}{\mathbb{L}}

\mathchardef\mhyphen="2D

\begin{document}

\pagenumbering{arabic}

\title{Log prismatic $F$-crystals and purity}

\author{Heng Du} 
\address{Yau Mathematical Sciences Center, Tsinghua University, Beijing 100084, China}
\email{hengdu@mail.tsinghua.edu.cn}
\author{Tong Liu}
\address{Department of Mathematics, Purdue University,
150 N. University Street, West Lafayette, IN 47907, U.S.A.}
\email{tongliu@purdue.edu}
\author{Yong Suk Moon} 
\address{Beijing Institute of Mathematical Sciences and Applications, Beijing 101408, China}
\email{ysmoon@bimsa.cn}
\author{Koji Shimizu}
\address{Yau Mathematical Sciences Center, Tsinghua University, Beijing 100084, China;
Beijing Institute of Mathematical Sciences and Applications, Beijing 101408, China}
\email{shimizu@tsinghua.edu.cn}

\begin{abstract} 
Our goal is to study $p$-adic local systems on a rigid-analytic variety with semistable formal model. We prove that such a local system is semistable if and only if so are its restrictions to the points corresponding to the irreducible components of the special fiber.
For this, the main body of the paper concerns analytic prismatic $F$-crystals on the absolute logarithmic prismatic site of a semistable $p$-adic log formal scheme.  Analyzing Breuil--Kisin log prisms, we obtain a prismatic purity theorem and deduce the above purity theorem for semistable local systems.
\end{abstract}

\maketitle

\tableofcontents
\section{Introduction} \label{sec:introduction}

Fix a prime $p$. Let $K$ be a complete discrete valuation field of mixed characteristic $(0,p)$ with perfect residue field $k$. Write $\calO_K$ for the ring of integers and let $\pi$ be a uniformizer.

In the study of $p$-adic Galois representations of $K$, Fontaine introduced the notion of crystalline, semistable, and de Rham representations. These classes are defined using period rings and reflect properties of Galois representations coming from geometry. For example, a deep result in $p$-adic Hodge theory shows that a Galois representation arising as a $p$-adic \'etale cohomology of an algebraic variety over $K$ is de Rham. Moreover, if the variety is smooth proper with good (resp.~ semistable) reduction, the resulting Galois representation is crystalline (resp.~semistable). 

If one is interested in a family of $p$-adic Galois representations parametrized by a variety over $K$, one is led to study $p$-adic representations of the \'etale fundamental group or $p$-adic \'etale local systems on the variety. Crystalline and de Rham local systems, as generalizations of crystalline and de Rham representations, have been well-studied: to name a few, see \cite{Faltings-CryscohandGalrep,faltings-almostetale, brinon-relative, Andreatta-Iovita-smooth, scholze-p-adic-hodge, Tan-Tong, liu-zhu-rigidity}.

We aim to study semistable representations in the relative setting: we work on a smooth adic space $\calX$ over $K$ with semistable reduction and consider semistable \'etale local systems on $\calX$.
There have been several works to study semistable local systems \cite{Faltings-CryscohandGalrep, faltings-almostetale, andreatta-iovita-semistable-relative, tsuji-cryst-shvs}. Here one needs to fix a $p$-adic formal scheme $X$ over $\calO_K$ whose generic fiber gives $\calX$ and use log geometry to understand semistable periods in geometric families, and as a result, the definitions of semistable local systems in these works are much more involved than those of semistable representations or crystalline and de Rham local systems. In particular, it is generally difficult to check whether a given $p$-adic \'etale local system is semistable. 

In this paper, we prove a \emph{purity theorem} that relates semistable local systems with semistable Galois representations: for each irreducible component of $X$, the completed local ring at the generic point is a complete discrete valuation ring (CDVR for short) and defines a rank-$1$ point of $\calX$, which we call an \emph{$X$-Shilov point}. Note that $p$-adic local systems on such a point correspond to $p$-adic Galois representations of the fraction field of the CDVR, and Fontaine's formulation of semistable representations is developed in \cite{fontaine-representations, Morita-imperfsemistable, Ohkubo}. Here is our main result in a crude form.

\begin{thm}[Purity for semistable local systems]\label{thm:intro purity for semistable local systems}
A $p$-adic local system on $\calX$ is semistable if and only if its restriction to each $X$-Shilov point corresponds to a semistable Galois representation.
\end{thm}

When $X$ is smooth over $\calO_K$, this type of a purity theorem for crystalline local systems was first obtained by Tsuji in \cite{tsuji-cryst-shvs}. Our approach is different from the aforementioned works: we employ logarithmic prismatic theory for a semistable $p$-adic formal scheme to study semistable local systems of the generic fiber and prove the purity theorem.

\medskip \noindent
\textbf{Logarithmic prismatic theory.}
Prismatic cohomology theory by Bhatt and Scholze 
\cite{bhatt-scholze-prismaticcohom} has brought new perspectives to integral $p$-adic cohomology theory. 
Here we focus on $F$-crystals on the absolute prismatic site and its logarithmic variant to study local systems.

For a bounded $p$-adic formal scheme $X$, the absolute prismatic site $X_\Prism$ consists of bounded prisms $(A,I)$ together with a map $\Spf (A/I)\rightarrow X$. 
Koshikawa \cite{koshikawa} develops logarithmic prismatic theory to apply to log formal schemes: a \emph{log prism} consists of a bounded log prism $(A,I)$ and a log structure $M_{\Spf A}$ on $\Spf A$ equipped with a so-called $\delta_\log$-structure. For a bounded $p$-adic log formal $(X,M_X)$ scheme with $M_X$ fine and saturated, one can define the \emph{absolute (saturated) prismatic site} $(X,M_X)_\Prism$ consisting of log prisms over $(X,M_X)$ with saturated log structures, as well as the full subcategory $(X,M_X)_\Prism^\str$ of \textit{strict} objects (see Definition~\ref{def:absolute-prismatic-site}); when $M_X$ is trivial, $(X,M_X)_\Prism^\str$ coincides with $X_\Prism$.

\begin{eg}[Breuil--Kisin log prism]\label{eg:intro-BK-log-prisms}
Suppose $X = \Spf R^0$ with 
\[
R^0 = \mathcal{O}_K \langle T_1, \ldots, T_m, T_{m+1}^{\pm 1}, \ldots, T_d^{\pm 1}\rangle / (T_1\cdots T_m - \pi).
\]
Let $M_X$ be the log structure on $X$ given by $\mathbf{N}^d \rightarrow R^0$, $e_i \mapsto T_i$ for $i = 1, \ldots, d$. Consider 
\[
\mathfrak{S} \coloneqq W(k) \langle T_1, \ldots, T_m, T_{m+1}^{\pm 1}, \ldots, T_d^{\pm 1}\rangle[\![u]\!] / (T_1\cdots T_m - u)
\]
equipped with Frobenius via $\varphi(T_i) = T_i^p$. Let $E(u) \in W(k)[u]$ be the monic irreducible polynomial for $\pi$. Then $(\mathfrak{S}, (E(u))) \in X_{\Prism}$ with the structure map $\mathfrak{S}/(E(u)) \cong R^0$. The prelog structure $\mathbf{N}^d \rightarrow \mathfrak{S}$, $e_i \mapsto T_i$, further gives rise to a log prism in $(X, M_X)_{\Prism}$, called the \emph{Breuil--Kisin log prism}. The same construction works when $X$ is \emph{small affine}, namely, when $X$ is affine and admits an \'etale map to $\Spf R_0$.
The Breuil--Kisin log prism will play a central role in our local study.   
\end{eg}

\medskip \noindent
\textbf{Prismatic $F$-crystals.}
The site $(X,M_X)_\Prism$ admits a structure sheaf $\calO_\Prism$ of rings together with a Frobenius endomorphism $\varphi$ and an ideal sheaf $\calI_\Prism$ by assigning $A$, $\varphi_A$, and $I$ to each log prism $(A,I, M_{\Spf A})$. These structures allow one to define various $F$-crystals on $X_\Prism$ and relate them to $\Z_p$-local systems on the generic fiber $X_\eta$.

A \emph{prismatic $F$-crystal} is a pair $(\calE,\varphi_\calE)$ of a vector bundle $\calE$ of $\calO_\Prism$-modules with an identification $\varphi_{\calE}\colon (\varphi^\ast\calE)[\calI_\Prism^{-1}]\xrightarrow{\cong}\calE[\calI_\Prism^{-1}]$. Let $\Vect^\varphi((X,M_X)_\Prism)$ denote the category of prismatic $F$-crystals.
When $X=\Spf \calO_K$ with trivial log structure, Bhatt and Scholze \cite{bhatt-scholze-prismaticFcrystal} prove that this category is equivalent to the category $\operatorname{Rep}_{\Z_p}^\cris(G_K)$ of finite free $\Z_p$-crystalline representations $T$ of $K$ (i.e., $T[p^{-1}]$ being crystalline), and Du and Liu \cite{du-liu-prismaticphiGhatmodule} prove a similar equivalence for the category $\operatorname{Rep}_{\Z_p}^\st(G_K)$ of finite free $\Z_p$-semistable representations of $K$.

\begin{thm}[\cite{bhatt-scholze-prismaticFcrystal, du-liu-prismaticphiGhatmodule}] \label{thm:intro-BhattScholze-DuLiu}
There are natural equivalence of categories:
\begin{enumerate}
    \item $\Vect^\varphi((\Spf\calO_K)_\Prism)\xrightarrow{\cong}
\operatorname{Rep}_{\Z_p}^\cris(G_K)$;
    \item $\Vect^\varphi((\Spf\calO_K, M_{\Spf \calO_K})_\Prism^{\str}) \xrightarrow{\cong}
\operatorname{Rep}_{\Z_p}^\st(G_K)$, where $M_{\Spf \calO_K}$ is given by $\calO_K\cap K^\times\rightarrow \calO_K$.
\end{enumerate}
\end{thm}

\noindent We remark that Yao \cite{yao-lattices-semistable-representations} proves the equivalence between $\operatorname{Rep}_{\Z_p}^\st(G_K)$ and the category of prismatic $F$-crystals on a variant of $(\Spf\calO_K, M_{\Spf \calO_K})_\Prism^{\str}$.

In the case where $X$ is smooth over $\calO_K$, the category of prismatic $F$-crystals is not large enough to capture all the $\Z_p$-crystalline local systems on $X_\eta$. The current authors \cite{du-liu-moon-shimizu-completed-prismatic-F-crystal-loc-system} and Guo and Reinecke \cite{GuoReinecke-Ccris} introduced certain categories which contain the category of prismatic $F$-crystals in vector bundles on $X$, and proved that they are equivalent to the category of crystalline $\Z_p$-local systems on the generic fiber of $X$.

Here we explain the formulation by \cite{GuoReinecke-Ccris}:
$\Vect^{\an,\varphi}(X_\Prism)$ denotes the category of \emph{analytic prismatic $F$-crystals} $(\calE,\varphi_\calE)$ on $X$, where $\calE$ assigns to each $(A,I)\in X_\Prism$ a vector bundle $\calE_A$ over $\Spec A \smallsetminus V(p,I)$, and $\varphi_\Prism$ gives an identification $(\varphi_A^\ast\calE_A)[I^{-1}]\xrightarrow{\cong}\calE_A[I^{-1}]$. 
Each prismatic $F$-crystal yields an analytic $F$-crystal by restricting the induced vector bundle on $\Spec A$ to $\Spec A \smallsetminus V(p,I)$. We remark that completed prismatic $F$-crystals introduced in \cite{du-liu-moon-shimizu-completed-prismatic-F-crystal-loc-system} consists of $(\calE,\varphi)$ with $\calE$ a sheaf of $\calO_\Prism$-modules satisfying the completed crystal property and certain finiteness conditions, and they correspond to effective analytic prismatic $F$-crystals. 

\begin{thm}[\cite{GuoReinecke-Ccris}, \cite{du-liu-moon-shimizu-completed-prismatic-F-crystal-loc-system}] 
There is a natural equivalence of categories 
\[
\Vect^{\an,\varphi}(X_\Prism)\xrightarrow{\cong}
\Loc_{\Z_p}^\cris(X_\eta).
\]
\end{thm}

By the above theorem, one can apply prismatic techniques to study crystalline local systems.
In fact, when $X$ is smooth over $\calO_K$, analytic prismatic $F$-crystals are governed, \'etale locally, by the evaluations over the Breuil--Kisin prism and its self-coproducts. This is an approach taken in \cite{du-liu-moon-shimizu-completed-prismatic-F-crystal-loc-system}, and as an application, \cite{moon-purity} gives a prismatic proof of the purity theorem for crystalline local systems in the smooth case.

\medskip \noindent
\textbf{Prismatic purity theorem.}
Motivated by the smooth case, this paper studies the category $\Vect^{\an,\varphi}((X,M_X)_\Prism)$ of analytic prismatic $F$-crystals in the semistable case. To state the prismatic purity theorem, let us recall another type of prismatic $F$-crystals:
a \emph{Laurent $F$-crystal} is a vector bundle $\calE$ of $\calO_\Prism[\calI_\Prism^{-1}]^\wedge_p$-modules equipped with a Frobenius isomorphism $\varphi_\calE\colon \varphi^\ast\calE\xrightarrow{\cong}\calE$. 
Let $\Vect((X,M_X)_\Prism,\calO_\Prism[\calI_\Prism^{-1}]^\wedge_p)^{\varphi=1}$ denote the category of Laurent $F$-crystals. By restricting a vector bundle on $\Spec A\smallsetminus V(p,I)$ to $\Spec A[I^{-1}]^\wedge_p$ for each log prism $(A,I,M_{\Spf A})$ gives the restriction functor
\[
\Vect^{\an,\varphi}((X,M_X)_\Prism)\rightarrow \Vect((X,M_X)_\Prism,\calO_\Prism[\calI_\Prism^{-1}]^\wedge_p)^{\varphi=1}.
\]
Bhatt and Scholze prove that the category of Laurent $F$-crystals is equivalent to the category of $\Z_p$-local systems on $X_\eta$ \cite[Cor.~3.8]{bhatt-scholze-prismaticFcrystal} in the non-log case (i.e., when $M_X$ is trivial): 
\[
\Vect(X_\Prism,\calO_\Prism[\calI_\Prism^{-1}]^\wedge_p)^{\varphi=1}\xrightarrow{\cong} \Loc_{\Z_p}(X_\eta).
\]
This equivalence continues to hold in the semistable case.

\begin{thm}[Theorem~\ref{thm:logetalerealizationofLaurentFcrystals-strict}, Remark~\ref{rem:strict vs saturated for Laurent F crystals}]
When $(X,M_X)$ is semistable, there is a natural equivalence of categories  
\[
\Vect((X,M_X)_\Prism,\calO_\Prism[\calI_\Prism^{-1}]^\wedge_p)^{\varphi=1}\xrightarrow{\cong} \Loc_{\Z_p}(X_\eta).
\]  
\end{thm}

\noindent We remark that when $(X, M_X)$ is a general bounded fs log $p$-adic formal scheme, Koshikawa and Yao study the category $\Loc_{\Z_p}((X,M_X)_\eta)$ of quasi-pro-Kummer-\'etale $\Z_p$-local systems on the log diamond generic fiber and define a functor $\Vect((X,M_X)_\Prism,\calO_\Prism[\calI_\Prism^{-1}]^\wedge_p)^{\varphi=1}\rightarrow \Loc_{\Z_p}((X,M_X)_\eta)$ (cf. \cite[Thm.~7.36]{koshikawa-yao} and Remark~\ref{rem:KY's equivalence for Laurent F-crystals}).

The \emph{\'etale realization functor} $T_X$ is defined as the composite of the functors 
\[
T_X\colon\Vect^{\an,\varphi}((X,M_X)_\Prism)\rightarrow \Vect((X,M_X)_\Prism,\calO_\Prism[\calI_\Prism^{-1}]^\wedge_p)^{\varphi=1} \xrightarrow{\cong} \Loc_{\Z_p}(X_\eta),
\]
which is fully faithful (Proposition~\ref{prop:etale-realization-fullyfaithful}). We say that a Laurent $F$-crystal \emph{extends} to an analytic $F$-crystal if it lies in the essential image of the restriction functor from $\Vect^{\an,\varphi}((X,M_X)_\Prism)$.

We turn to the prismatic purity theorem. Assume that $(X,M_X)$ is semistable over $\calO_K$ and let $\{\xi_1,\ldots,\xi_m\}\subset X$ be the set of the generic points of all the irreducible components of $X$. For each such point $\xi_j$, the completed local ring $\calO_{X,\xi_j}^\wedge$ of the local ring $\calO_{X,\xi_j}$ is a CDVR with $\pi$ a uniformizer. Set $\Delta_j\coloneqq \Spf \calO_{X,\xi_j}^\wedge$. The morphism $f_j\colon \Delta_j\rightarrow X$ of $p$-adic formal schemes induces a morphism of topoi $(f_j)_\Prism\colon \operatorname{Sh}((\Delta_j,M_{\Delta_j})\rightarrow \operatorname{Sh}((X,M_X)_\Prism)_\Prism)$.

\begin{thm}[Purity for analytic prismatic $F$-crystals: Theorem~\ref{thm:main-purity}] \label{thm:main-intro-purity-analytic-prismatic-F-crystals}
A Laurent $F$-crystal $\calE$ on $(X,M_X)$ extends to an analytic prismatic $F$-crystal on $(X,M_X)$ if and only if its pullback $(f_j)_\Prism^{-1}\mathscr{E}$ extends to an analytic prismatic $F$-crystal on $(\Delta_j,M_{\Delta_j})$ for each $1\leq j\leq m$.
\end{thm}

The proof of Theorem~\ref{thm:main-intro-purity-analytic-prismatic-F-crystals} uses the techniques introduced in \cite{du-liu-prismaticphiGhatmodule, du-liu-moon-shimizu-completed-prismatic-F-crystal-loc-system}, involving local studies using the Breuil--Kisin log prism in Example~\ref{eg:intro-BK-log-prisms} and its self-coproducts. More precisely, we use the observation that analytic prismatic $F$-crystals can be given (\'etale locally) by certain modules over the Breuil--Kisin log prism together with descent data in terms of its self-products, called \textit{Kisin descent data} (Definition~\ref{defn:BK-descent-data}, Lemma~\ref{lem:descentlemmamoduleversion}). 

To show the if-part of Theorem~\ref{thm:main-intro-purity-analytic-prismatic-F-crystals}, we need intersection arguments to produce a Kisin descent data. With the notation as in Example~\ref{eg:intro-BK-log-prisms}, the simplest example of this argument is the equality 
\[
\fkS=\bigcap_{j=1}^m(\fkS[E^{-1}]^\wedge_p\cap\fkS_j),
\]
where $\fkS_j$ is the Breuil--Kisin log prism for the CDVR $(R^0_{(T_j)})^\wedge_p$ (see Lemma~\ref{lem:intersection-basic-rings} for the precise formulation). With a more careful analysis, we are also able to prove a similar equality for the self-coproduct of the Breuil-Kisin log prism (Corollary`\ref{cor:purity-intersection}). For the proof of Theorem~\ref{thm:main-intro-purity-analytic-prismatic-F-crystals}, we work on the intersections of modules with Frobenius over these rings to construct a Kisin descent data, and the Frobenius structure plays a crucial role (\S\ref{sec:proof-purity}, especially, Proposition~\ref{prop:kisin-mod-invert-p-projectivity}).

\medskip \noindent
\textbf{Semistable local systems.}
Let us return to the notion of semistable local systems and Theorem~\ref{thm:intro purity for semistable local systems}. 
We say that a $\Z_p$-local system on a smooth adic space $\calX$ over $K$ with a semistable formal model $X$ is \emph{prismatic $X$-semistable} if it is in the essential image of the \'etale realization functor $T_X\colon \Vect^{\an,\varphi}((X,M_X)_\Prism)\rightarrow \Loc_{\Z_p}(\calX)$. 

Indeed, we also define the notion of prismatic $\Z_p$-local systems on $\Spa(L,\calO_L)$ using the formal model $\Spf \calO_L$ when $\calO_L$ is a CDVR with residue field admitting a finite $p$-basis. In this case, prismatic $\Z_p$-local systems correspond to semistable Galois representations of $L$ defined via the admissibility of the semistable period ring in \cite{fontaine-representations, Morita-imperfsemistable, Ohkubo}. The following generalizes Theorem~\ref{thm:intro-BhattScholze-DuLiu}(2): 

\begin{thm}[Part of Theorem~\ref{thm:CDVR-semistable-notions-equivalent}] \label{thm:intro-equivalence-CDVR}
A $\Z_p$-local system on $\Spa(L,\calO_L)$ is prismatic semistable, namely, arises from an analytic prismatic $F$-crystal if and only if $T[p^{-1}]$ is semistable, where $T$ is the corresponding $\Z_p$-Galois representation of $L$.    
\end{thm}

Combining the prismatic purity (Theorem~\ref{thm:main-intro-purity-analytic-prismatic-F-crystals}) with Theorem~\ref{thm:intro-equivalence-CDVR} and the pullback functoriality of analytic prismatic $F$-crystals, one immediately sees that the notion of prismatic $X$-semistable local systems on $\calX$ does not depend on the choice of a semistable formal model $X$ (Corollary~\ref{cor:independent-model}); we simply call these local systems \emph{prismatic semistable}. Moreover, Theorem~\ref{thm:intro purity for semistable local systems} is an obvious consequence of Theorem~\ref{thm:main-intro-purity-analytic-prismatic-F-crystals}, provided that semistable $\Z_p$-local systems on $\calX$ are interpreted as prismatic semistable local systems.

Let us now discuss other definitions of semistable local systems. Fix a semistable formal model $X$ of $\calX$. As in the work of Faltings \cite{Faltings-CryscohandGalrep}, one can define a notion of semistable local systems on $\calX$ in terms of association to (filtered) $F$-isocrystals on the crystalline site of the mod-$p$ fiber $(X_1,M_{X_1})$ of $(X,M_X)$ (or its convergent site). Unfortunately, it is hard to find a detailed exposition of such a definition in the semistable case. We provide the minimum foundation of the absolute crystalline site and $F$-isocrystals on $(X_1,M_{X_1})$ in Appendix \ref{sec:log-crystalline-site} and formulate the notion of \emph{association} of a $\Z_p$-local system on $\calX$ with an $F$-isocrystal on $(X_1,M_{X_1})_\CRIS$ in Definition~\ref{defn:association}. The resulting class of $\Z_p$-local systems coincides with that of prismatic semistable local systems:

\begin{thm}[{Corollary~\ref{cor:semistable-prismatic-Faltings-equivalent}}] 
A $\Z_p$-local system on $\calX$ is prismatic semistable if and only if it is associated to an $F$-isocrystal on $(X_1,M_{X_1})_\CRIS$.
\end{thm}

We define the \emph{crystalline realization functor} $D_\cris$ sending an analytic prismatic $F$-crystal on $(X,M_X)$ to an $F$-isocrystal on $(X_1,M_{X_1})_\CRIS$ in \S~\ref{sec:crystallinerealization} (Corollary~\ref{cor:crystalline-realization}). Comparing the \'etale and crystalline realizations, we see that a prismatic semistable local system is associated to an $F$-isocrystal on $(X_1,M_{X_1})_\CRIS$. The converse direction in the above theorem is much harder to obtain, and we again resort to the prismatic purity theorem: we treat the CDVR case in Theorem~\ref{thm:CDVR-semistable-notions-equivalent} and deduce it from the CDVR case by Theorem~\ref{thm:main-intro-purity-analytic-prismatic-F-crystals}. 

In fact, Theorem~\ref{thm:CDVR-semistable-notions-equivalent} (including Theorem~\ref{thm:intro-equivalence-CDVR}) becomes a powerful tool, as we have explained, when combined with the prismatic purity theorem. A main step of the proof is to construct Kisin descent data from $\mathbf{Z}_p$-lattices of semistable representations. Such descent data are obtained from the log connections on Breuil--Kisin modules associated to semistable representations (\S\ref{sec:quasi Kisin module and rational descent data}).  

Finally, yet another notion of semistable $\Z_p$-local systems on $\calX$ is given in \cite{andreatta-iovita-semistable-relative} if the fixed semistable formal model $X$ over $\calO_K$ admits a lift along $W[\![u]\!]\twoheadrightarrow \calO_K, u\mapsto \pi$.
For the relation of prismatic semistability with this definition, see Remark~\ref{rem:Andreatta-Iovita-semistability}.

\medskip \noindent
\textbf{Outline.} 
Section~\ref{sec:absolute-log-prismatic-site} explains several preliminary facts on the absolute logarithmic prismatic site: the definition of the site is given in \S\ref{subsec:defn-absolute-log-prismatic-site}. Sections~\ref{sec:Breuil--Kisin log prism} and \ref{sec:CDVR-case-BK-prism-selfprod} concern the Breuil--Kisin log prism and its self-products in the small affine case and CDVR case respectively, and \S\ref{sec:prelim-facts-rings} provides some commutative algebraic facts which will be used in studying the purity results.

In Section~\ref{sec:prismatic-F-crystal-etale-crystalline-realization}, we study analytic prismatic $F$-crystals and their \'etale realizations and crystalline realizations. We define analytic prismatic $F$-crystals and related objects in \S\ref{sec:analytic-prismatic-F-crystals}, and study a local description in terms of Kisin descent data in \S\ref{sec:analytic prismatic F-crystal in small affine case}. Laurent $F$-crystals and their relation to \'etale $\mathbf{Z}_p$-local systems are explained in \S\ref{subsec:Laurent-F-crystals}, where we define the \'etale realizations of analytic prismatic $F$-crystals and study first properties. In \S\ref{subsec:prismatic-semistable-local-system}, the notion of prismatic semistable $\mathbf{Z}_p$-local systems is defined. We construct the crystalline realizations of analytic prismatic $F$-crystals in \S\ref{sec:crystallinerealization}. Section~\ref{subsec: etale-cris-association} studies the notion of association between \'etale $\mathbf{Z}_p$-local systems and $F$-isocrystals on the absolute log crystalline site. In particular, we show that the \'etale and crystalline realizations of every analytic  prismatic $F$-crystal are associated (Proposition~\ref{prop:etaleandcrysassociated} and Corollary~\ref{cor:prismatic-semistable-deRham-semistable-case}).     

Section~\ref{sec:CDVR-semistable-representations} is devoted to proving Theorem~\ref{thm:CDVR-semistable-notions-equivalent} that states the equivalence among the various notions of semistable representations in the CDVR case. In \S\ref{sec:classically-semistable}, we consider the classical definition of semistable representations, and prove that if an \'etale $\mathbf{Z}_p$-local system is associated to an $F$-isocrystal as in \S\ref{subsec: etale-cris-association}, then it is classically semistable. In \S\ref{sec:quasi Kisin module and rational descent data}, we show that every classically semistable \'etale $\mathbf{Z}_p$-local system comes from an analytic prismatic $F$-crystal by constructing the corresponding Kisin descent datum.

Section \ref{sec:purity} proves the purity for analytic prismatic $F$-crystals (Theorem~\ref{thm:main-purity}). The statement and its corollaries are discussed in \S\ref{subsec:purity-statement-applications}, and the proof is given in \S\ref{sec:proof-purity}.

Some complements on log formal schemes and on the absolute log crystalline site are discussed in Appendices~\ref{sec:log-formal-scheme} and \ref{sec:log-crystalline-site} respectively. In Appendix~\ref{app: Functoriality}, we explain the functoriality of the absolute strict logarithmic prismatic sites.

\medskip
\noindent
\textbf{Notation and conventions}.
Throughout the paper, we fix a prime $p$ and let $k$ be a perfect field of characteristic $p$. Write $W=W(k)$ and let $K$ be a finite totally ramified extension of $K_0\coloneqq W[p^{-1}]$. Let $e \coloneqq [K : K_0]$ be the ramification index. Fix a uniformizer $\pi$ of $K$ and let $E(u)\in W[u]$ denote the monic minimal polynomial of $\pi$.
Fix an algebraic closure $\overline{K}$ of $K$ and choose a nontrivial compatible system $(\epsilon_n)_{n \geq 0}$ of $p$-power roots of unity such that $\epsilon_0 = 1$ and $\epsilon_{n+1}^p = \epsilon_n$. Write $\epsilon \coloneqq (\epsilon_n)_{n \geq 0}$ for the corresponding element in $\mathcal{O}_{\overline{K}}^{\flat}$.

We use the standard conventions and terminology regarding $p$-adic log formal schemes, some of which are reviewed in Appendix \ref{sec:log-formal-scheme}. When we work on the monoid $\N^d$, we write $e_i$ ($1\leq i\leq d$) for the element $(0,\ldots,0,1,0,\ldots,0)$ where $1$ is in the $i$-th entry.

A semistable $p$-adic formal scheme over $\calO_K$ refers to a $p$-adic formal scheme over $\calO_K$ that is, \'etale locally, \'etale over an affine formal scheme of the form $\Spf R^0$ with 
\[
R^0 = \mathcal{O}_K \langle T_1, \ldots, T_m, T_{m+1}^{\pm 1}, \ldots, T_d^{\pm 1}\rangle / (T_1\cdots T_m - \pi).
\]
For a semistable $p$-adic formal scheme $X$ over $\calO_K$, we always equip it with the log structure $M_X$ given by the subsheaf associated to the subpresheaf $\calO_{X,\et}\cap (\calO_{X,\et}[p^{-1}])^\times \hookrightarrow \calO_{X,\et}$ (cf.~\cite[1.6(2)]{Cesnavicius-Koshikawa}). Let $X_n$ denote the scheme $X\otimes_{\Z_p}\Z_p/p^n$ equipped with the pullback log structure $M_{X_n}$ of $M_X$.

We generally follow \cite[Notation~3.1]{bhatt-scholze-prismaticFcrystal} for the convention on the generic fiber and $\Z_p$-local systems: 
for a bounded $p$-adic formal scheme $X$, we let $X_\eta$ denote the generic fiber, regarded as a locally spatial diamond. For a locally spatial diamond $\calX$, let $\Loc_{\Z_p}(\calX)$ denote the category of $\Z_p$-local systems on the quasi-pro-\'etale site of $\calX$. However, when $X$ is a semistable $p$-adic formal scheme over $\calO_K$ (or more generally, when $X$ is a bounded $p$-adic formal scheme whose adic generic fiber is representable by a locally noetherian adic space over $\Spa(K,\calO_K)$), we also write $X_\eta$ for the adic generic fiber and let $X_{\eta,\proet}$ denote the pro-\'etale site in the sense of \cite{scholze-p-adic-hodge, Scholze-p-adicHodgeerrata}. In this case, $\Loc_{\Z_p}(X_\eta)$ is naturally identified with the category of $p$-torsion free lisse $\widehat{\Z}_p$-sheaves on $X_{\eta,\proet}$ in the sense of \cite[Def.~8.1]{scholze-p-adic-hodge} by \cite[Prop.~8.2]{scholze-p-adic-hodge}, \cite[Prop.~3.6, 3.7]{MannWerner-Localsystemsondiamonds}, and \cite[Lem.~15.6]{scholze-etalecohomologyofdiamonds}, and we simply call the objects \emph{$\Z_p$-local systems} on $X_\eta$ or $X_{\eta,\proet}$.

For a ring $A$, let $\Vect(A)\coloneqq \Vect(\Spec (A))$ denote the category of vector bundles over $\Spec(A)$. More generally, for an ideal $J$ on $A$, write $\Vect(\Spec(A)\smallsetminus V(J))$ for the category of vector bundles over $\Spec(A)\smallsetminus V(J)$.

For a site $\calC$, write $\operatorname{PSh}(\calC)$ (resp.~$\operatorname{Sh}(\calC)$) for the category of presheaves (resp.~sheaves) on $\calC$.

\medskip
\noindent
\textbf{Acknowledgments.}
We thank Teruhisa Koshikawa for answering questions about absolute saturated and strict logarithmic prismatic sites.  We are also grateful to Kentaro Inoue for valuable feedback, especially, for pointing out errors related to Remark~\ref{rem:log prismatic site when admitting finite free chart}(2) and  a gap in the proof of Lemma~\ref{lem:S-lattice of rational DD}, and helpful discussions on \cite[Thm.~7.36]{koshikawa-yao}.
The first author was partially supported by National Key R\&D Program of China No.~2023YFA1009703. 
The second author prepared this paper during his visit to the University of California, San Diego and the Institute for Advanced Study, and would like to thank both UCSD and IAS for their hospitality. He is supported by the Shiing-Shen Chern Membership during his stay at IAS. 
The third and the fourth authors were partially supported by AMS--Simons Travel Grant at the early stage of the project.

\section{Absolute logarithmic prismatic site} \label{sec:absolute-log-prismatic-site}

\subsection{Definition of the absolute logarithmic prismatic site} \label{subsec:defn-absolute-log-prismatic-site}

We work on log prisms and the logarithmic prismatic sites introduced in \cite{koshikawa, koshikawa-yao}. The notions of $\delta_\mathrm{log}$-rings and $\delta_\mathrm{log}$-triples are introduced in Definition~2.2 and page 12 of \textit{loc.~cit.}: a \emph{(bounded) prelog prism} is a $\delta_{\mathrm{log}}$-triple $(A,I,M)$ where $(A,I)$ is a (bounded) prism. A \emph{log prism} $(A,I,M_{\Spf A})$ consists of a bounded prism $(A, I)$ and a log structure $M_{\Spf A}$ on $(\Spf A)_\et$ together with $\delta_{\mathrm{log}}$-structure that comes from some prelog prism of the form $(A, I, M)$. Write $M_{\Spf (A/I)}$ for the log structure on $\Spf (A/I)$ obtained as the pullback of $M_{\Spf A}$.

For a bounded prelog prism $(A,I,M)$, we have the associated log prism $(A,I,M)^a\coloneqq (A,I,M_{\Spf A}^a)$ (see Definition~3.3 and the following paragraph of \textit{loc.~cit.} for details).

\begin{lem}\label{lem: chart of log prism}
Let $(A, I, M_{\Spf A})$ be a log prism with $M_{\Spf A}$ integral, and assume that $M_{\Spf (A/I)}$ admits a chart $P\rightarrow \Gamma(\Spf (A/I), M_{\Spf (A/I)})$.
Then the monoid homomorphism 
\[
P_A\coloneqq P\times_{\Gamma(\Spf (A/I), M_{\Spf (A/I)})}\Gamma(\Spf A, M_{\Spf A})\rightarrow  \Gamma(\Spf A, M_{\Spf A})
\]
 is a chart, and the map $(A, P_A)\rightarrow (A/I, P)$ is exact surjective with $P_A\rightarrow P$ a $(1+IA)$-torsor of monoids. Moreover, the $\delta_\log$-structure $\delta_\log\colon \Gamma(\Spf A, M_{\Spf A})\rightarrow A$ makes $(A,I, P_A)$ a bounded prelog prism such that $(A,I, P_A)^a \cong (A, I, M_{\Spf A})$.
\end{lem}

Note that a chart always exists for $M_{\Spf (A/I)}$ by the definition of a log prism.

\begin{proof}
By \cite[Lem.~3.8]{koshikawa}, $(A, \Gamma(\Spf A, M_{\Spf A}))\rightarrow (A/I, \Gamma(\Spf (A/I), M_{\Spf (A/I)}))$ is exact surjective, and the monoid map is a $(1+IA)$-torsor. Now the lemma follows from the proof of \cite[Prop.~3.7]{koshikawa}.
\end{proof}

\begin{cor}\label{cor: chart of lifts}
Let $(A,I,M_{\Spf A})$ and $(B,IB,M_{\Spf B})$ be log prisms with integral log structures. Let $\overline{f}\colon(\Spf(B/IB),M_{\Spf(B/IB)}) \to (\Spf(A/I),M_{\Spf(A/I)})$ be a morphism of log formal schemes. Suppose that $M_{\Spf(A/I)}$ admits a chart $P\to \Gamma(\Spf(A/I),M_{\Spf(A/I)})$ and that $\overline{f}$ lifts to a morphism  $f\colon (\Spf B,IB,M_{\Spf B}) \to (\Spf A,I,M_{\Spf A})$ of log prisms.
\begin{enumerate}
\item The monoid map
\[
P_{B,f}\coloneqq P\times_{\Gamma(\Spf(B/IB),\overline{f}^\ast M_{\Spf(A/IA)})}\Gamma(\Spf B,{f}^\ast M_{\Spf A}) \to \Gamma(\Spf B,{f}^\ast M_{\Spf A})
\]
defines a chart for $(\Spf B, f^\ast M_{\Spf A})$ and yields a log prism $(B,IB,f^\ast M_{\Spf A})$ equipped with the factorization of morphisms of log prisms
\[
f\colon (\Spf B,IB,M_{\Spf B}) \rightarrow (\Spf B,IB,f^\ast M_{\Spf A})\rightarrow (\Spf A,I,M_{\Spf A}). 
\]
\item Assume further that $P$ is free. For another prism $(A',I',M_{\Spf A'})$ with $(A/I,M_{\Spf A/I}) \simeq (A'/I',M_{\Spf A'/I'})$ and a morphism of log prisms $f'\colon (\Spf B,IB,M_{\Spf B}) \to (\Spf A',I',M_{\Spf A'})$ lifting $\overline{f}$, 
the pullback log structures $f^\ast M_{\Spf A}$ and $f'^\ast M_{\Spf A'}$ on $\Spf B$ coincide.
\end{enumerate}
\end{cor}

\begin{proof}
(1) Consider the natural monoid homomorphisms $P_A\rightarrow P_{B,f}\rightarrow \Gamma(\Spf B,{f}^\ast M_{\Spf A})\rightarrow B$. Observe that both $P_A$ and $P_{B,f}$ are charts for $f^\ast M_{\Spf A}$. Since $P_A\rightarrow A$ admits the $\delta_\log$-structure and $f$ is a morphism of log prisms, the second assertion also follows.

(2) Write $f_1=f$, $A_1=A$, $f_2=f'$, and $A_2=A'$. For $i=1,2$, we have the following commutative diagram of monoids
\[
\xymatrix{
P_{A_i}\ar[r]\ar@{->>}[d] & P_{B,f_i}\ar[r]\ar@{->>}[d] &B\ar@{->>}[d]\\
P\ar@{=}[r]& P\ar[r]& B/IB.
}
\]
Since $P$ is free, we can take a section $s$ of $P_{A_i}\rightarrow P$. Then $f^\ast M_{\Spf A_i}$ is associated to the monoid map $P\rightarrow P_{A_i}\rightarrow P_{B,f_i}\rightarrow B$ and it factors as $P\xrightarrow{\beta_i} \Gamma(\Spf B,M_{\Spf B})\rightarrow B$. Since $M_{\Spf B}$ is integral, $q\colon \Gamma(\Spf B,M_{\Spf B})\rightarrow \Gamma(\Spf (B/IB),M_{\Spf (B/IB)})$ is a $(1+IB)$-torsor (see \cite[p.3, Ex.~(iii)]{Beilinson-crystalline-period-map} or \cite[Prop.~IV.2.1.2.]{Ogus-log}). On the other hand, both $q\circ \beta_1$ and $q\circ \beta_2$ coincide with $P\rightarrow  \Gamma(\Spf(A/I),M_{\Spf(A/I)})\rightarrow\Gamma(\Spf (B/IB),M_{\Spf (B/IB)})$ induced by $\overline{f}$. Hence $\beta_1(x)$ and $\beta_2(x)$ differ by an element in $1+IB\subset B^\times$ for every $x\in P$, and we conclude that $f_1^\ast M_{\Spf A_1}$ and $f_2^\ast M_{\Spf A_2}$ are canonically identified.
\end{proof}

Let $(X,M_X)$ be a bounded fs $p$-adic log formal scheme.

\begin{defn}[{cf.~\cite[Rem.~4.6]{koshikawa}, \cite[Def.~7.34]{koshikawa-yao}}] \label{def:absolute-prismatic-site}
Let $(X,M_X)_\Prism$ denote the \emph{absolute (saturated) prismatic site} of $(X,M_X)$: an object is a diagram
\[
(\Spf A,M_{\Spf A})\hookleftarrow (\Spf (A/I),M_{\Spf (A/I)})\rightarrow (X,M_X),
\]
where $(A,I,M_{\Spf A})$ is a log prism, $M_{\Spf A}$ is saturated, and $(\Spf (A/I),M_{\Spf (A/I)})\rightarrow (X,M_X)$ is a morphism of log formal schemes that admits a saturated chart \'etale locally (recall that $M_{\Spf (A/I)}$ denotes the pullback log structure from $M_{\Spf A}$). To simplify the notation, we often write an object as $(\Spf A,I, M_{\Spf A})$, or even as $(A,I,M_{\Spf A})$ with the opposite arrows. Morphisms are the obvious commutative diagrams of the following form compatible with structures:
\begin{equation}\label{diag:morphism of prismatic site}
\xymatrix{
(\Spf A',M_{\Spf A'})\ar[d]^-{f} & (\Spf (A'/I'),M_{\Spf (A'/I')}) \ar@{_{(}->}[l] \ar[r] \ar[d]^-{\overline{f}}& (X, M_X)\ar@{=}[d] \\
(\Spf A,M_{\Spf A}) & (\Spf (A/I),M_{\Spf (A/I)})\ar@{_{(}->}[l]\ar[r] & (X, M_X).
}
\end{equation}
 We equip $(X,M_X)_\Prism$ with strict flat topology: a morphism $f\colon (\Spf A',I',M_{\Spf A'})\rightarrow (\Spf A,I,M_{\Spf A})$ is a cover if $\Spf A'\rightarrow \Spf A$ is a $(p,I)$-completely faithfully flat map of formal schemes and the induced map $f^\ast M_{\Spf A}\rightarrow M_{\Spf A'}$ is an isomorphism.

Write $(X,M_X)_\Prism^\str$ for the full subcategory of $(X,M_X)_\Prism$ consisting of objects $(A,I,M_{\Spf A})$ such that $(\Spf (A/I),M_{\Spf (A/I)})\rightarrow (X,M_X)$ is strict. We equip it with the induced topology. If the morphism in \eqref{diag:morphism of prismatic site} is in $(X,M_X)_\Prism^\str$, then $\overline{f}$ is strict; since all the log structures are integral, we conclude that $f$ is also strict by Lemma~\ref{lem: chart of log prism} in this case. In particular, when $M_X$ is trivial, $(X,M_X)_\Prism^\str$ agrees with the absolute prismatic site $X_\Prism$ in \cite[Def.~2.3]{bhatt-scholze-prismaticFcrystal}.
\end{defn}

\begin{rem}\label{rem:fiber product in the strict case}
Consider a diagram
\[
(\Spf A_1,I_1,M_{\Spf A_1})\xrightarrow{f_1} (\Spf A,I,M_{\Spf A}) \xleftarrow{f_2} (\Spf A_2,I_2,M_{\Spf A_2})
\]
in $(X,M_X)_\Prism$ with $f_2$ being strict and $(p,I)$-completely flat. In this case, one can construct the fiber product explicitly as follows: let $B$ denote the derived $(p,I)$-completion of $A_1\otimes^L_AA_2$. By the proofs of \cite[Lem.~3.7]{bhatt-scholze-prismaticcohom} and \cite[Lem.~3.3]{du-liu-moon-shimizu-completed-prismatic-F-crystal-loc-system}, $B$ agrees with the classical $(p,I)$-completion of $A_1\otimes_AA_2$, and $(B,IB)$ represents the pushout of the diagram $(A_1,I_1)\leftarrow (A,I) \rightarrow (A_2,I_2)$ of bounded prisms. We will show that $(\Spf B,IB)$ admits a log prism structure that represents the fiber product.

Take a chart $P\rightarrow \Gamma(\Spf (A_1/I_1), M_{\Spf (A_1/I_1)})$. With the notation as in Lemma~\ref{lem: chart of log prism}, $P_{A_1}$ gives a chart of $M_{\Spf A_1}$. The prelog structure $P_{A_1}\rightarrow A_1\rightarrow B$ comes with a $\delta_\log$-ring structure induced from a prelog prism $(A_1, I_1, P_{A_1})$. Let $(B,IB, M_{\Spf B})\coloneqq (B,IB,P_{A_1})^a$ be the associated log prism. Since $M_{\Spf A_1}$ is integral, so is $M_{\Spf B}$. Since $f_2$ is strict, we obtain a commutative diagram in $(X,M_X)_\Prism$:
\[
\xymatrix{
(\Spf B, IB, M_{\Spf B}) \ar[r]\ar[d] & (\Spf A_2,I_2,M_{\Spf A_2})\ar[d]_-{f_2}\\
(\Spf A_1,I_1,M_{\Spf A_1}) \ar[r]^-{f_1} & (\Spf A,I,M_{\Spf A}).
}
\]
It is now straightforward to see that $(\Spf B, IB, M_{\Spf B})$ represents the fiber product of $f_1$ and $f_2$. If $f_2$ is a cover, then $(\Spf B,IB,M_{\Spf B})\rightarrow (\Spf A_1,I_1,M_{\Spf A_1})$ is also a cover. 
If the underlying prism of each of the three log prisms is perfect, so is $(\Spf B,IB)$.
\end{rem}

\begin{rem}\label{rem:log prismatic site when admitting finite free chart}
    Assume that $M_X$ admits a finite free chart \'etale locally (e.g. $(X,M_X)$ is semistable). 
\begin{enumerate}
    \item 
For a log prism $(A,I,M_{\Spf A})$ such that $M_{\Spf (A/I)}$ admits saturated chart \'etale locally, every morphism $(\Spf (A/I),M_{\Spf (A/I)}) \rightarrow (X, M_X)$ admits a saturated chart \'etale locally by Lemma~\ref{lem:existence of saturated chart for finite free log structure}, making $(A,I,M_{\Spf A})$ an object of $(X,M_X)_\Prism$.

In particular, for every morphism $f\colon (Y,M_Y)\rightarrow (X,M_X)$ from another bounded fs $p$-adic log formal scheme $Y$, the association 
\[
    (Y,M_Y)_\Prism\ni
    (\Spf A, I, M_{\Spf A}, \iota\colon \Spf (A/I)\rightarrow Y)
    \mapsto (\Spf A, I, M_{\Spf A}, f\circ \iota)\in (X,M_X)_\Prism
    \]
    is well-defined and gives a cocontinuous functor. Hence we obtain a morphism of topoi
    \[
    f_\Prism=(f^{-1}_\Prism, f_{\Prism,\ast})\colon \operatorname{Sh}((Y,M_Y)_\Prism)\rightarrow\operatorname{Sh}((X,M_X)_\Prism).
 \]
\item 
Let $(X,M_X)_\Prism^{\str,\perf}$ denote the full subcategory of $(X,M_X)_\Prism^{\str}$ consisting of $(\Spf A,I,M_{Spf A})$ with $(A,I)$ being perfect.
If $(X,M_X)$ is semistable or the log formal spectrum of a CDVR (cf.~\S\ref{sec:CDVR-case-BK-prism-selfprod}), the forgetful functor 
\[
\iota\colon (X,M_X)_\Prism^{\mathrm{str},\perf}\rightarrow X_\Prism^\perf: (\Spf A, I, M_{\Spf A})\mapsto (\Spf A, I)
\]
 is an equivalence by \cite[Prop.~2.18]{min-wang-HT-crys-log-prism}: for every perfect prism $(A,I)\in X_\Prism$, there is a unique log structure $M_{\Spf A}$ on $\Spf A$ so that $(\Spf A,I,M_{\Spf A})$ is inside $(X,M_X)_\Prism^{\mathrm{str}}$.
In particular, it induces an equivalence $\operatorname{Sh}((X,M_X)_\Prism^{\mathrm{str},\perf})\xrightarrow{\cong}\operatorname{Sh}(X_\Prism^\perf)$ of topoi.
\end{enumerate}    
\end{rem}

\subsection{Breuil--Kisin log prism and its self-products} \label{sec:Breuil--Kisin log prism}

This subsection collects some preliminary facts on the Breuil--Kisin log prisms, which play a central role in studying prismatic $F$-crystals in the small affine case. 

For $d\geq m\geq 1$, set
\[
R^0\coloneqq \mathcal{O}_K \langle T_1, \ldots, T_m, T_{m+1}^{\pm 1}, \ldots, T_d^{\pm 1}\rangle / (T_1\cdots T_m - \pi),
\]
and consider the prelog structure $\mathbf{N}^d \rightarrow R^0$ sending $e_i \mapsto T_i$ for $1 \leq i \leq d$.
Let $R$ be a \emph{connected} $\mathcal{O}_K$-algebra equipped with a $p$-adically completed \'etale map 
\[
\square\colon R^0 \rightarrow R.
\]
We call a $p$-adic log formal scheme $(X,M_X)$ \emph{small affine with framing $\square$} if it is an affine log formal scheme of the form $(X,M_X)=(\Spf R, \N^d\rightarrow R^0\xrightarrow{\square} R)^a$.

Let $\fkS_{R^0}\coloneqq W(k)\langle T_1, \ldots, T_m, T_{m+1}^{\pm 1}, \ldots, T_d^{\pm 1}\rangle[\![u]\!] / (T_1\cdots T_m-u)$. Note that $\fkS_{R^0}$ is $(p,E(u))$-adically complete and $\fkS_{R^0}/(E(u)) \cong R^0$. We equip $\fkS_{R^0}$ with the prelog structure $\mathbf{N}^d \rightarrow \mathfrak{S}_{R^0}$, $e_i \mapsto T_i$ ($1 \leq i \leq d$) lifting the one on $R^0$. Since $\square$ is $p$-adically completed \'etale, it lifts uniquely to a $(p,E(u))$-adically completed \'etale map 
\[
\square_{\fkS}\colon \fkS_{R^0}\ra \fkS_{R,\square}.
\]
If there is no confusion, we simply write $\fkS$ for $\fkS_{R,\square}$.

Setting $\delta_{\mathrm{log}}(e_i) = 0 $ and $\delta(T_j) = 0$ makes $(\fkS_{R^0}, \delta, (E(u)), \N^d, \delta_{\mathrm{log}})$ a bounded prelog prism; we will also write it as $(\fkS_{R^0}, (E(u)), \N^d)$ for short. By \cite[Lem.~2.13]{koshikawa}, $(\fkS, (E(u)), \N^d\rightarrow \fkS_{R^0}\rightarrow \fkS)$ admits a bounded prelog prism structure.

\begin{defn}
The \emph{Breuil--Kisin log prism} $(\Spf \mathfrak{S}, (E(u)), M_{\Spf \mathfrak{S}})$ (with respect to $\pi$ and $\square$) is the log prism associated to $(\mathfrak{S}, (E(u)), \N^d)$. It is an object of $((\Spf R,\N^d)^a)_\Prism^\str$ via $R \cong \fkS/(E(u))$.
\end{defn}

\begin{lem}\label{lem:varphionfkSisfaithfullyflat}
The map $\varphi_{\fkS}\colon\fkS \to \fkS$ is faithfully flat, and $\fkS$ is finite free over $\varphi_{\fkS}(\fkS)$.
\end{lem}

\begin{proof}
Since $\mathfrak{S}$ is a $p$-torsion free $p$-adically complete Noetherian ring and $\mathfrak{S}/(p)$ is regular with a finite $p$-basis, the proof of \cite[Lem.~2.7]{du-liu-moon-shimizu-completed-prismatic-F-crystal-loc-system} works. 
\end{proof}

\begin{lem}\label{lem:BKcoversfinalobject}
Let $(X,M_X)=(\Spf R,\N^d)^a$ be small affine.
Then $(\Spf\mathfrak{S}, (E(u)), M_{\Spf \mathfrak{S}})$ covers the final object of the topos associated to $(X,M_X)_\Prism$. In particular, the same conclusion holds for $(X,M_X)_\Prism^{\mathrm{str}}$.
\end{lem}

We note that the statement for $(X,M_X)_\Prism^{\mathrm{str}}$ is proved in \cite[Lem.~3.16]{min-wang-HT-crys-log-prism}.

\begin{proof}
Let $R_{\infty}$ be $\mathfrak{S}[T_1^{1/p^{\infty}}, \ldots, T_d^{1/p^{\infty}}, u^{1/p^{\infty}}]\widehat{\otimes}_{\mathfrak{S}} R$, where $\widehat{\otimes}$ denotes $p$-completed tensor product. Then $R \rightarrow R_{\infty}$ is a quasi-syntomic cover since it is obtained as the $p$-completion of the colimit of extracting $p$-power roots. Let $A$ be the $u$-completion of $(\mathfrak{S}/(p))_{\mathrm{perf}}$ and $\xi \coloneqq E([u^{\flat}])$ where $u^{\flat} \in A^{\flat}$ is the element corresponding to $u$ under $A^\flat \xrightarrow{\cong} A$. We have $R_{\infty} \cong W(A^\flat)/(\xi)$. Thus, $R_{\infty}$ is a perfectoid ring by \cite[Lem.~3.6]{bhatt-iyengar-ma-regular-rings-perfectoid}. Write $\Prism_{R_\infty}\coloneqq W((R_\infty)^\flat)$ and let $(\Prism_{R_\infty}, I,\N^d)$ be the prelog prism where $\alpha\colon \N^d \to \Prism_{R_\infty}$ is given by $e_i \mapsto [T_i^\flat]$ with $T_i^\flat \in R^\flat_\infty$ defined by a fixed choice of compatible $p$-power roots of $T_i$ in $R_\infty$. 
The choice of $T_i^\flat$ defines a map of prelog prisms $(\fkS,(E),\N^d) \to (\Prism_{R_\infty},I,\N^d)$.

We are going to show that $(\Spf\Prism_{R_\infty}, I, \N^d)^a$ covers the final object in each topos in the statement. Take any $(\Spf B,J,M_{\Spf B})\in(X,M_X)_\Prism$.  By the proof of Lemma~\ref{lem:existence of saturated chart for finite free log structure}, we may further assume that the map $(\Spf(B/J), M_{\Spf(B/J)}) \to (X,M_X)$ admits a saturated chart which fits into the following commutative diagram as in \eqref{diag:chart}
\[
\begin{tikzcd}
\N^d \arrow[r]\arrow[d] & \Gamma(X,M_{X}) \arrow[d] \\
Q \arrow[r] & \Gamma(\Spf(B/J),M_{\Spf(B/J)})
\end{tikzcd}
\]
such that the composite $\N^d\rightarrow \Gamma(X,M_{X})\rightarrow R$ sends $e_i$ to $T_i$.
Consider the quasi-syntomic cover $B/J \to (B/J)\widehat{\otimes}_{R} R_{\infty}$. By \cite[Prop.~7.11(1)]{bhatt-scholze-prismaticcohom}, there exists a prism $(C,JC)$ over $X$ which covers $(B, JB)$ such that the structure map $R\to C/JC$ factors through $R \to (B/J)\widehat{\otimes}_{R} R_{\infty}$. In particular, $C/JC$ admits a map from $R_{\infty}$. Let $(C, JC, M_{\Spf C})$ be the log prism with the pullback log structure $M_{\Spf C}$ from $(B,J,M_{\Spf B})$. Then $(B,J,M_{\Spf B}) \to (C,JC,M_{\Spf C})$ is a cover by definition. We also have morphisms of log formal schemes over $(X,M_X)$
\[
(\Spf B/J,M_{\Spf B/J})\leftarrow (\Spf C/JC,M_{\Spf C/JC})\xrightarrow{\overline{f}}(\Spf R_\infty,(\N^d)^a).
\]
It suffices to show that there is a map $(\Prism_{R_\infty}, I, \N^d)^a \to (C,JC,M_{\Spf C})$ of log prisms lifting $\overline{f}$.
Since $Q\rightarrow B/J\rightarrow C/JC$ defines the log structure $M_{\Spf C/JC}$, one can consider the chart $Q_C$ of $M_{\Spf C}$ as in Lemma~\ref{lem: chart of log prism}. For each $D\in\{C,C/JC, \Prism_{R_\infty}, R_{\infty}\}$, set $M_D\coloneqq \Gamma(\Spf D, M_{\Spf D})$. Then $M_D$ is an integral monoid, and $(D,M_D\rightarrow D)$ is a log ring in the sense of \cite[Convention]{koshikawa-yao}. Moreover, the monoid map $M_C\rightarrow M_{C/JC}$ (resp.~$M_{\Prism_{R_\infty}}\rightarrow M_{R_\infty}$) is a torsor under $1+JC$ (resp.~$1+I\Prism_{R_\infty}$) by the proof of Lemma~\ref{lem: chart of log prism}. By \cite[Lem.~2.40]{koshikawa-yao}, the map $(R_\infty, M_{R_{\infty}})\rightarrow (C/JC,M_{C/JC})$ of log rings lifts uniquely to a map $(\Prism_{R_\infty}, I, M_{\Prism{R_\infty}}^\flat)\rightarrow (C,JC,M_C)$ of pre-log prisms, where $M_{\Prism_{R_\infty}}^\flat\coloneqq \varprojlim_{x\mapsto x^p}M_{\Prism_{R_\infty}}$. By definition of $\alpha\colon \N^d\rightarrow \Prism_{R_\infty}$, the induced map $\N^d\rightarrow M_{\Prism{R_\infty}}$ factors as $\N^d\rightarrow M_{\Prism{R_\infty}}^\flat\rightarrow M_{\Prism{R_\infty}}$ and yields a map $(\Prism_{R_\infty}, I, \N^d)\rightarrow (C,JC,M_C)$ of pre-log prisms. Now it remains to prove that the resulting monoid map $\N^d\rightarrow M_C$ factors through $Q_C\rightarrow M_C$. To see this, consider the following commutative diagram of monoids
\[
\xymatrix{
Q_C\ar[r]\ar[d]& M_{\Prism_{R_\infty}}\ar[d]&\N^d\ar[l]\ar@{=}[d]\\
Q\ar[r]& M_{R_\infty}&\N^d.\ar[l] \ar@/^1pc/[ll]
}
\]
Since the left and middle vertical morphisms are torsors under $1+JC$, it follows from the commutativity of the bottom diagram and the definition of $Q_C$ that $\N^d\rightarrow M_C$ factors through $Q_C\rightarrow M_C$. 
\end{proof}

We next discuss the existence of the product with the Breuil--Kisin log prism inside $(X,M_X)_\Prism^\mathrm{str}$. 

\begin{lem}\label{lem:coproduct with fkS is a covering}
Keep the small affine assumption and let $(A,I,M_{\Spf A})$ be any log prism in the opposite category $(X,M_X)_\Prism^{\mathrm{str},\mathrm{op}}$. Then the coproduct $(A,I,M_{\Spf A})\amalg (\mathfrak{S}, (E(u)), M_{\Spf \mathfrak{S}})$ of $(A,I,M_{\Spf A})$ and $(\mathfrak{S}, (E(u)), M_{\Spf \mathfrak{S}})$ exists in  $(X,M_X)_\Prism^{\mathrm{str},\mathrm{op}}$, and the structure map 
\[
(A,I,M_{\Spf A}) \to (A,I,M_{\Spf A})\amalg (\mathfrak{S}, (E(u)), M_{\Spf \mathfrak{S}}) 
\]
is a cover.
\end{lem}

\begin{proof}
We may assume that the prelog prism $(A, I, M_A)$ is orientable so that $I = (r)$ for some $r \in A$. Let 
\[
B_0 = A[S_1, \ldots, S_m][S_{m+1}^{\pm 1}, \ldots, S_{d}^{\pm 1}][v]/(S_1\cdots S_m -v)
\]
equipped with the $\delta$-structure given by $\delta(S_i) = 0 = \delta(v)$. Let $(B_0, M_B)$ be the associated prelog ring where $M_B\coloneqq M_A \times \mathbf{N}^d \rightarrow B_0$ is given by $M_A \rightarrow A \rightarrow B_0$ and $e_i \mapsto S_i$ ($1\leq i\leq d$). Consider the surjection $(B_0, M_B) \rightarrow (A/I, M_{A/I})$ induced by $S_i \mapsto T_i$ (so that $v \mapsto \pi$) and let $J\subset B_0$ be the kernel. Denote by $(B, M_B)$ the $(p, I)$-completed universal $\delta_{\mathrm{log}}$-ring over $(A, M_A)$ generated by $(B_0, M_B)$. Then by \cite[Prop.~3.9]{koshikawa}, the prelog prismatic envelope $(B', IB', M_{B'})$ of $(B, (JB)^{\wedge}_{(p, I)}, M_B)$ with exact surjection $(B', IB', M_{B'}) \rightarrow (B'/IB', M_{A/I})$ exists, and it is bounded and $(p, I)$-completely flat over $(A, I, M_A)$. Note that $(B', IB', M_{B'})^a$ gives the log prismatic envelope by \cite[Prop.~3.7]{koshikawa}. 

Note $J/I = (S_1-T_1, \ldots, S_d-T_d)$ and $E(v) \in J$. Since $r$ divides $E(v)$ in $B'$ and $B'$ is $(p, r)$-complete, we have a natural map of prelog prisms $(\mathfrak{S}, (E(u)), M_{\mathfrak{S}}) \rightarrow (B', IB', M_{B'})$ given by $T_i \mapsto S_i$ (so that $u \mapsto v$), and this induces a map of the associated log prisms in $(X, M_X)_{\Prism}$. By the universal property of log prismatic envelope, $(B', IB', M_{B'})^a$ is the coproduct of $(A,I,M_{\Spf A})$ and $(\mathfrak{S}, (E(u)), M_{\Spf \mathfrak{S}})$ in $(X, M_X)_{\Prism}^\mathrm{op}$. 

Finally, let us show that  $A \rightarrow B'$ is $(p, I)$-completely faithfully flat by following \cite[Prop.~2.4.9]{Bhatt-Lurie-absoluteprismticcohomology}. 
Take $x \in \Spec (A/(p, I)) \subset \Spec  A$, and let $k_1$ be a perfect field containing $\kappa(x)$. The map $A \rightarrow k_1$ lifts uniquely to a map of $\delta$-rings $A \rightarrow W(k_1)$, which gives a map of prisms $(A, I) \rightarrow (W(k_1), (p))$. The associated log prism $(W(k_1), (p), M_A\rightarrow A\rightarrow W(k_1))^a$ is an object of $(X,M_X)_\Prism^\mathrm{op}$. On the other hand, since $\Spec (\mathfrak{S}/(p, E(u))) \neq \emptyset$, we have a map of commutative rings $\mathfrak{S}/(p, E(u)) \rightarrow k_1$ after enlarging $k_1$ if necessary. Similarly, this induces a map $(\mathfrak{S}, (E(u)), M_{\mathfrak{S}})^a \rightarrow (W(k_1), (p), M_{\mathfrak{S}})^a$ in $(X,M_X)_\Prism^\mathrm{op}$. Since $(W(k_1), (p))$ is a perfect prism, we have $(W(k_1), (p), M_A)^a = (W(k_1), (p), M_{\mathfrak{S}})^a$ by \cite[Lem.~2.17]{min-wang-HT-crys-log-prism}, and thus, we obtain a map $(B', IB', M_{B'})^a \rightarrow (W(k_1), (p), M_A)^a$. The point in $\Spec  (B'/(p, I)B')$ given by $B'/(p, I)B' \rightarrow W(k_1)/(p) = k_1$ lifts $x$ by construction, so $A \rightarrow B'$ is $(p, I)$-completely faithfully flat. 
\end{proof}

For two framings $\square$ and $\square'$ of $R$, let $(\fkS_{\square,\square'},(E),M_{\square,\square'})$ be the coproduct of $(\fkS_{R,\square},(E),M_{\Spf\fkS_{R,\square}})$ and $(\fkS_{R,\square'},(E),M_{\Spf\fkS_{R,\square'}})$ in $(X,M_X)_\Prism^{\mathrm{str},\mathrm{op}}$. When the framing is fixed, we also write $(\mathfrak{S}^{(i)},(E),\N^d)^a$ for the $(i+1)$-st self-coproduct of $(\mathfrak{S},(E),\N^d)^a$ in $(X,M_X)_\Prism^{\mathrm{str},\mathrm{op}}$, where the prelog structure $\N^d\rightarrow \fkS^{(i)}$ is assumed to be given by $\N^d\rightarrow \fkS\xrightarrow{\text{1st}}\fkS^{(i)}$. 

The proof of Lemma~\ref{lem:coproduct with fkS is a covering} describes the self-product of the Breuil--Kisin log prism explicitly: set $B_0 = \mathfrak{S}[S_1, \ldots, S_m][S_{m+1}^{\pm 1}, \ldots, S_{d}^{\pm 1}]$ as in the proof. Let $B^{(1)} = \mathfrak{S}[\![1-\frac{T_{i, 2}}{T_{i, 1}}]\!]_{1\leq i \leq d}$; this is identified with the $(1-\frac{S_1}{T_1}, \ldots, 1-\frac{S_d}{T_d})$-completion of $B_0[\frac{S_1}{T_1}, \ldots, \frac{S_m}{T_m}]$ with $T_{i, 1} \coloneqq T_i$ and $T_{i, 2} \coloneqq S_i$. Since the kernel of the surjection $B^{(1)}\rightarrow R$ is given by $J^{(1)}\coloneqq (E(u), 1-\frac{T_{i, 2}}{T_{i, 1}})_{1\leq i\leq d}$, we have
\[
\mathfrak{S}^{(1)} = B^{(1)}\biggl\{\frac{1-\frac{T_{1, 2}}{T_{1, 1}}}{E(u)}, \ldots, \frac{1-\frac{T_{d, 2}}{T_{d, 1}}}{E(u)}\biggr\}^{\wedge}_{\delta}. 
\]
Set $w_i \coloneqq \frac{1-\frac{T_{i ,2}}{T_{i, 1}}}{E(u)}$. For $j = 1, 2$, define the map $p^1_j\colon \mathfrak{S} \rightarrow B^{(1)}$ by $T_i \mapsto T_{i, j}$. Then $p^1_j$ extends to the $j$-th projection map $\mathfrak{S} \rightarrow \mathfrak{S}^{(1)}$. More generally for $i \geq 1$ and $1 \leq j \leq i$, write $p_j^i\colon \mathfrak{S} \rightarrow \mathfrak{S}^{(i)}$ for the $j$-th projection map.  

\begin{lem}\label{lem:coproduct in two sites are the same} 
Let $(X,M_X)=(\Spf R,\N^d)^a$ be small affine. For two framings $\square$ and $\square'$ of $R$, $(\fkS_{\square,\square'},(E),M_{\square,\square'})$ represents the coproduct of $\fkS_{R,\square}$ and $\fkS_{R,\square'}$ in $(X,M_X)_\Prism^\mathrm{op}$.
Similarly, when the framing is fixed, $(\mathfrak{S}^{(i)},(E),\N^d)^a$ represents the $(i+1)$-st self-coproduct of $(\mathfrak{S},(E),\N^d)^a$ in $(X,M_X)_\Prism^\mathrm{op}$ for each $i\geq 1$.
\end{lem}

\begin{proof}
We only prove the first part of the statement in the lemma. The second part follows from a similar argument. Write $\fkS_1$ for $\fkS_{R,\square}$ and $\fkS_2$ for $\fkS_{R,\square'}$, and let $q_i\colon (\Spf\fkS_{\square,\square'},(E),M_{\square,\square'}) \to (\Spf\fkS_i,(E), M_{\Spf \fkS_i})$ be the projection in $(X,M_X)_\Prism^{\mathrm{str}}$. Take any $(\Spf A, I, M_{\Spf A}) \in (X,M_X)_\Prism$ and a map $f_i\colon (\Spf A, I, M_{\Spf A})\rightarrow(\Spf\mathfrak{S}_i,(E),M_{\Spf\mathfrak{S}_i})$ for $i=1, 2$. Since $(\Spf(\fkS_1/E),M_{\Spf(\fkS_1/E)})\simeq (X,M_X) \simeq (\Spf(\fkS_2/E),M_{\Spf(\fkS_2/E)})$ which admits a free chart, Corollary~\ref{cor: chart of lifts} implies that $f_1^\ast M_{\Spf \fkS_1}$ and $f_2^\ast M_{\Spf \fkS_2}$ define the same log structure on $\Spf A$, which we denote by $N_{\Spf A}$. Moreover, $f_1$ and $f_2$ can be written as $f_i=\tilde{f}_i\circ h$ with $h\colon (\Spf A, I, M_{\Spf A}) \to (\Spf A, I, N_{\Spf A})$ and  $\tilde{f}_i\colon (\Spf A, I, N_{\Spf A})\rightarrow(\Spf\mathfrak{S}_i,(E),M_{\Spf\mathfrak{S}_i})$ being strict. By the definition of $(\fkS_{\square,\square'},(E),M_{\square,\square'})$, there exists a unique map $b\colon (\Spf A, I, N_{\Spf A}) \to (\Spf \fkS_{\square,\square'},(E),M_{\square,\square'})$ such that $b \circ q_i = \widetilde{f_i}$ for $i=1, 2$. Thus, we have the following commutative diagram:
\[
\begin{tikzcd}
 & & (\Spf\mathfrak{S}_1,(E),M_{\Spf\mathfrak{S}_1})\\
(\Spf A, I, M_{\Spf A}) \arrow[r,"h"]\arrow[urr,"f_1",bend left=12] \arrow[drr,"f_2"',bend right=12] & (\Spf A, I, N_{\Spf A}) \arrow[r,"b"] \arrow[ur,"\widetilde{f_1}"] \arrow[dr,"\widetilde{f_2}"'] & (\Spf\fkS_{\square,\square'},(E),M_{\square,\square'}) \arrow[u,"q_1"']\arrow[d,"q_2"] \\
 & & (\Spf\mathfrak{S}_2,(E),M_{\Spf\mathfrak{S}_2})
\end{tikzcd} 
\]
This implies that $(\fkS_{\square,\square'},(E),M_{\square,\square'})$ is also the coproduct of $(\mathfrak{S}_1,(E),M_{\Spf \mathfrak{S}_1})$ and $(\mathfrak{S}_2,(E),M_{\Spf \mathfrak{S}_2})$ in $(X,M_X)_\Prism^\mathrm{op}$. 
\end{proof}

Recall that a prism $(A, I)$ is called \emph{transversal} if $A/I$ is $p$-torsion free \cite[Def.~2.1.3]{Bhatt-Lurie-absoluteprismticcohomology}; if $(A, I)$ is transversal, then $A$ is $p$-torsion free by \cite[Rem.~2.1.7]{Bhatt-Lurie-absoluteprismticcohomology}.

\begin{lem} \label{lem:transv-prism}
For an orientable transversal prism $(A, (d))$, we have 
\[
A[p^{-1}] \cap A[d^{-1}] = A,
\]
and $A[d^{-1}]$ is $p$-adically separated. Consequently, 
\[
A[p^{-1}] \cap A[d^{-1}]^{\wedge}_p = A.
\]
\end{lem}

\begin{proof}
The first statement follows from the injectivity of $A/(d) \rightarrow (A/(d))[p^{-1}]$. Since $A$ is classically $p$-complete, we deduce the second assertion from the first one. Finally, the  map $A[d^{-1}]/pA[d^{-1}] \rightarrow A[d^{-1}]^{\wedge}_p/pA[d^{-1}]^{\wedge}_p$ is an isomorphism (cf.~\cite[Rem.~3.15(i)]{du-liu-moon-shimizu-completed-prismatic-F-crystal-loc-system}), so the last statement follows. 
\end{proof}

\begin{cor} \label{cor:proj-map-faithful-flat}
The map $p^i_j\colon \mathfrak{S} \rightarrow \mathfrak{S}^{(i)}$ is classically faithfully flat for each $i \geq 1$ and $1 \leq j \leq i$. In particular, $(\mathfrak{S}^{(i)},(E))$ is transversal for all $i$, $\mathfrak{S}^{(i)}[E^{-1}]$ is $p$-adically separated, and 
\[
\mathfrak{S}^{(i)}[p^{-1}] \cap \mathfrak{S}^{(i)}[E^{-1}]^{\wedge}_p = \mathfrak{S}^{(i)}.
\]
\end{cor}

\begin{proof}
By inductively applying Lemma~\ref{lem:coproduct with fkS is a covering}, we see that $p^i_j\colon \mathfrak{S} \rightarrow \mathfrak{S}^{(i)}$ is $(p,E)$-completely faithfully flat and thus classically faithfully flat by \cite[Tag~0912]{stacks-project}. Hence each $\mathfrak{S}^{(i)}$ is transversal, and the last statement follows from Lemma~\ref{lem:transv-prism}.
\end{proof}

\begin{lem} \label{lem:modE-modp-complete}
The ring $\mathfrak{S}^{(1)}/(E)$ is classically $p$-complete, and $\mathfrak{S}^{(1)}/(p)$ is classically $E$-complete.   
\end{lem}

\begin{proof}
Consider the exact sequence
\begin{equation} \label{eq:SES-before-completion}
0 \rightarrow E\mathfrak{S}^{(1)} \rightarrow \mathfrak{S}^{(1)} \rightarrow \mathfrak{S}^{(1)}/(E) \rightarrow 0.   
\end{equation}
For each $n \geq 1$, $\mathfrak{S}^{(1)}/(E)$ is $p^n$-torsion free by Corollary~\ref{cor:proj-map-faithful-flat}, so the induced sequence
\[
0 \rightarrow E\mathfrak{S}^{(1)}/p^n E\mathfrak{S}^{(1)} \rightarrow \mathfrak{S}^{(1)}/(p^n) \rightarrow \mathfrak{S}^{(1)}/(E, p^n) \rightarrow 0
\]
is exact. By \cite[Tag~03CA]{stacks-project}, the above exact sequence (\ref{eq:SES-before-completion}) remains exact after classical $p$-completion. Since $E\mathfrak{S}^{(1)}$ is a free $\mathfrak{S}^{(1)}$-module of rank $1$, it is classically $p$-complete. Thus, $\mathfrak{S}^{(1)}/(E)$ is classically $p$-complete.
The statement for $\mathfrak{S}^{(1)}/(p)$ follows similarly.
\end{proof}

\begin{prop} \label{prop:frakS^(1)/E}
We have a natural isomorphism
\[
\mathfrak{S}^{(1)}/(E) \cong R \{w_1,\ldots,w_d\}^\wedge_\mathrm{PD},
\]
where the latter ring denotes the $p$-completed free PD-polynomial algebra over $R$ in variables $w_j$'s.
\end{prop}

\begin{proof}
This is \cite[Prop.~2.2.8(2)]{du-liu-prismaticphiGhatmodule}. For the convenience of the reader, we briefly recall the outline. Let $B^{(1)}[w_1^{(0)},\ldots,w_d^{(0)}]_{\delta}$ be the free $\delta$-ring over $B^{(1)}\coloneqq \mathfrak{S}[\![1-\frac{T_{i, 2}}{T_{i, 1}}]\!]_{1\leq i \leq d}$ in variables $w_j^{(0)}$'s. Set $J = (\delta^n(Ew_j^{(0)})-\delta^n(1-\frac{T_{j, 2}}{T_{j, 1}}))_{n \geq 0,\,1\leq j \leq d}$ and consider $\mathfrak{S}^{(1), \mathrm{nc}} \coloneqq B^{(1)}[w_1^{(0)},\ldots,w_d^{(0)}]_{\delta} / J$. One can check $\delta^n(1-\frac{T_{j, 2}}{T_{j, 1}}) \equiv 0 \pmod E$ and thus $\delta^n(Ew_j^{(0)}) \equiv 0 \pmod E$ in $\mathfrak{S}^{(1), \mathrm{nc}}$ for every $n\geq 0$. Using \cite[Lem.~2.2.5]{du-liu-prismaticphiGhatmodule}, one can define an isomorphism $\mathfrak{S}^{(1), \mathrm{nc}}/(E)\xrightarrow{\cong} R[\gamma_n(w_j)]_{n \geq 0, 1 \leq j \leq d}$, by inductively assigning the image of $\delta^n(w_j^{(0)})$ (see \cite[Prop.~2.2.8(1)]{du-liu-prismaticphiGhatmodule}).
Note that $\mathfrak{S}^{(1)}$ is naturally isomorphic to the classically $(p, E)$-completion of $\mathfrak{S}^{(1), \mathrm{nc}}$. Since $\mathfrak{S}^{(1)}/(E)$ is classically $p$-complete by Lemma~\ref{lem:modE-modp-complete}, it is the $p$-adic completion of $\mathfrak{S}^{(1), \mathrm{nc}}/(E)$ and thus we obtain
\[
\mathfrak{S}^{(1)}/(E) \cong 
 \varprojlim_{l} \mathfrak{S}^{(1), \mathrm{nc}}/(E, p^{l})\cong \varprojlim_{l} R[\gamma_n(w_1),\ldots, \gamma_n(w_d)]_{n \geq 0}/(p^{l}) = R \{w_1,\ldots, w_d \}^\wedge_\mathrm{PD}.    
\]
\end{proof}

\begin{eg}\label{eg:coproduct of Frob twists of BK prism}
One can also consider the self-products of the Frobenius twist of the Breuil--Kisin log prism, which will be used in the crystalline realization functor in \S\ref{sec:crystallinerealization}.

Consider $(\fkS,(\varphi(E)),\varphi^\ast\N^d)^a\in (X,M_X)_\Prism$ with the prelog structure 
$\varphi^\ast\N^d$ given by $\N^d\rightarrow \fkS: e_i \mapsto T_i^p$. Then $\varphi\colon (\fkS,(E),\N^d)^a \to (\fkS,(\varphi(E)),\varphi^\ast\N^d)^a$ defines a cover by Lemma~\ref{lem:varphionfkSisfaithfullyflat}. Let $\varphi_1^{(1)}\colon(\fkS^{(1)},(E),\N)^a \to (\widetilde{\fkS}^{(1)},(\varphi(E)),\varphi^\ast\N^d)^a$ be the pushout of $\varphi\colon(\fkS,(E),\N^d)^a \to (\fkS,(\varphi(E)),\varphi^\ast\N^d)^a$ along $p_1\colon (\fkS,(E),\N^d)^a \to (\fkS^{(1)},(E),\N)^a$ given by Remark~\ref{rem:fiber product in the strict case}, and consider 
\[
\widetilde{p}_2\colon (\fkS,(E),\N^d)^a \xrightarrow{p_2} (\fkS^{(1)},(E),\N)^a \xrightarrow{\varphi_1^{(1)}} (\widetilde{\fkS}^{(1)},(\varphi(E)),\varphi^\ast\N^d)^a.
\]
Write $(\widetilde{\fkS}^{(1)},(\varphi(E)),\varphi^\ast\N^d)^a \to (\fkS^{[1]},(\varphi(E)),\varphi^\ast\N^d)^a$ for the pushout of $\varphi\colon(\fkS,(E),\N^d)^a \to (\fkS,(\varphi(E)),\varphi^\ast\N^d)^a$ along $\widetilde{p}_2$. We summarize these constructions in the following commutative diagram with pushout squares:
\[
\xymatrix{
&(\fkS,(E),\N^d)^a\ar[d]_-{p_2}\ar[r]^-\varphi\ar@/^3.5pc/[dd]^-{\widetilde{p}_2}&(\fkS,(\varphi(E)),\varphi^\ast\N^d)^a\ar[dd]\\
(\fkS,(E),\N^d)^a\ar[d]_-\varphi\ar[r]^-{p_1}&(\fkS^{(1)},(E),\N)^a\ar[d]_-{\varphi_1^{(1)}}&\\
(\fkS,(\varphi(E)),\varphi^\ast\N^d)^a\ar[r]&(\widetilde{\fkS}^{(1)},(\varphi(E)),\varphi^\ast\N^d)^a\ar[r]&(\fkS^{[1]},(\varphi(E)),\varphi^\ast\N^d)^a.
}
\]
Then $(\fkS^{[1]},(\varphi(E)),\varphi^\ast\N^d)^a$ has two natural maps $(\fkS,(\varphi(E)),\varphi^\ast\N^d)^a \to (\fkS^{[1]},(\varphi(E)),\varphi^\ast\N^d)^a$, and one can check that it represents the self-coproduct $(\fkS,(\varphi(E)),\varphi^\ast\N^d)^a$ in both $(X,M_X)_\Prism^{\mathrm{str},\mathrm{op}}$ and $(X,M_X)_\Prism^{\mathrm{op}}$ using Lemma~\ref{lem:coproduct in two sites are the same}. Similarly, we construct the $(i+1)$-st self-coproduct of $(\fkS,(\varphi(E)),\varphi^\ast\N^d)^a$ in $(X,M_X)_\Prism^\mathrm{op}$ for each $i\in \N$ by induction, and denote it by $(\fkS^{[i]},(\varphi(E)),\varphi^\ast\N^d)^a$. Moreover, the projections are classically faithfully flat by Corollary~\ref{cor:proj-map-faithful-flat}.
\end{eg}

\subsection{Breuil--Kisin log prism in the CDVR case} \label{sec:CDVR-case-BK-prism-selfprod}

We discuss analogous results of the previous section when $R$ is a complete discrete valuation ring (CDVR) of mixed characteristic $(0, p)$. Here we also describe certain filtrations on the relevant rings in Lemmas~\ref{prop:CDVR-filtration}(3) and \ref{lem:filtration-h0-CDVR}. These are key to comparing different definitions of semistable representations via Kisin descent data in the CDVR case in \S\ref{sec:quasi Kisin module and rational descent data}.

Let $\mathcal{O}_{L_0}$ be a Cohen ring whose residue field has a finite $p$-basis given by the mod $p$ reduction of $\{X_1, \ldots, X_b\}\subset \calO_{L_0}$. Write $L_0 = \mathcal{O}_{L_0}[p^{-1}]$ and equip $\mathcal{O}_{L_0}$ with the Frobenius $\varphi\colon \mathcal{O}_{L_0} \rightarrow \mathcal{O}_{L_0}$ given by $\varphi(X_i) =X_i^p$ (see \cite[Cor.~1.2.7(ii)]{BM3}). Let $L$ be a finite totally ramified extension of $L_0$ and write $\mathcal{O}_L$ for its ring of integers. Let $\pi_L \in \mathcal{O}_L$ be a uniformizer, and let $E_L[u] \in \mathcal{O}_{L_0}[u]$ be a monic irreducible polynomial for $\pi_L$. When we study purity results in \S\ref{sec:prelim-facts-rings} and \S\ref{sec:purity}, we will restrict to the case where $L_0$ contains $K_0$, $L = L_0\otimes_{K_0} K$, and $E_L(u)$ is given by $E(u)$.   

Consider $Y = \Spf (\mathcal{O}_L)$ equipped with the log structure $M_Y$ given by the prelog structure $\mathbf{N} \rightarrow \mathcal{O}_L: n \mapsto \pi_L^n$. Denote by $(Y, M_Y)_{\Prism}$ the absolute prismatic site of $(Y, M_Y)$ as in Definition~\ref{def:absolute-prismatic-site}. Let $(\mathfrak{S}_L, (E_L(u)), M_{\Spf (\mathfrak{S}_L)}) \in (Y, M_Y)_{\Prism}$ be the Breuil--Kisin log prism where $\mathfrak{S}_L = \mathcal{O}_{L_0}[\![u]\!]$ with $\varphi(u) = u^p$ and $M_{\Spf (\mathfrak{S}_L)}$ given by the prelog structure $\mathbf{N} \rightarrow \mathfrak{S}_L: n \mapsto u^n$. 

\begin{rem}
The analogues of Lemmas~\ref{lem:BKcoversfinalobject} and \ref{lem:coproduct with fkS is a covering} hold for $(\mathfrak{S}_L, (E_L), M_{\Spf (\mathfrak{S}_L)})$: for the former, consider the $p$-completion $S_{\infty}$ of $\mathfrak{S}_L[X_1^{1/p^{\infty}}, \ldots, X_b^{1/p^{\infty}}, u^{1/p^{\infty}}]\otimes_{\mathfrak{S}} R$ and argue as in the proof of Lemma~\ref{lem:BKcoversfinalobject}; for the latter,  use $B_0 \coloneqq A[S_0, S_1^{\pm 1}, \ldots, S_b^{\pm 1}][v]/(S_0-v)\rightarrow A/I$ given by $S_0 = v \mapsto \pi_L$ and $S_i \mapsto X_i$ for $1 \leq i \leq b$ and follow the proof of Lemma~\ref{lem:coproduct with fkS is a covering}.   
\end{rem}

Let $(\mathfrak{S}_L^{(i)}, (E_L), \mathbf{N})^a$ denote the $(i+1)$-st self-coproduct of $(\mathfrak{S}_L, (E_L), \mathbf{N})^a$ in $(Y, M_Y)_{\Prism}^{\mathrm{str},\mathrm{op}}$. As in \S\ref{sec:Breuil--Kisin log prism}, it also represents the $(i+1)$-st self-coproduct in $(Y, M_Y)_{\Prism}^{\mathrm{str}}$ and has the following explicit description: let $B_0 = \mathfrak{S}_L[S_1^{\pm 1}, \ldots, S_b^{\pm 1}][v]$ as in the proof above. Let $B^{(1)}_L = \mathfrak{S}_L[\![1-\frac{u_2}{u_1}, 1-\frac{X_{i, 2}}{X_{i, 1}}]\!]_{1 \leq i \leq b}$, which is also the $(1-\frac{v}{u}, 1-\frac{S_1}{X_1}, \ldots, 1-\frac{S_b}{X_b})$-completion of $B_0[\frac{v}{u}]$ via $u_1\coloneqq u$, $u_2\coloneqq v$, $X_{i, 1}\coloneqq X_i$, $X_{i, 2} \coloneqq S_i$. The kernel $J_L^{(1)}$ of the surjection $B^{(1)}_L \rightarrow \mathcal{O}_L$ is given by $(E_L(u_1), 1-\frac{u_2}{u_1}, 1-\frac{X_{i, 2}}{X_{i, 1}})_{1 \leq i \leq b}$, and we have
\[
\mathfrak{S}_L^{(1)} = B_L^{(1)}\biggl\{\frac{1-\frac{u_2}{u_1}}{E_L(u_1)}, \frac{1-\frac{X_{1, 2}}{X_{1, 1}}}{E_L(u_1)}, \ldots, \frac{1-\frac{X_{b, 2}}{X_{b, 1}}}{E_L(u_1)}\biggr\}^{\wedge}_{\delta}.
\]
Set $z_0 = \frac{1-\frac{u_2}{u_1}}{E_L(u_1)}$ and $z_i = \frac{1-\frac{X_{1, 2}}{X_{1, 1}}}{E_L(u_1)}$ for $1 \leq i \leq b$, and write $E_L = E_L(u_1)$ for simplicity.

\begin{lem} \label{lem:CDVR-BK-self-products-properties}
The projection map $p^i_j\colon \mathfrak{S}_L \rightarrow \mathfrak{S}_L^{(i)}$ is classically faithfully flat for each $i \geq 1$ and $1 \leq j \leq i$. In particular, $(\mathfrak{S}_L^{(i)},(E_L))$ is transversal for all $i$, and so $\mathfrak{S}_L^{(i)}[E_L^{-1}]$ is $p$-adically separated and
\[
\mathfrak{S}_L^{(i)}[p^{-1}] \cap \mathfrak{S}_L^{(i)}[E_L^{-1}]^{\wedge}_p = \mathfrak{S}_L^{(i)}.
\]
\end{lem}

\begin{proof}
This follows from the same argument as in the proof of Corollary~\ref{cor:proj-map-faithful-flat}.
\end{proof}

We write 
\[
A^{(1)}_{L, \mathrm{max, log}}\coloneqq  \Bigl(\mathfrak{S}_L\Bigl[\frac{E_L(u_1)}{p}, \gamma_n (z_0),\ldots,\gamma_n (z_b) \;(n\geq 0)\Bigr]\Bigr)^{\wedge}_p,
\]
considered as a subring of $\mathfrak{S}_L[{\frac{1}{p}}]_{E_L}^\wedge [\![ z_0,\ldots,z_b]\!]$. 
When the residue field of $L$ is perfect, $A^{(1)}_{L, \mathrm{max, log}}$ is studied in \cite[\S2.3]{du-liu-prismaticphiGhatmodule} (and denoted instead by $A^{(2)}_{\mathrm{st}, \mathrm{max}}$). Let $\mathcal{O}_{\mathrm{max}}$ be the $p$-completion of $\mathfrak{S}_L$-subalgebra $\mathfrak{S}_L[\frac{E_L}{p}]$ of $\mathfrak{S}_L[p^{-1}]$, and $c \coloneqq \frac{\varphi(E_L)}{p} \in \mathcal{O}_{\mathrm{max}}$. By \cite[Lem.~2.2.2]{du-liu-prismaticphiGhatmodule}, we have $c^{-1} \in \mathcal{O}_{\mathrm{max}}$. Note that any element in $A^{(1)}_{L, \mathrm{max, log}}$ can be written uniquely as 
\[
\sum_{(i_0, \ldots, i_b)\in \N^{b+1}} a_{i_0, \ldots, i_b} \gamma_{i_0}(z_0)\cdots\gamma_{i_b}(z_b),
\quad a_{i_0, \ldots, i_b} \in \mathcal{O}_{\mathrm{max}} \quad\text{with}\quad \lim_{i_0+\cdots+i_b \to \infty}a_{i_0, \ldots, i_b}=0,
\]
where the last limit is taken with respect to $p$-adic topology.

\begin{lem} \label{lem:CDVR-A(1)maxlog-Frob}
We have an injective ring map $\varphi\colon A^{(1)}_{L, \mathrm{max, log}} \rightarrow A^{(1)}_{L, \mathrm{max, log}}$ extending $\varphi$ on $\mathfrak{S}_L$ and satisfying
\[
\varphi(z_0) = \frac{1-(\frac{u_2}{u_1})^p}{\varphi(E_L)} \qquad\text{and}\qquad \varphi(z_i) = \frac{1-(\frac{X_{i, 2}}{X_{i, 1}})^p}{\varphi(E_L)} \qquad(1 \leq i \leq b).
\]
\end{lem}

\begin{proof}
Since $z_0=(1-u_2/u_1)/E_L$, we compute
\[
\varphi(z_0)=\frac{1-(\frac{u_2}{u_1})^p}{\varphi(E_L)}=\frac{p}{\varphi(E_L)}\; \frac{E_L}{p}\; \frac{1-\frac{u_2}{u_1}}{E_L}\sum_{j=0}^{p-1}\biggl(\frac{u_2}{u_1}\biggr)^j=c^{-1}\frac{E_L}{p}z_0\sum_{j=0}^{p-1} (1-E_L z_0)^j,
\]
and so
\[
\varphi(\gamma_n(z_0)) = \gamma_n(z_0)(c^{-1}\frac{E_L}{p}\sum_{j=0}^{p-1} (1-E_L z_0)^j)^n.
\]
By the same computation, the last equality also holds with $z_0$ replaced by $z_i$ ($1\leq i\leq b$). It follows that $\varphi$ on $\mathfrak{S}_L$ extends to $\varphi\colon A^{(1)}_{L, \mathrm{max, log}}\rightarrow A^{(1)}_{L, \mathrm{max, log}}$. One can also deduce the injectivity from that of $\varphi$ on $\mathfrak{S}_L$.
\end{proof}

The map $\mathfrak{S}_L\rightarrow A^{(1)}_{L, \mathrm{max, log}}$ naturally extends to a ring map $\iota\colon B^{(1)}_L= \mathfrak{S}_L[\![1-\frac{u_2}{u_1}, 1-\frac{X_{i, 2}}{X_{i, 1}}]\!]_{1 \leq i \leq b} \rightarrow A^{(1)}_{L, \mathrm{max, log}}$ compatible with $\varphi$. Consider $\mathfrak{S}_L^{(1), \mathrm{nc}} = B^{(1)}_L[\delta^n(z_0),\ldots,\delta^n(z_b)]_{n\geq 0}$ as in the proof of Proposition~\ref{prop:frakS^(1)/E}, and let $\alpha\colon \mathfrak{S}_L^{(1), \mathrm{nc}} \rightarrow \mathfrak{S}_L^{(1), \mathrm{nc}}[p^{-1}]$ denote the natural map.  

\begin{lem} \label{lem:CDVR-gamma(z)}
For each $n \geq 0$ and $0 \leq i \leq b$, there exists $f_{ni}(X) \in \mathfrak{S}_L^{(1), \mathrm{nc}}[X]$ such that
\[
\gamma_n(z_i) = f_{ni}\Bigl(\frac{E_L}{p}\Bigr)
\]
as elements of $\mathfrak{S}_L^{(1), \mathrm{nc}}[p^{-1}]$ via $\alpha$. Furthermore, $\iota\colon B^{(1)}_L \rightarrow A^{(1)}_{L, \mathrm{max, log}}$ extends naturally to a ring map $\iota\colon \mathfrak{S}_L^{(1)} \rightarrow A^{(1)}_{L, \mathrm{max, log}}$ compatible with $\varphi$.
\end{lem}

\begin{proof}
We check by induction that $\delta^n(E_L z_i) \in (E_L z_i) \mathfrak{S}_L$ for all $n \geq 1$.  
By arguing as in the proof of \cite[Lem.~2.2.6]{du-liu-prismaticphiGhatmodule} and the paragraph after Rem.~2.2.7 in \textit{loc. cit.}, we obtain the first statement and a natural map $\iota'\colon \mathfrak{S}_L^{(1), \mathrm{nc}}[p^{-1}] \rightarrow  A^{(1)}_{L, \mathrm{max, log}}[p^{-1}]$ extending $\iota\colon B^{(1)}_L \rightarrow A^{(1)}_{L, \mathrm{max, log}}$ such that $(\iota'\circ \alpha) (\mathfrak{S}_L^{(1), \mathrm{nc}}) \subset A^{(1)}_{L, \mathrm{max, log}}$. Thus, we obtain a ring map $\iota\colon \mathfrak{S}_L^{(1), \mathrm{nc}} \rightarrow A^{(1)}_{L, \mathrm{max, log}}$ given by $\iota'\circ \alpha$. Note that $E_L = p\cdot \frac{E_L}{p} \in A^{(1)}_{L, \mathrm{max, log}}$, so $\iota$ is continuous with respect to the $(p, E_L)$-adic topology on $\mathfrak{S}_L^{(1), \mathrm{nc}}$ and $p$-adic topology on $A^{(1)}_{L, \mathrm{max, log}}$. Since $\mathfrak{S}_L^{(1)}$ is naturally isomorphic to the classical $(p, E_L)$-completion of $\mathfrak{S}_L^{(1), \mathrm{nc}}$, this induces
\[
\iota\colon \mathfrak{S}_L^{(1)} \rightarrow A^{(1)}_{L, \mathrm{max, log}},
\]
which is compatible with $\varphi$ by construction.
\end{proof}

For any subring $B \subset A^{(1)}_{L,\mathrm{max, log}}[p^{-1}]$ and $i \geq 0$, set 
\[
\mathrm{Fil}^i B \coloneqq B \cap E_L(u_1)^i A^{(1)}_{L, \mathrm{max, log}}[p^{-1}].
\]

\begin{prop} \label{prop:CDVR-filtration} \hfill
\begin{enumerate}
\item We have a natural isomorphism
\[
\mathfrak{S}_L^{(1)}/(E_L) \cong \mathcal{O}_L \{z_0,\ldots,z_b\}^\wedge_\mathrm{PD}.
\]

\item The map $\iota\colon \mathfrak{S}_L^{(1)} \rightarrow A^{(1)}_{L, \mathrm{max, log}}$ is injective.

\item We have $\mathrm{Fil}^i \mathfrak{S}_L^{(1)} = E_L^i \mathfrak{S}_L^{(1)}$ (where $\mathfrak{S}_L^{(1)}$ is a subring of $A^{(1)}_{L, \mathrm{max, log}}[p^{-1}]$ via $\iota$).
\end{enumerate}    
\end{prop}

Proposition~\ref{prop:CDVR-filtration} also implies that $\alpha\colon\mathfrak{S}_L^{(1), \mathrm{nc}} \rightarrow \mathfrak{S}_L^{(1), \mathrm{nc}}[p^{-1}]$ is injective (cf.~\cite[\S2.2]{du-liu-prismaticphiGhatmodule}).

\begin{proof}
(1) The proof of Proposition~\ref{prop:frakS^(1)/E} works with $\mathfrak{S}_L^{(1), \mathrm{nc}}$ in place of $\mathfrak{S}^{(1), \mathrm{nc}}$.

(2) We have a natural isomorphism $A^{(1)}_{L, \mathrm{max, log}}[p^{-1}]/(E_L) \cong \mathcal{O}_L \{z_0,\ldots,z_b\}^\wedge_\mathrm{PD}[p^{-1}]$. So $\iota\colon \mathfrak{S}_L^{(1)} \rightarrow A^{(1)}_{L, \mathrm{max, log}}[p^{-1}]$ induces the natural map
\[
\mathcal{O}_L \{z_0,\ldots,z_b \}^\wedge_\mathrm{PD} \rightarrow \mathcal{O}_L \{z_0,\ldots,z_b \}^\wedge_\mathrm{PD}[p^{-1}],
\]
which is injective. Since $A^{(1)}_{L, \mathrm{max, log}}[p^{-1}]$ is $E_L$-torsion free and $\mathfrak{S}_L^{(1)}$ is $E_L$-adically separated, $\iota$ is injective.

(3) This follows directly from (1).
\end{proof}

Next we introduce the ring $S_L^{(1)}$, which shows up when discussing the log crystalline site of $\calO_L/(p)$. Recall $E_L z_0=1-\frac{u_2}{u_1}$ and $E_L z_i=1-\frac{X_{i, 2}}{X_{i, 1}}$ ($1\leq i\leq b$). Let 
\[
S_{L, 1} \coloneqq B_L^{(1)}\Bigl[\gamma_n(E_L), \gamma_n\Bigl(1-\frac{u_2}{u_1}\Bigr), \gamma_n\Bigl(1-\frac{X_{1, 2}}{X_{1, 1}}\Bigr),\ldots,\gamma_n\Bigl(1-\frac{X_{b, 2}}{X_{b, 1}}\Bigr)\Bigr]_{n\geq 0} \subset B_L^{(1)}[p^{-1}].
\]
Note $S_{L, 1} \subset A^{(1)}_{L, \mathrm{max, log}}$ since $\gamma_n(1-\frac{u_2}{u_1}) = \gamma_n(z_0)E_L^n$ and similarly for $\gamma_n(1-\frac{X_{i, 2}}{X_{i, 1}})$'s. Since $E_L, 1-\frac{u_2}{u_1}, 1-\frac{X_{1, 2}}{X_{1, 1}}, \ldots, 1-\frac{X_{b, 2}}{X_{b, 1}}$ form a regular sequence in $B_L^{(1)}$, $S_{L, 1}$ is the PD-envelope of $B_L^{(1)}$ for $J_L^{(1)}$ by \cite[Cor.~2.39]{bhatt-scholze-prismaticcohom}. Write $S_L^{(1)}$ for the classical $p$-completion of $S_{L, 1}$, which has the $\delta$-structure compatible with that on $B_L^{(1)}$ by \textit{loc. cit.}. As a subring of $L_0[\![u, 1-\frac{u_2}{u_1}, 1-\frac{X_{i, 2}}{X_{i, 1}}]\!]_{1 \leq i \leq b}$, we have 
\begin{align*}
S_L^{(1)} = \Bigl\{\sum_{(i,i_0, \ldots, i_{b})\in \N^{b+2}} &a_{i,i_0, \ldots, i_b} \gamma_{i}(E_L)\gamma_{i_0}\Bigl(1-\frac{u_2}{u_1}\Bigr)\gamma_{i_{1}}\Bigl(1-\frac{X_{1, 2}}{X_{1, 1}}\Bigr)\cdots\gamma_{i_{b}}\Bigl(1-\frac{X_{b, 2}}{X_{b, 1}}\Bigr) \\
&\Bigm|~ a_{i,i_0, \ldots, i_b} \in  B_L^{(1)},\lim_{i+i_0+\cdots+ i_{b} \to \infty} a_{i,i_0, \ldots, i_b} = 0 \quad(\text{$p$-adically})\Bigr\}.
\end{align*}
Furthermore, the filtration on $S_L^{(1)}$ as a subring of $A^{(1)}_{L, \mathrm{max, log}}$ is compatible with the PD-filtration on $S_L^{(1)}$, i.e., with the above description, 
\[
\mathrm{Fil}^h S_L^{(1)} = \Bigl\{\sum a_{i,i_0, \ldots, i_{b+1}} \gamma_{i}(E_L)\gamma_{i_0}\Bigl(1-\frac{u_2}{u_1}\Bigr)\cdots\gamma_{i_{b}}\Bigl(1-\frac{X_{b, 2}}{X_{b, 1}}\Bigr)\in S_L^{(1)} \Bigm| i+i_0+\cdots+i_{b} \geq h\Bigr\}.
\]

\begin{lem} \label{lem:CDVR-frakS^(1)-to-S^(1)}
We have an embedding $\mathfrak{S}_L^{(1)} \stackrel{\varphi}{\lhook\joinrel\longrightarrow} S_L^{(1)}$ induced by $\varphi$.   
\end{lem}

\begin{proof}
Note $c^{-1} \in S_L^{(1)}$ where $c \coloneqq \varphi(E_L)/p$. Since $S_L^{(1)}$ is a $\delta$-ring, we have
\[
\varphi(z_0) = c^{-1}\frac{\varphi(1-\frac{u_2}{u_1})}{p} = c^{-1}\Bigl(\frac{(1-\frac{u_2}{u_1})^p}{p}+\delta\bigl(1-\frac{u_2}{u_1}\bigr)\Bigr) \in S_L^{(1)}
\]
and $\varphi(\delta^n(z_0)) = \delta^n(\varphi(z_0)) \in S_L^{(1)}$ for $n \geq 1$. Similarly, $\varphi(\delta^n(z_i)) \in S_L^{(1)}$ for $1 \leq i \leq b$ and $n \geq 0$. Thus, the Frobenius induces a map $\mathfrak{S}_L^{(1)} \stackrel{\varphi}{\longrightarrow} S_L^{(1)}$. This map is injective since $\varphi\colon A^{(1)}_{L, \mathrm{max, log}} \rightarrow A^{(1)}_{L, \mathrm{max, log}}$ is injective by Lemma~\ref{lem:CDVR-A(1)maxlog-Frob}.
\end{proof}

When $p = 2$, we consider auxiliary subrings $\widetilde{S}_L, \widehat{S}_L \subset A^{(1)}_{L, \mathrm{max, log}}$ defined by
\begin{align*}
\widetilde{S}_L &\coloneqq \mathfrak{S}_L^{(1)}\Bigl[\!\!\Bigl[\frac{E_L^2}{2}\Bigr]\!\!\Bigr] = \Bigl\{\sum_{i \geq 0} a_i\Bigl(\frac{E_L^2}{2}\Bigr)^i \Bigm| a_i \in \mathfrak{S}_L^{(1)} \Bigr\}\quad\text{and}\\    
\widehat{S}_L &\coloneqq \mathfrak{S}_L^{(1)}\Bigl[\!\!\Bigl[\frac{E_L^4}{2}\Bigr]\!\!\Bigr] = \Bigl\{\sum_{i \geq 0} a_i\Bigl(\frac{E_L^4}{2}\Bigr)^i \Bigm| a_i \in \mathfrak{S}_L^{(1)} \Bigr\}.
\end{align*}
Note that $\widetilde{S}_L$ and $\widehat{S}_L$ are stable under $\varphi$ on $A^{(1)}_{L, \mathrm{max, log}}$.

\begin{lem} \label{lem:CDVR-hatS}
Assume $p = 2$. We have $\varphi(A^{(1)}_{L, \mathrm{max, log}}) \subset \widetilde{S}_L$ and $\varphi(\widetilde{S}_L) \subset \widehat{S}_L$. Furthermore, for each $h \geq 1$,
\[
\mathrm{Fil}^h \widetilde{S}_L = \Bigl\{ \sum_{i \geq h} a_i \frac{E_L^i}{2^{\lfloor \frac{i}{2} \rfloor}} \Bigm| a_i \in \mathfrak{S}_L^{(1)} \Bigr\}\quad\text{and}\quad
\mathrm{Fil}^h \widehat{S}_L = \Bigl\{ \sum_{i \geq h} a_i \frac{E_L^i}{2^{\lfloor \frac{i}{4} \rfloor}} \Bigm| a_i \in \mathfrak{S}_L^{(1)} \Bigr\}.
\]
\end{lem}

\begin{proof}
Given Lemma~\ref{lem:CDVR-gamma(z)}, this follows from the same argument as in the proof of \cite[Lem.~2.2.12]{du-liu-prismaticphiGhatmodule}.    
\end{proof}

The following fact on filtrations will be used in Proposition~\ref{prop:rational Kisin descent datum in CDVR case} in \S\ref{sec:quasi Kisin module and rational descent data}. 

\begin{lem} \label{lem:filtration-h0-CDVR}
Fix a positive integer $r$. 
\begin{enumerate}
\item Suppose $p \geq 3$. There exists an integer $h_0 > r$ such that for any $h \geq h_0$ and $x \in S_L^{(1)}[E_L^{-1}]$ with $E_L^r x \in \mathrm{Fil}^h S_L^{(1)}$, we have $\varphi(x) = v+w$ for some $v \in \mathfrak{S}_L^{(1)}$ and $w \in \mathrm{Fil}^{h+1} S_L^{(1)}$ as elements in $A^{(1)}_{L, \mathrm{max}, \mathrm{log}}$.

\item Suppose $p = 2$. The statement analogous to (1) with $\widehat{S}_L$ in place of $S_L^{(1)}$ holds. 
\end{enumerate}
\end{lem}

\begin{proof}
Note that $E_L^r x = \sum_{i \geq h} c_i\gamma_i(E_L)$ for some $c_i \in \mathfrak{S}_L^{(1)}$. So (1) follows from the same argument as in the proof of \cite[Lem.~4.4]{du-liu-moon-shimizu-completed-prismatic-F-crystal-loc-system}. Given Lemma~\ref{lem:CDVR-hatS}, (2) follows from the same argument as in the proof of \cite[Lem.~4.8]{du-liu-moon-shimizu-completed-prismatic-F-crystal-loc-system}.   
\end{proof}

\subsection{Ring-theoretic techniques for purity} \label{sec:prelim-facts-rings}

In this subsection, we assume $X=\Spf(R)$ is a small semistable formal scheme over $\Spf(\calO_K)$ with a framing 
\[
\square\colon R^{0} = \mathcal{O}_K \langle T_1, \ldots, T_m, T_{m+1}^{\pm 1}, \ldots, T_d^{\pm 1}\rangle / (T_1\cdots T_m - \pi) \rightarrow R,
\]
as in \S\ref{sec:Breuil--Kisin log prism}. We call the generic points of $\lvert \Spf R\rvert=\lvert \Spec (R/\pi)\rvert$ \emph{$R$-Shilov points}. We want to study the relation between $R$ and the completed local rings of $R$ at Shilov points. 

For any Shilov point $x\in X$, the completed local ring $\calO_{X,x}^\wedge$ is a ($p$-adically) complete discrete valuation ring, and its residue field is the fraction field of the corresponding irreducible component of $\Spec (R/\pi)$. Since $R$ is regular and thus normal, the adic generic fiber $X_\eta$ of $X$ is the affinoid adic space $\Spa(R[p^{-1}], R)$. The affinoid field $(\calO_{X,x}^\wedge[p^{-1}],\calO_{X,x}^\wedge)$ gives a rank one point of $X_\eta$. The points arising in this way form the Shilov boundary of the $K$-affinoid $R[p^{-1}]$ in the sense of \cite[p.~36]{Berkovich}, which justifies our terminology (see also \cite[Def.~2.4]{bhatt_hansen-6-functor-rigid}).

In fact, we want to impose a further assumption on the framing $\square$ in this subsection and when we need explicit computations: by refining $X$ by an \'etale covering and localizing $R^0$ accordingly, we may and do assume that $\square$ induces a bijection on the sets of generic points of $\Spf R^0$ and $X=\Spf R$. Under this assumption, the Shilov points are precisely $\{(T_1),\ldots, (T_m)\}\subset X$ (note also $(\pi,T_j)=(T_j)$). For each $j = 1, \ldots, m$, let $\calO_{L_j}$ denote the $p$-adic completion of $R_{(T_j)}$. Choose a subring $\calO_{L_{0,j}}\subset \calO_{L_j}$ and a ring endomorphism $\varphi$ on $\calO_{L_{0,j}}$ such that $\calO_{L_{0,j}}$ is a Cohen ring of $\calO_{L_j}/(\pi)=\operatorname{Frac}(R/(\pi, T_j))$ with $T_i\in \calO_{L_{0,j}}$ ($1\leq i \leq d$, $i\neq j$) and such that $\varphi$ is a lift of Frobenius on $\calO_{L_{0,j}}/(p)$ with $\varphi(T_i)=T_i^p$. We see $\calO_{L_j}=\calO_{L_{0,j}}[\pi]$ and the pair $(\calO_{L_{0,j}},\varphi)$ satisfies the conditions at the beginning of Subsection~\ref{sec:CDVR-case-BK-prism-selfprod}.
It also follows that $R/(\pi) \rightarrow \prod_{j=1}^m \mathcal{O}_{L_j}/(\pi)$ is injective.

We are going to establish the ring-theoretic results for purity relating the constructions for $R$ and $\calO_{L_j}$'s.
Recall that the framing $\square\colon R^0 \to R$ uniquely lifts to a map 
\[
\mathfrak{S}_{R^{0}} = W(k)\langle T_1, \ldots, T_m, T_{m+1}^{\pm 1}, \ldots, T_d^{\pm 1}\rangle[\![u]\!] / (T_1\cdots T_m-u)\rightarrow \fkS
\] 
of $\delta_\log$-rings. Write $\mathcal{O}_{\mathcal{E}} = \mathfrak{S}[E(u)^{-1}]^{\wedge}_p$. For each $j = 1, \ldots, m$, consider 
\[
\mathfrak{S}_{L_j} \coloneqq \mathcal{O}_{L_{0, j}}[\![u]\!]\qquad\text{and}\qquad \mathcal{O}_{\mathcal{E}, L_j} \coloneqq \mathfrak{S}_{L_j}[E^{-1}]^{\wedge}_p.
\]
The ring $\mathfrak{S}_{L_j}$ gives a Breuil--Kisin log prism over $\mathcal{O}_{L_j}$, and the map $\mathfrak{S}_{R^{0}} \rightarrow \mathfrak{S}_{L_j}$ given by $T_j \mapsto u(T_1\cdots T_{j-1}T_{j+1}\cdots T_m)^{-1}$, $u \mapsto u$, and $T_i \mapsto T_i$ ($i \neq j$) lifts uniquely to $\mathfrak{S} \rightarrow \mathfrak{S}_{L_j}$, which is injective, flat, and compatible with $\varphi$. Consider the intersection $\mathcal{O}_{\mathcal{E}} \cap \mathfrak{S}_{L_j}$ inside $\mathcal{O}_{\mathcal{E}, L_j}$.  Then regarding this as a subring of $\calO_\calE$, we obtain a map
\[
f\colon \mathfrak{S} \rightarrow \bigcap_{j = 1}^m (\mathcal{O}_{\mathcal{E}} \cap \mathfrak{S}_{L_j}).
\]
For each $\mathbf{Z}$-module $A$ and each $n \geq 1$, write $A_n \coloneqq A/p^n A$.

\begin{lem} \label{lem:intersection-basic-rings}
The above map $f\colon \mathfrak{S} \rightarrow \bigcap_{j=1}^m (\mathcal{O}_{\mathcal{E}} \cap \mathfrak{S}_{L_j})$ is an isomorphism. Furthermore, the induced map $\mathfrak{S}_n \rightarrow \bigcap_{j=1}^m (\mathfrak{S}_n[E^{-1}] \cap \mathfrak{S}_{L_j, n})$ is also an isomorphism. 
\end{lem}

\begin{proof}
Consider the first statement. Write $A = \bigcap_{j=1}^m (\mathcal{O}_{\mathcal{E}} \cap \mathfrak{S}_{L_j})$. The natural injective map $\mathfrak{S} \rightarrow \mathcal{O}_{\mathcal{E}}$ factors through $f\colon \mathfrak{S} \rightarrow A$, so $f$ is injective. Since $\mathfrak{S}$ and $A$ are classically $p$-complete, we need to show that the induced map $\mathfrak{S}/p\mathfrak{S} \rightarrow A/pA$ is surjective. 
Consider the composite
\[
g\colon \mathfrak{S}/p\mathfrak{S} \longrightarrow A/pA \stackrel{h}{\longrightarrow} \bigcap_{j=1}^m (\mathcal{O}_{\mathcal{E}}/p\mathcal{O}_{\mathcal{E}} \cap \mathfrak{S}_{L_j}/p\mathfrak{S}_{L_i}).
\]
Since $\mathcal{O}_{\mathcal{E}}$ is $p$-torsion free, the second map $h$ is injective. So it suffices to show that the composite $g$ is surjective. The induced map
\[
\mathfrak{S}/(p, E)\mathfrak{S} = R/pR \rightarrow \prod_{j=1}^m \mathfrak{S}_{L_j}/(p, E)\mathfrak{S}_{L_j} = \prod_{j=1}^m \mathcal{O}_{L_j}/p\mathcal{O}_{L_j},
\]
is injective since $R/(\pi) \rightarrow \prod_{j=1}^m \mathcal{O}_{L_j}/(\pi)$ is injective. Thus, the map $g$ is surjective. The second statement also follows from this.
\end{proof}

\begin{lem} \label{lem:purity-prod-map-injective}
The natural map
\[
\mathfrak{S}^{(1)}/(p,E) \rightarrow \prod_{j=1}^m \mathfrak{S}_{L_j}^{(1)}/(p,E)
\]
is injective.
\end{lem}

\begin{proof}
We use the isomorphism $\mathfrak{S}^{(1)}/(E) \cong R\{w_1,\ldots,w_d \}^\wedge_\mathrm{PD}$ from Proposition~\ref{prop:frakS^(1)/E}, 
where $\displaystyle w_i = \frac{1-T_{i, 2}/T_{i, 1}}{E}$. Furthermore, since $u = T_1\cdots T_m$ in $R$, we deduce from Proposition~\ref{prop:CDVR-filtration} that $\mathfrak{S}_{L_j}^{(1)}/(E) \cong \mathcal{O}_{L_j}\{w_1,\ldots,w_d \}^\wedge_\mathrm{PD}$ for each $j = 1, \ldots, m$. Thus, it suffices to show that the natural map
\[
R/(p)\{w_1,\ldots,w_d \}_\mathrm{PD} \rightarrow \prod_{j=1}^m \mathcal{O}_{L_j}/(p)\{w_1,\ldots,w_d \}_\mathrm{PD}
\]
is injective, which follows from the injectivity of $R/(p) \rightarrow \prod_{j=1}^m \mathcal{O}_{L_j}/(p)$.
\end{proof}

Note that $\prod_{j=1}^m  \mathfrak{S}_{L_j}^{(1)}/(p)$ is $E$-torsion free by Lemma~\ref{lem:CDVR-BK-self-products-properties}, and $\mathfrak{S}^{(1)}/(p)$ is classically $E$-complete by Lemma~\ref{lem:modE-modp-complete}. So by Lemma~\ref{lem:purity-prod-map-injective}, the map $\mathfrak{S}^{(1)}/(p) \rightarrow \prod_{j=1}^m \mathfrak{S}_{L_j}^{(1)}/(p)$ is injective. Since $\mathfrak{S}_{L_j}[E^{-1}]^{\wedge}_p \rightarrow \mathfrak{S}_{L_j}^{(1)}[E^{-1}]^{\wedge}_p$ is classically flat by \cite[Tag~0912]{stacks-project}, $\mathfrak{S}_{L_j}^{(1)}[E^{-1}]^{\wedge}_p$ is $p$-torsion free. Thus, the map $\mathfrak{S}^{(1)}[E^{-1}]^{\wedge}_p \rightarrow \prod_{j=1}^m \mathfrak{S}_{L_j}^{(1)}[E^{-1}]^{\wedge}_p$ is injective.

\begin{cor} \label{cor:purity-intersection}
The following equality holds as subrings of $\prod_{j=1}^m \mathfrak{S}_{L_j}^{(1)}[E^{-1}]^{\wedge}_p$:
\[ 
\mathfrak{S}^{(1)} = \mathfrak{S}^{(1)} [E^{-1}]^{\wedge}_p \cap \prod_{j=1}^m \mathfrak{S}_{L_j}^{(1)}.
\]
\end{cor}

\begin{proof} 
It suffices to show that the natural map $f^{(1)}\colon \mathfrak{S}^{(1)} \to \mathfrak{S}^{(1)} [E^{-1}]^{\wedge}_p \cap  \prod_{j=1}^m  \mathfrak{S}_{L_j}^{(1)}$ is surjective after modulo $p$, as both sides are classically $p$-complete. Consider 
\[ 
g^{(1)}\colon \mathfrak{S}^{(1)}/(p) \xrightarrow{f^{(1)}\bmod{p}} 
\Bigl( \mathfrak{S}^{(1)} [E^{-1}]^{\wedge}_p \cap \prod_{j =1}^m  \mathfrak{S}_{L_j}^{(1)}\Bigr)/(p)
\xrightarrow{h^{(1)}} \mathfrak{S}^{(1)}[E^{-1}]/(p) \cap \prod_{j=1}^m \mathfrak{S}_{L_j}^{(1)}/(p).
\]
Since $\prod_{j =1}^m  \mathfrak{S}_{L_j}^{(1)}[E^{-1}]^{\wedge}_p$ is $p$-torsion free, the second map $h^{(1)}$ is injective. So the surjectivity of $f^{(1)}\bmod{p}$ is reduced to the surjectivity of $g^{(1)}$, which follows from Lemma~\ref{lem:purity-prod-map-injective}: take any $x\in \mathfrak{S}^{(1)} [E^{-1}]/(p)$ and assume $x\in \mathfrak{S}_{L_j}^{(1)}/(p)$ for each $j$. Also choose $n\geq 0$ such that $E^nx\in \mathfrak{S}^{(1)}/(p)$. If $n\geq 1$, then the lemma gives $E^{n-1}x\in \mathfrak{S}^{(1)}/(p)$. Repeating this completes the proof.
\end{proof}

The following general fact will be used later.

\begin{lem}[{\cite[Lem.~3.1]{du-liu-moon-shimizu-completed-prismatic-F-crystal-loc-system}}] \label{lem:intersection-modules-flat-base-change}
Let $A$ be a ring. 
\begin{enumerate}
    \item Let $M$ be a flat $A$-module, and $N_1$, $N_2$ submodules of an $A$-module $N$. Then as submodules of $M\otimes_A N$, we have
\[
M\otimes_A (N_1 \cap N_2) = (M\otimes_A N_1) \cap (M\otimes_A N_2).
\]

    \item Let $M$ be a finite projective $A$-module, and $N$ an $A$-module. Let $\mathcal{I}$ be a (possibly infinite) index set. Suppose for each $i \in \mathcal{I}$, we are given an $A$-submodule $N_i$ of $N$. Then as submodules $M\otimes_A N$, we have
\[
M\otimes_A \bigl(\,\bigcap_{i \in \mathcal{I}} N_i\bigr) = \bigcap_{i \in \mathcal{I}} (M\otimes_A N_i).
\]
\end{enumerate}
\end{lem}

This lemma will be used in particular when we consider the following localization:

\begin{notation} \label{notation:localization}
Define $\overline{R}$ to be the union of finite $R$-subalgebras $R'$ inside a fixed algebraic closure of $\Frac{R}$ such that $R'[p^{-1}]$ is \'etale over $R[p^{-1}]$. Let $\mathcal{P}$ denote the set of minimal prime ideals of $\overline{R}$ containing $p$. For each $\mathfrak{P} \in \mathcal{P}$, write $\mathfrak{p} \coloneqq R \cap \mathfrak{P}$. Then $ \mathcal{O}_L(\mathfrak{p}) \coloneqq (R_{\mathfrak{p}})^{\wedge}$ is a CDVR whose residue field has a $p$-basis given by $\{ T_1, \ldots, T_{j-1}, T_{j+1}, \ldots T_d \}$, where the omitted $T_j$ corresponds to $\mathfrak{p}$. Note that $\pi$ is a uniformizer of $\mathcal{O}_L(\mathfrak{p})$. Write $L(\mathfrak{p}) \coloneqq \mathcal{O}_L(\mathfrak{p})[p^{-1}]$. 

For each $\mathfrak{P} \in \mathcal{P}$, fix a continuous ring map $(\overline{R}_{\mathfrak{P}})^{\wedge} \rightarrow (\mathcal{O}_{\overline{L(\mathfrak{p}})})^{\wedge}$ which extends the map $R_{\mathfrak{p}} \rightarrow \mathcal{O}_{\overline{L(\mathfrak{p}})}$. 
Taking the product over $\mathfrak{P}$'s induces injective maps
\[
\overline{R}^{\wedge} \rightarrow \prod_{\mathfrak{P} \in \mathcal{P}} (\overline{R}_{\mathfrak{P}})^{\wedge} \rightarrow \prod_{\mathfrak{P} \in \mathcal{P}} (\mathcal{O}_{\overline{L(\mathfrak{p}})})^{\wedge} \quad \text{and} \quad \overline{R}^{\flat} \rightarrow \prod_{\mathfrak{P} \in \mathcal{P}} \mathcal{O}_{\overline{L(\mathfrak{p})}}^{\flat},
\]
where $A^\flat$ denotes the tilt $\varprojlim_{x\mapsto x^p} A/p$.
\end{notation}

\begin{lem} \label{lem:intersection-A_inf-localizations}
As subrings of $\prod_{\mathfrak{P} \in \mathcal{P}} W(\overline{L(\mathfrak{p})}^{\flat})$, 
\[
W(\overline{R}^\flat[(\pi^{\flat})^{-1}]) \cap \prod_{\mathfrak{P} \in \mathcal{P}} W(\mathcal{O}_{\overline{L(\mathfrak{p}})}^{\flat}) = \A_\mathrm{inf}(\overline{R}).
\]
Furthermore, we have
\[
W(\overline{R}^\flat[(\pi^{\flat})^{-1}])[p^{-1}] \cap \prod_{\mathfrak{P} \in \mathcal{P}} W(\mathcal{O}_{\overline{L(\mathfrak{p}})}^{\flat})[p^{-1}] = \A_\mathrm{inf}(\overline{R})[p^{-1}].
\]
\end{lem}

\begin{proof}
This follows from the same argument as in the proof of \cite[Lem.~3.32]{du-liu-moon-shimizu-completed-prismatic-F-crystal-loc-system}.    
\end{proof}

\section{Analytic prismatic \texorpdfstring{$F$}{F}-crystals, \'etale and crystalline realizations} \label{sec:prismatic-F-crystal-etale-crystalline-realization}

Let $(X, M_X)$ be a bounded $p$-adic log formal scheme such that $M_X$ is integral. 
In this section, we introduce analytic prismatic $F$-crystals and Laurent $F$-crystals on $(X,M_X)_\Prism$. We then study \'etale and crystalline realizations in the semistable case and show that they are \emph{associated} over the pro-\'etale site $X_{\eta,\proet}$ of the generic fiber. 

\subsection{Analytic prismatic \texorpdfstring{$F$}{F}-crystals} \label{sec:analytic-prismatic-F-crystals}

When $X$ is smooth, \cite{du-liu-moon-shimizu-completed-prismatic-F-crystal-loc-system} and \cite{GuoReinecke-Ccris} introduced certain categories which contain the category of prismatic $F$-crystals in vector bundles on $X$, and proved that they are equivalent to the category of crystalline $\Z_p$-local systems on the generic fiber of $X$. Here we introduce the category of \emph{analytic} prismatic $F$-crystals in the logarithmic case, following \cite{GuoReinecke-Ccris}. In \S\ref{sec:analytic prismatic F-crystal in small affine case}, we will use the Breuil--Kisin log prism and its products to give a local description of analytic prismatic $F$-crystals in the semistable case as in \cite{du-liu-moon-shimizu-completed-prismatic-F-crystal-loc-system}.

\begin{defn}
    Let $\calO_\Prism$ (resp.~$\calI_\Prism$) be the structure sheaf (resp.~the ideal sheaf of the Hodge--Tate divisor) on $(X,M_X)_\Prism$ defined by associating to $(\Spf A,I,M_{\Spf A})$ the ring $A$ (resp.~ the ideal $I$). Write $\Vect((X,M_X)_\Prism,\calO_\Prism)$ for the category of vector bundles on the ringed site $((X,M_X)_\Prism,\calO_\Prism)$ (see \cite[Notation~2.1]{bhatt-scholze-prismaticFcrystal}).
\end{defn}

\begin{prop}\label{prop:strict flat descent of vector bundles}
There is a natural equivalence of categories
\[
\Vect((X,M_X)_\Prism,\calO_\Prism)\cong \lim_{(\Spf A,I,M_{\Spf A})\in (X,M_X)_\Prism}\Vect(A).
\]
A similar assertion holds if $\Vect(-)$ is replaced by $D_\mathrm{perf}(-)$ and/or if $\calO_\Prism$ is replaced by $\calO_\Prism[p^{-1}]^\wedge_{\calI_\Prism}$ or $\calO_\Prism[\calI_\Prism^{-1}]^\wedge_p$.
Similarly, the assignment 
\[
(\Spf A,I,M_{\Spf A})\mapsto \Vect(\Spec A\smallsetminus V(p,I))\quad \text{(resp.~$\Vect(\Spec A\smallsetminus V(I))$)}
\]
is a sheaf of categories for the strict flat topology. 
\end{prop}

\begin{proof}
By Remark~\ref{rem:fiber product in the strict case}, this follows from the same argument as in the proof of \cite[Prop.~2.7]{bhatt-scholze-prismaticFcrystal}:
for example, to see the first assertion, it suffices to check the assignment $(\Spf A,I,M_{\Spf A})\mapsto \Vect(A)$ is a sheaf for the strict flat topology, which follows from the $(p,I)$-completely faithfully flat descent for vector bundles over prisms.
The remaining assertions can be deduced from \cite[Thm.~5.8]{mathew-descent}.
\end{proof}

For any $(A,I,M_{\Spf A}) \in (X,M_X)_\Prism$, the $\delta$-structure on $A$ induces the Frobenius $\varphi_A$ on $A$, and $\varphi_A$ acts on $\Spec(A)\smallsetminus V(p,I)$.

\begin{defn}\label{defn:analytic F crystal}\hfill
\begin{enumerate}
    \item Write $\Vect^{\varphi}(A,I)$ for the category of pairs $(\mathfrak{M}_A,\varphi_{\mathfrak{M}_A})$ where $\mathfrak{M}_A$ is a finite locally free $A$-module, and $\varphi_{\mathfrak{M}_A}$ is an isomorphism of $A$-modules
\[
    \varphi_{\mathfrak{M}_A}\colon \varphi_A^\ast(\mathfrak{M}_A)[I^{-1}] \simeq \mathfrak{M}_A[I^{-1}].
\]

    \item  Define $\Vect^{\mathrm{an},\varphi}(A,I)$ to be the category of pairs $(\mathcal{E}_A, \varphi_{\mathcal{E}_A})$ where $\mathcal{E}_A$ is a vector bundle over $\Spec(A)\smallsetminus V(p,I)$ and $\varphi_{\mathcal{E}_A}$ is an isomorphism of vector bundles
\[
    \varphi_{\mathcal{E}_A}\colon \varphi_A^\ast(\mathcal{E}_A)[I^{-1}] \simeq \mathcal{E}_A[I^{-1}].
\]
Such a pair $(\mathcal{E}_A,\varphi_{\mathcal{E}_A})$ is called \emph{effective} if $\varphi_{\mathcal{E}_A}$ comes from a morphism $\varphi_A^\ast(\mathcal{E}_A) \to \mathcal{E}_A$ of vector bundles over $\Spec(A)\smallsetminus V(p,I)$. 
  \item
A \emph{prismatic $F$-crystal} over $(X,M_X)$ is a vector bundle $\calE$ on $((X,M_X)_\Prism,\calO_\Prism)$ equipped with an identification $\varphi_\calE\colon \varphi^\ast(\calE)[\calI_\Prism^{-1}]\cong \calE[\calI_\Prism^{-1}]$. Write $\Vect^{\varphi}((X,M_X)_\Prism)$ for the category of prismatic $F$-crystals over $(X,M_X)$. Proposition~\ref{prop:strict flat descent of vector bundles} gives an identification
\[
\Vect^{\varphi}((X,M_X)_\Prism)=\lim_{(\Spf A,I,M_{\Spf A})\in (X,M_X)_\Prism} \Vect^{\varphi}(A,I).
 \]    
    \item Define the category $\Vect^{\mathrm{an},\varphi}((X,M_X)_\Prism)$ of \emph{analytic prismatic $F$-crystals}  over $(X,M_X)$ by
\[
\Vect^{\mathrm{an},\varphi}((X,M_X)_\Prism)\coloneqq\lim_{(\Spf A,I,M_{\Spf A})\in (X,M_X)_\Prism} \Vect^{\mathrm{an},\varphi}(A,I).
 \]
 An analytic prismatic $F$-crystal is usually denoted by $(\calE_\Prism,\varphi_{\calE_\Prism})$ or $\calE_\Prism$, and write $(\calE_{\Prism,A},\varphi_{\calE_{\Prism,A}})$ for the associated object in $\Vect^{\mathrm{an},\varphi}(A,I)$ for short. We also set
 \[
 \calE_{\Prism}(A)\coloneqq \calE_{\Prism}(A,I,M_{\Spf A})\coloneqq H^0(\Spec(A)\smallsetminus V(p,I),\calE_{\Prism,A}).
 \]
 
 An analytic prismatic $F$-crystal $(\calE_\Prism,\varphi_{\calE_\Prism})$ is called \emph{effective} if $(\calE_{\Prism,A},\varphi_{\calE_{\Prism,A}})$ is effective for each $(A, I, M_{\Spf A}) \in (X,M_X)_{\Prism}$. Write $\Vect^{\mathrm{an},\varphi}_{\mathrm{eff}}((X,M_X)_\Prism)$ for the full subcategory of effective analytic prismatic $F$-crystals over $(X,M_X)$.
\end{enumerate}
\end{defn}

We have the restriction functor $\Vect^{\varphi}((X,M_X)_\Prism)\rightarrow \Vect^{\mathrm{an},\varphi}((X,M_X)_\Prism)$. 
It follows from Proposition~\ref{prop:strict flat descent of vector bundles} that the assignment $(\Spf A,I,M_{\Spf A})\mapsto \Vect^{\an,\varphi}(A,I)$ is a sheaf of categories for strict flat topology.

Recall that for any prism $(A, I)$, its Breuil--Kisin twist $A\{1\}$ is defined in \cite[\S2.3, 2.5]{Bhatt-Lurie-absoluteprismticcohomology}. This yields an object $\calO_{\Prism}\{1\} \in \Vect^{\varphi}((X,M_X)_\Prism)$ whose underlying $\calO_\Prism$-module is invertible. For each $n \in \mathbf{Z}$, write $\calO_\Prism\{n\}\coloneqq  \bigl(\calO_\Prism\{1\}\bigr)^{\otimes n}$. We have a canonical isomorphism $\varphi^\ast\calO_\Prism\{n\}\cong \calI_\Prism^{-n}\calO_\Prism\{n\}$ by \textit{loc. cit.}. For any $\calE_\Prism \in \Vect^{\mathrm{an},\varphi}((X,M_X)_\Prism)$, we define $\mathcal{E}_{\Prism}\{n\} \coloneqq \calE_\Prism\otimes_{\calO_\Prism} \calO_{\Prism}\{n\}$. 

\begin{rem}\label{rem:corresponding notions for strict site}
All the definitions and results for $(X,M_X)_\Prism$ in this subsection admit analogues for $(X,M_X)_\Prism^\str$, which we use later: we will explain that the theories for $(X,M_X)_\Prism$ and $(X,M_X)_\Prism^\str$ are equivalent in the semistable or CDVR case. See Remark~\ref{rem:strict vs saturated for prismatic crystals}.
\end{rem}

\subsection{Analytic prismatic \texorpdfstring{$F$}{F}-crystals in the small affine case}\label{sec:analytic prismatic F-crystal in small affine case}
In this subsection, we assume that 
\begin{itemize}
    \item $(X, M_X)$ is small affine as in \S\ref{sec:Breuil--Kisin log prism}, or
    \item $(X, M_X)$ is the log formal spectrum of a CDVR as in \S\ref{sec:CDVR-case-BK-prism-selfprod}.
\end{itemize}
Write $(\mathfrak{S}, (E(u)), M_{\Spf \mathfrak{S}})$ for the associated Breuil--Kisin log prism. We will show that analytic prismatic $F$-crystals over $(X, M_X)$ in these cases can be described in terms of certain modules over the Breuil--Kisin log prism with descent data. 

Consider the two projection maps $p^1_1, p^1_2\colon \mathfrak{S} \rightarrow \mathfrak{S}^{(1)}$. By the rigidity of maps of prisms \cite[Lem.~3.5]{bhatt-scholze-prismaticcohom},  $p^1_1(E(u))\mathfrak{S}^{(1)} = p^1_2(E(u))\mathfrak{S}^{(1)}$, and the projection maps 
$p^1_1, p^1_2\colon \Spec (\mathfrak{S}^{(1)}) \rightarrow \Spec (\mathfrak{S})$
restrict to
\[
p^1_1, p^1_2\colon \Spec (\mathfrak{S}^{(1)})\smallsetminus V(p,E) \rightarrow \Spec (\mathfrak{S})\smallsetminus V(p,E).
\]
A similar assertion holds for the projection maps from the triple self-product.

\begin{lem} \label{lem:descentlemmavectorbundleversion}
The evaluation on the diagram 
\[
p^1_1, p^1_2\colon (\mathfrak{S},(E), M_{\Spf(\mathfrak{S})}) \to (\mathfrak{S}^{(1)}, (E), M_{\Spf(\mathfrak{S}^{(1)})})
\] gives an equivalence of categories from $\Vect^{\varphi}((X,M_X)_\Prism)$ (resp. $\Vect^{\mathrm{an},\varphi}((X,M_X)_\Prism)$) to the following category: 
\begin{itemize}
    \item the objects are triples $(\mathcal{E}_{\mathfrak{S}},\varphi_{\mathcal{E}_{\mathfrak{S}}}, g)$, where $(\mathcal{E}_{\mathfrak{S}},\varphi_{\mathcal{E}_{\mathfrak{S}}})$ is an object in $\Vect^{\varphi}(\mathfrak{S},E(u))$ (resp. $\Vect^{\mathrm{an},\varphi}(\mathfrak{S},E(u))$), and 
\[
g\colon (p^1_1)^\ast \mathcal{E}_{\mathfrak{S}} \cong (p^1_2)^\ast\mathcal{E}_{\mathfrak{S}}
\]
is an isomorphism of vector bundles over $\Spec(\mathfrak{S}^{(1)})$ (resp. $\Spec(\mathfrak{S}^{(1)})\smallsetminus V(p,E)$) that is compatible with the Frobenii $(p^1_1)^\ast\varphi_{\mathcal{E}_{\mathfrak{S}}}$ and $(p^1_2)^\ast\varphi_{\mathcal{E}_{\mathfrak{S}}}$, and satisfies the cocycle condition over $\Spec(\mathfrak{S}^{(2)})$ (resp. $\Spec(\mathfrak{S}^{(2)})\smallsetminus V(p,E)$);
\item the morphisms are maps of vector bundles compatible with all the structures.
\end{itemize}

Furthermore, a(n analytic) prismatic $F$-crystal is effective if and only if it corresponds to $(\mathcal{E}_{\mathfrak{S}},\varphi_{\mathcal{E}_{\mathfrak{S}}},g)$ with $(\mathcal{E}_{\mathfrak{S}},\varphi_{\mathcal{E}_{\mathfrak{S}}})$ being effective. 

The same results hold if $(X,M_X)_\Prism$ is replaced by $(X,M_X)_\Prism^{\mathrm{str}}$ under Remark~\ref{rem:corresponding notions for strict site}.
\end{lem}

\begin{proof}
It is obvious that the evaluation  defines a functor to the category in the statement. On the other hand, given a triple $(\mathcal{E}_{\mathfrak{S}},\varphi_{\mathcal{E}_{\mathfrak{S}}}, g)$ in the category, Lemma~\ref{lem:BKcoversfinalobject} and \cite[Thm.~7.8]{mathew-descent} yield a(n analytic) prismatic $F$-crystal over $(X,M_X)$ which gives a quasi-inverse to the above functor. Thus, $\mathcal{E}_\Prism \mapsto (\mathcal{E}_{\mathfrak{S}},\varphi_{\mathcal{E}_{\mathfrak{S}}}, g)$ is an equivalence. The last statement follows similarly from Lemma~\ref{lem:BKcoversfinalobject} and \cite[Thm.~7.8]{mathew-descent}. 
\end{proof}

\begin{rem}\label{rem:strict vs saturated for prismatic crystals}
An immediate consequence of this lemma is that the restriction functor
\[
 \Vect^{(\mathrm{an},)\varphi}((X,M_X)_\Prism) \to \Vect^{(\mathrm{an},)\varphi}((X,M_X)_\Prism^{\mathrm{str}}).
\]
is an equivalence. One can easily deduce from it that the same equivalence also holds in the semistable case. We also note that  a similar equivalence holds without Frobenius structure. 
Hence in these cases, the two absolute prismatic sites yield the same theory.
\end{rem}

\begin{defn} \label{defn:BK-descent-data}
Let $\mathrm{DD}_\fkS\coloneqq \mathrm{DD}_{\fkS_{R,\square}}$ denote the category of \emph{(integral) Kisin descent data}: objects are triples $(\mathfrak{M}, \varphi_\fkM,  f)$ where
\begin{itemize}
\item $\fkM$ is a torsion free finite $\mathfrak{S}$-module such that $\fkM[p^{-1}]$ (resp. $\fkM[E^{-1}]$) is projective over $\fkS[p^{-1}]$ (resp. $\fkS[E^{-1}]$), and $\fkM = \fkM[p^{-1}]\cap \fkM[E^{-1}]$;
\item $\varphi_\fkM\colon (\varphi_{\mathfrak{S}}^*\mathfrak{M})[E^{-1}] \cong \mathfrak{M}[E^{-1}]$ is an isomorphism of $\fkS[E^{-1}]$-modules;
\item $f\colon \mathfrak{S}^{(1)}\otimes_{p^1_1,\mathfrak{S}}\mathfrak{M}\xrightarrow{\cong}\mathfrak{S}^{(1)}\otimes_{p^1_2,\mathfrak{S}}\mathfrak{M}$ is an isomorphism of $\fkS^{(1)}$-modules that satisfies the cocycle condition over $\mathfrak{S}^{(2)}$ and is compatible with Frobenii after inverting $E$.
\end{itemize}
Morphisms of $\mathrm{DD}_\fkS$ are $\mathfrak{S}$-linear maps which are compatible with all the structures.
Write $\mathrm{DD}_\fkS^{\mathrm{vb}}$ for the full subcategory of  $\mathrm{DD}_\fkS$ consisting of triples $(\mathfrak{M}, \varphi_\fkM, f)$ where $\fkM$ is a finite projective $\fkS$-module.
\end{defn}

\begin{lem} \label{lem:descentlemmamoduleversion}
The restriction along $U\coloneqq\Spec(\fkS)\smallsetminus V(p,E)\hookrightarrow \Spec\fkS$ gives rise to an equivalence of categories
\[
\mathrm{DD}_\fkS
\xrightarrow{\cong} \Vect^{\mathrm{an},\varphi}((X,M_X)_\Prism),
\]
and the quasi-inverse $\calE_\Prism\mapsto(\mathfrak{M}, \varphi_\fkM, f)$ satisfies $\mathfrak{M} \cong H^0(U, \calE_{\Prism,\mathfrak{S}})$.
Similarly, we have $\mathrm{DD}_\fkS^{\mathrm{vb}}\xrightarrow{\cong}\Vect^{\varphi}((X,M_X)_\Prism)$. Moreover, the equivalences preserve effectivity.
\end{lem}

\begin{proof}
We will only show the case for analytic prismatic $F$-crystals, as the prismatic $F$-crystal case is easier. By Lemma~\ref{lem:descentlemmavectorbundleversion}, it suffices to have a natural equivalence between $\mathrm{DD}_{\mathfrak{S}}$ and the category of objects $(\mathcal{E}_{\mathfrak{S}}, \varphi_{\mathcal{E}_{\mathfrak{S}}}, g)$ as in Lemma~\ref{lem:descentlemmavectorbundleversion}. For $(\mathfrak{M}, \varphi_\fkM, f)\in \mathrm{DD}_\fkS$, its restriction to $U$ yields such a triple, which defines a functor. 

Now let $(\mathcal{E}_{\mathfrak{S}}, \varphi_{\mathcal{E}_{\mathfrak{S}}}, g)$ be a triple as in Lemma~\ref{lem:descentlemmavectorbundleversion}, and consider a $\fkS$-module $\mathfrak{M} = H^0(U, \mathcal{E}_{\mathfrak{S}})$. Let $j\colon U \hookrightarrow \Spec\mathfrak{S}$ be the open immersion. Since $(\mathfrak{S}, (E))$ is a transversal prism and $\mathfrak{S}$ is noetherian, $j_\ast \mathcal{E}_{\mathfrak{S}}$ is a coherent sheaf over $\mathfrak{S}$ and $\mathfrak{M} = H^0(\mathfrak{S}, j_\ast \mathcal{E}_{\mathfrak{S}})$ is a finitely generated $\mathfrak{S}$-module by \cite[Lem.~5.8]{GuoReinecke-Ccris}. Furthermore, $\mathcal{E}_{\mathfrak{S}}$ is the pullback to $U$ of the coherent sheaf on $\Spec \fkS$ attached to $\mathfrak{M}$ by \cite[Tag~0BK0]{stacks-project}. In particular, $\fkM[p^{-1}]$ (resp. $\fkM[E^{-1}]$) is projective over $\fkS[p^{-1}]$ (resp. $\fkS[E^{-1}]$), and $\fkM = \fkM[p^{-1}]\cap \fkM[E^{-1}]$. 

Since $\varphi\colon \mathfrak{S} \rightarrow \mathfrak{S}$ is flat by Lemma~\ref{lem:varphionfkSisfaithfullyflat}, $\varphi_{\mathcal{E}_{\mathfrak{S}}}$ induces an isomorphism $\varphi_\fkM\colon (\varphi_{\mathfrak{S}}^\ast\mathfrak{M})[E^{-1}] \cong \mathfrak{M}[E^{-1}]$ of $\fkS[E^{-1}]$-modules by Lemma~\ref{lem:intersection-modules-flat-base-change}. Similarly, since $p^1_1, p^1_2\colon \mathfrak{S} \rightarrow \mathfrak{S}^{(1)}$ is flat by Corollary~\ref{cor:proj-map-faithful-flat} and Lemma~\ref{lem:CDVR-BK-self-products-properties}, $f$ induces an isomorphism $f\colon \mathfrak{S}^{(1)}\otimes_{p^1_1,\mathfrak{S}}\mathfrak{M}\xrightarrow{\cong}\mathfrak{S}^{(1)}\otimes_{p^1_2,\mathfrak{S}}\mathfrak{M}$ of $\fkS^{(1)}$-modules that satisfies the cocycle condition over $\mathfrak{S}^{(2)}$ and is compatible with Frobenii after inverting $E$ (alternatively, use the flat base change \cite[Tag~02KH]{stacks-project}). Thus, we obtain a functor to $\mathrm{DD}_{\mathfrak{S}}$ given by $(\mathcal{E}_{\mathfrak{S}}, \varphi_{\mathcal{E}_{\mathfrak{S}}}, g)\mapsto(\mathfrak{M}, \varphi_\fkM, f)$.
By construction, this is quasi-inverse to the above functor, which gives the desired equivalence.

The last statement on effectivity follows from the corresponding one in Lemma~\ref{lem:descentlemmavectorbundleversion}.
\end{proof}

\begin{cor} \label{cor:CDVR-equivalence-descent-datum}
The restriction $\Vect^{\varphi}((X,M_X)_\Prism)\rightarrow\Vect^{\mathrm{an},\varphi}((X,M_X)_\Prism)$ is fully faithful. When $X = \Spf\calO_L$ for a CDVR $\calO_L$, it is an equivalence of categories. 
\end{cor}

\begin{proof}
The first part follows from Lemma~\ref{lem:descentlemmamoduleversion}. For the second part, note that $\mathfrak{S}$ is a regular local ring of dimension $2$ when $X = \Spf\calO_L$. So the global sections of any vector bundle over $\Spec{\mathfrak{S}}\smallsetminus V(p,E)$ is a finite free $\mathfrak{S}$-module (see also \cite[Rem.~3.18]{du-liu-moon-shimizu-completed-prismatic-F-crystal-loc-system}). 
\end{proof}

\begin{rem} \hfill
\begin{enumerate}
\item For any $\calE_\Prism \in \Vect^{\mathrm{an},\varphi}((X,M_X)_\Prism)$, $\calE_\Prism\{n\}$ is effective for $n \ll 0$.
\item We will show in Proposition~\ref{prop:kisin-mod-invert-p-projectivity} that for any finite $\mathfrak{S}$-module equipped with an isomorphism $\varphi_\fkM\colon (\varphi^*\mathfrak{M})[E^{-1}] \cong \mathfrak{M}[E^{-1}]$ of $\fkS[E^{-1}]$-modules, $\fkM[p^{-1}]$ is finite projective over $\fkS[p^{-1}]$.
\item The first statement of Corollary~\ref{cor:CDVR-equivalence-descent-datum} holds when $X$ is semistable over $\calO_K$: by taking an \'etale cover, we can reduce it to the small affine case. 
\end{enumerate}
\end{rem}

Finally, we discuss the localization set-up as in \S\ref{sec:prelim-facts-rings}: assume that $(X,M_X)$ is small affine for $X = \Spf R$ equipped with a framing 
\[
\square\colon R^{0} = \mathcal{O}_K \langle T_1, \ldots, T_m, T_{m+1}^{\pm 1}, \ldots, T_d^{\pm 1}\rangle / (T_1\cdots T_m - \pi) \rightarrow R \quad (m \geq 1)
\]
such that $\Spf (R/\pi R)$ has precisely $m$ irreducible components given by $R/(\pi, T_j)$ ($j = 1, \ldots, m$). Let $\Spf(\mathcal{O}_{L_j}) \rightarrow X$ be the Shilov point corresponding to $R/(\pi, T_j)$. We have maps of log prisms  $p_1^1,p_2^1\colon (\mathfrak{S}, (E), M_{\Spf \mathfrak{S}})\rightarrow (\mathfrak{S}^{(1)}, (E), M_{\Spf \mathfrak{S}^{(1)}})$ and $(\mathfrak{S}, (E), M_{\Spf \mathfrak{S}}) \rightarrow (\mathfrak{S}_{L_j}, (E), M_{\Spf \mathfrak{S}_{L_j}})$ in $(X, M_X)_{\Prism}^\mathrm{op}$.  

\begin{lem} \label{lem:properties-analy-prism-F-crystal}
For $\mathcal{E}_\Prism \in \Vect^{\mathrm{an},\varphi}((X,M_X)_\Prism)$, let $(\mathfrak{M}, \varphi_\fkM, f) \in \mathrm{DD}_{\mathfrak{S}}$ be the corresponding Kisin descent datum given in Lemma~\ref{lem:descentlemmamoduleversion} and set $\mathcal{M} \coloneqq \mathfrak{M}[E^{-1}]^{\wedge}_p$. The following properties hold:
\begin{enumerate}
    \item The natural map
    \[
    \mathfrak{S}^{(1)}\otimes_{p^1_i, \mathfrak{S}} \mathfrak{M} \rightarrow \mathfrak{S}^{(1)}[E^{-1}]^{\wedge}_p\otimes_{p^1_i, \mathfrak{S}} \mathfrak{M}
    \]
    is injective for $i = 1, 2$.
    
    \item For each $1 \leq j \leq m$, the space of global sections of $\mathcal{E}_\Prism(\mathfrak{S}_{L_j}, (E), M_{\Spf \mathfrak{S}_{L_j}})$ is isomorphic to $\mathfrak{M}_{L_j}\coloneqq\mathfrak{M}\otimes_{\mathfrak{S}} \mathfrak{S}_{L_j}$ as $\mathfrak{S}_{L_j}$-modules, and $\mathfrak{M}_{L_j}$ is free over $\mathfrak{S}_{L_j}$.  

    \item Consider $\mathcal{M} \cap \mathfrak{M}_{L_j}$ as an $\mathfrak{S}$-submodule of $\mathcal{M}\otimes_{\mathfrak{S}[E^{-1}]^{\wedge}_p} \mathfrak{S}_{L_j}[E^{-1}]^{\wedge}_p$. We have
    \[
    \mathfrak{M} = \bigcap_{j = 1}^m (\mathcal{M} \cap \mathfrak{M}_{L_j})
    \]
    as submodules of $\mathcal{M}$.
\end{enumerate}
\end{lem}

\begin{proof}
(1) The proof of \cite[Lem.~3.24]{du-liu-moon-shimizu-completed-prismatic-F-crystal-loc-system} also works here: $p^1_i$ is flat and $\mathfrak{S}^{(1)} \rightarrow \mathfrak{S}^{(1)}[E^{-1}]^{\wedge}_p$ is injective by Corollary~\ref{cor:proj-map-faithful-flat}. Moreover, $\mathfrak{M}[p^{-1}]$ is flat over $\mathfrak{S}$. Hence each of the two maps
\[
\mathfrak{S}^{(1)}\otimes_{p^1_i, \mathfrak{S}} \mathfrak{M}\rightarrow\mathfrak{S}^{(1)}\otimes_{p^1_i, \mathfrak{S}} \mathfrak{M}[p^{-1}] \rightarrow \mathfrak{S}^{(1)}[E^{-1}]^{\wedge}_p\otimes_{p^1_i, \mathfrak{S}} \mathfrak{M}[p^{-1}]
\]
is injective. Since the composite factors through $\mathfrak{S}^{(1)}\otimes_{p^1_i, \mathfrak{S}} \mathfrak{M} \rightarrow \mathfrak{S}^{(1)}[E^{-1}]^{\wedge}_p\otimes_{p^1_i, \mathfrak{S}} \mathfrak{M}$, the latter is also injective.

(2) Since $\mathfrak{S} \rightarrow \mathfrak{S}_{L_j}$ is flat, the space of global sections of $\mathcal{E}_\Prism(\mathfrak{S}_{L_j}, (E), M_{\Spf \mathfrak{S}_{L_j}})$ is isomorphic to $\mathfrak{M}_{L_j}=\mathfrak{M}\otimes_{\mathfrak{S}} \mathfrak{S}_{L_j}$ by \cite[Tag~02KH]{stacks-project} and the proof of Lemma~\ref{lem:descentlemmamoduleversion}. Moreover, we have
\[
\mathfrak{M}_{L_j}[p^{-1}] \cap \mathfrak{M}_{L_j}[E^{-1}] =(\mathfrak{M}[p^{-1}]\cap\mathfrak{M}[E^{-1}])\otimes_{\mathfrak{S}} \mathfrak{S}_{L_j} =\mathfrak{M}_{L_j}
\]
by Lemma~\ref{lem:intersection-modules-flat-base-change}(1). Since $\mathfrak{S}_{L_j}$ is a regular local ring of dimension $2$, $\mathfrak{M}_{L_j}$ is free over $\mathfrak{S}_{L_j}$.

(3) It suffices to show $\bigcap_{j = 1}^m (\mathcal{M} \cap \mathfrak{M}_{L_j}) \subset \mathfrak{M}$. Since $\mathfrak{M}[p^{-1}]$ is flat over $\mathfrak{S}$, we have
\[
\bigcap_{j = 1}^m (\mathcal{M} \cap \mathfrak{M}_{L_j})[p^{-1}] = \mathfrak{M}[p^{-1}]
\]
by Lemmas~\ref{lem:intersection-basic-rings} and \ref{lem:intersection-modules-flat-base-change}. Similarly, $\bigcap_{j = 1}^m (\mathcal{M} \cap \mathfrak{M}_{L_j})[E^{-1}] = \mathfrak{M}[E^{-1}]$. So
\[
\bigcap_{j = 1}^m (\mathcal{M} \cap \mathfrak{M}_{L_j}) \subset \mathfrak{M}[p^{-1}] \cap \mathfrak{M}[E^{-1}] = \mathfrak{M}.
\]
\end{proof}

\subsection{Laurent \texorpdfstring{$F$}{F}-crystals and \'etale realizations} \label{subsec:Laurent-F-crystals}

Let $(X,M_X)$ be a bounded fs $p$-adic log formal scheme. We recall the definition of Laurent $F$-crystals over the absolute logarithmic prismatic site and show that they are related to $\Z_p$-local systems on the generic fiber $X_\eta$ of $X$ as in \cite[\S3]{bhatt-scholze-prismaticFcrystal} and \cite[\S7.5]{koshikawa-yao}.

\begin{defn} \label{def:Laurent-F-crystals}
A \emph{Laurent $F$-crystal} on $(X,M_X)_\Prism$  consists of a pair $(\mathscr{E},\varphi_{\mathscr{E}})$ where $\mathscr{E}$ is a crystal in vector bundles on the ringed site $((X,M_X)_\Prism, \calO_\Prism[\calI_\Prism^{-1}]^\wedge_p)$ and $\varphi_{\mathscr{E}}$ is an isomorphism $\varphi_{\mathscr{E}}\colon \varphi^\ast {\mathscr{E}}\xrightarrow{\cong}{\mathscr{E}}$.
Write $\Vect((X,M_X)_\Prism,\calO_\Prism[\calI_\Prism^{-1}]^\wedge_p)^{\varphi=1}$ for the category of Laurent $F$-crystals on $(X,M_X)_\Prism$.
Similarly, we define the category $\Vect((X,M_X)_\Prism^{\mathrm{str}(, \perf)},\calO_\Prism[\calI_\Prism^{-1}]^\wedge_p)^{\varphi=1}$ of Laurent $F$-crystals on the absolute (perfect) strict logarithmic prismatic site $(X,M_X)_\Prism^{\mathrm{str}(, \perf)}$.
\end{defn}

Let $\Vect(A[I^{-1}]^\wedge_p)^{\varphi=1}$ denote the category of \'etale $\varphi$-modules over $A[I^{-1}]^\wedge_p$, i.e., pairs $(\mathcal{M},\varphi_{\mathcal{M}})$ where $\mathcal{M}$ is a finite locally free $A[I^{-1}]^\wedge_p$-module, and $\varphi_{\mathcal{M}}\colon \varphi^\ast \mathcal{M} \xrightarrow{\cong} \mathcal{M}$ is an isomorphism of $A[I^{-1}]^\wedge_p$-modules.

\begin{lem} \label{lem:Laurent log F-crystals as limit}
There is a natural equivalence
\[
\Vect((X,M_X)_\Prism,\calO_\Prism[\calI_\Prism^{-1}]^\wedge_p)^{\varphi=1}\cong
\lim_{(\Spf A,I,M_{\Spf A})\in (X,M_X)_\Prism}\Vect(A[I^{-1}]^\wedge_p)^{\varphi=1}.
\]
A similar equivalence holds with $(X,M_X)_\Prism^\mathrm{str}$ or $(X,M_X)_\Prism^{\mathrm{str},\perf}$ in place of $(X,M_X)_\Prism$.
\end{lem}

\begin{proof}
This follows from the same argument as in the proof of \cite[Prop.~2.7]{bhatt-scholze-prismaticFcrystal} based on \cite[Thm.~5.8]{mathew-descent}.
\end{proof}

The category of Laurent $F$-crystals on the absolute strict prismatic site $(X,M_X)_\Prism^{\mathrm{str}}$ corresponds to the category $\Loc_{\Z_p}(X_\eta)$ of (quasi-pro-)\'etale $\Z_p$-local systems on the generic fiber $X_\eta$ forgetting the log structure in the following set-up:

\begin{thm} \label{thm:logetalerealizationofLaurentFcrystals-strict}
Assume that $(X,M_X)$ is semistable or the log formal spectrum of a CDVR. There is a natural equivalence of categories
\[
\Vect((X,M_X)_\Prism^\mathrm{str},\calO_\Prism[\calI_\Prism^{-1}]^\wedge_p)^{\varphi=1}\cong \Loc_{\Z_p}(X_\eta).
\]
\end{thm}

\begin{proof}
We will obtain the desired equivalence as follows:
\begin{align}
\Vect((X,M_X)_\Prism^\mathrm{str},\calO_\Prism[\calI_\Prism^{-1}]^\wedge_p)^{\varphi=1}
&\xrightarrow{\cong} \Vect((X,M_X)_\Prism^{\mathrm{str},\perf},\calO_\Prism[\calI_\Prism^{-1}]^\wedge_p)^{\varphi=1}
\label{eq:first eq of etale realization in strict case}\\
&\xleftarrow{\cong} \Vect(X_\Prism^\perf,\calO_\Prism[\calI_\Prism^{-1}]^\wedge_p)^{\varphi=1}
\label{eq:second eq of etale realization in strict case}
\\
&\cong \Loc_{\Z_p}(X_\eta).
\label{eq:third eq of etale realization in strict case}
\end{align}
For \eqref{eq:first eq of etale realization in strict case}, consider the equivalences in Lemma~\ref{lem:Laurent log F-crystals as limit}:
\[
\Vect((X,M_X)_\Prism^{\mathrm{str}},\calO_\Prism[\calI_\Prism^{-1}]^\wedge_p)^{\varphi=1}\cong
\lim_{(\Spf A,I,M_{\Spf A})\in (X,M_X)_\Prism^\mathrm{str}}\Vect(A[I^{-1}]^\wedge_p)^{\varphi=1}
\]
and
\[
\Vect((X,M_X)_\Prism^{\mathrm{str},\perf},\calO_\Prism[\calI_\Prism^{-1}]^\wedge_p)^{\varphi=1}\cong
\lim_{(\Spf A,I,M_{\Spf A})\in (X,M_X)_\Prism^{\mathrm{str},\perf}}\Vect(A[I^{-1}]^\wedge_p)^{\varphi=1}.
\]
So it suffices to compare the target categories. For $(\Spf A,I,M_{\Spf A})\in (X,M_X)_\Prism^\mathrm{str}$, let $(A',I')$ denote the perfection of the prism $(A,I)$. As in \cite[Rem.~2.6]{min-wang-HT-crys-log-prism}, the pair $(A',I')$ gives rise to a log prism $(\Spf A',I',M_{\Spf A'})$ with perfect underlying prism that is initial in the category of such log prisms over $(\Spf A,I,M_{\Spf A})$. By the proof of \cite[Cor.~3.7]{bhatt-scholze-prismaticFcrystal},  the pullback
\[
\Vect(A[I^{-1}]^\wedge_p)^{\varphi=1}\rightarrow \Vect(A'[I'^{-1}]^\wedge_p)^{\varphi=1}
\]
is an equivalence of categories, which gives equivalence \eqref{eq:first eq of etale realization in strict case}.

Equivalence \eqref{eq:second eq of etale realization in strict case} follows from Remark~\ref{rem:log prismatic site when admitting finite free chart}(2) (here we use the assumption on $(X,M_X)$) and Proposition~\ref{prop:strict flat descent of vector bundles}, and \eqref{eq:third eq of etale realization in strict case} is obtained in the proof of \cite[Cor.~3.7, 3.8]{bhatt-scholze-prismaticFcrystal}. 
\end{proof}

\begin{rem}\label{rem:KY's equivalence for Laurent F-crystals}
For a general bounded fs log $p$-adic formal scheme $(X,M_X)$, 
Koshikawa and Yao study the category $\Loc_{\Z_p}((X,M_X)_\eta)$ of quasi-pro-Kummer-\'etale $\Z_p$-local systems on the log diamond generic fiber and defines an \'etale realization functor 
\[
\Vect((X,M_X)_\Prism,\calO_\Prism[\calI_\Prism^{-1}]^\wedge_p)^{\varphi=1}\rightarrow \Loc_{\Z_p}((X,M_X)_\eta).
\]
See \cite[Thm.~7.35]{koshikawa-yao} for the construction.
In the set-up of Theorem~\ref{thm:logetalerealizationofLaurentFcrystals-strict}, one can show that the following diagram is $2$-commutative:
\[
\xymatrix{
\Vect((X,M_X)_\Prism^\mathrm{str},\calO_\Prism[\calI_\Prism^{-1}]^\wedge_p)^{\varphi=1}
\ar[d]_-\alpha\ar[r]^-{\cong}_-{\text{\ref{thm:logetalerealizationofLaurentFcrystals-strict}}}
&
\Loc_{\Z_p}(X_\eta) \ar[d]^-\beta
\\
\Vect((X,M_X)_\Prism,\calO_\Prism[\calI_\Prism^{-1}]^\wedge_p)^{\varphi=1}
\ar[r]
& 
\Loc_{\Z_p}((X,M_X)_\eta),
}
\]
where $\alpha$ is induced by the pullback of Laurent $F$-crystals on the strict sites along the structure map $\Spf (A/I)\rightarrow X$ in 
Proposition~\ref{prop: prismatic pullback of crystals}(2) for each $(\Spf A,I, M_{\Spf A})\in (X,M_X)_{\Prism}$, and $\beta$ is defined via a morphism of sites $(X,M_X)_{\eta,\mathrm{qpk\acute{e}t}}\rightarrow X_{\eta,\mathrm{qpro\acute{e}t}}$ given by  \cite[Prop.~7.16(1)]{koshikawa-yao}.
Note that $\alpha$ is an equivalence with quasi-inverse given by the restriction functor (cf.~ Remark~\ref{rem:strict vs saturated for Laurent F crystals}).
\end{rem}

\begin{defn} \label{defn:EprismtoEet}
For any log prism $(\Spf A,I,M_{\Spf A})\in (X,M_X)_\Prism$, restricting vector bundles over $\Spec(A)\smallsetminus V(p,I)$ to $\Spec(A)\smallsetminus V(I)$ and then extending the scalar along $A[I^{-1}]\to A[I^{-1}]^\wedge_p$ defines a functor $\Vect^{\an,\varphi}(A,I) \to \Vect(A[I^{-1}]^\wedge_p)^{\varphi=1}$. By Lemma~\ref{lem:Laurent log F-crystals as limit}, we obtain the scalar extension functor
\[
\Vect^{\an,\varphi}((X,M_X)_\Prism^{(\str)}) \to \Vect((X,M_X)_\Prism^{(\str)}, \calO_\Prism[\calI_\Prism^{-1}]^\wedge_p)^{\varphi=1}.
\]
We say that a Laurent $F$-crystal \emph{can be extended} to an analytic prismatic $F$-crystal if it lies in the essential image of the above functor. For any analytic prismatic $F$-crystal $\calE_\Prism$ over $(X,M_X)_{\Prism}^{(\str)}$, write $\calE_{\mathrm{\acute{e}t}}$ for the associated Laurent $F$-crystal. 
\end{defn}

\begin{defn}[\'Etale realization of analytic prismatic $F$-crystals]\label{defn:etale realization}

Assume $(X, M_X)$ is semistable or the log formal spectrum of a CDVR. The \emph{\'etale realization} functor is defined to be the composite
\[
T=T_X\colon \Vect^{\an,\varphi}((X,M_X)_\Prism^{\str})\rightarrow \Vect((X,M_X)_\Prism^{\str},\calO_\Prism[\calI_\Prism^{-1}]^\wedge_p)^{\varphi=1}\underset{\text{Thm.~\ref{thm:logetalerealizationofLaurentFcrystals-strict}}}{\xrightarrow{\cong}} \Loc_{\Z_p}(X_\eta).
\]
By abuse of notation, we will also use $T$ to denote the equivalence functor defined in Theorem~\ref{thm:logetalerealizationofLaurentFcrystals-strict} and call it the \emph{\'etale realization functor for Laurent $F$-crystals.}

Recall from Remark~\ref{rem:strict vs saturated for prismatic crystals} that the restriction functor $\Vect^{(\mathrm{an},)\varphi}((X,M_X)_\Prism) \to \Vect^{(\mathrm{an},)\varphi}((X,M_X)_\Prism^{\mathrm{str}})$ is an equivalence. We also define the \'etale realization 
\[
T\colon \Vect^{(\mathrm{an},)\varphi}((X,M_X)_\Prism) \rightarrow \Loc_{\Z_p}(X_\eta)
\]
to be the composite (cf. Remarks~\ref{rem:KY's equivalence for Laurent F-crystals} and \ref{rem:strict vs saturated for Laurent F crystals}).
\end{defn}

\begin{prop}\label{prop:pullback commutes with etale realization}
Let $f\colon (X,M_X) \rightarrow (Y,M_Y)$ be a morphism of $p$-adic fs log formal schemes that are semistable or the log formal spectrum of a CDVR, and let $f_\eta\colon X_{\eta} \rightarrow Y_{\eta}$ be the induced map on the generic fibers. Then the following is a $2$-commutative diagram of categories:
\[
\xymatrix{
\Vect^{\mathrm{an},\varphi}((Y,M_Y)_\Prism^{\str})\ar[r]\ar[d]^-{f_\Prism^{-1}}
& \Vect((Y,M_Y)_\Prism^{\str},\calO_\Prism[\calI_\Prism^{-1}]^\wedge_p)^{\varphi=1}
\ar[d]^-{f_\Prism^{-1}}\ar[r]^-{\cong}
&
\Loc_{\Z_p}(Y_\eta) \ar[d]^-{f_\eta^{-1}}
\\
\Vect^{\mathrm{an},\varphi}((X,M_X)_\Prism^{\str})\ar[r]
&
\Vect((X,M_X)_\Prism^{\str},\calO_\Prism[\calI_\Prism^{-1}]^\wedge_p)^{\varphi=1}
\ar[r]^-{\cong}
& 
\Loc_{\Z_p}(X_\eta).
}
\]
\end{prop}

See Propositions~\ref{prop:functoriality of log prismatic site}, \ref{prop: prismatic pullback of crystals} and Definition~\ref{defn:pullback-analytic-prismatic-F-crystals} for the pullback $f_\Prism^{-1}$ for the strict prismatic sites.

\begin{proof}
The $2$-commutativity of the first square directly follows from the definition. 
 By the proof of Theorem~\ref{thm:logetalerealizationofLaurentFcrystals-strict}, the $2$-commutativity of the second square follows from Lemma~\ref{lem:pullback functoriality for etale realization in the non-log case} and Remark~\ref{rem:functoriality of perfect log prismatic site} below.
\end{proof}

\begin{lem}\label{lem:pullback functoriality for etale realization in the non-log case}
Let $f\colon X\rightarrow Y$ be a morphism of bounded $p$-adic formal schemes whose adic generic fiber $f_\eta\colon \colon X_\eta\rightarrow Y_\eta$ exists as a morphism of locally noetherian adic space over $\Spa(K,\calO_K)$.
Then the following is a $2$-commutative diagram of categories:
\[
\xymatrix{
\Vect(Y_\Prism^{\perf},\calO_\Prism[\calI_\Prism^{-1}]^\wedge_p)^{\varphi=1}
\ar[r]^-{\cong}\ar[d]^-{(f_\Prism^\perf)^{-1}}
&
\Loc_{\Z_p}(Y_\eta)\ar[d]^-{f_\eta^{-1}} \\
\Vect(X_\Prism^{\perf},\calO_\Prism[\calI_\Prism^{-1}]^\wedge_p)^{\varphi=1}
\ar[r]^-{\cong}
& 
\Loc_{\Z_p}(X_\eta).
}
\]
\end{lem}

The assumption on $f$ can surely be relaxed, but we content ourselves with the above one to avoid finding precise references for a more general setting in the literature.

\begin{proof}
We use the v-descent for $\Z_p$-local systems: 
$\Loc_{\Z_p}(X_\eta) \xrightarrow{\cong}\lim_{\Spa(S,S^+)\rightarrow X_\eta}\Loc_{\Z_p}(\Spa(S,S^+))$, where $(S,S^+)$ runs over all affinoid perfectoid rings over $(K,\calO_K)$ (cf. \cite[Prop.~3.5, 3.7, 3.9]{MannWerner-Localsystemsondiamonds}).

For any such $(S,S^+)$, a map $\Spa(S,S^+)\rightarrow X_\eta$ yields a map $\Spf S^+\rightarrow X$ of bounded $p$-adic formal schemes, which in turn defines a perfect prism $(W(S^{+,\flat}),\Ker(W(S^{+,\flat})\rightarrow S^+))$ of $X_\Prism^\perf$.
We then have a $2$-commutative diagram
\[
\xymatrix{
\Vect(Y^\perf_\Prism,\calO_\Prism[\calI_\Prism^{-1}]^\wedge_p)^{\varphi=1}
\ar[d]^-{(f^\perf_\Prism)^{-1}}\ar[r]&\Vect(W(S^{\flat}))^{\varphi = 1}  \ar@{=}[d]\ar[r]^-\cong& \Loc_{\Z_p}(\Spa(S,S^+))\ar@{=}[d]\\
\Vect(X^\perf_\Prism,\calO_\Prism[\calI_\Prism^{-1}]^\wedge_p)^{\varphi=1}
\ar[r]&\Vect(W(S^{\flat}))^{\varphi = 1}  \ar[r]^-\cong& \Loc_{\Z_p}(\Spa(S,S^+)).
}
\]
Here the second horizontal equivalence follows from \cite[Prop.~3.4]{bhatt-scholze-prismaticFcrystal}, or
\cite[Thm.~8.5.5(a)$\iff$(d)]{kedlaya-liu-relative-padichodge} and 
\cite[Prop.~3.7]{MannWerner-Localsystemsondiamonds}.
By varying $\Spa(S,S^+)$ over $X_\eta$, we obtain the desired $2$-commutativity since $f_\eta^{-1}\colon \Loc_{\Z_p}(Y_\eta)\rightarrow \Loc_{\Z_p}(X_\eta)$ is identified as the composite 
\[
\lim_{\Spa(S,S^+)\rightarrow Y_\eta}\mspace{-20mu}\Loc_{\Z_p}(\Spa(S,S^+))\rightarrow \lim_{\Spa(S,S^+)\rightarrow  X_\eta\rightarrow Y_\eta}\mspace{-30mu}\Loc_{\Z_p}(\Spa(S,S^+)) \rightarrow \lim_{\Spa(S,S^+)\rightarrow X_\eta}\mspace{-20mu}\Loc_{\Z_p}(\Spa(S,S^+)).
\]
\end{proof}

\bigskip

\noindent
\textbf{More on Laurent $F$-crystals.}
In the remaining, we assume that $(X, M_X)$ is semistable or the log formal spectrum of a CDVR, and explain additional properties of Laurent $F$-crystals.

\begin{prop} \label{prop:etale-realization-fullyfaithful}
The scalar extension functor 
\[
\Vect^{\an,\varphi}((X,M_X)_\Prism^{(\str)}) \to \Vect((X,M_X)_\Prism^{(\str)},\calO_\Prism[\calI_\Prism^{-1}]^\wedge_p)^{\varphi=1}
\]
is fully faithful. 
In particular, if a Laurent $F$-crystal can be extended to an analytic prismatic $F$-crystal, then such an extension is unique. Moreover, extendability can be checked \'etale locally on $(X, M_X)$.
\end{prop}

\begin{proof}
We only treat the semistable case as the CDVR case is essentially included in the proof below.
The faithfulness follows from the construction. For the fullness, we can reduce the general case to the small affine case in \S\ref{sec:prelim-facts-rings} by taking an \'etale cover and using \cite[Rem.~4.2]{koshikawa}. So we assume $X = \Spf R$ where $R$ is small affine equipped with a $p$-completed \'etale map
\[
\square\colon R^{0} = \mathcal{O}_K \langle T_1, \ldots, T_m, T_{m+1}^{\pm 1}, \ldots, T_d^{\pm 1}\rangle / (T_1\cdots T_m - \pi) \rightarrow R \quad (m \geq 1)
\]
such that $\Spf (R/\pi R)$ has precisely $m$ irreducible components given by $R/(\pi, T_j)$ for $j = 1, \ldots, m$. Let $\Spf (\mathcal{O}_{L_j}) \rightarrow X$ be the corresponding Shilov points.   

We adopt a localization argument in the proof of 
\cite[Thm.~3.29(i)]{du-liu-moon-shimizu-completed-prismatic-F-crystal-loc-system}. Let $\mathcal{E}_{\Prism, 1}, \mathcal{E}_{\Prism, 2} \in \Vect^{\an,\varphi}((X,M_X)_\Prism^{(\str)})$, and suppose we have a map $h\colon \mathcal{E}_{\et, 1} \rightarrow \mathcal{E}_{\et, 2}$ of Laurent $F$-crystals. For each $i = 1, 2$, let $(\mathfrak{M}_i, \varphi_{\mathfrak{M}_i}, f_i) \in \mathrm{DD}_{\mathfrak{S}}$ be the object corresponding to $\mathcal{E}_{\Prism, i}$ via Lemma~\ref{lem:descentlemmamoduleversion}. Write $\mathcal{M}_i \coloneqq \mathfrak{M}_i[E^{-1}]^{\wedge}_p$. By Lemma~\ref{lem:equiv-Laurent-F-cryst-descent-datum} below, $h$ induces a map $h\colon \mathcal{M}_1 \rightarrow \mathcal{M}_2$
of \'etale $\varphi$-modules over $\mathcal{O}_{\mathcal{E}}$, which is compatible with the descent data
\[
\mathfrak{S}^{(1)}[E^{-1}]^{\wedge}_p \otimes_{p^1_1, \mathcal{O}_{\mathcal{E}}} \mathcal{M}_i \xrightarrow{\cong} \mathfrak{S}^{(1)}[E^{-1}]^{\wedge}_p\otimes_{p^1_2, \mathcal{O}_{\mathcal{E}}} \mathcal{M}_i
\]
induced from $f_i$. For each $1 \leq j \leq m$, write $\mathfrak{M}_{i, L_j} \coloneqq \mathfrak{M}_i\otimes_{\mathfrak{S}} \mathfrak{S}_{L_j}$ and $\mathcal{M}_{i, L_j} \coloneqq \mathcal{M}_i\otimes_{\mathcal{O}_{\mathcal{E}}} \mathcal{O}_{\mathcal{E}, L_j}$. By Lemma~\ref{lem:properties-analy-prism-F-crystal}(2) and \cite[Prop.~4.2.7]{gao-integral-padic-hodge-imperfect}, the map
\[
h\otimes_{\mathcal{O}_{\mathcal{E}}}\mathcal{O}_{\mathcal{E}, L_j} \colon \mathcal{M}_{1, L_j} \rightarrow \mathcal{M}_{2, L_j}
\]
is refined to a map $g_j\colon \mathfrak{M}_{1, L_j} \rightarrow \mathfrak{M}_{2, L_j}$ such that $g_j\otimes_{\mathfrak{S}_{L_j}} \mathcal{O}_{\mathcal{E}, L_j} = h\otimes_{\mathcal{O}_{\mathcal{E}}}\mathcal{O}_{\mathcal{E}, L_j}$. Then Lemma~\ref{lem:properties-analy-prism-F-crystal}(3) induces a map $g\colon \mathfrak{M}_1 \rightarrow \mathfrak{M}_2$
such that $g\otimes_{\mathfrak{S}} \mathcal{O}_{\mathcal{E}} = h \colon \mathcal{M}_1 \rightarrow \mathcal{M}_2$. Note that $g$ is compatible with the descent data $f_i$ by Lemma~\ref{lem:properties-analy-prism-F-crystal}(1). Hence, the fullness follows from Lemma~\ref{lem:descentlemmamoduleversion}. 

The last assertion follows from the fully faithfulness and \cite[Thm.~7.8]{mathew-descent}. 
\end{proof}

\begin{lem} \label{lem:equiv-Laurent-F-cryst-descent-datum}
Assume that $(X, M_X)$ is small affine as in \S\ref{sec:Breuil--Kisin log prism} or the log formal spectrum of a CDVR as in \S\ref{sec:CDVR-case-BK-prism-selfprod}.
The evaluation on the diagram
\[
p^1_1, p^1_2\colon (\mathfrak{S},(E), M_{\Spf(\mathfrak{S})}) \to (\mathfrak{S}^{(1)}, (E), M_{\Spf(\mathfrak{S}^{(1)})})
\]
gives an equivalence from $\Vect((X,M_X)_\Prism,\calO_\Prism[\calI_\Prism^{-1}]^\wedge_p)^{\varphi=1}$ to the category of triples $(\mathcal{M}, \varphi_{\mathcal{M}}, f_{\mathrm{\acute{e}t}})$ where $(\mathcal{M}, \varphi_{\mathcal{M}}) \in \Vect(\mathfrak{S}[E^{-1}]^\wedge_p)^{\varphi=1}$, and 
\[
f_{\mathrm{\acute{e}t}}\colon \mathfrak{S}^{(1)}[E^{-1}]^{\wedge}_p \otimes_{p^1_1, \mathcal{O}_{\mathcal{E}}} \mathcal{M} \stackrel{\cong}{\rightarrow} \mathfrak{S}^{(1)}[E^{-1}]^{\wedge}_p \otimes_{p^1_2, \mathcal{O}_{\mathcal{E}}} \mathcal{M}
\]
is an isomorphism of $\mathfrak{S}^{(1)}[E^{-1}]^{\wedge}_p$-modules compatible with Frobenii and satisfying the cocycle condition over $\mathfrak{S}^{(2)}[E^{-1}]^{\wedge}_p$.
The same result holds if $(X,M_X)_\Prism$ is replaced by $(X,M_X)_\Prism^{\mathrm{str}}$.
\end{lem}

\begin{proof}
This follows from Lemma~\ref{lem:BKcoversfinalobject} and \cite[Thm.~7.8]{mathew-descent} as in the proof of Lemma~\ref{lem:descentlemmavectorbundleversion}.
\end{proof}

\begin{rem} \label{rem:strict vs saturated for Laurent F crystals}
As in Remark~\ref{rem:strict vs saturated for prismatic crystals}, the above lemma directly implies that the restriction functor  
\[
\Vect((X,M_X)_\Prism,\calO_\Prism[\calI_\Prism^{-1}]^\wedge_p)^{\varphi=1} \rightarrow \Vect((X,M_X)_\Prism^{\str},\calO_\Prism[\calI_\Prism^{-1}]^\wedge_p)^{\varphi=1}
\]
is an equivalence, which also holds when $(X, M_X)$ is semistable.
\end{rem}

\begin{notation}\label{notation: Ainf and GR}
Assume further $(X,M_X)=(\Spf R, (\N^d)^a)$: either (1) small affine as in \S\ref{sec:Breuil--Kisin log prism}, or (2) $R=\calO_L$ is a CDVR as in \S\ref{sec:CDVR-case-BK-prism-selfprod}. 

Let $\overline{R}$ denote the union of finite $R$-subalgebras $R'$ inside a fixed algebraic closure of $\Frac{R}$ such that $R'[p^{-1}]$ is \'etale over $R[p^{-1}]$. Set 
\[
\calG_{R}\coloneqq \Gal(\overline{R}[p^{-1}]/R[p^{-1}]).
\]
In the CDVR case, $\calG_{R}$ is the absolute Galois group $G_L$ of $L = \calO_L[p^{-1}]$. Let $\mathrm{Rep}_{\Z_p}(\calG_R)$ (resp. $\mathrm{Rep}_{\Q_p}(\calG_R)$) denote the category of finite free $\Z_p$-modules (resp. finite-dimensional $\Q_p$-vector spaces) with continuous $\calG_R$-action. Note that $\mathrm{Rep}_{\Z_p}(\calG_R)$ (resp. $\mathrm{Rep}_{\Q_p}(\calG_R)$) is equivalent to the category (resp. the $\Q$-isogeny category) of $\Z_p$-local systems on $X_\eta$ (cf.~\cite[\S1.4, 8.2, 8.4]{kedlaya-liu-relative-padichodge}, \cite[\S15]{scholze-etalecohomologyofdiamonds}).    

In the small affine (resp. CDVR) case, choose a compatible system $(T_{i,n})$ (resp.~$(X_{i,n})$) of $p$-power roots of $T_i$ (resp.~$X_i$) in $\overline{R}$ for each $i$. This choice defines a unique $W(k)$-linear map $\mathfrak{S} \to \A_\mathrm{inf}(\overline{R})$ sending $T_i$ to $[(T_{i,n})]$ (resp. $X_i$ to $[(X_{i,n})]$). The prelog structure on $\mathfrak{S}$ induces a prelog structure on $\A_\mathrm{inf}(\overline{R})$, giving a log prism $(\A_\mathrm{inf}(\overline{R}), (E), \N^d)^a \in (X,M_X)_\Prism^{\str, \mathrm{op}}$. In fact, the induced log structure on $\A_\mathrm{inf}(\overline{R})$ is independent of the choice of $p$-power roots of coordinates by Remark~\ref{rem:log prismatic site when admitting finite free chart}. In particular, $\mathcal{G}_R$ acts on $(\A_\mathrm{inf}(\overline{R}), (E), \N^d)^a$ as well as the $\Z_p$-module $\calE(\A_\mathrm{inf}(\overline{R}))^{\varphi_\calE=1}\coloneqq\calE((\A_\mathrm{inf}(\overline{R}), (E), \N^d)^a)^{\varphi_\calE=1}$ for any Laurent $F$-crystal $(\calE,\varphi_{\calE})$ on $(X,M_X)_\Prism$. 
\end{notation}

\begin{prop} \label{prop: explicit description of etale realization for Laurent F-crystal}
With the notation as above, the equivalence of categories
\[
\Vect((X,M_X)_\Prism^{(\str)},\calO_\Prism[\calI_\Prism^{-1}]^\wedge_p)^{\varphi=1}\underset{\text{Thm.~\ref{thm:logetalerealizationofLaurentFcrystals-strict}}}{\xrightarrow{\cong}}\Loc_{\Z_p}(X_\eta)\xrightarrow{\cong} \mathrm{Rep}_{\Z_p}(\calG_R)
\]
sends a Laurent $F$-crystal $(\calE,\varphi_\calE)$ to the $\Z_p$-module $\calE(\A_\mathrm{inf}(\overline{R}))^{\varphi_\calE=1}$ equipped with a continuous $\calG_R$-action.  
\end{prop}

\begin{proof}
Since $(\Spf \A_\mathrm{inf}(\overline{R}),(E))\in X_\Prism^\perf$, the assertion follows from the proof of Theorem~\ref{thm:logetalerealizationofLaurentFcrystals-strict} and \cite[Prop.~4.7]{min-wang-rel-phi-gamma-prism-F-crys} (one of the assumptions in \emph{ibid.} is that $X$ is smooth over $\calO_K$, but the proof only requires that $X_\eta$ is smooth over $K$). 
\end{proof}

Recall that $\epsilon = (\epsilon_n)_{n \geq 0} \in \mathcal{O}_{\overline{K}}^{\flat}$ is defined in Notation and conventions in \S\ref{sec:introduction}. Write $\mu \coloneqq [\epsilon]-1$.

\begin{lem} \label{lem:isomafterinvertingmu}
Let $X=\Spf R$ be small affine with a framing as in \S\ref{sec:prelim-facts-rings}. Let $\calE_\Prism$ be an analytic prismatic $F$-crystal over $(X,M_X)$ and set $V\coloneqq T(\calE_\Prism)[p^{-1}]$. There is a natural $\calG_{R}$-equivariant isomorphism
\[
\mathcal{E}_{\Prism,\Ainf(\overline{R})}(\Spec\Ainf(\overline{R})[p^{-1}])[\mu^{-1}] \cong V\otimes_{\mathbf{Q}_p} \Ainf(\overline{R})[p^{-1}, \mu^{-1}]
\]
as $\Ainf(\overline{R})$-modules.
Here $\mathcal{E}_{\Prism,\Ainf(\overline{R})}$ denotes the associated vector bundle on $\Spec \Ainf(\overline{R})\smallsetminus V(p,E)$.
\end{lem}

\begin{proof}
With the notation as in \S\ref{sec:analytic prismatic F-crystal in small affine case}, let $\mathfrak{M} = H^0(U, \mathcal{E}_{\Prism,\mathfrak{S}})$ where $U \coloneqq \Spec\mathfrak{S} \smallsetminus V(p, E)$. Without loss of generality, we may assume that $\mathcal{E}_{\Prism,\mathfrak{S}}$ is effective, and let $r$ be a positive integer such that $(\mathfrak{M}, \varphi_{\mathfrak{M}})$ has $E$-height $\leq r$. Consider the $\Ainf(\overline{R})$-module $M \coloneqq \mathfrak{M}\otimes_{\mathfrak{S}} \Ainf(\overline{R})$ with the induced tensor product Frobenius. Since $\mathfrak{M}[p^{-1}] = \mathcal{E}_{\Prism,\mathfrak{S}}(\Spec \mathfrak{S}[p^{-1}])$, we have
\[
M[p^{-1}] \cong \mathcal{E}_{\Prism,\Ainf(\overline{R})}(\Spec\Ainf(\overline{R})[p^{-1}])
\]
as $\Ainf(\overline{R})$-modules compatibly with Frobenius (see the proof of Lemma~\ref{lem:descentlemmamoduleversion}). Furthermore, by Proposition~\ref{prop: explicit description of etale realization for Laurent F-crystal} and the theory of \'etale $\varphi$-modules, we have a natural identification
\[
V = (M[p^{-1}]\otimes_{\Ainf(\overline{R})} W(\overline{R}^\flat[(\pi^{\flat})^{-1}]))^{\varphi = 1}
\]
and
\[
M[p^{-1}]\otimes_{\Ainf(\overline{R})} W(\overline{R}^\flat[(\pi^{\flat})^{-1}])) = V\otimes_{\mathbf{Z}_p} W(\overline{R}^\flat[(\pi^{\flat})^{-1}]).
\]
We need to show that the above equality descends to that over $\Ainf(\overline{R})[p^{-1}, \mu^{-1}]$ as in the statement of the lemma.

We follow a similar argument as in the proof of \cite[Thm.~3.29(i)]{du-liu-moon-shimizu-completed-prismatic-F-crystal-loc-system} using the localization in the end of Section~\ref{sec:prelim-facts-rings}. As in Notation~\ref{notation:localization}, for each $\mathfrak{P} \in \mathcal{P}$, consider the map $R \rightarrow \mathcal{O}_{L(\mathfrak{p})}$ and a choice of map $\overline{R}[p^{-1}] \rightarrow \overline{L(\mathfrak{p})}$. Then $\mathcal{E}_\Prism|_{(\Spf \mathcal{O}_{L(\mathfrak{p})}, M_{\Spf \mathcal{O}_{L(\mathfrak{p}})})_{\Prism}}$ is an analytic prismatic $F$-crystal over $(\Spf \mathcal{O}_{L(\mathfrak{p})}, M_{\Spf \mathcal{O}_{L(\mathfrak{p})}})$. By Corollary~\ref{cor:CDVR-equivalence-descent-datum}, this extends uniquely to a prismatic $F$-crystal over $(\Spf \mathcal{O}_{L(\mathfrak{p})}, M_{\Spf \mathcal{O}_{L(\mathfrak{p})}})$, which we denote by $\mathcal{E}_\Prism|_{\mathcal{O}_{L(\mathfrak{p})}}$ abusing notation. For the \'etale realization $T(\mathcal{E}_\Prism|_{\mathcal{O}_{L(\mathfrak{p})}})$ of $\mathcal{E} _\Prism|_{\mathcal{O}_{L(\mathfrak{p})}}$, we have $T(\mathcal{E}_\Prism|_{\mathcal{O}_{L(\mathfrak{p})}})[p^{-1}] = V|_{G_{\overline{L(\mathfrak{p})}}}$ as $\mathbf{Q}_p$-representations of $G_{\overline{L(\mathfrak{p})}}$ by Proposition~\ref{prop:pullback commutes with etale realization}. Write $M_{\mathfrak{p}} \coloneqq H^0(\Spec(\Ainf(\overline{\mathcal{O}_{L(\mathfrak{p})}})),(\mathcal{E}_\Prism|_{\mathcal{O}_{L(\mathfrak{p})}})_{\Ainf(\overline{\mathcal{O}_{L(\mathfrak{p})}})})$, which is a finite free module over $\Ainf(\overline{\mathcal{O}_{L(\mathfrak{p})}})$. Note $M[p^{-1}]\otimes_{\Ainf(\overline{R})} \Ainf(\overline{\mathcal{O}_{L(\mathfrak{p})}}) = M_{\mathfrak{p}}[p^{-1}]$.

By the same argument as in the proof of \cite[Thm.~3.29(i)]{du-liu-moon-shimizu-completed-prismatic-F-crystal-loc-system}, we deduce
\[
M_{\mathfrak{p}} \subset T(\mathcal{E}_\Prism |_{\mathcal{O}_{L(\mathfrak{p})}})\otimes_{\mathbf{Z}_p} \Ainf(\overline{\mathcal{O}_{L(\mathfrak{p})}}) \subset \frac{1}{\mu^r} M_{\mathfrak{p}}
\]
as $\Ainf(\overline{\mathcal{O}_{L(\mathfrak{p})}})$-modules. Thus, by Lemmas~\ref{lem:intersection-modules-flat-base-change}(2) and \ref{lem:intersection-A_inf-localizations}, we obtain 
\[
M[p^{-1}, \mu^{-1}] \cong V\otimes_{\mathbf{Q}_p}\Ainf(\overline{R})[p^{-1}, \mu^{-1}].
\]
\end{proof}

\subsection{Semistable \texorpdfstring{$\Z_p$}{Zp}-local systems} \label{subsec:prismatic-semistable-local-system}

In this subsection, we assume that $(X,M_X)$ is semistable or the log formal spectrum of a CDVR.

\begin{defn} \label{defn:semistablity-via-analy-prismatic-F-crystal}
The \'etale realization functor $T_X$ for analytic prismatic $F$-crystals is fully faithful by Proposition~\ref{prop:etale-realization-fullyfaithful}. A $\Z_p$-local system over $X_\eta$ is said to be \emph{(($X$-)prismatic) semistable} if it is in the essential image of $T_X$.
\end{defn}

When $X_\eta$ is an affinoid, an affine semistable formal model is unique if it exists. 

\begin{rem} \label{rem:first-properties-of-prismatic-semistability}
By Proposition~\ref{prop:pullback commutes with etale realization}, prismatic semistable $\Z_p$-local systems are stable under pullbacks. Furthermore, the following results will be proved in this paper: 
\begin{enumerate}
\item When a rigid-analytic variety $\calX$ over $K$ admits two semistable models $X$ and $X'$ over $\mathcal{O}_K$, a $\Z_p$-local system on $\calX$ is $X$-semistable if and only if it is $X'$-semistable (Corollary~\ref{cor:independent-model}).
\item A $\Z_p$-local system on $X_{\eta}$ is prismatic semistable if and only if it is associated to an $F$-isocrystal on $(X_1, M_{X_1})_{\CRIS}$ as in Definition~\ref{defn:association} (Corollary~\ref{cor:semistable-prismatic-Faltings-equivalent}).
\item In the CDVR case where $X=\Spf \mathcal{O}_L$, prismatic semistable $\Z_p$-local systems correspond to lattices in semistable $\Q_p$-representations of $G_L$ (Theorem~\ref{thm:CDVR-semistable-notions-equivalent}).
\item Prismatic semistable $\Z_p$-local systems are de Rham: when $(X, M_X)$ is semistable, this is Corollary~\ref{cor:prismatic-semistable-deRham-semistable-case}. In the CDVR case, it follows from Theorem~\ref{thm:CDVR-semistable-notions-equivalent}. 
\end{enumerate}
\end{rem}

\subsection{Crystalline realization}\label{sec:crystallinerealization}

In this subsection, we assume that $(X, M_X)$ is semistable or the log formal spectrum of a CDVR. Set $(X_1,M_{X_1})\coloneqq (X,M_X)\otimes_{\Z_p}\Z_p/p$ and let $(X_1,M_{X_1})_\CRIS$ denote the absolute logarithmic crystalline site of $(X_1,M_{X_1})$ introduced in Definition~\ref{def:big absolute log crystallien site}. We will attach to an analytic prismatic $F$-crystal on $X$ an $F$-isocrystal on $(X_1,M_{X_1})_\CRIS$: first, we give a construction in the small affine case and then show that this local construction glues and globalizes. Since $\Vect^{\an,\varphi}((X,M_X)_\Prism)\xrightarrow{\cong}\Vect^{\an,\varphi}((X,M_X)_\Prism^\mathrm{str})$ under our assumption on $X$ by Remark~\ref{rem:strict vs saturated for prismatic crystals}, we work on the strict prismatic site $(X,M_X)_\Prism^\mathrm{str}$ below.

\smallskip

\noindent
\textbf{Small affine case.} Assume that $(X, M_X)=(\Spf(R),\N^d)^a$ is small affine with a framing 
$\square$ as in the beginning of \S\ref{sec:Breuil--Kisin log prism}.
The framing defines the Breuil--Kisin log prism $(\fkS_{R,\square},(E),M_{\Spf\fkS})$. 

We define $S_{R,\square}$ to be the $p$-adic completion of the log PD-envelope of $\fkS_{R,\square}$ with respect to the kernel of the exact surjection $\fkS_{R,\square} \to R/p$ (see Example~\ref{eg:BreuilSR}). 
Set $D_X \coloneqq \Spf(S_{R,\square})$ and write $M_{D_X}$ for the log structure on $D_X$ induced from the prelog structure $\N^d$ with $e_i \mapsto T_i$. Then the exact closed immersion $X_1\hookrightarrow D_X$ is a log PD-thickening, and $(X_1, D_X, M_{D_X})$ defines an ind-object $\varinjlim_n (X_1,\Spec (S_{R,\square}/p^n))$ of $(X_1,M_{X_1})_{\CRIS}^{\mathrm{aff}}$.

\begin{lem}\label{lem:SRweaklyfinal}
The pair $(X_1,D_X,M_{D_X})$ is a weakly final ind-object in $(X_1,M_{X_1})_{\CRIS}^{\mathrm{aff}}$.
\end{lem}

\begin{proof}
This follows from Lemma~\ref{lem:weakly-final} and Example~\ref{eg:BreuilSR}.
\end{proof}

For each $i\in\N$, let $(X_1,\Spf S_{R,\square}^{(i)},M_{\Spf S_{R,\square}^{(i)}})$ be the $(i+1)$-st self-product of $(X_1, D_X, M_{D_X})$ as an ind-object of $(X_1,M_{X_1})_{\CRIS}^{\aff}$ (see Lemma~\ref{lem:nonempty finite limit in abs log cristalline site} and the paragraph before Lemma~\ref{lem:weakly-final}). 

By \cite[Cor.~2.39]{bhatt-scholze-prismaticcohom}, we have $S_{R, \square} \simeq \fkS_{R,\square}\{\varphi(E)/p\}_\delta^\wedge$, where $\{-\}_\delta$ denotes adjoining elements as $\delta$-ring. This gives a well-defined map $\varphi\colon(\fkS_{R,\square},(E)) \to (S_{R, \square},(p))$. The pullback of the prelog structure on $\fkS_{R,\square}$ defines a prelog structure $\varphi^\ast\N^d$ on $S_{R, \square}$, and we obtain $(S_{R, \square},(p),\varphi^\ast \N^d)^a \in (X,M_X)_\Prism^\mathrm{op}$. Similarly, a prelog structure $\varphi^\ast \N^d\rightarrow S_{R,\square}\xrightarrow{} S_{R,\square}^{(i)}$ defines $(S_{R,\square}^{(i)},(p), \varphi^\ast \N^d)^a\in (X,M_X)_\Prism^\mathrm{op}$ with a map 
\[
\varphi\colon (\fkS_{R,\square}^{(i)},(E), \N^d)^a\rightarrow(S_{R,\square}^{(i)},(p), \varphi^\ast \N^d)^a, 
\]
and the induced log structure is independent of the choice of the projection $S_{R,\square}\xrightarrow{} S_{R,\square}^{(i)}$. We remark that the $\delta$-structures of $S_{R,\square}$ and $S_{R,\square}^{(i)}$ depends on the choice of $\square$ while the underlying ring structure of $S_{R,\square}$ and the induced log structure on $\Spf (S_{R,\square})$ are independent of the choice of $\square$. 

\begin{construction} \label{construction:crystal-small-affine}
Keep the small affine assumption and the notation as above.
Let $(\calE_\Prism,\varphi_{\calE_\Prism})$ be an analytic prismatic $F$-crystal over $(X,M_X)_\Prism$. Let $(\mathfrak{M}, \varphi_{\mathfrak{M}}, f)$ be the object in $\mathrm{DD}_{\mathfrak{S}_{R,\square}}$ associated to $(\mathcal{E}_{\Prism},\varphi_{\calE_\Prism})$ by Lemma~\ref{lem:descentlemmamoduleversion}. For each $n \geq 1$, set 
\[
\mathcal{M}_n \coloneqq S_{R, \square}\otimes_{\varphi, \mathfrak{S}_{R, \square}} \mathfrak{M}/p^n\mathfrak{M}. 
\]
Via the map $\varphi\colon \mathfrak{S}_{R, \square}^{(i)} \rightarrow S_{R, \square}^{(i)}$ for $i = 1, 2$, the isomorphism
\[
f\colon \mathfrak{S}_{R, \square}^{(1)}\otimes_{p^1_1,\mathfrak{S}_{R, \square}}\mathfrak{M}\xrightarrow{\cong}\mathfrak{S}_{R, \square}^{(1)}\otimes_{p^1_2,\mathfrak{S}_{R, \square}}\mathfrak{M}
\]
induces an isomorphism $\eta_n\colon S_{R, \square}^{(1)}\otimes_{p^1_1,S_{R, \square}}\mathcal{M}_n \xrightarrow{\cong} S_{R, \square}^{(1)}\otimes_{p^1_2, S_{R, \square}}\mathcal{M}_n$ for each $n$, satisfying the cocycle condition on $S_{R, \square}^{(2)}$. By Proposition~\ref{prop: crystals and quasi-nilpotent connections}, the pair $(\mathcal{M},\eta)\coloneqq \varprojlim_n(\calM_n,\eta_n)$ corresponds to a finitely generated quasi-coherent crystal $\calE_{\cris,S_{R,\square}} \in \CR((X_1,M_{X_1})_\CRIS)$.
We have an isomorphism $\mathcal{M}[p^{-1}] \cong S_{R, \square}[p^{-1}]\otimes_{\varphi, \mathfrak{S}_{R, \square}[p^{-1}]} \mathfrak{M}[p^{-1}]$ of finite projective ${S}_{R, \square}[p^{-1}]$-modules (see \cite[Lem.~3.24(iv)]{du-liu-moon-shimizu-completed-prismatic-F-crystal-loc-system}). Furthermore, $\varphi_{\mathfrak{M}}$ induces a Frobenius 
\[
\varphi_{\mathcal{M}}\colon \varphi^* \mathcal{M}[p^{-1}] \stackrel{\cong}{\rightarrow} \mathcal{M}[p^{-1}]
\]
compatible with $\eta$. This yields an $F$-isocrystal on $(X_1,M_{X_1})_\CRIS$. 
\end{construction}

\smallskip
\noindent
\textbf{Globalization.} To go beyond the small affine case, we glue $\calE_{\cris,R,\square}$ to obtain a crystal on $(X_1,M_{X_1})_\CRIS$. For this, we need to show that $\calE_{R,\square}$ is independent of the choice of framing $\square$.

\begin{lem}\label{lem:independence of framing for associated crystalline crystal}
For two framings $\square$ and $\square'$ of $R$, there is a canonical identification $\calE_{\cris,R,\square}\cong \calE_{\cris,R,\square'}$ of crystals on $(X_1,M_{X_1})_\CRIS$, which is compatible with the Frobenii on the associated $F$-isocrystals.
\end{lem}

\begin{proof}
Let $(\fkS_{\square,\square'},(E),M_{\square,\square'})$ be the coproduct $(\fkS_{R,\square},(E),M_{\Spf \fkS_{R,\square}})\amalg(\fkS_{R,\square'},(E),M_{\Spf \fkS_{R,\square'}})$ in $(X,M_X)_\Prism^{\mathrm{str},\mathrm{op}}$, given by Lemma~\ref{lem:coproduct with fkS is a covering}. Let $S_{\square,\square'}$ denote the ring corresponding to the product of ind-objects $\Spf S_{R,\square}$ and $\Spf S_{R,\square'}$ in $(X_1,M_{X_1})_\CRIS$; the product is represented by an affine $p$-adic log PD-formal scheme by Lemma~\ref{lem:nonempty finite limit in abs log cristalline site}(1). By construction, $(S_{\square, \square'}, (p), M_{\Spf S_{\square, \square'}}) \in (X,M_X)_\Prism^{\mathrm{str},\mathrm{op}}$, and we have a commutative diagram
\[
\xymatrix{
S_{R,\square}\ar[r]_-{\pr_{\square}}
& S_{\square,\square'}
& S_{R,\square'}\ar[l]^-{\pr_{\square'}}\\
\fkS_{R,\square}\ar[u]^-{\varphi} \ar[r]_-{\pr_{\square}}
& \fkS_{\square,\square'}\ar[u]
& \fkS_{R,\square'}\ar[u]^-{\varphi}\ar[l]^-{\pr_{\square'}}.
}
\]
Thus, we obtain a canonical isomorphism 
\[
\pr_{\square}^\ast \calE_{\cris,S_R,\square}
\cong \pr_{\square'}^\ast \calE_{\cris,S_R,\square'}
\]
of $S_{\square, \square'}$-modules. 

By considering the self-product of $\fkS_{\square,\square'}$ (resp.~$S_{\square,\square'}$) in the absolute log prismatic site (resp.~absolute log crystalline site) and arguing similarly, we see that the identification is upgraded to the identification of the pullbacks of the pairs $(\calE_{\cris,S_{R,\square}}, \eta)$ and $(\calE_{\cris,S_{R,\square'}}, \eta)$. By Proposition~\ref{prop: crystals and quasi-nilpotent connections}, we obtain a canonical identification
$\calE_{\cris,R,\square}\cong \calE_{\cris,R,\square'}$ of crystals on $(X_1,M_{X_1})_\CRIS$. This is compatible with the Frobenii on the associated isocrystals, since all the relevant maps above are compatible with Frobenii.
\end{proof}

\begin{cor}[Crystalline realization of analytic prismatic $F$-crystals]  \label{cor:crystalline-realization}
Let $(X,M_X)$ be a semistable $p$-adic formal scheme over $\calO_K$.
For $(\calE_\Prism,\varphi_{\calE_\Prism})\in \Vect^{\an,\varphi}((X,M_X)_\Prism)$, there exists a unique finitely generated quasi-coherent crystal $\calE_\cris$ on $(X_1,M_{X_1})_\CRIS$ such that for each small affine open $\Spf R$ of $X$ with framing $\square$, the restriction of $\calE_{\cris}$ to $((\Spf R)_1,M_{(\Spf R)_1})_\CRIS$ coincides with $\calE_{\cris,R,\square}$. 

Moreover, $\calE_{\cris}$ gives rise to an $F$-isocrystal $(\calE_{\cris, \mathbf{Q}},\varphi_{\calE_{\cris, \mathbf{Q}}})$ on $(X_1,M_{X_1})_\CRIS$ such that 
\begin{equation}\label{eq:crystalline realization}
\calE_\Prism((S_{R,\square}^{(i)},(p),\varphi^{\ast}\N^d)^a) \cong \calE_{\cris, \mathbf{Q}}(X_1,\Spf S_{R,\square}^{(i)},M_{\Spf S_{R,\square}^{(i)}})
\end{equation}
for each $i \geq 0$, and the assignment $\calE_\Prism \mapsto \calE_{\cris, \mathbf{Q}}$ defines a functor 
\[
D_{\cris}\colon \Vect^{\an,\varphi}((X,M_X)_\Prism) \to \Vect_\Q^\varphi((X_1,M_{X_1})_\CRIS).
\]
\end{cor}

We call $D_\cris$ the \emph{crystalline realization functor}.

\begin{proof}
For the first statement, take a $p$-completely \'etale covering $\{f_\lambda\colon X_\lambda\rightarrow X\}_{\lambda\in \Lambda}$ such that each $X_\lambda$ is small affine with framing $\square_\lambda$. Construction~\ref{construction:crystal-small-affine} gives a finitely generated quasi-coherent crystal $\calE_{\cris,\square_\lambda}$ on $(X_{\lambda,1},M_{X_{\lambda,1}})_\CRIS$. By applying  Lemma~\ref{lem:independence of framing for associated crystalline crystal} to an affine open covering of $X_{\lambda\lambda'}\coloneqq X_\lambda\times_XX_{\lambda'}$, we see that there is a canonical identification of the pullbacks $\calE_{\cris,\square_\lambda}|_{X_{\lambda\lambda',1}}\cong \calE_{\cris,\square_{\lambda'}}|_{X_{\lambda\lambda',1}}$ on $(X_{\lambda\lambda',1},M_{X_{\lambda\lambda',1}})_\CRIS$ that satisfies the cocycle condition over $(X_{\lambda\lambda'\lambda'',1},M_{X_{\lambda\lambda'\lambda'',1}})_\CRIS$ for $X_{\lambda\lambda'\lambda''}\coloneqq X_\lambda\times_XX_{\lambda'}\times_XX_{\lambda''}$ ($\lambda,\lambda',\lambda''\in \Lambda$). Since \'etale descent holds for quasi-coherent crystals, we obtain the desired crystal $\calE_\cris$ on $(X_1,M_{X_1})_\CRIS$. The remaining statements follow from Construction~\ref{construction:crystal-small-affine} and \cite[Lem.~3.24(iv)]{du-liu-moon-shimizu-completed-prismatic-F-crystal-loc-system}. 
\end{proof}

\smallskip
\noindent
\textbf{Crystalline realization in terms of $(X_1,M_{X_1})_\Prism^{\mathrm{str}}$.}
In the rest of this subsection, we will give another description of the crystalline realization functor\footnote{We remark that this alternative description (Lemma~\ref{lem:S-lattice of rational DD} and Propositions~\ref{prop: commutative-diagram-Dcris} and \ref{prop: alternative-crystalline-realization}) is not necessary for the main results in this paper, and is included here mainly to explain the relation with \cite[Ex.~4.7, Const.~4.12]{bhatt-scholze-prismaticFcrystal}.} similar to \cite[Ex.~4.7, Const.~4.12]{bhatt-scholze-prismaticFcrystal} and \cite[Const.~3.9, Thm.~6.4]{GuoReinecke-Ccris}:
let $i_X\colon(X_1,M_{X_1}) \to (X,M_X)$ denote the exact closed immersion. We will explain that the restriction
$\Vect^{\an,\varphi}((X_1,M_{X_1})_\Prism)\rightarrow\Vect^{\an,\varphi}((X_1,M_{X_1})_\Prism^{\mathrm{str}})$ is an equivalence, and
there is a fully faithful functor $\Vect^{\an,\varphi}((X_1,M_{X_1})_\Prism^{(\mathrm{str})}) \xhookrightarrow{} \Vect_\Q^\varphi((X_1,M_{X_1})_\CRIS)$ which fits into the following 2-commutative diagram:
\begin{equation}\label{diagram:crystalline realization}
\xymatrix{
\Vect^{\an,\varphi}((X,M_X)_\Prism^{(\mathrm{str})})\ar[rd]_-{D_\cris}\ar[r]^-{i_{X,\Prism}^{-1}}
& \Vect^{\an,\varphi}((X_1,M_{X_1})_\Prism^{(\mathrm{str})})\ar@{^{(}->}[d]^-{D_{\cris,X_1}}\\
& \Vect_\Q^\varphi((X_1,M_{X_1})_\CRIS).
}
\end{equation}
Note that every prism $(A,I,M_{\Spf A})$ in $ (X_1,M_{X_1})_\Prism$ is crystalline, i.e., $I=(p)$, and vector bundles over $\Spec A\smallsetminus V(p,I)$ correspond to finite projective $A[p^{-1}]$-modules. 

First assume that $X=\Spf R$ is small affine. Set $\fkS\coloneqq \fkS_{R,\square}$ and $S^{(i)}\coloneqq S^{(i)}_{R,\square}$, and keep the notation at the beginning of this subsection.

\begin{lem}\label{lem: coproduct of S in prismatic site}
For each $i \geq 0$, the $(i+1)$-st self-coproducts of $(S,(p), \varphi^\ast \N^d)^a$ in 
$(X_1,M_{X_1})_\Prism^{\mathrm{str},\mathrm{op}}$ and $(X_1,M_{X_1})_\Prism^{\mathrm{op}}$ exist and coincide. Denote the common log prism by $(S^{[i]},(p),\varphi^\ast \N^d)^a$. Then all the projections are $p$-complete faithfully flat. Moreover, the map $(S^{[i]},(p),\varphi^\ast \N^d)^a \to (S^{(i)},(p),\varphi^\ast \N^d)^a$ induced by the universal property of $(S^{[i]},(p),\varphi^\ast \N^d)^a$ is an isomorphism.
\end{lem}

\begin{proof}
In Example~\ref{eg:coproduct of Frob twists of BK prism}, for each $i \geq 0$, we constructed the $(i+1)$-st self-coproduct $(\fkS^{[i]},(\varphi(E)),\varphi^\ast\N^d)^a$ of $(\fkS,(\varphi(E)),\varphi^\ast\N^d)^a$ inside $(X,M_{X})_\Prism^{\mathrm{str},\mathrm{op}}$, which also represents the $(i+1)$-st self-coproduct in $(X,M_{X})_\Prism^{\mathrm{op}}$. We claim $S^{[i]} \simeq \fkS^{[i]}\{\varphi(E)/p\}^\wedge_\delta$. Once the claim is proved, \cite[Prop. 3.13(3)]{bhatt-scholze-prismaticcohom} shows $\fkS^{[i]}\{\varphi(E)/p\}^\wedge_\delta \simeq \fkS^{[i]}\widehat{\otimes}^\mathbb{L}_{\fkS}\fkS\{\varphi(E)/p\}^\wedge_\delta$. By Example~\ref{eg:coproduct of Frob twists of BK prism},  $\fkS^{[i]}\widehat{\otimes}^\mathbb{L}_{\fkS}\fkS\{\varphi(E)/p\}^\wedge_\delta \xrightarrow{\cong} \fkS^{[i]}\widehat{\otimes}_{\fkS}\fkS\{\varphi(E)/p\}^\wedge_\delta$, and the projections are $p$-completely faithfully flat. 

To see the claim, observe that the self-coproducts of $(S,(p), \varphi^\ast \N^d)^a$ in $(X_1,M_{X_1})_\Prism^{(\mathrm{str},)\mathrm{op}}$ and $(X,M_X)_\Prism^{(\mathrm{str},)\mathrm{op}}$ agree and that the map $(\fkS,(\varphi(E)),\varphi^\ast\N^d) \to (S,(p),\varphi^\ast\N^d)$ is initial among all maps $(\fkS,(\varphi(E)),\varphi^\ast\N^d) \to (A,I,M_{\Spf A})$ with $I=(p)$. This implies $S^{[i]} \simeq \fkS^{[i]}\{\varphi(E)/p\}^\wedge_\delta$ by the universal property of $(\fkS^{[i]},(\varphi(E)),\varphi^\ast\N^d)^a$.

To show that the induced map $S^{[i]} \to S^{(i)}$ is an isomorphism, we will construct its inverse. For simplicity, we only deal with the case $i=1$. By the claim and \cite[Cor.~2.39]{bhatt-scholze-prismaticcohom}, $S^{[1]}$ is the $p$-adic completion of PD-envelope of $\fkS^{[1]}$ with respect to the ideal generated by $E(u\otimes 1)$, which is $\Z_p$-flat. By the construction of $\fkS^{[1]}$, we have $(T_j^p\otimes 1)/(1\otimes T_j^p) \in (\fkS^{[1]})^\times$ for all $1\leq j \leq d$. It follows that the two prelog structures $\N^d \to S^{[1]}$ defined by $e_j \mapsto T_j\otimes 1$ and $e_j \mapsto 1\otimes T_j$ induce the same log structure $N_{\Spf S^{[1]}}$. Let $I$ be the $p$-adic closure of the PD-ideal of $S^{[1]}$ generated by $p$ and $(E(u\otimes1))$. Then $(\Spf(S^{[1]}/I), \Spf S^{[1]},N_{\Spf S^{[1]}})$ defines an ind-object in $(X_1,M_{X_1})_{\CRIS}^{\mathrm{aff}}$, which admits two maps from $(X_1, D_X, M_{D_X})$ where $D_X=\Spf S$. By the universal property, these give a map $S^{(1)} \to S^{[1]}$, which can be checked to be the inverse of $S^{[1]} \to S^{(1)}$.
\end{proof}

\begin{lem}\label{lem:S covers final object over X_1}
    The log prism $(S,(p),\varphi^\ast\N^d)^a$ covers the final object in $\operatorname{Sh}((X_1,M_{X_1})_\Prism)$ and $\operatorname{Sh}((X_1,M_{X_1})_\Prism^{\mathrm{str}})$.
\end{lem}

\begin{proof}
Observe that, for every prism $(C,(\varphi(E)),\varphi^\ast\N^d)^a$ that is strict over $(\mathfrak{S},(\varphi(E)),\varphi^\ast\N^d)^a$ and every $(B,(p),M_{\Spf B})\in (X_1,M_{X_1})_\Prism^\mathrm{op}$, any map $(C,(\varphi(E)),\varphi^\ast\N^d)^a \to (B,(p),M_{\Spf B})$ in $(X, M_X)_{\Prism}^{\mathrm{op}}$ uniquely factors as
\[
(C,(\varphi(E)),\varphi^\ast
\N^d)^a \rightarrow (C\{\varphi(E)/p\}_\delta^\wedge,(p),\varphi^\ast\N^d)^a \rightarrow(B,(p),M_{\Spf B}).
\]
To see this, use the rigidity of distinguished elements \cite[Lem.~2.24]{bhatt-scholze-prismaticcohom}. 

Let $(\Prism_{R_\infty},I,\N^d)$ be the perfectoid prelog prism over $(\mathfrak{S},(E),
\N^d)$ constructed in Lemma~\ref{lem:BKcoversfinalobject}, and write $\Prism_{R_\infty/p}\coloneqq \Prism_{R_\infty}\{\varphi(E)/p\}^\wedge_\delta$. Then $\fkS \to \Prism_{R_\infty}$ induces a map $(S,(p),\varphi^\ast\N^d)^a \to (\Prism_{R_\infty/p},(p),\varphi^\ast\N^d)^a$ via the identification of $S$ with $\mathfrak{S}\{\varphi(E)/p\}_\delta^\wedge$. It is enough to show $(\Prism_{R_\infty/p},(p),\varphi^\ast\N^d)^a$ covers the final object in $\operatorname{Sh}((X_1,M_{X_1})_\Prism)$. 

Take any $(A,(p),M_{\Spf A})\in (X_1,M_{X_1})_\Prism^\mathrm{op}$. By Lemma~\ref{lem:BKcoversfinalobject}, there exists a cover $f\colon (A,(p),M_{\Spf A})\rightarrow (P,(p),M_{\Spf P})$ whose target admits a map from $(\Prism_{R_\infty},(E),\N^d)^a$. Consider the following commutative diagram 
\[
\begin{tikzcd}
(\Prism_{R_\infty},(E),\N^d)^a \arrow[r] \arrow[d, "\varphi"] & (P,(p),M_{\Spf P}) \arrow[d]  & (A,(p),M_{\Spf A}), \arrow[l,"f"] \arrow[dl, "g"] \\
(\Prism_{R_\infty},(\varphi(E)),\varphi^\ast\N^d)^a \arrow[r] \arrow[d] & (B,(p),M_{\Spf B}) &  \\
(\Prism_{R_\infty/p},(p),\varphi^\ast\N^d)^a \arrow[ru,dashrightarrow,bend right=10] & &
\end{tikzcd}
\]
where $(B,(p),M_{\Spf B})$ is defined so that the left upper corner is a pushout diagram, and the dotted arrow from $(\Prism_{R_\infty/p},(p),\varphi^\ast\N^d)^a$ is given by the above observation. Since both $\varphi$ and $f$ are covers, so is $g$. Hence $(\Prism_{R_\infty/p},(p),\varphi^\ast\N^d)^a$ covers the final object in $\operatorname{Sh}((X_1,M_{X_1})_\Prism)$.
\end{proof}

Let $\mathrm{DD}_{S[p^{-1}]}$ denote the category consisting of triples $(\calM,\varphi_{\calM},f_{\calM})$ where $\calM$ is a finitely projective $S[p^{-1}]$-module, $\varphi_{\calM}\colon \varphi^\ast \calM \to \calM$ is an $S[p^{-1}]$-linear isomorphism, and $f_{\calM}\colon S^{(1)}[p^{-1}]\otimes_{p_1^1,S[p^{-1}]}\calM\xrightarrow{\cong}S^{(1)}[p^{-1}]\otimes_{p_2^1,S[p^{-1}]}\calM$ is a descent isomorphism of $\calM$ over $S^{(1)}[p^{-1}]$ that is $\varphi$-equivariant and satisfies the cocycle condition over $S^{(2)}[p^{-1}]$.

Every object of $\mathrm{DD}_{S[p^{-1}]}$ admits a finitely generated $p$-adically complete $S$-lattice stable under the descent isomorphism.

\begin{lem}\label{lem:S-lattice of rational DD}
For every $(\calM,\varphi_{\calM},f_{\calM})\in \mathrm{DD}_{S[p^{-1}]}$, there exists a finitely generated $p$-adically complete $S$-submodule $\calM^+\subset\calM$ such that $\calM^+[p^{-1}]=\calM$ and $f_\calM$ restricts to an HPD-stratification on $\calM^+$. Via Proposition~\ref{prop: crystals and quasi-nilpotent connections}, this construction defines a fully faithful functor
\[
\mathcal{F} \colon\mathrm{DD}_{S[p^{-1}]} \hookrightarrow\Vect_\Q^\varphi((X_1,M_{X_1})_\CRIS). 
\]
\end{lem}

\begin{proof}
The first assertion follows from \cite[Lem.~3.21]{Inoue-log-prismatic-F-crystals}. By Proposition~\ref{prop: crystals and quasi-nilpotent connections}, such an $S$-submodule $\calM^+$ given by the first statement defines a coherent crystal $\calF_{\calM^+}$ on $(X_1,M_{X_1})_\CRIS$, and the associated isocrystal $\calF_{\calM^+,\Q}$ is independent of the choice of $\calM^+$ and equipped with the Frobenius structure. Note that $\calF_{\calM^+,\Q}$ is finite locally free by Lemma~\ref{lem:DLMS3.24(iv)} and is functorial in $(\calM,\varphi_{\calM},f_{\calM})\in \mathrm{DD}_{S[p^{-1}]}$. Since $(\calM,\varphi_{\calM},f_{\calM})$ can be recovered from the $F$-isocrystal $\calF_{\calM^+,\Q}$ by evaluating on the diagram 
    \[
    (X_1, \Spf S,M_{\Spf S}) \rightarrow (X_1, \Spf S^{(1)},M_{\Spf S^{(1)}})\leftarrow (X_1, \Spf S,M_{\Spf S}),
    \]
the second assertion follows.
\end{proof}

\begin{prop} \label{prop: commutative-diagram-Dcris}
Keep the assumption and notation as above.
\begin{enumerate}
\item Evaluating on the diagram 
    \[
    (S,(p),\varphi^\ast \N)^a \rightarrow (S^{(1)},(p),\varphi^\ast \N)^a \leftarrow (S,(p),\varphi^\ast \N)^a
    \]
    defines a commutative diagram 
    \[
\xymatrix{
\Vect^{\an,\varphi}((X_1,M_{X_1})_\Prism)\ar[r]\ar[rd]_-{\cong}
&\Vect^{\an,\varphi}((X_1,M_{X_1})^{\mathrm{str}}_\Prism)\ar[d]^-{\cong}\\
&\mathrm{DD}_{S[p^{-1}]},}
\]
where the evaluation functors are equivalences. In particular, the horizontal restriction functor is also an equivalence.
    \item The functors in (1) and Lemma~\ref{lem:S-lattice of rational DD} define a fully faithful functor  
\[
D_{\cris,X_1}\colon\Vect^{\an,\varphi}((X_1,M_{X_1})_\Prism^{(\mathrm{str})})\xrightarrow{\cong} \mathrm{DD}_{S[p^{-1}]} \hookrightarrow \Vect_\Q^\varphi((X_1,M_{X_1})_\CRIS),
\]
and $D_{\cris}$ is identified with $D_{\cris,X_1}\circ i_\Prism^{-1}$. 
\end{enumerate}
\end{prop}

\begin{proof}
(1) follows from Lemmas~\ref{lem: coproduct of S in prismatic site} and \ref{lem:S covers final object over X_1} (cf.~proof of Lemma~\ref{lem:descentlemmavectorbundleversion}).

For (2), take $(\calE_\Prism,\varphi_{\calE_\Prism})\in \Vect^{\mathrm{an},\varphi}((X,M_X)_\Prism)$. Since $i\colon (X_1,M_{X_1})\rightarrow (X,M_X)$ is strict, we have
\[
(i_\Prism^{-1}\calE_\Prism)((S,(p),\varphi^\ast \N)^a)=\calE_\Prism((S,(p),\varphi^\ast \N)^a).
\]
The compatibility follows from this and the isomorphism \eqref{eq:crystalline realization}.
\end{proof}

Arguing as in Remark~\ref{rem:strict vs saturated for prismatic crystals} and Lemma~\ref{lem:independence of framing for associated crystalline crystal}, we obtain the global results:

\begin{prop} \label{prop: alternative-crystalline-realization}
If $(X,M_X)$ is semistable, then there is a fully faithful functor
\[
D_{\cris,X_1}\colon
\Vect^{\an,\varphi}((X_1,M_{X_1})_\Prism)\xrightarrow{\cong} \Vect^{\an,\varphi}((X_1,M_{X_1})_\Prism^{\mathrm{str}})\hookrightarrow
\Vect_\Q^\varphi((X_1,M_{X_1})_\CRIS),
\]
which makes the diagram \eqref{diagram:crystalline realization} $2$-commutative.
\end{prop}

\subsection{Association of \'etale and crystalline realizations}\label{subsec: etale-cris-association}
Keep the notation as before and assume that $(X,M_X)$ is semistable or the log formal spectrum of a CDVR. To an analytic prismatic $F$-crystal $(\calE_\Prism,\varphi_{\calE_\Prism})$ on $(X,M_X)$, we have attached
\begin{itemize}
    \item a $\Z_p$-\'etale local system $T_X(\calE_\Prism,\varphi_{\calE_{\cris},\Q})$ on $X_\eta$ (Definition~\ref{defn:etale realization}), and
    \item an $F$-isocrystal $D_\cris(\calE_\Prism,\varphi_{\calE_\Prism})=(\calE_{\cris,\Q},\varphi_{\calE_{\cris},\Q})$ on $(X_1,M_{X_1})_{\CRIS}$ (Corollary~\ref{cor:crystalline-realization}).
\end{itemize}
We will show that these \'etale and crystalline realizations are associated over the crystalline period ring. This notion of association for $\Q_p$-local systems on $X_\eta$ originates from the work of Faltings \cite[V.f)]{Faltings-CryscohandGalrep}. Here we use the pro-\'etale site $X_{\eta,\proet}$ of the adic generic fiber $X_\eta$ in \cite{scholze-p-adic-hodge} to define the association. 

Let $\Perfd/X_{\eta,\proet}$ denote the full subcategory of $X_{\eta,\proet}$ consisting of \emph{affinoid} perfectoids \emph{that contain all $p$-power roots of unity}. We equip $\Perfd/X_{\eta,\proet}$ with the induced topology. Then the inclusion induces an equivalence of topoi $\operatorname{Sh}(X_{\eta,\proet})\xrightarrow{\cong}\operatorname{Sh}(\Perfd/X_{\eta,\proet})$.
Let $\bA_\cris$, $\bB_\cris^+$, and $\bB_\cris$ denote the crystalline period sheaves on $X_{\eta,\proet}$ introduced in \cite[\S2A]{Tan-Tong}. See below for a description of $\bB_\cris$ over $\Perfd/X_{\eta,\proet}$.

\begin{construction}
For $U\in \Perfd/X_{\eta,\proet}$, let $\widehat{U}=\Spa(A,A^+)$ denote the associated affinoid perfectoid. It comes with a morphism $\Spf A^+\rightarrow X$ of $p$-adic formal schemes, so equip $\Spf A^+$ with the pullback log structure $M_{\Spf A^+}$ from $M_X$.
Since $A^+$ is (integral) perfectoid, the surjection $\theta\colon \A_{\mathrm{inf}}(A^+)\coloneqq W((A^+)^\flat)\ra A^+$ gives a perfect prism over $X$, which in turn yields a log prism 
\[
(\Spf \Ainf(A^+),\Ker\theta, M_{\Spf \Ainf(A^+)}) \in (X,M_X)_\Prism 
\]
by Remark~\ref{rem:log prismatic site when admitting finite free chart}(2).

Let $\A_{\cris}(A^+)$ denote the $p$-adically completed (log) PD-envelope of the exact surjection $\theta$. 
Set $\B^+_{\cris}(A^+)\coloneqq \A_{\cris}(A^+)[p^{-1}]$ and $\B_{\cris}(A^+)\coloneqq \B_{\cris}^+(A^+)[t^{-1}]$ with $t \coloneqq \operatorname{log}[\epsilon]\in \A_{\cris}(A^+)$. We obtain an ind-object $(\Spec (A^+/p), \Spf(\A_{\cris}(A^+)), M_{\Spf(\A_{\cris}(A^+))})$ in $(X_1,M_{X_1})_{\CRIS}^{\mathrm{aff}}$. We also have a map $\A_{\cris}(A^+)\ra \bA_{\cris}(U)$, which induces an isomorphism $\B_{\cris}(A^+) \xrightarrow{\cong} \bB_{\cris}(U)$ by \cite[Cor.~2.8(2)]{Tan-Tong}. 
\end{construction}

\begin{construction}
For a finite locally free isocrystal $\calE_\Q$ on $(X_1,M_{X_1})_{\CRIS}$, we define a sheaf $\calE_\Q(\bB_\cris^+)$ of $\bB^+_\cris$-modules on $X_{\eta,\proet}$ by
\[
\calE_\Q(\bB_\cris^+)(U)\coloneqq \calE_\Q(\A_{\cris}(A^+))\otimes_{\B_{\cris}^+(A^+)}\bB_\cris^+(U)
\]
where $U\in \Perfd/X_{\eta,\proet}$ is associated to the affinoid perfectoid $\Spa (A,A^+)$ and $\calE_\Q(\A_{\cris}(A^+))$ denotes the evaluation of $\calE_\Q$ on the ind-object $(\Spec (A^+/p), \Spf(\A_{\cris}(A^+)), M_{\Spf(\A_{\cris}(A^+))})$ (see Construction~\ref{construction: eval of isocrystal at log PD-formal scheme}), which is a finite projective $\B_{\cris}^+(A^+)$-module. Since $\bB_\cris^+$ is a sheaf, it follows from the isocrystal property and the finite projectivity that $\calE_\Q(\bB_\cris^+)$ satisfies the sheaf condition on $\Perfd/X_{\eta,\proet}$.
Finally, we define
\[
\calE_\Q(\bB_\cris)\coloneqq \calE_\Q(\bB_\cris^+)\otimes_{\bB_\cris^+}\bB_\cris,
\]
which is a sheaf of $\bB_\cris$-modules on $X_{\eta,\proet}$.
If $\calE_\Q$ is an $F$-isocrystal, then $\calE_\Q(\bB_{\cris})$ is equipped with the induced Frobenius $F$.
\end{construction}

\begin{defn}[{cf.~\cite[Def.~2.31]{GuoReinecke-Ccris}}] \label{defn:association}
Let $\bL$ be a $\mathbf{Z}_p$-local system on $X_{\eta,\proet}$ and let $(\calE_\Q,\varphi_{\calE_\Q})$ be an $F$-isocrystal on $(X_1,M_{X_1})_{\CRIS}$. We say that $\bL$ and $(\calE_\Q,\varphi_{\calE_\Q})$ are \emph{associated} if there exists an isomorphism 
\[
\alpha_{\cris,\bL,\calE_\Q}\colon \calE_\Q(\bB_{\cris})\xrightarrow{\cong}\bL\otimes_{\Z_p}\bB_\cris
\]
as sheaves of $\bB_\cris$-modules on $X_{\eta, \proet}$ such that $\alpha_{\cris,\bL,\calE_\Q}$ is compatible with Frobenii.
\end{defn}

\begin{rem}\label{rem:association over general base} \hfill
\begin{enumerate}
\item By Remark~\ref{rem:F-isocrystalmoduloprestricttomodulopiissurjective}, the category of $F$-isocrystals on $(X_1,M_{X_1})_{\CRIS}$ is equivalent to the category of $F$-isocrystal on $(X_k,M_{X_k})_{\CRIS}$, where $X_k$ is the reduced special fiber of $X$.  
In particular, Definition~\ref{defn:association} is equivalent to that $\mathbb{L}$ is associated with an $F$-isocrystal over $(X_k,M_{X_k})$.
\item Unlike the work of Faltings \cite{Faltings-CryscohandGalrep}, we do not require $\calE_\Q$ to be a \emph{filtered} $F$-isocrystal. See Proposition~\ref{prop:association-filtration-compatibility} for the discussion on filtrations.
\end{enumerate}
\end{rem}

\begin{rem}\label{rem:GaloisactiononcrysllinerealizationoverAcris}
Suppose that $(X,M_X)=(\Spf(R),\N^d)^a$ is small and follow Notation~\ref{notation: Ainf and GR}.
\begin{enumerate}
    \item 
Let $T$ denote the continuous $\Z_p$-representation of $\calG_R$ corresponding to $\bL$. Then Definition~\ref{defn:association} of $T$ and $\calE_\Q$ being associated is equivalent to the following: there is an isomorphism
\[
\alpha_{\cris}\colon \calE_\Q(\bB_{\cris})(\Spa(\overline{R}[p^{-1}], \overline{R})) \xrightarrow{\cong} T\otimes_{\Z_p}\B_\cris(\overline{R})
\]
compatible with $\varphi$ and $\calG_R$-actions.
\item When the $F$-isocrystal $\calE_\Q$ arises as 
 $\calE_\Q=\calE_{\cris,\Q}\coloneqq D_\cris(\calE_\Prism)$ from an analytic prismatic $F$-crystal $\calE_\Prism$ on $(X,M_X)$ via the crystalline realization functor, one can describe the $\calG_R$-action on $\calE_\Q(\bB_{\cris})(\Spa(\overline{R}[p^{-1}], \overline{R}))$ as follows:
for any map $\fkS \to \Ainf(\overline{R})$ defined as in Notation~\ref{notation: Ainf and GR}, it induces a map $(S_R,(p),\varphi^\ast\N^d)^a \to (\A_{\cris}(\overline{R}),(p),\varphi^\ast\N^d)^a$. The isomorphism \eqref{eq:crystalline realization} for $i=0$ in Corollary~\ref{cor:crystalline-realization} extends to a $\varphi$-equivariant isomorphism
\[
\calE_\Prism((\A_{\cris}(\overline{R}),(p),\varphi^\ast\N^d)^a)\simeq \calE_{\cris,\Q}(\A_{\cris}(\overline{R})))
\]
of $\B_{\cris}^+(\overline{R})$-modules.
Moreover, using the universal property of $S_R^{(1)}$ and the cocycle condition of the isomorphism \eqref{eq:crystalline realization} for $i=1$, one can check that this isomorphism is also $\calG_R$-equivariant.
The base change of the isomorphism along $\B_{\cris}^+(\overline{R})\rightarrow\B_{\cris}(\overline{R})$ gives a $\calG_R$-equivariant description of $\calE_\Q(\bB_{\cris})(\Spa(\overline{R}[p^{-1}], \overline{R}))$.
For later use, we note that there is a natural isomorphism 
\[
\calE_\Prism((\A_{\cris}(\overline{R}),(p),\varphi^\ast\N^d)^a)\simeq \calE_\Prism((\A_{\inf}(\overline{R}),(E),\N^d)^a)[p^{-1}] \otimes_{\A_{\inf}(\overline{R})[p^{-1}],\varphi} \B_{\cris}^+(\overline{R})
\]
by applying the crystal property of $\calE_\Prism$ to $(\A_{\inf}(\overline{R}),(E),\N^d)^a \xrightarrow{\varphi}(\A_{\cris}(\overline{R}),(p),\varphi^\ast\N^d)^a$.
\end{enumerate}
\end{rem}

\begin{prop} \label{prop:etaleandcrysassociated}
For an analytic prismatic $F$-crystal $\calE_{\Prism}$ over $(X,M_X)$, its \'etale realization $T(\calE_{\Prism})$ (Definition~\ref{defn:etale realization}) and its crystalline realization $\calE_{\cris, \mathbf{Q}}\coloneqq D_\cris(\calE_\Prism)$  (Corollary~\ref{cor:crystalline-realization}) are associated in the sense of Definition~\ref{defn:association}. 
\end{prop}

\begin{proof}
Let $\{\Spf (R_{\lambda})\}_{\lambda\in\Lambda}$ be a $p$-complete \'etale cover of $X$ consisting of affine formal schemes $\Spf (R_{\lambda})$ such that each $(\Spf (R_{\lambda}), M_X|_{\Spf (R_{\lambda})})$ is a small affine log formal scheme as in \S~\ref{sec:Breuil--Kisin log prism}. We may further assume that each $R_{\lambda}$ admits a framing as in \S~\ref{sec:prelim-facts-rings}. Let $U_{\lambda}$ denote the object of $\Perfd/X_{\eta, \proet}$ corresponding to $\Spa((\overline{R_{\lambda}})^{\wedge}_p)[p^{-1}], (\overline{R_{\lambda}})^{\wedge}_p))$. By Remark~\ref{rem:GaloisactiononcrysllinerealizationoverAcris}(2), there is a $\calG_{R_\lambda}$-equivariant isomorphism
\[
\mathcal{E}_{\Prism}((\mathbf{A}_{\mathrm{cris}}((\overline{R_{\lambda}})^{\wedge}_p), (p),\varphi^\ast\N^d)^a) \cong \mathcal{E}_{\mathrm{cris}, \mathbf{Q}}(\mathbf{A}_{\mathrm{cris}}((\overline{R_{\lambda}})^{\wedge}_p)
\]
of $\mathbf{B}_{\mathrm{cris}}^+((\overline{R_{\lambda}})^{\wedge}_p)$-modules compatible with $\varphi$ for each $\lambda$. On the other hand, taking the base change of the isomorphism in Lemma~\ref{lem:isomafterinvertingmu} along the map $(\mathbf{A}_{\mathrm{inf}}(\overline{R}_{\lambda}), \Ker \theta) \stackrel{\varphi}{\rightarrow} (\mathbf{A}_{\mathrm{cris}}((\overline{R_{\lambda}})^{\wedge}_p), (p))$ gives a $\calG_{R_\lambda}$ and Frobenius equivariant isomorphism
\[
\mathcal{E}_{\Prism}((\mathbf{A}_{\mathrm{cris}}((\overline{R_{\lambda}})^{\wedge}_p), (p),\varphi^\ast\N^d)^a)\otimes_{\mathbf{A}_{\mathrm{cris}}((\overline{R_{\lambda}})^{\wedge}_p)} \mathbf{B}_{\mathrm{cris}}((\overline{R_{\lambda}})^{\wedge}_p) \cong T(\mathcal{E}_{\Prism})\otimes_{\mathbf{Z}_p} \mathbf{B}_{\mathrm{cris}}((\overline{R_{\lambda}})^{\wedge}_p).
\]
Thus, we obtain a canonical $\varphi$-equivariant isomorphism
\[
\mathcal{E}_{\mathrm{cris}, \mathbf{Q}}(\mathbb{B}_{\mathrm{cris}})|_{U_{\lambda}} \cong T(\mathcal{E}_{\Prism})\otimes_{\mathbf{Z}_p}\mathbb{B}_{\mathrm{cris}} |_{U_{\lambda}}.
\]

Note that $\{U_{\lambda}\}_{\lambda\in \Lambda}$ defines a cover of $X_{\eta}$ in $X_{\eta, \proet}$. It follows directly from the construction of the above isomorphisms for $U_{\lambda}$'s that they descend to an isomorphism
\[
\mathcal{E}_{\mathrm{cris}, \mathbf{Q}}(\mathbb{B}_{\mathrm{cris}}) \cong T(\mathcal{E}_{\Prism})
\otimes_{\mathbf{Z}_p} \mathbb{B}_{\mathrm{cris}}
\]
of sheaves of $\mathbb{B}_{\mathrm{cris}}$-modules on $X_{\eta, \proet}$ compatible with $\varphi$.
\end{proof}

In the rest of this section, we explain Remark~\ref{rem:association over general base}(2) on why we do not impose any condition on filtrations in Definition~\ref{defn:association}. 
In the CDVR case, we will show in the next section that Definition~\ref{defn:association} is equivalent to the classical semistability (Theorem~\ref{thm:CDVR-semistable-notions-equivalent}). Hence from now on we assume that $(X, M_X)$ is semistable. Recall that de Rham period sheaves $\bB_{\dR}^+$,  $\bB_{\dR}$, $\mathcal{O}\bB_{\dR}^+$ and $\mathcal{O}\bB_{\dR}$ on $X_{\eta, \proet}$ are defined in  \cite[Def.~6.1, 6.8]{scholze-p-adic-hodge} and come with filtrations. We also set $\bB_{\dR}^{+,n}=\bB_{\dR}^+/(\Ker\theta)^{n+1}$ so that $\bB_{\dR}^+=\varprojlim_n\bB_{\dR}^{+,n}$.

\begin{construction} \label{const: filtration-on-BdR}
Let $((\calE_\Q,\varphi_{\calE_\Q}),(\mathbf{E},\nabla_{\mathbf{E}},\Fil^\bullet\mathbf{E}))$ be a filtered $F$-isocrystal on $(X,M_X)$ as in Definition~\ref{def:filtered F-isocrystals}: 
$(\calE_\Q,\varphi_{\calE_\Q})$ is an $F$-isocrystal on $(X_1,M_{X_1})_{\CRIS}$; $(\mathbf{E},\nabla_{\mathbf{E}})$ is the associated vector bundle with integrable connection  on $X_{\eta}$ by Remark~\ref{rem:integrable connection associated to F-isocrystal}; $\Fil^\bullet\mathbf{E}$ is a $\Z$-indexed decreasing separated and exhaustive filtration of $\mathcal{O}_{X_{\eta}}$-submodules of $\mathbf{E}$ satisfying the Griffiths transversality such that $\mathrm{Fil}^i \mathbf{E} / \mathrm{Fil}^{i+1} \mathbf{E}$ is locally free over $\mathcal{O}_{X_{\eta}}$ for each $i$.

For $\bB\in\{\bB^+_\dR, \bB^{+,n}_\dR, \bB_\dR\}$, consider the sheaf 
\[
\mathcal{E}_{\mathbf{Q}}(\mathbb{B}) \coloneqq \calE_\Q(\bB_{\cris}^+)\otimes_{\bB_{\cris}^+} \mathbb{B}
\]
on $X_{\eta,\proet}$. We are going to define a filtration on $\mathcal{E}_{\mathbf{Q}}(\mathbb{B}_{\dR})$ as follows. 

First consider the small affine case $(X, M_X) = (\Spf R, \mathbf{N}^d)^a$ with framing 
\[
\square\colon R^0 = \mathcal{O}_K \langle T_1, \ldots, T_m, T_{m+1}^{\pm 1}, \ldots, T_d^{\pm 1}\rangle / (T_1\cdots T_m - \pi) \rightarrow R
\]
as in \S\ref{sec:Breuil--Kisin log prism}. 
Let $U\in \Perfd/X_{\eta,\proet}$ with the associated affinoid perfectoid $\widehat{U}=\Spa(A,A^+)$, and suppose that $A^{+}$ admits compatible systems $T_i^{\flat}$ of $p$-power roots of $T_i$ for each $i = 1, \ldots, d$. Note that $\pi^{\flat} \coloneqq T_1^\flat \cdots T_m^{\flat}$ gives a compatible system of $p$-power roots of $\pi$ in $A^{+}$. 
We write $\B_\dR^{+}(A^+)$ for $\B_\dR^+(A^+)$ defined in \cite[p.49]{scholze-p-adic-hodge}, which is isomorphic to $\bB_\dR^{+}(U)$ by \cite[Thm.~6.5(ii)]{scholze-p-adic-hodge}. We write $\B_\dR^{+,n}(A^+)=\B_\dR^+(A^+)/(\Ker\theta)^{n+1}$; note $\B_\dR^+(A^+)\xrightarrow{\cong}\varprojlim_n\B_\dR^{+,n}(A^+)$ and the elements $[\pi^\flat]$ and $[T_i^\flat]$ ($1\leq i\leq m$) are invertible in these rings as $\pi$ is invertible and $[\pi^\flat]/\pi-1\in \Ker\theta$. In the discussion below, we will also use the notation in Example~\ref{eg:BreuilSR} and Remark~\ref{rem:integrable connection associated to F-isocrystal}, including $\widetilde{R}$, $S=S_\square$, $S^{(1)}$, and $R^{(1),n}$.

For a choice of $1\leq j\leq m$, consider the $K$-algebra maps 
\[
s_n\coloneqq s_{\square,j,n}\colon R[p^{-1}]\rightarrow \B_\dR^{+,n}(A^+) \quad\text{and}\quad
s\coloneqq s_{\square,j}\coloneqq \varprojlim_n s_{\square,j,n}\colon R[p^{-1}]\rightarrow \B_\dR^{+}(A^+)
\]
given by $T_i \mapsto [T_i^{\flat}]$ for $i\neq j$ and $T_j\mapsto \pi([T_1^\flat]\cdots[T_{j-1}^\flat][T_{j+1}^\flat]\cdots [T_{m}^\flat])^{-1}$. Note that these maps exist and are unique: ${s_n}|_{R^0}$ exists and takes values in some open bounded subring of $\B_\dR^{+,n}(A^+)$ (equipped with $p$-adic topology), and thus it uniquely extends to $R$ and $R[p^{-1}]$. Moreover, the composite $R[p^{-1}] \xrightarrow{s} \B_{\dR}(A^+) \xrightarrow{\theta} A$ is the structure map $R[p^{-1}] \rightarrow A$.

Write $\mathbf{E}_R$ for the finite projective $R[p^{-1}]$-module corresponding to $\mathbf{E}$.
We will construct isomorphisms of $\B_{\dR}^{+}(A^+)$-modules
\[
\mathcal{E}_{\mathbf{Q}}(\mathbb{B}_{\dR}^{+,n})(U) \underset{\cong}{\xrightarrow{\alpha_{s,n}}}\mathbf{E}_R\otimes_{R[p^{-1}],s} \B_{\dR}^{+,n}(A^+)
\quad\text{and}\quad
\mathcal{E}_{\mathbf{Q}}(\mathbb{B}_{\dR}^{+})(U) \underset{\cong}{\xrightarrow{\alpha_{s}}} \mathbf{E}_R\otimes_{R[p^{-1}],s} \B_{\dR}^{+}(A^+).
\]
It is enough to construct $\alpha_{s,n}$, and for this, consider the following diagram
\[
\xymatrix{
S\ar[r]^-{p_{S,1}}\ar[d] & S^{(1)}\ar@{..>}[d]^-{f} &S \ar[l]_-{p_{S,2}} \ar[d]\\
\A_{\cris}(A^+)\ar[r]&\B_{\dR}^{+,n}(A^+)&R[p^{-1}].\ar[l]_-s
}
\]
Here the first and third vertical maps are induced by $\square$, and there is a unique map $f$ that makes the diagram commutative: since $\pi$, $[\pi^\flat]$, $s(T_j)$, and $[T_j^\flat]$ are all invertible in $\B_{\dR}^{+,n}(A^+)$, the two maps $\widetilde{R}\rightarrow S\rightrightarrows\B_{\dR}^{+,n}(A^+)$ yields uniquely to a map $\widetilde{R}^{(1)}\rightarrow\B_{\dR}^{+,n}(A^+)$, which further extends to a map $S^{(1),\mathrm{nc}}\rightarrow\B_{\dR}^{+,n}(A^+)$ with values in some open and bounded subring as the image of $\widetilde{I}^{(1)}$ is nilpotent in $\B_{\dR}^{+,n}(A^+)$. This further extends to $f$. Now the isomorphism $\alpha_{s,n}$ is defined to be the base change of the HPD-stratification on $\calE_{\Q}(S)$ along $f$ by Remark~\ref{rem:integrable connection associated to F-isocrystal}.

Equip $\mathcal{E}_{\mathbf{Q}}(\mathbb{B}_{\dR}^+)(U)$ with the tensor product filtration $\sum_{a+b=\bullet}\mathrm{Fil}^{a}(\mathbf{E}_R)\otimes_{R[p^{-1}], s} \mathrm{Fil}^{b}\B_{\dR}^+(A^+)$ via $\alpha_s$. We claim that the filtration is independent of the choice of $\square$ and $j$: to see this, take another $K$-algebra homomorphism $s'\coloneqq s_{\square',j'}\colon R[p^{-1}]\rightarrow\B_{\dR}^+(A^+)$. Since the log structure on $\Spf R$ is given by $(R[p^{-1}]^\times)\cap R$, it is straightforward to see that $s$ and $s'$ yield maps $g^n\colon R^{(1),n}[p^{-1}]\rightarrow \B_{\dR}^{+,n}(A^+)$ and $g\coloneqq \varprojlim_ng^n\colon R^{(1)}[p^{-1}]^\wedge\coloneqq\varprojlim_n R^{(1),n}[p^{-1}]\rightarrow \B_{\dR}^{+}(A^+)$ making the following diagrams commutative:
\[
\xymatrix{
R[p^{-1}]\ar[r]^-{p_{1}^n}\ar[dr]_-{s_n} & R^{(1),n}[p^{-1}]\ar[d]^-{g^n} &R[p^{-1}] \ar[l]_-{p_{2}^n} \ar[dl]^-{s_n'}\\
&\B_{\dR}^{+,n}(A^+)&
}
\quad\text{and}\quad
\xymatrix{
R[p^{-1}]\ar[r]^-{p_{1}}\ar[dr]_-s & R^{(1)}[p^{-1}]^\wedge\ar[d]^-{g} &R[p^{-1}] \ar[l]_-{p_{2}} \ar[dl]^-{s'}\\
&\B_{\dR}^{+}(A^+).&
}
\]
The pullback along $g$ of the infinitesimal log stratification $\eta_{\mathrm{inf}}\colon p_1^\ast\mathbf{E}_R\xrightarrow{\cong}p_2^\ast\mathbf{E}_R$ on $\mathbf{E}_R$ corresponding to $(\mathbf{E},\nabla_\mathbf{E})$ defines the canonical $\B_{\dR}^+(A^+)$-linear isomorphism
\[
\eta_{s,s'}\colon
\mathbf{E}_R\otimes_{R[p^{-1}],s} \B_{\dR}^+(A^+)\xrightarrow{\cong}\mathbf{E}_R\otimes_{R[p^{-1}],s'} \B_{\dR}^+(A^+).
\]
Furthermore, we can check $\eta_{s,s''}=\eta_{s',s''}\circ \eta_{s,s'}$ for another $s''$ by considering the cocycle condition of $\eta_{\mathrm{inf}}$ over $\varprojlim_n R^{(2),n}[p^{-1}]$.
We also note $\alpha_{s'}=\eta_{s,s'}\circ\alpha_s$: for this, consider the following \emph{non-commutative} diagram
\[
\xymatrix{
S=S_{\square}\ar[r]\ar[d] & \A_{\cris}(A^+)\ar[d]^-{\mathrm{can}} &S_{\square'} \ar[l] \ar[d]\\
R[p^{-1}]\ar[r]^-s&\B_{\dR}^{+}(A^+)&R[p^{-1}].\ar[l]_-{s'}
}
\]
Using the diagram $S_{\square}\xrightarrow{\pr_\square} S_{\square,\square'}\xleftarrow{\pr_{\square'}} S_{\square'}$ as in the proof of Lemma~\ref{lem:independence of framing for associated crystalline crystal} and repeating a similar argument as above, Example~\ref{eg:BreuilSR}, and Remark~\ref{rem:integrable connection associated to F-isocrystal}, we can verify $\alpha_{s'}=\eta_{s,s'}\circ\alpha_s$.

Now, to see the independence of the filtration on $\mathcal{E}_{\mathbf{Q}}(\mathbb{B}_{\dR}^+)(U)$, it is enough to show that $\eta_{s',s}=\eta_{s',s}^{-1}$ is a map of filtered modules: applying the same argument to $\eta_{s,s'}$ will imply that they are filtered isomorphisms.
To see that $\eta_{s',s}$ preserves the filtrations, 
we compute $\eta_{s',s}(e\otimes 1)$ for a section $e\in\mathbf{E}_R$. For this, fix a local coordinate $(x_2,\dots,x_d)$ of $R$ (e.g. the one coming from $\square$), and write $\nabla_i\coloneqq \nabla_{\mathbf{E},\partial/\partial x_{i}}\in \End(\mathbf{E}_R)$. With the convention \eqref{eq:connection and stratification}, a standard computation gives
\[
\eta_{s',s}(e\otimes 1)
=\sum_{(n_2,\ldots,n_d)\in\N^{d-1}}\frac{(\nabla_2^{n_2}\circ\cdots\circ \nabla_d^{n_d})(e)}{n_2!\cdots n_d!}\otimes (s'(x_2)-s(x_2))^{n_2}\cdots (s'(x_d)-s(x_d))^{n_d}.
\]
Since $s'(x_i)-s(x_i)\in \Ker\theta=\Fil^1\B_{\dR}^+(A^+)$, we deduce from the Griffiths transversality of $\nabla$ that $\eta_{s',s}$ is a filtered map. 

For a general semistable $p$-adic formal scheme $X$, we can glue local constructions as above to obtain a filtration on $\mathcal{E}_{\mathbf{Q}}(\mathbb{B}^+_{\dR})$. Finally, we equip $\calE_\Q(\mathbb{B}_\dR)=\mathcal{E}_{\mathbf{Q}}(\mathbb{B}^+_{\dR})\otimes_{\mathbb{B}^+_{\dR}}\mathbb{B}_{\dR}$ with the tensor product filtration.
\end{construction}

In what follows, $\calO_{X_\eta}$ and $\mathbf{E}$ are also regarded as sheaves on $X_{\eta,\proet}$ via the pullback along the projection $X_{\eta,\proet}\rightarrow X_{\mathrm{an}}$.

\begin{prop} \label{prop:connection-filtration-compatibility}
Let $\mathcal{E}_{\mathbf{Q}}$ be an $F$-isocrystal on $(X_1,M_{X_1})_{\CRIS}$, and write $(\mathbf{E}, \nabla_{\mathbf{E}})$ for the induced vector bundle with integrable connection on $X_{\eta}$ given by Remark~\ref{rem:integrable connection associated to F-isocrystal}. 
Let $\star\in \{+,\emptyset\}$.
\begin{enumerate}
    \item There is a natural isomorphism 
\[
(\mathbf{E}\otimes_{\mathcal{O}_{X_{\eta}}} \mathcal{O}\bB^{\star}_{\dR}, \nabla_{\mathbf{E}}\otimes\mathrm{id}+\mathrm{id}\otimes\nabla_{\mathcal{O}\bB^{\star}_{\dR}}) \cong (\mathcal{E}_{\mathbf{Q}}(\bB^{\star}_{\dR})\otimes_{\bB^{\star}_{\dR}} \mathcal{O}\bB^{\star}_{\dR}, \mathrm{id}\otimes \nabla_{\mathcal{O}\bB^{\star}_{\dR}})
\]
of $\mathcal{O}\bB^{\star}_{\dR}$-vector bundles with connections.
    \item Suppose that $(\mathbf{E},\nabla_{\mathbf{E}})$ is equipped with a Griffiths-transversal filtration and let $\mathrm{Fil}^{\bullet}\mathcal{E}_{\mathbf{Q}}(\bB^{\star}_{\dR})$ be the induced filtration on $\mathcal{E}_{\mathbf{Q}}(\bB^{\star}_{\dR})$ as in Construction~\ref{const: filtration-on-BdR}. Then the tensor product filtrations $\mathrm{Fil}^{\bullet} \mathbf{E}\otimes_{\mathcal{O}_{X_{\eta}}} \mathrm{Fil}^{\bullet}\mathcal{O}\bB^{\star}_{\dR}$ and $\mathrm{Fil}^{\bullet}\mathcal{E}_{\mathbf{Q}}(\bB^{\star}_{\dR})\otimes_{\bB^{\star}_{\dR}} \mathrm{Fil}^{\bullet}\mathcal{O}\bB^{\star}_{\dR}$ agree under the isomorphism $\mathbf{E}\otimes_{\mathcal{O}_{X_{\eta}}} \mathcal{O}\bB^{\star}_{\dR}\cong \mathcal{E}_{\mathbf{Q}}(\mathbb{B}^{\star}_{\dR})\otimes_{\bB^{\star}_{\dR}} \mathcal{O}\bB^{\star}_{\dR}$ in (1).
\end{enumerate}
\end{prop}

\begin{proof}
We follow a similar argument as in the proof of \cite[Prop.~2.36]{GuoReinecke-Ccris}. 
It suffices to treat the $\mathcal{O}\bB^{+}_{\dR}$-case.

First consider the small affine case $(X, M_X) = (\Spf R, \mathbf{N}^d)^a$ with a framing $\square\colon R^0\rightarrow R$ and use the notation in Construction~\ref{const: filtration-on-BdR}. Take any $U\in \Perfd/X_{\eta,\proet}$ with the associated affinoid perfectoid $\widehat{U}=\Spa(A,A^+)$ such that $A^+$ admits compatible systems $T_i^{\flat}$ of $p$-power roots of $T_i$ for each $i = 1, \ldots, d$. Write $U=\varprojlim_\lambda\Spa(A_\lambda,A_\lambda^+)$ so that $A^+$ is the completed direct limit of $A_\lambda^+$. Set $\OB_{\dR}^{+}(A^+)\coloneqq \calOB_{\dR}^+(U)$ and consider the ring
\[
\OB_{\dR}^{(1)+}(A^+)\coloneqq \varinjlim_\lambda\varprojlim_n (A_\lambda^+\widehat{\otimes}_{\calO_K}A_\lambda^+\widehat{\otimes}_{W(k)}\A_{\mathrm{inf}}(A^+))[p^{-1}]/(\Ker\theta_\lambda)^{n+1},
\]
where $\widehat{\otimes}$ denotes the $p$-completed tensor product and $\theta_\lambda$ is the natural surjection to $A$.
By an argument similar to \cite[(3)]{Scholze-p-adicHodgeerrata}, we see that $\varprojlim_n (A_\lambda^+\widehat{\otimes}_{\calO_K}A_\lambda^+\widehat{\otimes}_{W(k)}\A_{\mathrm{inf}}(A^+))[p^{-1}]/(\Ker\theta_\lambda)^{n+1} \rightarrow \OB_{\dR}^{(1)+}(A^+)$ is an isomorphism for a sufficiently large $\lambda$, and we obtain a commutative diagram
\[
\xymatrix{
R^{(1)}[p^{-1}]^\wedge\ar[r]^-{h^{(1)}} & \OB_{\dR}^{(1)+}(A^+)& \B_{\dR}^{+}(A^+)[\![Z_{2,1},\ldots,Z_{d,1},Z_{2,2},\ldots,Z_{d,2}]\!]\ar[l]_-\cong\\
R[p^{-1}]\ar[r]^-h \ar@<.5ex>[u]^-{p_1}\ar@<-.5ex>[u]_-{p_2} \ar@<.5ex>[ur]^-{p_{h,1}}\ar@<-.5ex>[ur]_-{p_{h,2}} & \OB_{\dR}^{+}(A^+)\ar@<.5ex>[u]^-{q_1}\ar@<-.5ex>[u]_-{q_2}& \B_{\dR}^{+}(A^+)[\![Z_2,\ldots,Z_d]\!]\ar@<.5ex>[u]^-{q'_1}\ar@<-.5ex>[u]_-{q'_2}\ar[l]_-\cong\\
 & \B_{\dR}^{+}(A^+),\ar[u]^-u\ar[ur] & 
}
\]
where $h$ and $h^{(1)}$ are the canonical structure maps (note $(R\widehat{\otimes}_{\calO_K}R)[p^{-1}]\xrightarrow{\cong}R^{(1)}[p^{-1}]$), the horizontal isomorphisms are $\B_{\dR}^{+}(A^+)$-algebra isomorphisms given by $Z_i\mapsto T_i\otimes 1-1\otimes[T_i^\flat]$, $Z_{i,1}\mapsto T_i\otimes 1\otimes 1-1\otimes 1\otimes[T_i^\flat]$, $Z_{i,2}\mapsto 1\otimes T_i\otimes 1-1\otimes 1\otimes[T_i^\flat]$, and the vertical and diagonal maps are the obvious ones. We also set $u^{(1)}\coloneqq q_1\circ u=q_2\circ u$.

Write $(\mathbf{E}_R,\eta_{\mathrm{inf}}\colon p_1^\ast\mathbf{E}_R\xrightarrow{\cong}p_2^\ast\mathbf{E}_R)$ for the finite projective $R[p^{-1}]$-module together with infinitesimal (log) stratification corresponding to $(\mathbf{E},\nabla_{\mathbf{E}})$.
By construction, $(\mathbf{E}\otimes_{\mathcal{O}_{X_{\eta}}} \mathcal{O}\bB^+_{\dR})(U)=\mathbf{E}_R\otimes_{R[p^{-1}],h}\OB_{\dR}^{+}(A^+)$, and $\nabla_{\mathbf{E}}\otimes\mathrm{id}+\mathrm{id}\otimes\nabla_{\mathcal{O}\bB^{+}_{\dR}}$ corresponds to an infinitesimal stratification 
\[
h^{(1),\ast}\eta_{\mathrm{inf}}\colon q_1^\ast(\mathbf{E}\otimes_{\mathcal{O}_{X_{\eta}}} \mathcal{O}\bB^+_{\dR})(U)\cong p_{h,1}^\ast\mathbf{E}_R\xrightarrow{\cong}p_{h,2}^\ast\mathbf{E}_R \cong q_2^\ast(\mathbf{E}\otimes_{\mathcal{O}_{X_{\eta}}} \mathcal{O}\bB^+_{\dR})(U).
\]
To see this, observe that $I_{\dR}/I_{\dR}^2$, where $I_{\dR}\coloneqq \Ker(\OB_{\dR}^{(1)+}( A^+)\rightarrow\OB_{\dR}^{+}(A^+))$, is canonically identified with $\Omega_{X_\eta/K}(X_\eta)\otimes_{R[p^{-1}]}\OB^+_{\dR}(A^+)$ and recall the formula \eqref{eq:connection and stratification}.

Consider the map $s=s_{\square,j}\colon R[p^{-1}]\rightarrow \B_{\dR}^{+}(A^+)$ in Construction~\ref{const: filtration-on-BdR}. The argument therein gives an isomorphism $\alpha_{s}\colon\mathcal{E}_{\mathbf{Q}}(\mathbb{B}_{\dR}^{+})(U) \xrightarrow{\cong} \mathbf{E}_R\otimes_{R[p^{-1}],s} \B_{\dR}^{+}(A^+)$, which in turn yields the identification 
$(\mathcal{E}_{\mathbf{Q}}(\bB^+_{\dR})\otimes_{\bB^+_{\dR}} \calOB^+_{\dR})(U)\cong\mathbf{E}_R\otimes_{R[p^{-1}],s}\OB_{\dR}^{+}(A^+)$.
With this identification, the connection $\mathrm{id}_{\mathcal{E}_{\mathbf{Q}}(\bB^+_{\dR})}\otimes \nabla_{\calOB^+_{\dR}}$ corresponds to the infinitesimal stratification $(u^{(1)}\circ g)^\ast \eta_{\mathrm{inf}}$, where $g\colon R^{(1)}[p^{-1}]^\wedge\rightarrow \B_{\dR}^{+}(A^+)$ is the map characterized by $g\circ p_1=s=g\circ p_2$. 

Now, to show (1) and (2) in the small affine case, consider a commutative diagram of the form
\[
\xymatrix{
R[p^{-1}]\ar[r]^-{p_{1}}\ar[dr]_-h & R^{(1)}[p^{-1}]^\wedge\ar[d]^-{H} &R[p^{-1}] \ar[l]_-{p_{2}} \ar[dl]^-{u\circ s}\\
&\OB_{\dR}^{+}(A^+).&
}
\]
The existence of such a diagram is proved as in Construction~\ref{const: filtration-on-BdR}, and this gives an $\OB_{\dR}^{+}(A^+)$-linear isomorphism
\[
H^\ast\eta_{\mathrm{inf}}\colon (\mathbf{E}\otimes_{\mathcal{O}_{X_{\eta}}} \mathcal{O}\bB^+_{\dR})(U)\xrightarrow{\cong}(\mathcal{E}_{\mathbf{Q}}(\bB^+_{\dR})\otimes_{\bB^+_{\dR}} \calOB^+_{\dR})(U).
\]
Moreover, the cocycle condition for $\eta_{\mathrm{inf}}$ implies $(u^{(1)}\circ s)^\ast \eta_{\mathrm{inf}}\circ q_1^\ast H^\ast\eta_{\mathrm{inf}}=q_2^\ast H^\ast\eta_{\mathrm{inf}}\circ h^{(1),\ast}\eta_{\mathrm{inf}}$.  
This is equivalent to the desired isomorphism of $\mathcal{O}\bB_{\dR}$-vector bundles with connections in (1) in this case. For (2), we need to show that $H^\ast\eta_{\mathrm{inf}}$ and $(H^\ast\eta_{\mathrm{inf}})^{-1}$ are maps of filtered modules, which follows from the Griffiths transversality as in the second to last paragraph of Construction~\ref{const: filtration-on-BdR} by using the local coordinates $Z_2,\ldots,Z_d$ of $\OB_{\dR}^{+}(A^+)$.

To show the general case, we need to compare the local constructions for two framings $\square$ and $\square'$, which is routine and left as an exercise to the reader.
\end{proof}

Consider the projection $\nu\colon X_{\eta,\proet}\rightarrow X_{\eta,\et}$. For a $\mathbf{Z}_p$-local system $\mathbb{L}$ on $X_{\eta,\proet}$, we set
\[
(D_\dR(\mathbb{L}),\nabla_{D_\dR(\mathbb{L})},\Fil^\bullet D_\dR(\mathbb{L}))\coloneqq 
\nu_\ast(\mathbb{L}\otimes_{\Z_p} \calOB_{\dR},\mathrm{id}_{\mathbb{L}}\otimes\nabla_{\calOB_{\dR}}, 
 \mathbb{L}\otimes_{\Z_p} \Fil^\bullet\calOB_{\dR}).
\]
By \cite[Thm.~3.9]{liu-zhu-rigidity}, this is a filtered vector bundle with integrable connection on $X_\eta$ satisfying the Griffiths transversality (here and below, we identify the category of vector bundles on $X_{\eta,\et}$ and the category of vector bundles on $X_\eta$ via \cite[Prop.~9.2(i)]{scholze-p-adic-hodge}); $\mathbb{L}$ is de Rham if and only if $\rank_{\Z_p}\mathbb{L}=\rank_{\calO_{X_\eta}}D_\dR(\mathbb{L})$, in which case the counit map for $(\nu^\ast,\nu_\ast)$ yields a filtered isomorphism of $\calOB_{\dR}$-vector bundles with integrable connections
\[
(D_\dR(\mathbb{L}),\nabla_{D_\dR(\mathbb{L})},\Fil^\bullet D_\dR(\mathbb{L}))\otimes_{\calO_{X_{\eta,\et}}} \calOB_\dR\underset{\cong}{\xrightarrow{\alpha_{\dR,\bL}}}(\mathbb{L}\otimes_{\Z_p} \calOB_{\dR},\mathrm{id}_{\mathbb{L}}\otimes\nabla_{\calOB_{\dR}}, 
 \mathbb{L}\otimes_{\Z_p} \Fil^\bullet\calOB_{\dR}),
\]
where the domain is equipped with the tensor product connection and filtration.

\begin{prop} \label{prop:association-filtration-compatibility}
Let $\mathbb{L}$ be a $\mathbf{Z}_p$-local system on $X_{\eta,\proet}$, and assume that it is associated with an $F$-isocrystal $(\calE_\Q,\varphi_{\calE_\Q})$ on $(X_1,M_{X_1})_{\CRIS}$ via  $\alpha_{\cris,\bL,\calE_\Q}\colon \calE_\Q(\bB_{\cris})\xrightarrow{\cong}\bL\otimes_{\Z_p}\bB_\cris$. 
\begin{enumerate}
    \item The isomorphism $\alpha_{\cris,\bL,\calE_\Q}$ induces an isomorphism
\[
(\mathbf{E},\nabla_{\mathbf{E}})\xrightarrow{\cong}(D_\dR(\mathbb{L}),\nabla_{D_\dR(\mathbb{L})})
\]
of vector bundles with connections on $X_\eta$.
Moreover, $\mathbb{L}$ is de Rham, and the triple $(\bL\otimes_{\Z_p}\bB_{\dR}, \bL\otimes_{\Z_p}\Fil^\bullet\bB_{\dR})$ gives a filtered $F$-isocrystal structure on $(\calE_\Q,\varphi_{\calE_\Q})$ such that 
    \[
\alpha_{\cris,\bL,\calE_\Q}\otimes_{\bB_{\cris}} \bB_{\dR}\colon (\mathcal{E}_{\mathbf{Q}}(\bB_{\dR}),\Fil^\bullet \mathcal{E}_{\mathbf{Q}}(\bB_{\dR})) \cong  (\bL\otimes_{\Z_p}\bB_{\dR}, \bL\otimes_{\Z_p}\Fil^\bullet\bB_{\dR})
\]
is a filtered isomorphism on $X_{\eta,\proet}$. 
    \item If $(\mathbf{E},\nabla_{\mathbf{E}},\Fil^\bullet\mathbf{E})$ is a filtered $F$-isocrystal structure on $(\calE_\Q,\varphi_{\calE_\Q})$ such that $\alpha_{\cris,\bL,\calE_\Q}\otimes_{\bB_{\cris}} \bB_{\dR}$ is a filtered isomorphism, 
then the isomorphism in (1) becomes a filtered isomorphism
\[
(\mathbf{E},\nabla_{\mathbf{E}},\Fil^\bullet\mathbf{E})\xrightarrow{\cong}(D_\dR(\mathbb{L}),\nabla_{D_\dR(\mathbb{L})},\Fil^\bullet D_\dR(\mathbb{L})).
\]
In particular, such a filtration $\Fil^\bullet\mathbf{E}$ on $\mathbf{E}$ is unique and given by $\Fil^\bullet D_\dR(\mathbb{L})$.
\item In the situation of (2), $\alpha_{\cris,\bL,\calE_\Q}$ induces an isomorphism
\[
\bL[p^{-1}]\cong \Ker\bigl(\calE_\Q(\bB_{\cris})^{\varphi=1}\rightarrow \Fil^0(\mathbf{E}\otimes_{\calO_{X_\eta}}\calOB_\dR)^{\nabla=0}\bigl).
\]
\end{enumerate}
\end{prop}

\begin{proof}
(1) By Proposition~\ref{prop:connection-filtration-compatibility}(1), $\alpha_{\cris,\bL,\calE_\Q}\otimes_{\bB_{\cris}} \calOB_{\dR}$ gives an isomorphism
\[
\alpha_{\calOB_{\dR}}\colon (\mathbf{E} \otimes_{\calO_{X_\eta}}\calOB_{\dR},\nabla_{\mathbf{E}}\otimes\mathrm{id}+\mathrm{id}\otimes\nabla_{\mathcal{O}\bB_{\dR}})\xrightarrow{\cong}  (\bL\otimes_{\Z_p}\calOB_{\dR},\mathrm{id}\otimes\nabla_{\mathcal{O}\bB_{\dR}})
\]
of $\calOB_\dR$-vector bundles with connections.
It follows from the projection formula and $\nu_\ast\calOB_\dR=\calO_{X_{\eta,\et}}$ \cite[Cor.~6.19]{scholze-p-adic-hodge} that the pushforward $\nu_\ast$ of the above isomorphism gives an isomorphism
\[
\alpha_{\mathbf{E}}\colon (\mathbf{E},\nabla_{\mathbf{E}})\xrightarrow{\cong}(D_\dR(\mathbb{L}),\nabla_{D_\dR(\mathbb{L})})
\]
of vector bundles with connections on $X_\eta$. In particular, $\rank_{\Z_p}\mathbb{L}=\rank_{\calO_{X_\eta}}D_\dR(\mathbb{L})$, and thus $\bL$ is de Rham.
Now consider the isomorphisms of $\calOB_{\dR}$-vector bundles with integrable connections
\[
\xymatrix{
(\mathbf{E} \otimes_{\calO_{X_\eta}}\calOB_{\dR},\nabla_{\mathbf{E}}\otimes\mathrm{id}+\mathrm{id}\otimes\nabla_{\mathcal{O}\bB_{\dR}})\ar[d]_-{\alpha_{\mathbf{E}}\otimes\calOB_\dR}^-\cong\ar[dr]_-\cong^-{\alpha_{\calOB_\dR}}\ar[r]^-\cong
& (\mathcal{E}_{\mathbf{Q}}(\bB_{\dR})\otimes_{\bB_\dR}\calOB_\dR,\mathrm{id}\otimes\nabla_{\mathcal{O}\bB_{\dR}})\ar[d]^-{\alpha_{\cris,\bL,\calE_\Q}\otimes_{\bB_{\cris}} \calOB_{\dR}}_-\cong\\
(D_\dR(\mathbb{L})\otimes_{\calO_{X_\eta}}\calOB_{\dR},\nabla_{D_\dR(\mathbb{L})}\otimes\mathrm{id}+\mathrm{id}\otimes\nabla_{\mathcal{O}\bB_{\dR}})\ar[r]^-{\alpha_{\dR,\bL}}_-\cong&(\mathbb{L}\otimes_{\Z_p} \calOB_{\dR},\mathrm{id}_{\mathbb{L}}\otimes\nabla_{\calOB_{\dR}}),
 }
\]
where the top horizontal isomorphism is given by Proposition~\ref{prop:connection-filtration-compatibility}(1). By construction, this diagram is commutative.
Equip $\mathbf{E}$ and $\calE_\Q(\bB_\dR)$ with the filtrations induced from $\Fil^\bullet D_\dR(\mathbb{L})$. Then the horizontal maps and $\alpha_{\mathbf{E}}\otimes\calOB_\dR$ in the above commutative diagram are filtered isomorphisms (use Proposition~\ref{prop:connection-filtration-compatibility}(2) for the top horizontal map). Hence so is $\alpha_{\cris,\bL,\calE_\Q}\otimes_{\bB_{\cris}} \calOB_{\dR}$. By taking the horizontal sections, we conclude that $\alpha_{\cris,\bL,\calE_\Q}\otimes_{\bB_{\cris}} \bB_{\dR}$ is a filtered isomorphism on $X_{\eta,\proet}$.

(2) 
By the assumption and Proposition~\ref{prop:connection-filtration-compatibility}(2), 
the tensor product filtrations on $\mathbf{E} \otimes_{\calO_{X_\eta}}\calOB_{\dR}$ and $\bL\otimes_{\Z_p}\calOB_{\dR}$ coincide under the above isomorphism $\alpha_{ \calOB_{\dR}}$. Using the projection formula and the isomorphisms $\calO_{X_{\eta},\et}\xrightarrow{\cong}\nu_\ast\Fil^0\calOB_\dR\xrightarrow{\cong}\calOB_\dR$ and $\nu_\ast\Fil^n\calOB_\dR=0$ ($n\geq 1$), we see that the map $\Fil^n\mathbf{E}\rightarrow \nu_\ast\Fil^n(\mathbf{E} \otimes_{\calO_{X_\eta}}\calOB_{\dR})$ is an isomorphism for each $n\in\Z$, which gives the first assertion. The second assertion is obvious.

(3) The assumption and Proposition~\ref{prop:connection-filtration-compatibility} give natural isomorphisms
\begin{align*}
\Ker\bigl(\calE_\Q(\bB_{\cris})^{\varphi=1}\rightarrow \Fil^0(\mathbf{E}\otimes_{\calO_{X_\eta}}\calOB_\dR)^{\nabla=0}\bigl)
   &\cong \Ker\bigl((\bL\otimes_{\Z_p}\bB_\cris)^{\varphi=1}\rightarrow \Fil^0(\bL\otimes_{\Z_p}\calOB_{\dR})^{\nabla=0}\bigl)\\
&\cong \bL\otimes_{\Z_p}\Ker\bigl(\bB_\cris^{\varphi=1}\rightarrow \Fil^0(\calOB_{\dR})^{\nabla=0}\bigl)\\
&\cong \bL\otimes_{\Z_p}\Ker\bigl(\bB_\cris^{\varphi=1}\rightarrow \Fil^0\bB_{\dR}\bigl)=\bL\otimes_{\Z_p}\Q_p=\bL[p^{-1}],
\end{align*}
where the last isomorphism follows from the fundamental exact sequence in $p$-adic Hodge theory (cf.~\cite[Prop.~6.2.24]{brinon-relative}, \cite[Lem.~2.22(ii)]{GuoReinecke-Ccris}). 
\end{proof}

\begin{cor} \label{cor:prismatic-semistable-deRham-semistable-case}
Prismatic semistable $\Z_p$-local systems are de Rham. 
More precisely, for every analytic prismatic $F$-crystal $\calE_\Prism$ on a semistable $p$-adic log formal scheme $(X,M_X)$ over $\calO_K$,
\begin{enumerate}
    \item the $\Z_p$-local system $T_X(\calE_\Prism)$ and the $F$-isocrystal $D_\cris(\calE_\Prism)$ are associated;
    \item $T_X(\calE_\Prism)$ is de Rham and $(D_\cris(\calE_\Prism), D_\dR(T_X(\calE_\Prism)))$ is a filtered $F$-isocrystal on $(X,M_X)$.
 \end{enumerate}
\end{cor}

\begin{proof}
(1) is Proposition~\ref{prop:etaleandcrysassociated}, and (2) follows from (1) and Proposition~\ref{prop:association-filtration-compatibility}(1).
Recall that a $\Z_p$-local system is called prismatic semistable if it is in the essential image of $T_X$.
\end{proof}

\section{Semistable representations in the CDVR case} \label{sec:CDVR-semistable-representations}

In this section, we compare the three different notions of semistable representations in the CDVR case. Together with the purity result in the next section (Theorem~\ref{thm:main-purity}), Theorem~\ref{thm:CDVR-semistable-notions-equivalent}  below will imply that the notion of prismatic semistability in Definition~\ref{defn:semistablity-via-analy-prismatic-F-crystal} agrees with that of being associated in Definition~\ref{defn:association} for  (Corollary~\ref{cor:semistable-prismatic-Faltings-equivalent}) for $\mathbf{Z}_p$-local systems on the generic fiber of a semistable $p$-adic log formal scheme.

Let $\mathcal{O}_{L_0}$ be a Cohen ring whose residue field has a finite $p$-basis given by $\{X_1, \ldots, X_b\}$, and write $L_0 = \mathcal{O}_{L_0}[p^{-1}]$. We equip $\mathcal{O}_{L_0}$ with the Frobenius $\varphi\colon \mathcal{O}_{L_0} \rightarrow \mathcal{O}_{L_0}$ given by $\varphi(X_i) =X_i^p$. Let $\widehat{\Omega}_{\mathcal{O}_{L_0}}$ denote the $p$-adic completion of the K\"{a}hler differentials $\Omega_{\mathcal{O}_{L_0}/\mathbf{Z}}$: note 
$\widehat{\Omega}_{\mathcal{O}_{L_0}} \cong \bigoplus_{i=1}^b \mathcal{O}_{L_0}\cdot d\operatorname{log} X_i$.

Let $L$ be a finite totally ramified extension of $L_0$. Fix a uniformizer $\pi_L \in \mathcal{O}_L$, and let $E_L[u] \in \mathcal{O}_{L_0}[u]$ be the monic minimal polynomial for $\pi_L$. Set $G_L \coloneqq \Gal (\overline L / L)$. Let $Y=\Spf(\mathcal{O}_L)$ and equip it with the log structure $M_Y$ associated to the prelog structure $\mathbf{N} \rightarrow \mathcal{O}_L: n \mapsto \pi_L^n$.

In \cite{Morita-imperfsemistable}, the notion of semistable representations of $G_L$ is given within Fontaine's framework of admissibility using the period ring $\OB_{\st}$, which we recall in \S\ref{sec:classically-semistable}. The following is the main theorem of this section.

\begin{thm} \label{thm:CDVR-semistable-notions-equivalent}
Let $\mathbb{L}$ be a $\mathbf{Z}_p$-local system on the generic fiber $Y_{\eta}$. The following properties are equivalent:
\begin{enumerate}
    \item $\mathbb{L}$ is semistable in the sense of Definition~\ref{defn:semistablity-via-analy-prismatic-F-crystal};
    \item $\mathbb{L}$ is associated to an $F$-isocrystal on $(X_1,M_{X_1})_{\CRIS}$ as in Definition~\ref{defn:association};
    \item the $\mathbf{Q}_p$-representation of $G_L$ associated with $\mathbb{L}$ is classically semistable (i.e., $\OB_{\st}$-admissible).
\end{enumerate}
In particular, the \'etale realization functor induces an equivalence between the category of prismatic $F$-crystals on $(\Spf \calO_L, M_{\Spf \calO_L})_\Prism$ and the category of lattices in classically semistable representations of $G_L$.
\end{thm}

\begin{proof}
Proposition~\ref{prop:etaleandcrysassociated} shows that (1) implies (2). The implications (2)$\implies$(3) and (3)$\implies$(1) will be given in Propositions~\ref{prop:CDVR-semistable-implies-Bst-admissible} and \ref{prop:classically-semistable-implies-semistable}, respectively.
\end{proof}

\subsection{Classically semistable representations} \label{sec:classically-semistable}

We first recall the construction of period rings. They are given in \cite{brinon-crys-rep-imperfect-residue} and \cite{Morita-imperfsemistable}, but we follow the notation of \cite{scholze-p-adic-hodge}.

Let $\Ainf \coloneqq \Ainf(\mathcal{O}_{\overline{L}})$; it comes with the surjection $\theta\colon \Ainf \rightarrow \mathcal{O}_{\overline{L}}^{\wedge} \coloneqq (\mathcal{O}_{\overline{L}})^{\wedge}_p$. Define $\B_{\mathrm{dR}}^{+}$ to be the $\Ker(\theta[p^{-1}])$-adic completion of $\Ainf[p^{-1}]$. 
Write $\OA_\mathrm{inf}$ for the $(p,\Ker\theta_L)$-adic completion of $\mathcal{O}_{L}\otimes_{\mathbf{Z}_p} \Ainf$, where $\theta_L\colon \mathcal{O}_{L}\otimes_{\mathbf{Z}_p} \Ainf\rightarrow \mathcal{O}_{\overline{L}}^{\wedge}$. Then $\theta_L$ extends to $\theta_L[p^{-1}]\colon \OA_\mathrm{inf}[p^{-1}]\rightarrow \mathcal{O}_{\overline{L}}^{\wedge}[p^{-1}]$. Define $\OB_{\mathrm{dR}}^{+}$ to be the $\Ker (\theta_L[p^{-1}])$-adic completion of $\OA_\mathrm{inf}[p^{-1}]$. 
We have $t \coloneqq \operatorname{log}[\epsilon] \in \B_{\mathrm{dR}}^{+} \subset \OB_{\mathrm{dR}}^{+}$. Define the de Rham period ring by $\OB_{\mathrm{dR}} \coloneqq \OB_{\mathrm{dR}}^{+}[t^{-1}]$, which is equipped with a natural $G_L$-action and decreasing filtration.

For the semistable period ring, consider the $\mathcal{O}_{L_0}$-linear extension $\theta_{L_0}\colon \mathcal{O}_{L_0}\otimes_{\mathbf{Z}_p} \Ainf \rightarrow \mathcal{O}_{\overline{L}}^{\wedge}$ of $\theta$, and write $\OA_{\cris}$ for the $p$-adically completed PD-envelope of $\mathcal{O}_{L_0}\otimes_{\mathbf{Z}_p} \Ainf$ with respect to $\Ker \theta_{L_0}$. The natural map $\OA_{\cris} \rightarrow \OB_{\mathrm{dR}}^{+}$ is injective. 
Define $\OB_{\cris}^{+} \coloneqq \OA_{\st}[p^{-1}]$ and $\OB_{\cris} \coloneqq \OB_{\cris}^{+}[t^{-1}]$ (note $t \in \OA_{\cris}$).
Choose a compatible system $p^{\flat} = (p, p^{\frac{1}{p}}, \ldots) \in \mathcal{O}_{\overline{L}}^{\flat}$ of $p$-power roots of $p$. Then $\operatorname{log}([p^{\flat}]/p)\in\B_{\dR}^{+}$ is transcendental over $\OB_{\cris}$. We define the semistable period ring by  $\OB_{\st} \coloneqq \OB_{\cris}[\operatorname{log}([p^{\flat}]/p)] \subset \OB_{\mathrm{dR}}$; it is a polynomial ring over $\OB_{\cris}$ in $\operatorname{log}([p^{\flat}]/p)$ and is equipped with a natural Frobenius and $G_L$-action extending those on $\mathcal{O}_{L_0}\otimes_{\mathbf{Z}_p} \Ainf$ and satisfying $\varphi(\operatorname{log}([p^{\flat}]/p))=p\operatorname{log}([p^{\flat}]/p)$. Furthermore, the map $L\otimes_{L_0} \OB_{\st} \rightarrow \OB_{\dR}$ is injective and compatible with $G_L$-actions.

Let $\pi_L^{\flat} = (\pi_{L, 0}, \pi_{L, 1}, \ldots)$ (resp. $X_i^{\flat} = (X_{i, 0}, X_{i, 1}, \ldots)$ for $i = 1, \ldots, b$) be the element in $\mathcal{O}_{\overline{L}}^{\flat}$ given by choosing a compatible system of $p$-power roots of $\pi_L$ (resp. $X_i$). We consider two $\varphi$-equivariant embeddings $\lambda_1, \lambda_2\colon L_0 \rightarrow \OB_{\st}$ given by $\lambda_1(X_i) = X_i$ and $\lambda_2(X_i) = [X_i^{\flat}]$. We have an $\Ainf$-linear topologically quasi-nilpotent integrable connection 
\[
\nabla\colon \OB_{\st} \rightarrow \OB_{\st}\otimes_{\lambda_1, \mathcal{O}_{L_0}} \widehat{\Omega}_{\mathcal{O}_{L_0}}
\]
which is compatible with the universal derivation $\mathcal{O}_{L_0} \rightarrow \widehat{\Omega}_{\mathcal{O}_{L_0}}$ via $\lambda_1$. Let $N_u\colon \OB_{\st} \rightarrow \OB_{\st}$ be the $\OB_{\cris}$-linear derivation given by $N_u(\operatorname{log}([p^{\flat}]/p)) = 1$. 

For any finite-dimensional $\mathbf{Q}_p$-representation $V$ of $G_L$, let
\[
D_{\st}(V) \coloneqq (\OB_{\st}\otimes_{\mathbf{Q}_p} V)^{G_L}.
\]
It is equipped with an $L_0$-semilinear Frobenius $\varphi\colon D_{\st}(V) \rightarrow D_{\st}(V)$, an integrable connection
\[
\nabla\colon D_{\st}(V) \rightarrow D_{\st}(V)\otimes_{\mathcal{O}_{L_0}} \widehat{\Omega}_{\mathcal{O}_{L_0}}, 
\]
and a derivation $N_u\colon D_{\st}(V) \rightarrow D_{\st}(V)$ induced from those on $\OB_{\st}$. Furthermore, the embedding $L\otimes_{L_0} \OB_{\st} \rightarrow \OB_{\dR}$ induces a decreasing filtration on $L\otimes_{L_0} D_{\st}(V)$. By \cite{Morita-imperfsemistable}, $D_{\st}(V)$ is a finite-dimensional $L_0$-vector space with $\dim_{L_0} D_{\st}(V) \leq \dim_{\mathbf{Q}_p}(V)$, and the equality holds if and only if the natural map
\[
\OB_{\st}\otimes_{\lambda_1, L_0} D_{\st}(V)  \rightarrow \OB_{\st}\otimes_{\mathbf{Q}_p} V
\]
is an isomorphism. We say that $V$ is \textit{classically semistable} if $\dim_{L_0} D_{\st}(V) = \dim_{\mathbf{Q}_p}(V)$.

Consider the $\varphi$-compatible embedding $S_L \rightarrow \A_{\cris}(\mathcal{O}_{\overline{L}})$ given by $u \mapsto [\pi_L^{\flat}]$ and $X_i \mapsto [X_i^{\flat}]$. Write $N_u\colon S_L \rightarrow S_L$ for the $\mathcal{O}_{L_0}$-linear derivation given by $N_u(u) = u$\footnote{This is different from the convention $N_u (u) = -u$ in \cite{breuil-representations}, but it does not affect anything up to signs.}, and let $\nabla\colon S_L \rightarrow S_L\otimes_{\mathcal{O}_{L_0}} \widehat{\Omega}_{\mathcal{O}_{L_0}}$ be the integrable connection given by the universal derivation on $\mathcal{O}_{L_0}$, which commutes with $N_u$. 

\begin{prop} \label{prop:CDVR-semistable-implies-Bst-admissible}
Let $T$ be a finite free $\mathbf{Z}_p$-representation of $G_L$, and denote by $\mathbb{L}$ the corresponding $\mathbf{Z}_p$-local system on $Y_{\eta}$. Suppose $\mathbb{L}$ is associated with an $F$-isocrystal $(\calE_\Q,\varphi_{\calE_\Q})$ on $(Y_1,M_{Y_1})_{\CRIS}$ as in Definition~\ref{defn:association}. Then $V = T[p^{-1}]$ is classically semistable, namely, $\dim_{L_0} D_{\st}(V) = \dim_{\mathbf{Q}_p} V$.
\end{prop}

\begin{proof}
Consider the evaluation $\mathscr{D} \coloneqq \mathcal{E}_\Q(S_L)$ on the ind-object $(Y_1,\Spf S_L,M_{\Spf S_L})$, which is a finite projective $S_L[p^{-1}]$-module equipped with Frobenius $\varphi_{\mathscr{D}}$ and quasi-nilpotent integrable log connection 
\[
\nabla\colon \mathscr{D} \rightarrow \mathscr{D}\otimes_{\mathfrak{S}_L} \omega^1_{(Z,M_Z) / \mathbf{Z}_p} = \mathscr{D}\cdot d\operatorname{log} u \oplus \bigoplus_{i=1}^b \mathscr{D}\cdot d\operatorname{log} X_i
\]
as in Corollary~\ref{cor:iso-crystals and quasi-nilpotent connections}, where $Z=\Spf \calO_{L_0}\langle u\rangle$ with log structure induced from $\N\ni 1\mapsto u$. Let $N_u\colon \mathscr{D} \rightarrow \mathscr{D}$ (resp. $N_{X_i}\colon \mathscr{D} \rightarrow \mathscr{D}$ for each $i = 1, \ldots, b$) be the derivation given by $\nabla$ composed with the projection to the factor of $d\operatorname{log} u$ (resp. $d\operatorname{log} X_i$). Write $\partial_{X_i}\colon \mathscr{D} \rightarrow \mathscr{D}$ for the derivation $\partial_{X_i} = X_i^{-1}N_{X_i}$.

Let $q\colon S_L \rightarrow \mathcal{O}_{L_0}$ be the projection given by $u \mapsto 0$, and let $D\coloneqq \mathscr{D}\otimes_{S_L[p^{-1}], q} L_0$ equipped with the induced $\varphi$, $N_u$, and $N_{X_i}$-actions. Since $1\otimes\varphi_{\mathscr{D}} \colon S_L\otimes_{\varphi, S_L} \mathscr{D} \rightarrow \mathscr{D}$ is an isomorphism, so is the induced map $1\otimes\varphi_D \colon L_0\otimes_{\varphi, L_0} D \rightarrow D$. We have $N_u\varphi = p\varphi N_u$ and $N_{X_i}\varphi = p\varphi N_{X_i}$, so $N_u$ and $N_{X_i}$ on $D$ are nilpotent. 

We first claim that there exists a unique $\varphi$-compatible section $s\colon D \rightarrow \mathscr{D}$, and such section is compatible with $N_u$ and $N_{X_i}$'s. We follow \cite[Lem.~1.2.6, Pf.]{kisin-crystalline} and \cite[Lem.~3.14]{kim-groupscheme-relative}: choose a finitely generated $S_L$-submodule $\mathscr{M} \subset \mathscr{D}$ such that $\mathscr{M}[p^{-1}] = \mathscr{D}$, and let $M \coloneqq q(\mathscr{M}) \subset D$. Then $M$ is an $\mathcal{O}_{L_0}$-lattice of $D$ such that $M[p^{-1}] = D$. Take  $m_1, m_2\in\N$ such that $\varphi_{\mathscr{D}}(\mathscr{M}) \subset p^{-m_1}\mathscr{M}$ and $(1\otimes\varphi_D)^{-1} (M) \subset p^{-m_2}( \mathcal{O}_{L_0}\otimes_{\varphi, \mathcal{O}_{L_0}} M) \subset L_0\otimes_{\varphi, L_0} D$. Choose any $\mathcal{O}_{L_0}$-linear section $s_0\colon M \rightarrow \mathscr{M}$ of $q$. For each $n \geq 0$, write $\varphi^{* n} s_0 \coloneqq 1 \otimes_{\varphi^n, \mathcal{O}_{L_0}} s_0 \colon \mathcal{O}_{L_0}\otimes_{\varphi^n, \mathcal{O}_{L_0}} D \rightarrow S_L\otimes_{\varphi^n, S_L} \mathscr{D}$ and
\[
h_{n+1}\coloneqq (1\otimes\varphi_{\mathscr{D}}^{n+1})\circ (\varphi^{* n+1}s_0)\circ (1\otimes\varphi_D)^{-n-1} - (1\otimes\varphi_{\mathscr{D}}^{n})\circ (\varphi^{* n}s_0)\circ (1\otimes\varphi_D)^{-n}.
\]
Then $h_{n+1}(M) \subset p^{-m(n+1)}u^{p^n}(\mathscr{M})$ where $m = m_1+m_2$. Since the $p$-adic norm $\displaystyle |p^{-m(n+1)} (\lfloor p^n/e_L \rfloor!) |_p$ where $e_L \coloneqq [L:L_0]$ goes to $0$ as $n \rightarrow \infty$, the map $s \colon D \rightarrow \mathscr{D}$ given by $s\coloneqq s_0+\sum_{n=0}^{\infty} h_n$ is well-defined. By construction, $s$ is compatible with $\varphi$. The induced map $1\otimes s\colon S_L[p^{-1}]\otimes_{L} D \rightarrow \mathscr{D}$ is a map of finite projective $S_L[p^{-1}]$-modules of the same rank. So as in \cite[Lem.~4.19, Pf.]{moon-strly-div-latt-cryst-cohom-CDVF}, we see that $1\otimes s$ is an isomorphism. Furthermore, by a similar argument as in \cite[Lem.~3.14, Pf.]{kim-groupscheme-relative}, $s$ is the unique $\varphi$-compatible section and is compatible with $N_u$ and $N_{X_i}$'s. 

Since $\mathcal{E}_\Q$ is an isocrystal, we have $\mathcal{E}_\Q(\A_{\cris}(\mathcal{O}_{\overline{L}})) \cong \A_{\cris}(\mathcal{O}_{\overline{L}})\otimes_{S_L} \mathscr{D}$, and $G_L$-action on the left hand side is given via this identification by
\begin{equation} \label{eq:Galois-action-isocrystal}
g(a\otimes x) = g(a)\sum_{(j_0, \ldots, j_b)\in \N^{b+1}} \gamma_{j_0}\Bigl(\operatorname{log}\frac{g([\pi_L^{\flat}])}{[\pi_L^{\flat}]}\Bigr)\prod_{i=1}^b \gamma_{j_i}\Bigl(\operatorname{log}\frac{g([X_i^{\flat}])}{[X_i^{\flat}]}\Bigr) \cdot N_u^{j_0}N_{X_1}^{j_1}\cdots N_{X_b}^{j_b}(x)    
\end{equation}
for $g \in G_L$, $a \in \A_{\cris}(\mathcal{O}_{\overline{L}})$ and $x \in \mathscr{D}$. For any $x \in s(D)$, set
\[
\overline{x} \coloneqq \sum_{(j_0, \ldots, j_b)\in \N^{b+1}} \gamma_{j_0}\Bigl(-\operatorname{log}\frac{[\pi_L^{\flat}]}{\pi_L}\Bigr)\prod_{i=1}^b \gamma_{j_i}\Bigl(-\operatorname{log}\frac{[X_i^{\flat}]}{X_i}\Bigr) \cdot N_u^{j_0}N_{X_1}^{j_1}\cdots N_{X_b}^{j_b}(x).
\]
Note that only finitely many terms in the above sum are non-zero since $N_u$ and $N_{X_i}$'s are nilpotent on $s(D)$, and so $\overline{x}$ is a well-defined element in $\OB_{\st}(\mathcal{O}_{\overline{L}})\otimes_{S_L} \mathscr{D}$. We have $\overline{\lambda_2(a)\cdot x} = \lambda_1(a)\cdot \overline{x} \in \OB_{\st}(\mathcal{O}_{\overline{L}})\otimes_{S_L} \mathscr{D}$ for any $a \in L_0$.

Now, consider the $L_0$-vector space $\overline{D} \coloneqq \{\overline{x} \mid x \in s(D)\}$. Since $s$ is injective, $\dim_{L_0} D = \dim_{L_0} \overline{D}$. Furthermore, since the $G_L$-action is given by Equation~(\ref{eq:Galois-action-isocrystal}), a similar computation as in \cite[\S7.2]{liu-fontaineconjecture} shows that $g(\overline{x}) = \overline{x}$ for any $g \in G_L$ and $\overline{x} \in \overline{D}$. Consider the $G_L$-equivariant isomorphism
\[
\OB_{\st}(\mathcal{O}_{\overline{L}})\otimes_{S_L} \mathscr{D} \cong \OB_{\st}(\mathcal{O}_{\overline{L}})\otimes_{\mathbf{Q}_p} V
\]
obtained from $\A_{\cris}(\mathcal{O}_{\overline{L}})\otimes_{S_L} \mathscr{D} \cong \A_{\cris}(\mathcal{O}_{\overline{L}})\otimes_{\mathbf{Z}_p} V$ by the base change along $\A_{\cris}(\mathcal{O}_{\overline{L}}) \rightarrow \OB_{\st}(\mathcal{O}_{\overline{L}})$. We then have 
\[
\overline{D} \subset (\OB_{\st}(\mathcal{O}_{\overline{L}})\otimes_{\mathbf{Q}_p} V)^{G_L},
\]
so $\dim_{L_0} D_{\st}(V) \geq \dim_{L_0} \overline{D} = \dim_{\mathbf{Q}_p} V$.
\end{proof}

\begin{rem}
The log crystalline site $(Y_1, M_\mathrm{triv})_\CRIS$ for the trivial log structure $M_\mathrm{triv}$ agrees with the usual big absolute crystalline site on $Y_1$.  With the notation as in the proof of Proposition~\ref{prop:CDVR-semistable-implies-Bst-admissible}, every $F$-isocrystal on $(Y_1, M_\mathrm{triv})_\CRIS$ defines an $F$-isocrystal on $(Y_1, M_{Y_1})_\CRIS$ such that the associated $N_u\colon \mathscr{D}\rightarrow\mathscr{D}$ is the zero map. In particular, if $\bL$ is associated with an $F$-isocrystal on $(Y_1, M_\mathrm{triv})_\CRIS$ in Proposition~\ref{prop:CDVR-semistable-implies-Bst-admissible}, then $V$ is classically crystalline, namely, $\OB_\cris$-admissible.
\end{rem}

\subsection{Kisin descent data} \label{sec:quasi Kisin module and rational descent data} 

Let $T$ be a finite free $\mathbf{Z}_p$-representation of $G_L$ such that $V \coloneqq T[p^{-1}]$ is classically semistable. In this subsection, we show that the Laurent $F$-crystal on $(Y, M_Y)_{\Prism}$ corresponding to $T$ extends to a prismatic $F$-crystal, and so it is semistable in the sense of Definition~\ref{defn:semistablity-via-analy-prismatic-F-crystal}. We follow \cite[\S4]{du-liu-moon-shimizu-completed-prismatic-F-crystal-loc-system} with suitable modifications.

Let $(\mathfrak{S}_L, (E_L), M_{\operatorname{Spf}\mathfrak{S}_L})$ (resp. $(S_L, (p), M_{\operatorname{Spf}S_L})$) be the Breuil--Kisin prism (resp. Breuil prism) in $(Y, M_Y)_{\Prism}$ as before. With fixed choices $\pi_L^{\flat} = (\pi_{L, 0}, \pi_{L, 1}, \ldots)$ and  $X_i^{\flat} = (X_{i, 0}, X_{i, 1}, \ldots)$ of $p$-power roots, let
\[
\widetilde{L}_{\infty} \coloneqq \bigcup_{n \geq 0} L(\pi_{L,n}, X_{1, n}, \ldots, X_{b, n}) \subset \overline{L}
\]
and $G_{\widetilde{L}_{\infty}} \coloneqq \mathrm{Gal}(\overline{L} / \widetilde{L}_{\infty})$. Write $\mathcal{O}_{\mathcal{E}, L} \coloneqq \mathfrak{S}_L[u^{-1}]^{\wedge}_p$. The $\varphi$-equivariant embedding $\mathfrak{S}_L \rightarrow W(\mathcal{O}_{\overline{L}}^{\flat})$ given by $u \mapsto [\pi_L^{\flat}]$ and $X_i \mapsto [X_i^{\flat}]$ extends uniquely to $\mathcal{O}_{\mathcal{E}, L} \rightarrow W(\overline{L}^{\flat})$. Let $\widehat{\mathcal{O}}_{\mathcal{E}, L}^{\mathrm{ur}} \subset W(\overline{L}^{\flat})$ be the $p$-adic completion of the ring of integers of the maximal unramified extension of $\mathcal{O}_{\mathcal{E}, L}[p^{-1}]$ inside $W(\overline{L}^{\flat})[p^{-1}]$. 

Without loss of generality, we may assume that the Hodge-Tate weights of $T[p^{-1}]$ lie in $[0, r]$ for some integer $r \geq 1$. By \cite[Thm.~3.2.3, Lem.~4.2.9]{gao-integral-padic-hodge-imperfect} (see \cite[\S4.3]{moon-strly-div-latt-cryst-cohom-CDVF}), we can naturally associate to $T$ an $\mathfrak{S}_L$-module $\mathfrak{M}$ satisfying the following properties:
\begin{itemize}
    \item $\mathfrak{M}$ is finite free over $\mathfrak{S}_L$;
    \item $\mathfrak{M}$ is equipped with a $\varphi$-semilinear endomorphism $\varphi_{\mathfrak{M}}\colon \mathfrak{M} \rightarrow \mathfrak{M}$ with $E_L$-height $\leq r$;
    \item Set 
    \[
    M \coloneqq \mathcal{O}_{L_0}\otimes_{\varphi,\mathcal{O}_{L_0}}\mathfrak{M}/u\mathfrak{M}
    \]
    equipped with the induced tensor-product Frobenius. We have a natural isomorphism of $L_0$-modules $M[p^{-1}] \cong D_{\mathrm{st}}(V)$ compatible with Frobenii;
    \item $\mathfrak{M}\otimes_{\mathfrak{S}_L} \mathcal{O}_{\mathcal{E}, L}$ is the \'etale $\varphi$-module over $\mathcal{O}_{\mathcal{E}, L}$ such that 
    \[
    ((\mathfrak{M}\otimes_{\mathfrak{S}_L} \mathcal{O}_{\mathcal{E}, L})\otimes_{\mathcal{O}_{\mathcal{E}, L}} \widehat{\mathcal{O}}_{\mathcal{E}, L}^{\mathrm{ur}})^{\varphi = 1} =  T|_{G_{\widetilde{L}_{\infty}}}.
    \]
\end{itemize}

Consider
\[
\mathscr{M} \coloneqq S_L\otimes_{\varphi, \mathfrak{S}_L} \mathfrak{M}
\]
equipped with the tensor product Frobenius. The projection $q\colon S_L \twoheadrightarrow \mathcal{O}_{L_0}$ given by $u \mapsto 0$ induces the projection $q\colon \mathscr{M} \twoheadrightarrow M$ compatible with $\varphi$. By \cite[Lem.~4.2]{du-liu-moon-shimizu-completed-prismatic-F-crystal-loc-system}, $q$ admits a unique $\varphi$-compatible section $s\colon M[p^{-1}] \rightarrow \mathscr{M}[p^{-1}]$, and $1\otimes s\colon S_L[p^{-1}]\otimes_{L_0}M[p^{-1}] \rightarrow \mathscr{M}[p^{-1}]$ is an isomorphism. Set 
\[
\mathscr{D} \coloneqq S_L[p^{-1}]\otimes_{L_0} M[p^{-1}] \cong S_L[p^{-1}]\otimes_{L_0} D_{\mathrm{st}}(V).
\]
Then $\varphi$, $N_u$, and $\nabla$ on $D_{\mathrm{st}}(V)$ naturally extend to $\mathscr{D}$ by the tensor product of corresponding structures on $S_L$. We remark that \cite[\S4.3]{du-liu-moon-shimizu-completed-prismatic-F-crystal-loc-system} studies the crystalline case, in which we have an extra condition that $N_u=0$ on $D_\st(V)$, so $N_u$ on $\mathscr{D}$ in \textit{loc. cit.} is defined as $N_u=N_{u,S_L}+1$.

For each integer $i \geq 0$, let $\mathrm{Fil}^i S_L$ be the PD-filtration of $S_L$. Define a decreasing filtration on $\mathscr{D}$ by setting $\mathrm{Fil}^0 \mathscr{D} = \mathscr{D}$ and inductively for $i \geq 1$ 
\[
\mathrm{Fil}^i \mathscr{D} = \{x \in \mathscr{D} ~|~ N_u(x) \in \mathrm{Fil}^{i-1} \mathscr{D}, ~q_{\pi_L}(x) \in \mathrm{Fil}^i(L\otimes_{L_0} D_{\mathrm{st}}(V))\}
\]
where $q_{\pi}\colon \mathscr{D} \rightarrow L\otimes_{L_0} D_{\mathrm{st}}(V)$ is the projection induced by $u\mapsto \pi_L$. The following fact is shown in \cite{moon-strly-div-latt-cryst-cohom-CDVF}.

\begin{lem}[{\cite[Lem.~4.2]{moon-strly-div-latt-cryst-cohom-CDVF}}] \label{lem:Griffiths-transversality-CDVR}
The connection $\nabla$ on $\mathscr{D}$ satisfies Griffiths transversality:
\[
\nabla(\mathrm{Fil}^{i+1} \mathscr{D}) \subset \mathrm{Fil}^i \mathscr{D}\otimes_{\mathcal{O}_{L_0}} \widehat{\Omega}_{\mathcal{O}_{L_0}}.
\]
\end{lem}

Using the $\varphi$-compatible isomorphism $\mathscr{D} \cong \mathscr{M}[p^{-1}]$, we consider another decreasing filtration on $\mathscr{D}$ given by 
\[
\mathrm{F}^i \mathscr{M}[p^{-1}] = \{x \in \mathscr{M}[p^{-1}] ~|~ (1\otimes \varphi_{\mathfrak{M}})(x) \in \mathrm{Fil}^i S_L[p^{-1}]\otimes_{\mathfrak{S}_L} \mathfrak{M}\}. 
\]

\begin{lem} \label{lem:filtrations-compatible}
$\mathrm{Fil}^i \mathscr{D} = \mathrm{F}^i \mathscr{M}[p^{-1}]$.
\end{lem}

\begin{proof}
As in \cite[Lem.~4.15, Pf.]{du-liu-moon-shimizu-completed-prismatic-F-crystal-loc-system}, we can reduce the general case to the perfect residue case, which is treated in \cite[Cor.~3.2.3, Pf.]{liu-semistable-lattice-breuil}. 
\end{proof}

Now we construct the Kisin descent datum associated with $T$: for $i = 1, \ldots, b$, let $N_{X_i}\colon \mathscr{D} \rightarrow \mathscr{D}$ be the derivation given as the composite 
\[
\nabla\colon \mathscr{D} \xrightarrow{\nabla} \mathscr{D}\otimes_{\mathcal{O}_{L_0}} \widehat{\Omega}_{\mathcal{O}_{L_0}} = \bigoplus\limits_{i=1}^b \mathscr{D} \cdot d\operatorname{log} X_i\xrightarrow{\pr_i}\mathscr{D}.
\]
 Define $f_{S_L}\colon S_L^{(1)}\otimes_{p^1_1, S_L} \mathscr{D} \rightarrow S_L^{(1)}\otimes_{p^1_2, S_L} \mathscr{D}$ by
\begin{equation} \label{eq:descent isom on S_L^1}
f_{S_L}(x) = \sum_{(j_0, \ldots, j_b)\in\N^{b+1}} \gamma_{j_0}\Bigl( \operatorname{log}\frac{p^1_2(u)}{p^1_1(u)}\Bigr) \prod_{i=1}^n \gamma_{j_i}\Bigl( \operatorname{log}\frac{p^1_2(X_i)}{p^1_1(X_i)}\Bigr) \cdot N_u^{j_0}N_{X_1}^{j_1}\cdots N_{X_b}^{j_b}(x).  
\end{equation} 
Here the sum converges since $N_u$ and $N_{X_i}$'s are nilpotent on $D_{\st}(V)$. By a standard computation, $f_{S_L}$ is a $\varphi$-compatible isomorphism of $S_L^{(1)}[p^{-1}]$-modules satisfying the cocycle condition over $S_L^{(2)}[p^{-1}]$. Furthermore, $f_{S_L}$ is compatible with filtrations by Lemma~\ref{lem:Griffiths-transversality-CDVR}.  

\begin{prop}\label{prop:rational Kisin descent datum in CDVR case}
There exists a unique isomorphism
\[
f\colon \mathfrak{S}_L^{(1)}[p^{-1}]\otimes_{p^1_1, \mathfrak{S}_L} \mathfrak{M} \stackrel{\cong}{\rightarrow} \mathfrak{S}_L^{(1)}[p^{-1}]\otimes_{p^1_2, \mathfrak{S}_L} \mathfrak{M}
\]
of $\mathfrak{S}_L^{(1)}$-modules satisfying the cocycle condition over $\mathfrak{S}_L^{(2)}$ such that $\mathrm{id}_{S_L^{(1)}}\otimes_{\varphi, \mathfrak{S}_L^{(1)}} f = f_{S_L}$.
\end{prop}

\begin{proof}
Given Lemmas~\ref{prop:CDVR-filtration}(3), \ref{lem:filtration-h0-CDVR} and \ref{lem:filtrations-compatible}, this follows from the same argument as in the proof of \cite[Lem.~4.5, Prop.~4.6, 4.9]{du-liu-moon-shimizu-completed-prismatic-F-crystal-loc-system}.
\end{proof}

Let $(\mathcal{M}, \varphi_{\mathcal{M}}, f_{\mathrm{\acute{e}t}})$ be the triple associated with the $G_L$-representation $T$ via Theorem~\ref{thm:logetalerealizationofLaurentFcrystals-strict} and Lemma~\ref{lem:equiv-Laurent-F-cryst-descent-datum}. Note that $\mathcal{O}_{\mathcal{E}, L}\otimes_{\mathfrak{S}_L} \mathfrak{M} \cong \mathcal{M}$ as \'etale $\varphi$-modules over $\mathcal{O}_{\mathcal{E}, L}$ since $((\mathfrak{M}\otimes_{\mathfrak{S}_L} \mathcal{O}_{\mathcal{E}, L})\otimes_{\mathcal{O}_{\mathcal{E}, L}} \widehat{\mathcal{O}}_{\mathcal{E}, L}^{\mathrm{ur}})^{\varphi = 1} =  T|_{G_{\widetilde{L}_{\infty}}}$. So $f_{\mathrm{\acute{e}t}}$ corresponds to a $\mathfrak{S}_L^{(1)}$-linear isomorphism
\[
f_{\mathrm{\acute{e}t}}\colon \mathfrak{S}_L^{(1)}[E_L^{-1}]^{\wedge}_{p}\otimes_{p^1_1, \mathfrak{S}_L} \mathfrak{M} \stackrel{\cong}{\rightarrow} \mathfrak{S}_L^{(1)}[E_L^{-1}]^{\wedge}_p\otimes_{p^1_2, \mathfrak{S}_L} \mathfrak{M}
\]
compatible with Frobenii and satisfying the cocycle condition over $\mathfrak{S}_L^{(2)}[E_L^{-1}]^{\wedge}_{p}$.

\begin{prop} \label{prop:CDVR-Galois-reps-compatible}
As maps of $\mathfrak{S}_L^{(1)}[E_L^{-1}]^{\wedge}_p [p^{-1}]$-modules, we have 
\[
f_{\mathrm{\acute{e}t}}[p^{-1}]=\mathrm{id}_{\mathfrak{S}_L^{(1)}[E_L^{-1}]^{\wedge}_p[p^{-1}]}\otimes_{\mathfrak{S}_L^{(1)}[p^{-1}]} f.
\]
\end{prop} 

\begin{proof}
We follow \cite[Prop.~4.33, Pf.]{du-liu-moon-shimizu-completed-prismatic-F-crystal-loc-system}. By \cite[Prop.~4.34]{du-liu-moon-shimizu-completed-prismatic-F-crystal-loc-system}, there exists a $\mathfrak{S}_L$-submodule $\mathfrak{N} \subset \mathfrak{M}$ stable under $\varphi_{\mathfrak{M}}$ such that $\mathfrak{N}$ is free over $\mathfrak{S}_L$ of $E_L$-height $\leq r$ satisfying $\mathfrak{N}[p^{-1}] = \mathfrak{M}[p^{-1}]$ and $f$ induces an isomorphism of $\mathfrak{S}_L^{(1)}$-modules
\[
f\colon \mathfrak{S}_L^{(1)}\otimes_{p^1_1, \mathfrak{S}_L} \mathfrak{N} \stackrel{\cong}{\rightarrow} \mathfrak{S}_L^{(1)}\otimes_{p^1_2, \mathfrak{S}_L} \mathfrak{N}.
\] 
So $\mathcal{N}\coloneqq \mathcal{O}_{\mathcal{E}, L}\otimes_{\mathfrak{S}_L}\mathfrak{N}$ is a finite free \'etale $\varphi$-module over $\mathcal{O}_{\mathcal{E}, L}$, and $f$ induces an isomorphism 
\[
f\colon \mathfrak{S}_L^{(1)}[E_L^{-1}]^{\wedge}_{p}\otimes_{p^1_1, \mathcal{O}_{\mathcal{E}, L}} \mathcal{N} \stackrel{\cong}{\rightarrow} \mathfrak{S}_L^{(1)}[E_L^{-1}]^{\wedge}_p\otimes_{p^1_2, \mathcal{O}_{\mathcal{E}, L}} \mathcal{N}
\]
compatible with Frobenii and satisfying the cocycle condition over $\mathfrak{S}_L^{(2)}[E_L^{-1}]^{\wedge}_{p}$. By Theorem~\ref{thm:logetalerealizationofLaurentFcrystals-strict} and Lemma~\ref{lem:equiv-Laurent-F-cryst-descent-datum}, this corresponds to a finite free $\mathbf{Z}_p$-representation $T'$ of $G_L$, and $T'$ is determined by the $G_L$-action on $W(\mathcal{O}_{\overline{L}}^{\flat}[(\pi_L^{\flat})^{-1}])\otimes_{\mathcal{O}_{\mathcal{E}, L}}\mathcal{N}$. Furthermore, since $\mathfrak{N}$ has a finite $E_L$-height, arguing as in the proof of Lemma~\ref{lem:isomafterinvertingmu} leads to the identification
\[
\A_{\inf}(\mathcal{O}_{\overline{L}})[p^{-1}, \mu^{-1}]\otimes_{\mathfrak{S}_L} \mathfrak{N} = \A_{\inf}(\mathcal{O}_{\overline{L}})[p^{-1}, \mu^{-1}]\otimes_{\mathbf{Q}_p} T'[p^{-1}]. 
\]
Here, the map $\mathfrak{S}_L \rightarrow \A_{\inf}(\mathcal{O}_{\overline{L}})$ is given by $u \mapsto [\pi_L^{\flat}]$ and $X_i \mapsto [X_i^{\flat}]$. It suffices to show that $T'[p^{-1}] = V$ as representations of $G_L$.

Consider the map of log prisms $(S_L, (p), M_{\operatorname{Spf}(S_L)}) \rightarrow (\A_{\cris}(\mathcal{O}_{\overline{L}}), (p), M_{\operatorname{Spf}(\A_{\cris}(\mathcal{O}_{\overline{L}}))})$ in $(Y, M_Y)_{\Prism}$ given by $u \mapsto [\pi_L^{\flat}]$ and $X_i \mapsto [X_i^{\flat}]$ as in \S\ref{sec:classically-semistable}. Note that $\mathscr{D} = S_L[p^{-1}]\otimes_{\varphi, \mathfrak{S}_L} \mathfrak{M}$. By the construction of the descent datum $f$ from $f_{S_L}$ in Equation~(\ref{eq:descent isom on S_L^1}), the $G_L$-action on $\A_{\cris}(\mathcal{O}_{\overline{L}}))[p^{-1}]\otimes_{\varphi, \mathfrak{S}_L} \mathfrak{N} = \A_{\cris}(\mathcal{O}_{\overline{L}}))\otimes_{S_L} \mathscr{D}$ is given by Equation~(\ref{eq:Galois-action-isocrystal}). On the other hand, by \cite[\S4.1]{moon-strly-div-latt-cryst-cohom-CDVF}, we have a natural $G_L$-equivariant isomorphism
\[
V \cong \mathrm{Fil}^0(\A_{\cris}(\mathcal{O}_{\overline{L}}))\otimes_{S_L} \mathscr{D})^{\varphi = 1}
\]
where the right hand side is equipped with the $G_L$-action via Equation~(\ref{eq:Galois-action-isocrystal}). Note that the relevant functors in \cite{moon-strly-div-latt-cryst-cohom-CDVF} are contravariant, and it is straightforward to check that the result also hold for the covariant version by the same argument. Thus, we have $T'[p^{-1}] = V$ as representations of $G_L$.
\end{proof}

\begin{prop} \label{prop:classically-semistable-implies-semistable}
The Laurent $F$-crystal on $(Y, M_Y)_{\Prism}$ corresponding to $T$ via Theorem~\ref{thm:logetalerealizationofLaurentFcrystals-strict} extends to a prismatic $F$-crystal.    
\end{prop}

\begin{proof}
We have $\mathfrak{S}_L^{(1)}[p^{-1}] \cap \mathfrak{S}_L^{(1)}[E_L^{-1}]^{\wedge}_p = \mathfrak{S}_L^{(1)}$ by Lemma~\ref{lem:CDVR-BK-self-products-properties}. So by Proposition~\ref{prop:CDVR-Galois-reps-compatible}, we obtain from $f$ and $f_{\mathrm{\acute{e}t}}$ a $\mathfrak{S}_L^{(1)}$-linear isomorphism
\[
\mathfrak{S}_L^{(1)}\otimes_{p^1_1, \mathfrak{S}_L} \mathfrak{M} \stackrel{\cong}{\rightarrow} \mathfrak{S}_L^{(1)}\otimes_{p^1_2, \mathfrak{S}_L} \mathfrak{M}
\]
compatible with Frobenii and satisfying the cocycle condition over $\mathfrak{S}_L^{(2)}$. By Lemma~\ref{lem:descentlemmamoduleversion}, this gives a prismatic $F$-crystal extending the Laurent $F$-crystal associated with $T$.  
\end{proof}

\section{Purity for logarithmic prismatic \texorpdfstring{$F$}{F}-crystals} \label{sec:purity}

\subsection{Statement and applications} \label{subsec:purity-statement-applications}

Let $(X, M_X)$ be a semistable $p$-adic log formal scheme over $\calO_K$. Let $\{\xi_1,\ldots,\xi_m\}\subset X$ be the set of the generic points of all the irreducible components of  $X$. For each such point $\xi_j$, the $p$-adic completion $\calO_{X,\xi_j}^\wedge$ of the local ring $\calO_{X,\xi_j}$ is a complete discrete valuation ring with $\pi$ a uniformizer. Set $\Delta_j\coloneqq \Spf \calO_{X,\xi_j}^\wedge$ and equip it with the pullback log structure $M_{\Delta_j}$ from $M_X$. Note that $M_{\Delta_j}$ is associated to the prelog structure $\N\ni 1\mapsto \pi\in \calO_{X,\xi_j}^\wedge$. The morphism $f_j\colon(\Delta_j,M_{\Delta_j})\rightarrow (X,M_X)$ induces a morphism of topoi $(f_j)_{\Prism}\colon \operatorname{Sh}(X,M_X)_\Prism\rightarrow \operatorname{Sh}(\Delta_j,M_{\Delta_j})_\Prism$. 

In this section, we prove the following purity result:

\begin{thm}[Prismatic purity] \label{thm:main-purity}
A Laurent $F$-crystal $\mathscr{E}$ on $(X,M_X)_\Prism$ extends to an analytic prismatic $F$-crystal on $(X,M_X)$ if and only if $(f_j)_{\Prism}^{-1}\mathscr{E}$ extends to an analytic prismatic $F$-crystal on $(\Delta_j,M_{\Delta_j})$ for each $1\leq j\leq m$.    
\end{thm}

Let us first explain the applications; the if-part of the purity theorem is proved in the next subsection.

\begin{cor} \label{cor:semistable-prismatic-Faltings-equivalent}
Let $\mathbb{L}$ be a $\mathbf{Z}_p$-local system on $X_{\eta}$. Then $\mathbb{L}$ is semistable in the sense of Definition~\ref{defn:semistablity-via-analy-prismatic-F-crystal} if and only if it is associated to an $F$-isocrystal on $(X_1,M_{X_1})_{\CRIS}$ as in Definition~\ref{defn:association}.    
\end{cor}

\begin{proof}
This follows from Proposition~\ref{prop:etaleandcrysassociated}, Theorems~\ref{thm:CDVR-semistable-notions-equivalent} and \ref{thm:main-purity}.    
\end{proof}

In fact, the above conditions are also equivalent to being semistable point-wise on $X_\eta$ and thus independent of the choice of the semistable model $X$. To explain this, let us recall the following definition in \cite{bhatt_hansen-6-functor-rigid}.

\begin{defn}[{\cite[Def.~2.5]{bhatt_hansen-6-functor-rigid}}]
Let $\mathcal{X}$ be an adic space that is topologically of finite type over $K$. A rank-$1$ point $x\in \mathcal{X}$ is called \emph{weakly Shilov} if there is an open affinoid subset $\Spa(A, A^\circ)\subset \mathcal{X}$ containing $x$ such that the image of $x$ under the specialization map $\mathrm{sp}_{A^\circ}\colon \lvert \Spa(A, A^\circ) \rvert \to \lvert \Spec((A^\circ/\pi A^\circ)_{\mathrm{red}})\rvert$ is a generic point. 
\end{defn}

\begin{rem} \label{rem:Sholiv pts and generic pts}
For any affine open $U=\Spf(R)$ of $X$ containing $\xi_j$ as in the beginning of this section, the map $\Delta_j \to U$ defines a continuous valuation of $R[p^{-1}]$. Moreover, by \cite[Prop.~2.2]{bhatt_hansen-6-functor-rigid}, this valuation defines the unique rank-$1$ point of $\Spa(R[p^{-1}],R)$ mapping to $\xi_j$ under $\mathrm{sp}_{R}$. The rank-$1$ points of $X_\eta$ corresponding to $\{\xi_1,\ldots, \xi_m\}$ are called \emph{$X$-Shilov points}.
\end{rem}

\begin{cor} \label{cor:independent-model}
Let $X$ be a semistable formal scheme, and let $\mathbb{L}$ be a $\mathbf{Z}_p$-local system on $X_{\eta}$. Then the following properties are equivalent:
\begin{enumerate}
    \item $\mathbb{L}$ is semistable in the sense of Definition~\ref{defn:semistablity-via-analy-prismatic-F-crystal};
    \item $x^\ast\mathbb{L}$ is a semistable representation of $G_{\kappa(x)}$ for all rank-$1$ points $x$ of $X_\eta$;
    \item $x^\ast\mathbb{L}$ is a semistable representation of $G_{\kappa(x)}$ for all rank-$1$ weakly Shilov points $x$ of $X_\eta$;
    \item $x^\ast\mathbb{L}$ is a semistable representation of $G_{\kappa(x)}$ for all $X$-Shilov points $x$ of $X_\eta$.
\end{enumerate}
In particular, $\mathbb{L}$ being semistable is independent of the choice of semistable integral models $X$.
\end{cor}

\begin{proof}
The implications $(2)\implies(3)\implies(4)$ are obvious.
By Proposition~\ref{prop:pullback commutes with etale realization} and Theorem~\ref{thm:CDVR-semistable-notions-equivalent}, (1) implies (2). Finally, (1) follows from (4) by Theorems~\ref{thm:CDVR-semistable-notions-equivalent} and \ref{thm:main-purity}.
\end{proof}

\begin{rem} \label{rem:Andreatta-Iovita-semistability}
Suppose that $(X, M_X)$ is small affine and write $X = \Spf R$. By fixing a lift $\widetilde{R}$ of $R$ over $W[\![u]\!]$, Andreatta and Iovita construct the period ring $B_{\mathrm{log}}^{\cris}(\widetilde{R})$ and define the notion of semistable representations of $\mathcal{G}_R$ via $B_{\mathrm{log}}^{\cris}(\widetilde{R})$-admissibility, that is, representations that satisfy the equivalent conditions in \cite[Prop.~3.60]{andreatta-iovita-semistable-relative}. For a $\mathbf{Z}_p$-local system $\mathbb{L}$ on $X_{\eta}$ that is semistable in the sense of Definition~\ref{defn:semistablity-via-analy-prismatic-F-crystal}, then by $(1)$ in Remark~\ref{rem:GaloisactiononcrysllinerealizationoverAcris} and the fact that there is a $\mathcal{G}_R$-equivariant map from 
$\mathrm{B}_{\cris}(\overline{R})$ to $B_{\mathrm{log}}^{\cris}(\widetilde{R})$, $\mathbb{L}$ is $B_{\mathrm{log}}^{\cris}(\widetilde{R})$-admissible. Conversely, if $\mathbb{L}$ is $B_{\mathrm{log}}^{\cris}(\widetilde{R})$-admissible, then its pullback at every $X$-Shilov point is semistable (we remark that \cite[Prop.~3.62(i)]{andreatta-iovita-semistable-relative} and its proof contain an error, but an analogue of Theorem~\ref{thm:CDVR-semistable-notions-equivalent} holds for their $B_{\mathrm{log}}^{\cris}(\widetilde{R})$-admissibility). It follows from Corollary~\ref{cor:independent-model} that $\mathbb{L}$ is semistable in the sense of Definition~\ref{defn:semistablity-via-analy-prismatic-F-crystal}.  
\end{rem}

\subsection{Proof of Theorem~\ref{thm:main-purity}} \label{sec:proof-purity}

To prove the if-part of Theorem~\ref{thm:main-purity}, we may assume that $(X, M_X) = (\operatorname{Spf}R, M_{\operatorname{Spf}R})$ is small affine by Proposition~\ref{prop:etale-realization-fullyfaithful}. Furthermore, we may assume that $R$ admits a framing
\[
\square\colon R^{0} = \mathcal{O}_K \langle T_1, \ldots, T_m, T_{m+1}^{\pm 1}, \ldots, T_d^{\pm 1}\rangle / (T_1\cdots T_m - \pi) \rightarrow R,
\]
which induces a bijection between the Shilov points. For each $j = 1, \ldots, m$, let $R \rightarrow \mathcal{O}_{L_j}$ be the map corresponding to the Shilov point $x_j$ of $X_{\eta}$ as in \S~\ref{sec:prelim-facts-rings}. 

Let $(\mathscr{E}, \varphi_{\mathscr{E}})$ be a Laurent $F$-crystal on $(X, M_X)_{\Prism}$, and let $(\mathcal{M}, \varphi_{\mathcal{M}}, f_{\mathrm{\acute{e}t}})$ be the triple corresponding to $(\mathscr{E}, \varphi_{\mathscr{E}})$ by Lemma~\ref{lem:equiv-Laurent-F-cryst-descent-datum}. For each $1 \leq j \leq m$, suppose that the pullback $\mathscr{E}_j$ of $\mathscr{E}$ along $\operatorname{Spf} \mathcal{O}_{L_j} \rightarrow \operatorname{Spf} R$ extends to a prismatic $F$-crystal $\mathscr{E}_{\Prism, L_j}$ on $(\operatorname{Spf} \mathcal{O}_{L_j}, M_{\operatorname{Spf} \mathcal{O}_{L_j}})_{\Prism}$. We are going to show that $\mathscr{E}$ extends to an analytic prismatic $F$-crystal on $(X, M_X)_{\Prism}$ by constructing a Kisin descent datum $(\fkM,\varphi_\fkM,f)$ in Definition~\ref{defn:BK-descent-data}. 

Without loss of generality, we may assume that each $\mathscr{E}_{\Prism, L_j}$ is effective. Let 
\[
(\mathfrak{M}_{L_j}, \varphi_{\mathfrak{M}_{L_j}}, f_{L_j}) \in \mathrm{DD}_{\mathfrak{S}_{L_j}}^{\mathrm{vb}}
\]
be the Kisin descent datum corresponding to $\mathscr{E}_{\Prism, L_j}$ via Lemma~\ref{lem:descentlemmamoduleversion}. Note that 
\[
\mathfrak{S}_{L_j}[E^{-1}]^{\wedge}_p\otimes_{\mathfrak{S}_{L_j}} \mathfrak{M}_{L_j} = \mathfrak{S}_{L_j}[E^{-1}]^{\wedge}_p\otimes_{\mathcal{O}_{\mathcal{E}}} \mathcal{M} 
\]  
compatibly with Frobenii, and $\mathrm{id}_{\mathfrak{S}_{L_j}^{(1)}[E^{-1}]^{\wedge}_p}\otimes_{\mathfrak{S}^{(1)}[E^{-1}]^{\wedge}_p} f_{\mathrm{\acute{e}t}} = \mathrm{id}_{\mathfrak{S}_{L_j}^{(1)}[E^{-1}]^{\wedge}_p}\otimes_{\mathfrak{S}_{L_j}^{(1)}} f_{L_j}$ as maps of $\mathfrak{S}_{L_j}^{(1)}[E^{-1}]^{\wedge}_p$-modules. Choose a positive integer $r$ such that $\mathfrak{M}_{L_j}$ has $E$-height $\leq r$ for all $j = 1, \ldots, m$.

Consider the $\mathfrak{S}$-module
\[
\mathfrak{M} \coloneqq \bigcap_{j=1}^m (\mathcal{M}\cap \mathfrak{M}_{L_j})
\]
equipped with the induced Frobenius endomorphism $\varphi_\fkM$. Note that $\mathfrak{M}$ is $p$-adically complete since $\mathcal{M}$ and $\mathfrak{M}_{L_j}$ are $p$-adically complete.

\begin{prop} \label{prop:E-height}
The $\fkS$-module $\mathfrak{M}$ with Frobenius has $E$-height $\leq r$.
\end{prop}

\begin{proof}
The map $1\otimes\varphi\colon \varphi^*\mathfrak{M} \rightarrow \mathfrak{M}$ is injective since $1\otimes\varphi\colon \varphi^*\mathcal{M} \rightarrow \mathcal{M}$ is an isomorphism. Let $z \in \mathfrak{M}$. There exists a unique element $w \in \varphi^*\mathcal{M} \cong \mathfrak{S}\otimes_{\varphi, \mathfrak{S}}\mathcal{M}$ such that $(1\otimes\varphi)(w) = E(u)^r z$. For each $1 \leq j \leq m$, since $\mathfrak{M}_{L_j}$ has $E$-height $\leq r$, there exists a unique element $w_j \in \mathfrak{S}_{L_j}\otimes_{\varphi, \mathfrak{S}_{L_j}} \mathfrak{M}_{L_j}$ such that $(1\otimes\varphi)(w_j) = E(u)^r z$. 

Note that $\mathfrak{S}/p\mathfrak{S}$ has a $p$-basis given by $\{T_1, \ldots, T_d\}$, which also gives a $p$-basis for $\mathfrak{S}_{L_j}/p\mathfrak{S}_{L_j}$. Thus, we have a natural isomorphism $\mathfrak{S}_{L_j}\otimes_{\varphi, \mathfrak{S}_{L_j}} \mathfrak{M}_{L_j} \cong \mathfrak{S}\otimes_{\varphi, \mathfrak{S}} \mathfrak{M}_{L_j}$. Since $\varphi\colon \mathfrak{S} \rightarrow \mathfrak{S}$ is flat, we have $w = w_j$ and
\[
w \in \mathfrak{S}\otimes_{\varphi, \mathfrak{S}} \bigcap_{j=1}^m (\mathcal{M}\cap \mathfrak{M}_{L_j}) = \varphi^*\mathfrak{M}
\]
by uniqueness and Lemma~\ref{lem:intersection-modules-flat-base-change}(1).
\end{proof}

For each $n \geq 1$, set $\mathfrak{S}_{n} \coloneqq \mathfrak{S}/p^n \mathfrak{S}$, $\mathfrak{S}_{L_j, n} \coloneqq \mathfrak{S}_{L_j}/p^n \mathfrak{S}_{L_j}$, $\mathcal{M}_n \coloneqq \mathcal{M}/p^n \mathcal{M}$, and $\mathfrak{M}_{L_j, n} \coloneqq \mathfrak{M}_{L_j}/p^n \mathfrak{M}_{L_j}$. Consider the $\mathfrak{S}_n$-module
\[
\mathfrak{M}_{(n)} \coloneqq \bigcap_{j=1}^m (\mathcal{M}_n \cap \mathfrak{M}_{L_j, n})
\]
equipped with the induced $\varphi$.

\begin{lem} \label{lem:torsion-finitely-generated}
The module $\mathfrak{M}_{(n)}$ is finitely generated over $\mathfrak{S}_{n}$.
\end{lem}

\begin{proof}
Let $l$ be the $\mathbf{Z}_p$-rank of $\mathbb{L}$. Then $\mathcal{M}_n$ is projective over $\mathfrak{S}_n [u^{-1}]$ of rank $l$, so there exists a non-zero divisor $g \in \mathfrak{S}_n$ such that $\mathcal{M}_n [g^{-1}]$ is free of rank $l$ over $\mathfrak{S}_n [u^{-1}][g^{-1}]$. Choose a basis $\{e_1, \ldots, e_{l}\}$ of $\mathcal{M}_n[g^{-1}]$ over $\mathfrak{S}_n[u^{-1}][g^{-1}]$ such that the $\mathfrak{S}_n$-submodule $\mathfrak{N}$ of $\mathcal{M}_n[g^{-1}]$ generated by $e_1, \ldots, e_{l}$ satisfies $\mathcal{M}_n \subset \mathfrak{N}[u^{-1}]$ as $\mathfrak{S}_n[u^{-1}]$-modules. It suffices to show $\mathfrak{M}_{(n)} \subset \frac{1}{(ug)^h}\cdot \mathfrak{N}$ as $\mathfrak{S}_n$-modules for some $h\in\N$.

For $1 \leq j \leq m$, let $\{f_{j1}, \ldots, f_{jl}\}$ be a basis of $\mathfrak{M}_{L_j, n}$ over $\mathfrak{S}_{L_j, n}$. Then
\[
(e_1, \ldots, e_{l})^t = A_j\cdot (f_{j1}, \ldots, f_{jl})^t
\]
for some invertible matrix $A_j$ with entries in $\mathfrak{S}_{L_j, n}[u^{-1}][g^{-1}]$. Let $x = b_1 e_1+\cdots+b_{l} e_{l} \in \mathfrak{N}[u^{-1}]$ with $b_1, \ldots, b_{l} \in \mathfrak{S}_n[u^{-1}]$. We have $x \in \mathfrak{M}_{L_j, n}$ if and only if $(b_1, \ldots, b_{l}) = (c_{j1}, \ldots, c_{jl})A_j^{-1}$ for some $c_{j1}, \ldots, c_{jl} \in \mathfrak{S}_{L_j, n}$. Thus, by Lemma~\ref{lem:intersection-basic-rings}, 
\[
\bigcap_{j=1}^m (\mathfrak{N}[u^{-1}] \cap \mathfrak{M}_{L_j, n}) \subset \frac{1}{(ug)^h}\cdot \mathfrak{N}
\]
as $\mathfrak{S}_n$-modules for any integer $h$ such that $(ug)^hA_j^{-1}$ has entries in $\mathfrak{S}_{L_j, n}$.
\end{proof}

\begin{prop} \label{prop:kisin-mod-finitely-generated-saturated}
The $\mathfrak{S}$-module $\mathfrak{M}$ is finitely generated and satisfies 
\[
\mathfrak{M}[p^{-1}] \cap \mathfrak{M}[E^{-1}]^{\wedge}_p = \mathfrak{M}.
\]
\end{prop}

\begin{proof}
Consider the first statement. Since $\mathfrak{M}$ is $p$-adically complete, it suffices to show that $\mathfrak{M}/p\mathfrak{M}$ is finitely generated over $\mathfrak{S}$. For $1 \leq j \leq m$, we have
\[
p\mathcal{M} \cap \mathfrak{M}_{L_j} \subset p(\mathfrak{M}_{L_j}[u^{-1}]^{\wedge}_p) \cap \mathfrak{M}_{L_j} = p\mathfrak{M}_{L_j},
\]
since $p(\mathfrak{S}_{L_j}[u^{-1}]^{\wedge}_p) \cap \mathfrak{S}_{L_j} = p\mathfrak{S}_{L_j}$ and $\mathfrak{M}_{L_j}$ is flat over $\mathfrak{S}_{L_j}$. Thus,
\[
p\mathcal{M}\cap \mathfrak{M}  = \bigcap_{j=1}^m (p\mathcal{M} \cap \mathfrak{M}_{L_j}) = \bigcap_{j=1}^m (p\mathcal{M} \cap p\mathfrak{M}_{L_j}) = p\mathfrak{M},    
\]
where the last equality follows from the fact that $\mathfrak{M}_{L_j}[u^{-1}]^{\wedge}_p$ is $p$-torsion free. So the natural map $\mathfrak{M}/p\mathfrak{M} \rightarrow \mathcal{M}/p\mathcal{M}$ is injective. This map factors through the natural map $\mathfrak{M}/p\mathfrak{M} \rightarrow \bigcap_{j=1}^m (\mathcal{M}/p \mathcal{M}) \cap (\mathfrak{M}_{L_j}/p\mathfrak{M}_{L_j})$, which therefore is injective. Thus, $\mathfrak{M}$ is finitely generated over $\mathfrak{S}$ by Lemma~\ref{lem:torsion-finitely-generated}.

For the second statement, note 
\[
\mathfrak{M}[p^{-1}] \cap \mathfrak{M}[E^{-1}] \subset \mathcal{M}[p^{-1}] \cap \mathcal{M}[E^{-1}] = \mathcal{M}
\]
and
\[
\mathfrak{M}[p^{-1}] \cap \mathfrak{M}[E^{-1}] \subset \mathfrak{M}_{L_j}[p^{-1}] \cap \mathfrak{M}_{L_j}[E^{-1}] = \mathfrak{M}_{L_j}
\]
for each $j = 1, \ldots, m$. Hence, 
\[
\mathfrak{M}[p^{-1}] \cap \mathfrak{M}[E^{-1}] \subset \bigcap_{j=1}^m (\mathcal{M}\cap \mathfrak{M}_{L_j}) = \mathfrak{M}.
\]
By \cite[Rem.~3.15(i)]{du-liu-moon-shimizu-completed-prismatic-F-crystal-loc-system}, this implies $\mathfrak{M}[p^{-1}] \cap \mathfrak{M}[E^{-1}]^{\wedge}_p = \mathfrak{M}$.
\end{proof}

\begin{lem}
The natural map $\mathfrak{M} \rightarrow \varprojlim_n \mathfrak{M}_{(n)}$ is an isomorphism.
\end{lem}

\begin{proof}
Similarly as in the proof of Proposition~\ref{prop:kisin-mod-finitely-generated-saturated}, we see that the map $\mathfrak{M}/p^n\mathfrak{M} \rightarrow \mathfrak{M}_{(n)}$ is injective for each $n \geq 1$. Thus, the assertion follows from the same argument as in \cite[Lem.~4.22, Pf.]{du-liu-moon-shimizu-completed-prismatic-F-crystal-loc-system}.
\end{proof}

\begin{lem}
We have natural isomorphisms $\mathcal{O}_{\mathcal{E}}\otimes_{\mathfrak{S}} \mathfrak{M}_{(n)}  \cong \mathcal{M}_n$ and $\mathfrak{S}_{L_j}\otimes_{\mathfrak{S}} \mathfrak{M}_{(n)} \cong \mathfrak{M}_{L_j, n}$ for $j = 1, \ldots, m$. Furthermore, $\mathfrak{M}_{(n)}$ has $E$-height $\leq r$.
\end{lem}

\begin{proof}
The first isomorphism is given by
\[
\mathcal{O}_{\mathcal{E}}\otimes_{\mathfrak{S}}\mathfrak{M}_{(n)} \cong \mathfrak{M}_{(n)}[u^{-1}] \cong \bigcap_{i=1}^m (\mathcal{M}_n \cap \mathfrak{M}_{L_i, n}[u^{-1}]) = \mathcal{M}_n.
\]
For the second isomorphism, note that $u = T_1\cdots T_m$ is invertible in $\mathfrak{S}_{L_i}\otimes_{\mathfrak{S}} \mathfrak{S}_{L_j} / (p^n)$ if $i \neq j$. So for any $i = 2, \ldots, m$, we have natural isomorphisms
\[
\mathfrak{S}_{L_1}\otimes_{\mathfrak{S}} \mathfrak{M}_{L_i, n} \cong (\mathfrak{S}_{L_1}\otimes_{\mathfrak{S}} \mathfrak{S}_{L_i})\otimes_{\mathfrak{S}} \mathcal{M}_n \cong \mathfrak{S}_{L_i}\otimes_{\mathfrak{S}} \mathfrak{M}_{L_1, n}.
\]
Since $\mathfrak{S} \rightarrow \mathfrak{S}_{L_1}$ is flat and $\mathfrak{M}_{L_1}$ is free over $\mathfrak{S}_{L_1}$, we deduce
\[
\mathfrak{S}_{L_1}\otimes_{\mathfrak{S}} \mathfrak{M}_{(n)} \cong \bigcap_{i=1}^m (\mathfrak{S}_n[u^{-1}] \cap \mathfrak{S}_{L_i, n}) \otimes_{\mathfrak{S}_n} \mathfrak{M}_{L_1, n} \cong \mathfrak{M}_{L_1, n}
\]
by Lemmas~\ref{lem:intersection-basic-rings} and \ref{lem:intersection-modules-flat-base-change}(1). A similar argument yields the second isomorphism for any $1 \leq j \leq m$.

Lastly, $\mathfrak{M}_{(n)}$ has $E$-height $\leq r$ essentially by the same argument as in the proof of Proposition~\ref{prop:E-height}.
\end{proof}

Let $\mathscr{A}_n$ denote the set of $\mathfrak{S}$-submodules $\mathfrak{N}$ of $\mathfrak{M}_{(n)}$ which are stable under $\varphi$, have $E$-height $\leq r$, and satisfy $\mathfrak{N}[u^{-1}] = \mathcal{M}_n$. Consider
\[
\mathfrak{M}_{(n)}^{\circ} \coloneqq \bigcap_{\mathfrak{N} \in \mathscr{A}_n} \mathfrak{N} \subset \mathfrak{M}_{(n)}.
\]
For any positive integers $n' > n$, let $q_{n',n}\colon \mathfrak{M}_{(n')} \rightarrow \mathfrak{M} _{(n)}$ denote the restriction of the reduction-modulo-$p^n$ map $\mathcal{M}_{n'} \rightarrow \mathcal{M}_n$.

\begin{lem} \label{lem:intersection-A_n}
We have $\mathfrak{M}_{(n)}^{\circ} \in \mathscr{A}_n$ and $\mathfrak{M}_{(n)}^{\circ} \subset q_{n+1, n}(\mathfrak{M}_{(n+1)}^{\circ})$.
\end{lem}

\begin{proof}
This follows from the same argument as in \cite[Lem.~4.25, Pf.]{du-liu-moon-shimizu-completed-prismatic-F-crystal-loc-system}.
\end{proof}

\begin{prop} \label{prop:kisin-mod-base-change-compatible}
The natural maps $\mathcal{O}_{\mathcal{E}}\otimes_{\mathfrak{S}} \mathfrak{M} \rightarrow \mathcal{M}$ and $\mathfrak{S}_{L_j}\otimes_{\mathfrak{S}} \mathfrak{M} \rightarrow \mathfrak{M}_{L_j}$ for $1 \leq j \leq m$ are isomorphisms. 
\end{prop}

\begin{proof}
Given Lemma~\ref{lem:intersection-A_n} and Proposition~\ref{prop:E-height}, this follows from essentially the same argument as in the proof of \cite[Prop.~4.26]{du-liu-moon-shimizu-completed-prismatic-F-crystal-loc-system}.
\end{proof}

\begin{prop} \label{prop:kisin-mod-invert-p-projectivity}
The module $\mathfrak{M}[p^{-1}]$ is projective over $\mathfrak{S}[p^{-1}]$.
\end{prop}

\begin{proof}
Note that $\mathfrak{S}$ is a regular ring. Let $\mathfrak{n} \subset \mathfrak{S}$ be any maximal ideal. We have $p \in \mathfrak{n} \smallsetminus \mathfrak{n}^2$, and $u = T_1\cdots T_m \in \mathfrak{n}$. Denote by $\mathfrak{S}_{\mathfrak{n}}^{\wedge}$ the $\mathfrak{n}$-adic completion of $\mathfrak{S}_{\mathfrak{n}}$ and equip $\mathfrak{M}_{\mathfrak{n}} \coloneqq \mathfrak{S}_{\mathfrak{n}}^{\wedge}\otimes_{\mathfrak{S}} \mathfrak{M}$ with the tensor product Frobenius. By the same argument as in the last two paragraphs of \cite[Prop.~4.13, Pf.]{du-liu-moon-shimizu-completed-prismatic-F-crystal-loc-system}, it suffices to show that $\mathfrak{M}_{\mathfrak{n}}[p^{-1}]$ is projective over $\mathfrak{S}_{\mathfrak{n}}^{\wedge}[p^{-1}]$. 

Suppose first that there exists precisely one $i \in \{1, \ldots, m\}$ such that $T_i \in \mathfrak{n}$. Then $u \notin \mathfrak{n}^2$, and $\{p, u\}$ can be extended to a minimal set of generators of $\mathfrak{n}$. So
\[
\mathfrak{S}_{\mathfrak{n}}^{\wedge} \cong W(\kappa_{\mathfrak{n}})[\![s_1, \ldots, s_{d-1}]\!][\![u]\!]
\]
for some $s_1, \ldots, s_{d-1} \in \mathfrak{n}$ and a finite extension $\kappa_{\mathfrak{n}}$ of $k$. The Frobenius endomorphism on $\mathfrak{S}$ naturally extends to that on $\mathfrak{S}_{\mathfrak{n}}^{\wedge}$, and $W(\kappa_{\mathfrak{n}}) \subset \mathfrak{S}_{\mathfrak{n}}^{\wedge} \cong W(\kappa_{\mathfrak{n}})[\![s_1, \ldots, s_{d-1}]\!][\![u]\!]$ is stable under $\varphi_{\mathfrak{S}_{\mathfrak{n}}^{\wedge}}$ by the universal property of $p$-typical Witt vectors. Furthermore, we can choose $s_1, \ldots, s_{d-1}$ such that $W(\kappa_{\mathfrak{n}})[\![s_1, \ldots, s_{d-1}]\!] \subset W(\kappa_{\mathfrak{n}})[\![s_1, \ldots, s_{d-1}]\!][\![u]\!]$ is stable under $\varphi_{\mathfrak{S}_{\mathfrak{n}}^{\wedge}}$: first prove this directly for $\fkS_{R^0}$ by writing $\mathfrak{S}_{\mathfrak{n}}^{\wedge}$ of the form $(W(k)\langle T_1,\ldots,T_d^{\pm 1}\rangle[\![ T_i]\!])^\wedge_{\fkn'}$ and then use the fact that $\mathfrak{S}_{R^0} \rightarrow \mathfrak{S}$ is an \'etale map of the same dimension $d+1$. By Proposition~\ref{prop:E-height} and \cite[Lem.~4.11, 4.12]{du-liu-moon-shimizu-completed-prismatic-F-crystal-loc-system}, $\mathfrak{M}_{\mathfrak{n}}[p^{-1}]$ is projective over $\mathfrak{S}_{\mathfrak{n}}^{\wedge}[p^{-1}]$.

Next we deal with the remaining case. Without loss of generality, we may assume $T_1, \ldots, T_a \in \mathfrak{n}$ and $T_{a+1}, \ldots, T_m \notin \mathfrak{n}$ for some $2 \leq a \leq m$. Then $T_1, \ldots, T_a \in \mathfrak{n} \smallsetminus \mathfrak{n}^2$, and $\{p, T_1, \ldots, T_a\}$ can be extended to a minimal set of generators of $\mathfrak{n}$. Thus,
\[
\mathfrak{S}_{\mathfrak{n}}^{\wedge} \cong W(\kappa_{\mathfrak{n}})[\![s_1, \ldots, s_{d-a}]\!][\![T_1, \ldots, T_a]\!]
\]
for some $s_1, \ldots, s_{d-a} \in \mathfrak{n}$ and a finite extension $\kappa_{\mathfrak{n}}$ of $k$. Similarly as above, $W(\kappa_n) \subset \mathfrak{S}_{\mathfrak{n}}^{\wedge} \cong W(\kappa_{\mathfrak{n}})[\![s_1, \ldots, s_{d-a}]\!][\![T_1, \ldots, T_a]\!]$ is stable under $\varphi_{\mathfrak{S}_{\mathfrak{n}}^{\wedge}}$, and we can choose $s_1, \ldots s_{d-a}$ such that $W(\kappa_{\mathfrak{n}})[\![s_1, \ldots, s_{d-a}]\!] \subset W(\kappa_{\mathfrak{n}})[\![s_1, \ldots, s_{d-a}]\!][\![T_1, \ldots, T_a]\!]$ is stable under $\varphi_{\mathfrak{S}_{\mathfrak{n}}^{\wedge}}$. By \cite[Lem.~4.11]{du-liu-moon-shimizu-completed-prismatic-F-crystal-loc-system}, we may assume that $\varphi_{\mathfrak{S}_{\mathfrak{n}}^{\wedge}}(s_i)$ has zero constant term for each $i = 1, \ldots, d-a$.

Let $J$ be the non-zero Fitting ideal of $\mathfrak{M}_{\mathfrak{n}}$ of the smallest index, and suppose $J\mathfrak{S}_{\mathfrak{n}}^{\wedge}[p^{-1}] \neq \mathfrak{S}_{\mathfrak{n}}^{\wedge}[p^{-1}]$. Since Fitting ideals are compatible under base change, Proposition~\ref{prop:E-height} gives
\begin{equation} \label{eq:Fitting-ideal}
J\mathfrak{S}_{\mathfrak{n}}^{\wedge}[E(u)^{-1}] = \varphi_{\mathfrak{S}_{\mathfrak{n}}^{\wedge}}(J)\mathfrak{S}_{\mathfrak{n}}^{\wedge}[E(u)^{-1}],
\end{equation}
and so
\begin{equation} \label{eq:Fitting-ideal-2}
(\mathfrak{S}_{\mathfrak{n}}^{\wedge}/J)[E(u)^{-1}] = (\mathfrak{S}_{\mathfrak{n}}^{\wedge}/\varphi_{\mathfrak{S}_{\mathfrak{n}}^{\wedge}}(J))[E(u)^{-1}].    
\end{equation}

Write $K_1 = W(\kappa_{\mathfrak{n}})[p^{-1}]$ and let $\lvert \cdot\rvert$ denote the $p$-adic norm with $|p| = p^{-1}$. Consider the rigid analytic open unit ball $B$ in coordinates $(s_1, \ldots, s_{d-a}, T_1, \ldots, T_a)$: the set of $\overline{K_1}$-valued points of $B$ is given by
\[
\{(s_1, \ldots, s_{d-a}, T_1, \ldots, T_a) \in \overline{K_1}^d ~|~ 0 \leq |s_1|, \ldots, |s_{d-a}|, |T_1|, \ldots, |T_a| < 1 \}.
\]
We have a natural map $\mathfrak{S}_{\mathfrak{n}}^{\wedge}[p^{-1}] \rightarrow \mathcal{O}_B(B)$ whose image is dense, and a functorial bijection between the set of maximal ideals of $\mathfrak{S}_{\mathfrak{n}}^{\wedge}[p^{-1}]$ and the points of $B$ by \cite[Lem.~7.1.9]{deJong-dieudonnemodule}. The Frobenius $\varphi_{\mathfrak{S}_{\mathfrak{n}}^{\wedge}}$ on $\mathfrak{S}_{\mathfrak{n}}^{\wedge}$ induces an endomorphism on $B$. Set
\[
V_{|\pi|} \coloneqq \{(x_1, \ldots, x_{d-a}, y_1, \ldots, y_a) \in \mathbf{R}^d ~|~ 0 \leq x_i, y_j < 1, ~~y_1\cdots y_a = |\pi| \}.
\]

Consider the $\overline{K_1}$-valued points in $\operatorname{Spec}(\mathfrak{S}_{\mathfrak{n}}^{\wedge}[p^{-1}]/J)$, and let $Z = \{(|s_1|, \ldots, |s_{d-a}|, |T_1|, \ldots, |T_a|)\}$ be the set of corresponding $d$-tuple norms. Write
\[
Z' \coloneqq \{ (|s_1|, \ldots, |s_{d-a}|, |T_1|, \ldots, |T_a|) \in \mathbf{R}^d ~|~ 0 \leq |s_i|, |T_j| < 1, ~~(|\varphi_{\mathfrak{S}_{\mathfrak{n}}^{\wedge}}(s_i)|, |T_j|^p) \in Z \}.
\]
By \eqref{eq:Fitting-ideal-2}, we have $Z-V_{|\pi|} = Z'-V_{|\pi|}$. Furthermore, for $s_i, s_i' \in \overline{K}$ ($1 \leq i \leq d-a$) with $0 \leq |s_i|, |s_i'| < 1$ and $\varphi_{\mathfrak{S}_{\mathfrak{n}}^{\wedge}}(s_i) = s_i'$, we have
\[
\max_{1 \leq i \leq d-a} \{|s_i'|\} \leq \max_{1 \leq i \leq d-a} \{|s_i|^p, p^{-1}|s_i| \}   
\]
since $\varphi_{\mathfrak{S}_{\mathfrak{n}}^{\wedge}}(s_i)$ has zero constant term.

Now, observe $|u| = |T_1\cdots T_a|$ since $|T_i| = 1$ for any $a+1 \leq i \leq m$. Thus, using $Z-V_{|\pi|} = Z'-V_{|\pi|}$ and the above inequality, we deduce from the same argument as in \cite[Lem.~4.12, Pf.]{du-liu-moon-shimizu-completed-prismatic-F-crystal-loc-system} that $Z$ contains a point with $|T_1 \cdots T_a| < |\pi|$ and thus $(0, \ldots, 0) \in Z$, i.e., 
\[
J\mathfrak{S}_{\mathfrak{n}}^{\wedge}[p^{-1}] \subset (s_1, \ldots, s_{d-a}, T_1, \ldots, T_a)\mathfrak{S}_{\mathfrak{n}}^{\wedge}[p^{-1}].
\]

We claim $J\mathfrak{S}_{\mathfrak{n}}^{\wedge}[p^{-1}] \subset I\mathfrak{S}_{\mathfrak{n}}^{\wedge}[p^{-1}]$ for $I = (s_1, \ldots, s_{d-a})$. To see this, consider the $\varphi$-equivariant projection $\mathfrak{S}_{\mathfrak{n}}^{\wedge} \rightarrow \mathfrak{S}_1 \coloneqq \mathfrak{S}_{\mathfrak{n}}^{\wedge} / I \cong W(\kappa_n)[\![T_1, \ldots, T_a]\!]$, and let $J_1 \subset \mathfrak{S}_1$ be the image of $J$. Note $J_1\mathfrak{S}_1[p^{-1}] \subset (T_1, \ldots, T_a)\mathfrak{S}_1[p^{-1}]$. By \eqref{eq:Fitting-ideal}, we have
\[
E(u)^s J_1\mathfrak{S}_1[p^{-1}] \subset (T_1^p, \ldots, T_a^p)\mathfrak{S}_1[p^{-1}]
\]
for some $s \geq 0$. Since $\mathfrak{S}_1[p^{-1}]/(T_1^p, \ldots, T_a^p)$ is $E(u)$-torsion free, we obtain
\[
J_1\mathfrak{S}_1[p^{-1}] \subset (T_1^p, \ldots, T_a^p)\mathfrak{S}_1[p^{-1}].
\]
By repeating a similar argument, we deduce $J_1\mathfrak{S}_1[p^{-1}] \subset (T_1^{p^m}, \ldots, T_a^{p^m})\mathfrak{S}_1[p^{-1}]$ for each $m \geq 1$. Since $\mathfrak{S}_1[p^{-1}]$ is $(T_1, \ldots, T_a)$-adically separated, we get $J_1 = (0)$, i.e., $J\mathfrak{S}_{\mathfrak{n}}^{\wedge}[p^{-1}] \subset I\mathfrak{S}_{\mathfrak{n}}^{\wedge}[p^{-1}]$. 

Finally, consider the $\varphi$-equivariant map $W(\kappa_n)[\![T_1, \ldots, T_a]\!] \rightarrow W(\kappa')[\![T_1]\!] = W(\kappa')[\![u]\!]$ where $\kappa' \coloneqq \mathrm{Frac}(\kappa_n[\![T_2, \ldots, T_a]\!])^{\mathrm{perf}}$, induced by the natural map $\kappa_n[\![T_2, \ldots, T_a]\!] \rightarrow \kappa'$ (note $a \geq 2$). This induces a $\varphi$-equivariant map $\mathfrak{S}_{\mathfrak{n}}^{\wedge} \rightarrow \mathfrak{S}' \coloneqq W(\kappa')[\![s_1, \ldots, s_{d-a}]\!][\![u]\!]$. Note $J\mathfrak{S}'[p^{-1}] \neq 0, ~\mathfrak{S}'[p^{-1}]$. Since $J\mathfrak{S}'[p^{-1}]$ is a Fitting ideal of $\fkS'[p^{-1}]\otimes_{\fkS_{\fkn}^\wedge}\fkM_{\fkn}$, this contradicts the proof of \cite[Lem.~4.12]{du-liu-moon-shimizu-completed-prismatic-F-crystal-loc-system} applied to $\fkS'$. 
\end{proof}

Now, by the Beauville--Laszlo theorem and Propositions~\ref{prop:kisin-mod-finitely-generated-saturated}, \ref{prop:kisin-mod-base-change-compatible} and \ref{prop:kisin-mod-invert-p-projectivity},$\mathfrak{M}[E^{-1}]$ is projective over $\mathfrak{S}[E^{-1}]$. Furthermore, by a similar argument as in \cite[Lem.~4.10, Pf.]{du-liu-moon-shimizu-completed-prismatic-F-crystal-loc-system} using Corollary~\ref{cor:purity-intersection} and Proposition~\ref{prop:kisin-mod-finitely-generated-saturated}, we obtain from $f_{\mathrm{\acute{e}t}}$ and $f_{L_j}$'s an isomorphism of $\mathfrak{S}^{(1)}$-modules
\[
f\colon \mathfrak{S}^{(1)}\otimes_{p^1_1,\mathfrak{S}}\mathfrak{M}\xrightarrow{\cong}\mathfrak{S}^{(1)}\otimes_{p^1_2,\mathfrak{S}}\mathfrak{M}
\]
compatible with Frobenii and satisfying the cocycle condition over $\mathfrak{S}^{(2)}$. Hence, $(\mathfrak{M}, \varphi_{\mathfrak{M}}, f) \in \mathrm{DD}_{\mathfrak{S}}$, which by Lemma~\ref{lem:descentlemmamoduleversion} gives an analytic prismatic $F$-crystal over $(X, M_X)_{\Prism}$ extending $\mathscr{E}$. This completes the proof of Theorem~\ref{thm:main-purity}.

\appendix

\section{Complements on log formal schemes} \label{sec:log-formal-scheme}

For the convenience of the reader, we recall several definitions about log formal schemes. See \cite{Kato-log,Ogus-log} for log schemes and \cite[\S2]{Shiho-I} for $p$-adic log formal schemes.

Let $M_X$ be a log structure on $X$. A \emph{chart} of $M_X$ is a monoid homomorphism $M\rightarrow\Gamma(X,M_X)$ which induces an isomorphism $\underline{M}^a\xrightarrow{\cong}M_X$, where $\underline{M}$ is the constant sheaf of monoids on $X_\et$ defined by $M$ and $\underline{M}^a$ is the associated log structure (cf.~\cite[Def.~2.1.2]{Shiho-I}); for a property $P$ of monoids, if $M$ satisfies $P$, we say that the chart satisfies $P$. We say that $M_X$ is \emph{quasi-coherent} (resp.~\emph{coherent}) if it admits a chart (resp.~ a finitely generated chart) \'etale locally.

Recall that a monoid $Q$ is called integral if $q, q',q''\in Q$ and $q+q'=q+q''$ implies $q'=q''$; $Q$ is saturated if it is integral and if $q\in Q^\mathrm{gp}$ with $nq\in Q$ for some $n\geq 1$ implies $q\in Q$, where $Q^\mathrm{gp}$ denotes the group associated to $Q$.
A log structure $M_X$ on $X$ is called \emph{integral} (resp.~\emph{saturated}) if it is a sheaf of integral (resp.~saturated) monoids, or equivalently, $M_{X,\overline{x}}$ is integral (resp.~saturated) for every geometric point $x$ (by a formal scheme version of \cite[Prop.~II.2.1.3]{Ogus-log}). 
We say that $M_X$ is \emph{fine} if it is integral and coherent; $M_X$ is called \emph{fs} if it is fine and saturated. Note that $M_X$ is fine (resp.~fs) if and only if it admits a fine (resp.~fs) chart \'etale locally (cf.~\cite[Cor.~II.2.3.6]{Ogus-log}). 

A \emph{chart} of a morphism of log formal schemes $f\colon (X',M_{X'})\rightarrow (X,M_X)$ is a monoid homomorphism $\theta\colon Q\rightarrow Q'$ together with a commutative diagram
\begin{equation}\label{diag:chart}
\xymatrix{
Q'\ar[r]^-{\beta'}&\Gamma(X',M_{X'})\\
Q \ar[r]^-\beta\ar[u]^-\theta& \Gamma(X,M_X)\ar[u]_-{f^\flat},\\
}
\end{equation}
where the horizontal maps are charts. A chart exists \'etale locally if $M_X$ and $M_{X'}$ are both coherent (cf.~\cite[Prop.~II.2.4.2]{Ogus-log}). When $Q$ and $Q'$ are both saturated (resp. fine, resp. fs), we call the chart \eqref{diag:chart} \emph{saturated (resp. fine, resp. fs)}, which appears in Definition~\ref{def:absolute-prismatic-site}. 

\begin{lem}\label{lem:existence of saturated chart for finite free log structure}
    Let $f\colon (X',M_{X'})\rightarrow (X,M_X)$ be a morphism of $p$-adic log formal schemes. Assume that $M_X$ admits a finite free chart \'etale locally and $M_{X'}$ admits a saturated chart \'etale locally. Then $f$ admits a saturated chart of the form  \eqref{diag:chart} with $Q$ being finite free \'etale locally.
\end{lem}

\begin{proof}
    Since the assertion is of \'etale local nature, we may assume by shrinking $X$ and $X'$ that there exist charts $\beta\colon Q\rightarrow \Gamma(X,M_X)$ and $\alpha\colon P\rightarrow \Gamma(X',M_{X'})$ such that $Q$ is finite free with basis $e_1,\ldots,e_d$ for some $d\geq 0$ and such that $P$ is saturated. Consider $f^\flat(\beta(e_i))\in \Gamma(X',M_{X'})$. By localizing $X'$ (around each geometric point), we may further assume that for each $i$, there exist $q'_i\in P$ and $u'_i\in \Gamma(X',\calO_{X'}^\times)$ such that $f^\flat(\beta(e_i))=\alpha(q'_i) u'_i$ in $\Gamma(X',M_{X'})$. Observe that $\beta'\colon Q'\coloneqq P\oplus \Gamma(X',\calO_{X'}^\times)\rightarrow \Gamma(X',M_{X'})$ is a saturated chart of $M_X'$. Define a monoid homomorphism $\theta\colon Q\rightarrow Q'$ by $\theta(e_i)=(q_i',u_i')$. Now the triple $(\theta, \beta,\beta')$ forms a commutative diagram of the form \eqref{diag:chart} and gives a saturated chart for $f$.
\end{proof}

\section{Big absolute logarithmic crystalline site} \label{sec:log-crystalline-site}

This appendix briefly explains the absolute logarithmic crystalline site used in our work. 
We refer the reader to \cite{du-moon-shimizu-cris-pushforward} for the details as well as the relative logarithmic crystalline site.

In this appendix, we always consider the trivial log structure on $\Z_p$ and $\Z_p/p^n$.

\subsection{Generalities}
Let $(Y,M_Y)$ be a log scheme over $\Z_p$ such that $p^{n_0}\calO_Y=0$ for some $n_0\geq 1$ and $M_Y$ is integral and quasi-coherent. In the following, $n$ denotes an integer with $n\geq n_0$.

\begin{defn}
The \emph{big absolute logarithmic crystalline site} $((Y,M_Y)/(\Z_p/p^n))_\CRIS$ is defined as follows:
\begin{itemize}
    \item an object is a tuple $(U,T,M_T, f, i,\gamma)$ where $f\colon U\ra Y$ is a morphism of schemes, $(T,M_T)$ is a $\Z_p/p^n$-scheme with an integral and quasi-coherent log structure in which $p$ is nilpotent, $i$ is an exact closed immersion $(U,M_U\coloneqq f^\ast M_Y)\ra (T,M_T)$, and $\gamma$ is a PD-structure on the ideal of $\calO_T$ defining $U$ compatible with the standard PD-structure on the ideal $(p)$ of  $\Z_p/p^n$. Such a tuple is called \emph{affine} if $T$ (hence $U$) is affine;
    \item a morphism from $(U,T,M_T,i,\gamma)$ to $(U',T',M_{T'},i',\gamma')$ is a pair of morphisms of log schemes $g_T\colon (T,M_T)\rightarrow (T', M_{T'})$ and $g_U\colon (U,M_U)\rightarrow (U',M_{U'})$ such that $g_T\circ i=i'\circ g_U$, $f=f'\circ g_U$, and $g_T$ is compatible with the PD-structures $\gamma,\gamma'$. Note that $g_T$ is strict since $g_U$ is strict and all the log structures are integral. Such a morphism is called \emph{Cartesian} if the diagram
    \[
    \xymatrix{
    U \ar@{^(->}[r]^{i}\ar[d]_{g_U}&T\ar[d]^{g_T}\\ U'\ar@{^(->}[r]^{i'}&T'
    }
    \]
    is Cartesian, and it is called \emph{\'etale} if $g_T$ and $g_U$ are \'etale as morphisms of underlying schemes;
    \item a cover is a family of \'etale Cartesian morphisms that is jointly surjective.
\end{itemize}
By Lemma~\ref{lem:nonempty finite limit in abs log cristalline site}(2) below, $((Y,M_Y)/(\Z_p/p^n))_\CRIS$ is a site. 
To simplify the notation, we often write $(U,T)$ or $(U,T,M_T)$ for $(U,T,M_T,i,\gamma)$. Write $\calJ_T$ for the nil ideal sheaf of $T$ that defines $U$. 

Let $((Y,M_Y)/(\Z_p/p^n))_{\CRIS}^{\mathrm{aff}}\subset ((Y,M_Y)/(\Z_p/p^n))_\CRIS$ denote the full subcategory of affine objects with induced topology. By Lemma~\ref{lem:nonempty finite limit in abs log cristalline site}(2), it is also a site. The inclusion is a special cocontinuous functor in the sense of \cite[Tag~03CG]{stacks-project} and thus induces an equivalence of topoi 
\[
\operatorname{Sh}(((Y,M_Y)/(\Z_p/p^n))_{\CRIS}^{\mathrm{aff}})\xrightarrow{\cong} \operatorname{Sh}(((Y,M_Y)/(\Z_p/p^n))_\CRIS).
\]
\end{defn}

\begin{rem}\label{rem:fine log structure in absolute log crystalline site}
Assume that $M_Y$ is fine. Then for every $(U,T)\in ((Y,M_Y)/(\Z_p/p^n))_\CRIS$, $(T,M_T)$ is a fine log scheme. In fact, since $(U,M_U)\hookrightarrow (T,M_T)$ is an exact closed immersion with $M_T$ integral, it is a log thickening in the sense of \cite[Def.~IV.2.1.1]{Ogus-log}. Since $M_U$ is the pullback of $M_Y$ and thus fine, we conclude that $M_T$ is also fine by \cite[Prop.~IV.2.1.3.1]{Ogus-log}.
\end{rem}

\begin{lem}\label{lem:nonempty finite limit in abs log cristalline site}
\hfill
\begin{enumerate}
    \item The category $((Y,M_Y)/(\Z_p/p^n))_\CRIS$ has non-empty finite products. If $Y$ is affine, then the non-empty finite products of affine objects are also affine.
    \item For a diagram $(U_1,T_1)\rightarrow (V,Q) \leftarrow (U_2,T_2)$ in $((Y,M_Y)/(\Z_p/p^n))_\CRIS$, the fiber product exists. 
    Moreover, the fiber product is an affine object if $T_1$, $T_2$, and $Q$ are all affine.
\end{enumerate}
\end{lem}

\begin{proof}
See \cite[Lem.~3.3]{du-moon-shimizu-cris-pushforward}.
\end{proof}

\begin{defn}
For $(U,T)\in ((Y,M_Y)/(\Z_p/p^n))_\CRIS$, set
\[
\calO_{Y/(\Z_p/p^n)}(U,T)=\Gamma(T,\calO_T) \quad\text{and}\quad
\calJ_{Y/(\Z_p/p^n)}(U,T)=\Gamma(T,\calJ_T).
\]
These presheaves are indeed sheaves on $((Y,M_Y)/(\Z_p/p^n))_\CRIS$.
\end{defn}

\begin{defn}\label{def:big absolute log crystallien site}
Define \emph{the big absolute logarithmic crystalline site} 
\[
(Y,M_Y)_\CRIS\coloneqq ((Y,M_Y)/\Z_p)_\CRIS
\]
to be the colimit over $n$ of the sites $((Y,M_{Y})/(\Z_p/p^n))_\CRIS$ (see \cite[Tag~0EXI]{stacks-project} and \cite[1.12]{Beilinson-crystalline-period-map}): an object of $(Y,M_Y)_\CRIS$ is an object of $((Y,M_Y)/(\Z_p/p^n))_\CRIS$ for some $n\geq n_0$.
It has the structure sheaf $\calO_{Y/\Z_p}$ and an ideal sheaf $\calJ_{Y/\Z_p}$ defined in an obvious way. Similarly, we define $(Y,M_Y)_{\CRIS}^{\mathrm{aff}}$ when $Y$ is affine.
\end{defn}

Big absolute logarithmic crystalline sites enjoy functoriality:
let $f\colon (Y',M_{Y'})\ra (Y,M_Y)$ be a morphism of integral and quasi-coherent log schemes (over $\Z_p/p^{n_0}$ for some $n_0$).

\begin{construction}\label{construction:strict log structure for functoriality}
For $(U',T')\in(Y',M_{Y'})_\CRIS$, let $M_{(T',f)}$ be the \'etale sheaf of monoids on $T'$ defined as the fiber product $(f\circ g')^\ast M_Y \times_{M_{U'}}M_{T'}$ appearing in the diagram
\[
\xymatrix{
M_{(T',f)} \ar[r]\ar[d]   & M_{T'}\ar[d] \ar[r]& \calO_{T'} \ar[d]\\
(f\circ g')^\ast M_Y \ar[r]  &M_{U'} \ar[r]&\calO_{U'}
}
\]
under the equivalence $T'_\et\xrightarrow{\cong}U'_\et$ of \'etale sites.
The pair $(T',M_{(T',f)})$ is an integral and quasi-coherent log scheme with exact closed immersion $i'\colon (U',(f\circ g')^\ast M_Y) \hookrightarrow (T',M_{(T',f)})$ and defines an object $(U',T',i',M_{(T',f)},\gamma')\in (Y,M_Y)_\CRIS$. 
\end{construction}

\begin{prop}\label{prop:functoriality of log crystalline topoi}
There exists a unique morphism of topoi
\[
f_\CRIS=(f_\CRIS^\ast,f_{\CRIS,\ast})\colon \operatorname{Sh}((Y',M_{Y'})_\CRIS)\rightarrow \operatorname{Sh}((Y,M_Y)_\CRIS)
\]
such that $(f_\CRIS^\ast\calF)(U',T',i',M_{T'},\gamma')=\calF(U',T',i',M_{(T',f)},\gamma')$.
\end{prop}

\begin{proof}
See \cite[Prop.~4.4]{du-moon-shimizu-cris-pushforward}.
\end{proof}

Note $f^\ast_\CRIS\calO_{Y/\Z_p}=\calO_{Y'/\Z_p}$. Hence the inverse image functor for the morphism of ringed topoi $(\operatorname{Sh}((Y',M_{Y'})_\CRIS),\calO_{Y'/\Z_p})\rightarrow (\operatorname{Sh}((Y,M_Y)_\CRIS),\calO_{Y/\Z_p})$ is also computed by the above formula, which will be used implicitly when we consider the inverse image of crystals and isocrystals.

\begin{eg}\label{eg:Frobenius pullback on crystalline site}
Assume that $(Y,M_Y)$ is a log scheme over $\F_p$. Then it admits the absolute Frobenius morphism $F=(F_Y,F_{M_Y})\colon (Y,M_Y)\ra (Y,M_Y)$ where $F_Y$ is the absolute $p$-th Frobenius on the scheme $Y$ and, $F_{M_Y}$ is the multiplication-by-$p$ map. Hence we obtain a morphism of topoi $F_\CRIS\colon\operatorname{Sh}((Y,M_Y)_\CRIS)\rightarrow \operatorname{Sh}((Y,M_Y)_\CRIS)$.
\end{eg}

We now turn to isocrystals on the absolute logarithmic crystalline site. Analogous results in different set-ups are discussed in \cite[\S6]{Kato-log}, \cite[\S1.7]{Beilinson-crystalline-period-map}, \cite[\S7]{Ogus-convergenttopos}, \cite[\S4.3]{Shiho-II}, and \cite{drinfeld-stackycrystal}. 
For our purpose, we work on locally finite free objects in the isogeny category of \emph{finitely generated quasi-coherent} crystals. 

\begin{defn}[{cf.~\cite[\S3.2]{drinfeld-stackycrystal}}]
A \emph{coherent crystal}\footnote{We also call it a \emph{quasi-coherent $\calO_{Y/\Z_p}$-module of finite type} (cf.~\cite[Prop.~3.9]{du-moon-shimizu-cris-pushforward}).} on $(Y,M_Y)_\CRIS$ is a sheaf $\calF$ of $\calO_{Y/\Z_p}$-modules such that 
\begin{enumerate}
    \item for each $(U,T)\in (Y,M_Y)_\CRIS$, the induced Zariski sheaf $\calF_T$ on $T$ is a \emph{finitely generated quasi-coherent} $\calO_T$-module, and such that
    \item for each morphism $g\colon (U',T')\ra (U,T)$ in $(Y,M_Y)_\CRIS$, the induced map
    \[
    g^\ast\calF_T\coloneqq \calO_{T'}\otimes_{g^{-1}\calO_T}\calF_T\ra \calF_{T'}
    \]
    is an isomorphism.
\end{enumerate}
We write $\CR((Y,M_Y)_\CRIS)$ for the category of coherent crystals on $(Y,M_Y)_\CRIS$. 
Note that giving a coherent crystal on $(Y,M_Y)_\CRIS$ is equivalent to giving a sheaf of $\calO_{Y/\Z_p}$-modules on $(Y,M_Y)_{\CRIS}^{\mathrm{aff}}$ satisfying the same conditions (1) and (2) for objects $(Y,M_Y)_{\CRIS}^{\mathrm{aff}}$. It is also straightforward to see that the association $Y\mapsto \CR((Y,M_Y)_\CRIS)$ is a stack in the \'etale topology. 
\end{defn}

\begin{defn}\label{defn:crystallineisocrystal}
Define the category $\CR_\Q((Y,M_Y)_\CRIS)$ of \emph{isocrystals} on $(Y,M_Y)_\CRIS$ to be the isogeny category of $\CR((Y,M_Y)_\CRIS)$. Namely, objects are coherent crystals on $(Y,M_Y)_\CRIS$; for such $\calF\in \CR((Y,M_Y)_\CRIS)$, we write $\calF_\Q$ for the corresponding object in $\CR_\Q((Y,M_Y)_\CRIS)$. Morphisms are defined as
\[
\Mor_{\CR_\Q((Y,M_Y)_\CRIS)}(\calF_\Q,\calF'_\Q)\coloneqq 
\Mor_{\CR((Y,M_Y)_\CRIS)}(\calF,\calF')\otimes_\Z \Q.
\]
\end{defn}

\begin{rem}
Set $Y_1\coloneqq Y\otimes_{\Z_p} \Z_p/p$ and equip it with the induced log structure $M_{Y_1}$ from $M_Y$. Then the inverse image functor along $Y_1\hookrightarrow Y$ gives equivalences of categories
\[
\CR((Y,M_Y)_\CRIS)\xrightarrow{\cong}\CR((Y_1,M_{Y_1})_\CRIS)
\quad\text{and}\quad 
\CR_\Q((Y,M_Y)_\CRIS)\xrightarrow{\cong}\CR_\Q((Y_1,M_{Y_1})_\CRIS).
\]
The first equivalence is proved by a similar argument to \cite[Thm.~IV.1.4.1]{Berthelot-book}: the quasi-inverse functor is given by sending $\calF_1\in \CR((Y_1,M_{Y_1})_\CRIS)$ to the crystal $(U,T)\in (Y,M_Y)_\CRIS\mapsto \calF_1(U_1,T)$. The second equivalence follows from the first.
Combining this with $F_{Y_1,\CRIS}^\ast$ for $Y_1$ in Example~\ref{eg:Frobenius pullback on crystalline site}, one obtains the Frobenius pullback $F_\CRIS^\ast=F_{Y,\CRIS}^\ast$ on $\CR((Y,M_Y)_\CRIS)$ and $\CR_\Q((Y,M_Y)_\CRIS)$.
\end{rem}

We consider the evaluation of crystals and isocrystals at a more general $p$-adic log formal scheme: let $T$ be a $p$-adic formal scheme and set $T_n\coloneqq T\otimes_{\Z_p}\Z_p/p^n$ for $n\geq 1$. 

\begin{defn}[{cf.~\cite[\S2.3, 3.6]{drinfeld-stackycrystal}}]
A \emph{quasi-coherent $\calO_T$-module} is a sheaf $\calF$ of $\calO_T$-modules on the Zariski site of $T$ such that $\calF_n\coloneqq \calF/p^n\calF$ is a quasi-coherent $\calO_{T_n}$-module for each $n\geq 1$ and such that the natural map $\calF\rightarrow \varprojlim_n\calF_n$ is an isomorphism. We say that a quasi-coherent $\calO_T$-module $\calF$ is \emph{finitely generated} if each $\calF_n$ is finitely generated.
Let $\QCoh(T)$ denote the category of quasi-coherent $\calO_T$-modules and write $\QCoh^{\mathrm{fin}}(T)$ for the full subcategory consisting of finitely generated objects.

Let $\QCoh_\Q(T)$ (resp.~$\QCoh^\mathrm{fin}_{\Q}(T)$) denote the isogeny category of $\QCoh(T)$ (resp.~$\QCoh^{\mathrm{fin}}(T)$). An object of $\QCoh_\Q(T)$ can and will be regarded as a sheaf of $\calO_T[p^{-1}]$-modules on $T$ where $\calO_T[p^{-1}]\coloneqq \calO_T\otimes_\Z\Q$. We also call objects of $\QCoh^{\mathrm{fin}}_{\Q}(T)$ \emph{isocoherent sheaves} on $T$.
\end{defn}

\begin{rem}
A (finitely generated) quasi-coherent $\calO_T$-module in our sense corresponds to an adically quasi-coherent sheaf (of finite type) in the sense of \cite[Def.~I.3.1.3]{fujiwara-kato}. The category of finitely generated $\calO_T$-modules on the Zariski site of $T$ is equivalent to the one on the small \'etale site of $T$ (similarly defined), and thus we do not distinguish them (cf.~\cite[Prop.~I.6.2.12]{fujiwara-kato}). When $T=\Spf A$ is affine, the global section functor defines an equivalence of categories between the category of (finitely generated) quasi-coherent $\calO_T$-modules and that of $p$-adically complete (finitely generated) $A$-modules by \cite[Thm.~I.3.2.8]{fujiwara-kato}. Note that in both categories, morphisms are automatically continuous. 
\end{rem}

\begin{construction}\label{construction: eval of isocrystal at log PD-formal scheme}
Assume that $T$ is equipped with a $p$-adic log PD-formal scheme structure $(T,M_{T}, \calJ_{T})$ such that $M_T$ is integral and quasi-coherent. Set $(\overline{T},M_{\overline{T}})\coloneqq (V(\calJ_{T}),M_{T}|_{\overline{T}})$ and write $(T_n,M_{T_n},\calJ_{T_n})$ for the base change of $(T,M_{T}, \calJ_{T})$ to $\Z_p/p^n$.

Consider the big absolute logarithmic crystalline site $(Y,M_Y)_\CRIS$ as before and let $f\colon (\overline{T},M_{\overline{T}})\rightarrow (Y,M_Y)$ be a morphism of log schemes such that $M_{\overline{T}}\cong f^\ast M_Y$. Then $(\overline{T},T_n,M_{T_n})\in (Y,M_Y)_\CRIS$. By abuse of notation, we regard $(\overline{T},T,M_{T})$ as an ind-object $\varinjlim_n (\overline{T},T_n,M_{T_n})$ of $(Y,M_Y)_\CRIS$ and say that $T$ is a \emph{$p$-adic log PD-formal scheme in $(Y,M_Y)_\CRIS$}.

For $\calF\in \CR((Y,M_Y)_\CRIS)$, $\calF_{T_n}$ is a finitely generated quasi-coherent $\calO_{T_n}$-module with $\calF_{T_{n+1}}\otimes \Z_p/p^n\cong\calF_{T_n}$. Hence $\calF_T\coloneqq \varprojlim_n\calF_n\in \QCoh^{\mathrm{fin}}(T)$ with $\calF/p^n\calF\cong \calF_n$ by \cite[Prop.~I.3.4.1]{fujiwara-kato}.

For $\calF_\Q\in \CR_\Q((Y,M_Y)_\CRIS)$, define a sheaf $\calF_{\Q,T}$ of $\calO_{T}[p^{-1}]$-modules on $T$ by $\calF_{\Q,T}\coloneqq \calF_{T}\otimes_{\calO_{T}}\calO_{T}[p^{-1}]$; this is independent of the choice of a representative $\calF$ of $\calF_\Q$. We also set $\calF(T)\coloneqq \calF_{T}(T)$ and $\calF_{\Q}(T)\coloneqq \calF_{\Q,T}(T)$. If $T=\Spf A$ is affine, we often write $\calF_\Q(A)$ for $\calF_\Q(T)$.

Moreover, this construction is functorial: suppose there is another $p$-adic integral quasi-coherent log PD-formal scheme $(T',M_{T'},\calJ_{T'})$ together with a morphism $g\colon (T',M_{T'},\calJ_{T'})\ra (T,M_{T},\calJ_{T})$. Then the natural morphism
$\varprojlim_n (\calO_{T_n'}\otimes_{g_n^{-1}\calO_{T_n}}g_n^{-1}\calF_{T_n})\rightarrow \calF_{T'}$
is an isomorphism by the crystal property of $\calF$, where $g_n$ denotes the induced map $T_n'\rightarrow T_n$. In particular, if $T=\Spf A$ and $T'=\Spf B$ are affine, we have $(\calF(T)\widehat{\otimes}_AB)\otimes_BB[p^{-1}]\xrightarrow{\cong}\calF(T')\otimes_BB[p^{-1}]=\calF_\Q(T')$, where $\widehat{\otimes}$ denotes $p$-adically completed tensor product.
\end{construction}

We use the above construction to define a notion of finite locally free isocrystals, which works well for our purpose.

\begin{defn}
 An isocrystal $\calF_\Q$ on $(Y,M_Y)_\CRIS$ is said to be \emph{finite locally free} if for every $p$-adic log PD-formal scheme $T$ in $((Y,M_Y)/\Z_p)_\CRIS$ such that $T=\Spf A$ is affine and flat over $\Z_p$,  $\calF_{\Q}(T)$ is a finite projective $A[p^{-1}]$-module.
\end{defn}

\begin{rem}
\hfill
\begin{enumerate}
    \item 
Let $\calF_\Q$ be a finite locally free isocrystal on $(Y,M_Y)_\CRIS$. Then for every morphism $T'=\Spf B\rightarrow T=\Spf A$ between affine $\Z_p$-flat $p$-adic log PD-formal schemes in $(Y,M_Y)_\CRIS$, we have an isomorphism $\calF_\Q(T)\otimes_{A}B\xrightarrow{\cong}\calF_{\Q}(T')$ by Lemma~\ref{lem:DLMS3.24(iv)} below. 
    \item Finite local freeness is of \'etale local nature: let $\{Y_i\rightarrow Y\}$ be an \'etale cover. Then an isocrystal $\calF_\Q$ on $(Y,M_Y)_\CRIS$ is finite locally free if and only if  $\calF_\Q|_{Y_i}$ is finite locally free on $(Y_i,(M_Y)|_{Y_i})_\CRIS$ for every $i$. This follows from Lemma~\ref{lem:p fflat descent for f proj 1/p} below.
\end{enumerate}
\end{rem}

\begin{lem}\label{lem:DLMS3.24(iv)}
Let $A\rightarrow B$ be a map of $p$-adically complete rings with $A$ $p$-torsion free and let $M$ be a finitely generated $p$-adically complete $A$-module. If $M[p^{-1}]$ is a finite projective $A[p^{-1}]$-module, then the map
\[
M[p^{-1}]\otimes_AB\rightarrow (M\widehat{\otimes}_AB)\otimes_BB[p^{-1}]
\]
is an isomorphism. In particular, the latter module is finite projective over $B[p^{-1}]$.
\end{lem}

\begin{proof}
The proof of \cite[Lem.~3.24(iv)]{du-liu-moon-shimizu-completed-prismatic-F-crystal-loc-system} works in the current setup with obvious modification  (without assuming that $A$ is noetherian nor $M$ is $p$-torsion free).
\end{proof}

\begin{lem}[{\cite[Lem.~2.12]{drinfeld-stackycrystal}}]\label{lem:p fflat descent for f proj 1/p}
Let $A\rightarrow B$ be a $p$-adically complete faithfully flat map between $p$-adically complete $\Z_p$-flat rings. Let $M$ be a finitely generated $p$-adically complete $A$-module and set $M_B \coloneqq M\widehat{\otimes}_A B$. If $M_{B}[p^{-1}]$ is finite projective over $B[p^{-1}]$, then $M[p^{-1}]$ is finite projective over $A[p^{-1}]$.
\end{lem}

\begin{proof}
Choose a presentation $0 \to N \to A^{\oplus m} \to M \to 0$ of $M$; note that $N$ is a $p$-adically complete $p$-torsion free $A$-module. Since $A/p^n\rightarrow B/p^n$ is faithfully flat, applying $-\widehat{\otimes}_A B$ to this exact sequence gives an exact sequence of $B$-modules $0 \to N\widehat{\otimes}_A B \to B^{\oplus m} \to M_B \to 0$.

It also follows that $N\widehat{\otimes}_A B$ is $p$-adically complete and $p$-torsion free. After inverting $p$, the resulting exact sequence $0 \to (N\widehat{\otimes}_A B)[p^{-1}] \to B[p^{-1}]^{\oplus m} \to M_B[p^{-1}] \to 0$ splits since $M_B[p^{-1}]$ is projective. In particular, $(N\widehat{\otimes}_A B)[p^{-1}]$ is finite projective. Hence there is a $p$-adically complete finitely generated $B$-submodule $N'_B$ of $(N\widehat{\otimes}_A B)[p^{-1}]$ with $N'_B[p^{-1}] = (N\widehat{\otimes}_A B)[p^{-1}]$. By \cite[Prop.~2.7]{drinfeld-stackycrystal}, we see that $N[p^{-1}]$ is finite generated over $A[p^{-1}]$ and obtain a finite presentation
$A[p^{-1}]^{\oplus n}\xrightarrow{\varphi} A[p^{-1}]^{\oplus m} \to M[p^{-1}] \to 0$
of $M[p^{-1}]$, whose base change to $B[p^{-1}]$ gives the corresponding presentation $B[p^{-1}]^{\oplus n}\xrightarrow{\varphi_B} B[p^{-1}]^{\oplus m} \to M_B[p^{-1}] \to 0$
of the finite projective module $M_B[p^{-1}]$. Finally, we conclude that $M[p^{-1}]$ is finite projective over $A[p^{-1}]$ by \cite[Cor.~2.6]{drinfeld-stackycrystal}.
\end{proof}

\begin{defn}\label{def:F-isocrystals}
An \emph{$F$-isocrystal} on $(Y,M_Y)_\CRIS$ is a pair $(\calF_\Q, \varphi_{\calF_\Q})$ where $\calF_\Q$ is a \emph{finite locally free} isocrystal on $(Y,M_Y)_\CRIS$ arising from a coherent crystal and $\varphi_{\calF_\Q}\colon F_{\CRIS}^\ast \calF_{\Q}\xrightarrow{\cong}\calF_{\Q}$ is an isomorphism of isocrystals (see the above remark for $F_\CRIS^\ast$).  
Write $\Vect_\Q^\varphi((Y,M_{Y})_\CRIS)$ for the category of $F$-isocrystals over $(Y,M_Y)_\CRIS$.
\end{defn}

\begin{rem}\label{rem:F-isocrystalmoduloprestricttomodulopiissurjective}
Assume that $(Y, M_Y)$ is a quasi-compact log scheme over $\F_p$, and let $i\colon (Y_0,M_{Y_0})\hookrightarrow (Y,M_Y)$ be an exact closed immersion defined by a nilpotent ideal sheaf. Then the restriction $i_\CRIS^\ast\colon\CR_{\Q}((Y,M_Y)_\CRIS)\rightarrow\CR_{\Q}((Y_0,M_{Y_0})_\CRIS)$ is not an equivalence in general. However, it gives an equivalence for \emph{$F$-isocrystals} by Dwork’s trick:  choose $N\gg 0$ such the $N$-th iterated Frobenius $F_{Y}^N$ on $Y$ factor through $Y_0$. Namely, there is a morphism $\rho\colon (Y,M_Y) \to (Y_0,M_{Y_0})$ such that $F_{Y}^N=i\circ \rho$ and $F_{Y_0}^N=\rho\circ i$. Hence we obtain the commutative diagram
\[
\xymatrix{
\CR_{\Q}((Y,M_Y)_\CRIS)\ar[r]^-{F_{Y,\CRIS}^{N,\ast}}\ar[d]_-{i_\CRIS^\ast}&\CR_{\Q}((Y,M_Y)_\CRIS)\ar[d]^-{i_\CRIS^\ast} \\
\CR_{\Q}((Y_0,M_{Y_0})_\CRIS)\ar[r]^-{F_{Y_0,\CRIS}^{N,\ast}}\ar[ru]^-{\rho^{\ast}_\CRIS} &\CR_{\Q}((Y_0,M_{Y_0})_\CRIS).
}
\]
In particular, given an $F$-isocrystal $(\calF_\Q,\varphi_{\calF_\Q})$ on $(Y,M_Y)_\CRIS$, we have an isomorphism 
\[
\varphi_{\calF_\Q}^N\coloneqq\varphi_{\calF_\Q}\circ (F_{Y,\CRIS}^\ast\varphi_{\calF_\Q})\circ \cdots\circ(F_{Y,\CRIS}^{N-1,\ast}\varphi_{\calF_\Q})\colon \rho_\CRIS^{\ast}i_\CRIS^\ast(\calF_\Q,\varphi_{\calF_\Q})\xrightarrow{\cong}(\calF_\Q,\varphi_{\calF_\Q})
\]
of $F$-isocrystals.
This shows that  $i_\CRIS^\ast\colon\Vect_\Q^\varphi((Y,M_{Y})_\CRIS)\rightarrow\Vect_\Q^\varphi((Y_0,M_{Y_0})_\CRIS)$ is fully faithful. The essential surjectivity is proved similarly since the finite local freeness is preserved under $\rho_\CRIS^{\ast}$ by Construction~\ref{construction:strict log structure for functoriality}.
\end{rem}

\subsection{More on crystals and isocrystals}
We describe crystals and isocrystals on $Y$ when $(Y,M_Y)$ is embedded into a log smooth $p$-adic formal scheme: let $k$ be a perfect field of characteristic $p$. Set $W\coloneqq W(k)$ and $W_n\coloneqq W_n(k)$ for $n\in \N$. We always equip these rings with the trivial log structure. Assume that $(Y,M_Y)$ is a log scheme over $W_n$ and $M_Y$ is \emph{fine}.

Let $(Z,M_Z)$ be a $p$-adic fine log formal scheme that is log smooth over $W$, namely, the morphism of log schemes $(Z_n,M_{Z_n})\coloneqq (Z,M_Z)\times_WW_n\rightarrow (\Spec W_n,M_{\Spec W_n})$ is smooth for every $n$. 
Since $M_{\Spec W_n}$ is trivial, $(Z_n,M_{Z_n})\rightarrow (\Spec W_n,M_{\Spec W_n})$ is integral, and thus the underlying morphism of schemes is flat by \cite[Cor.~4.5]{Kato-log}.
We let 
\[
\omega^1_{(Z,M_Z)/\Z_p}\coloneqq \varprojlim_n \omega_{(Z_n,M_{Z_n})/(\Spec (\Z_p/p^n),M_{\Spec(\Z_p/p^n)})}^1
\]
denote the sheaf of formal log differentials (see \cite[Def.~2.1.14(2)]{Shiho-I}) and define $\omega^1_{(Z,M_Z)/W}$ similarly. Since $k$ is perfect, the maps $(Z,M_Z)\rightarrow (\Spf W,M_{\Spf W})\rightarrow (\Spf \Z_p, M_{\Spf \Z_p})$ induce an isomorphism $\omega_{(Z,M_Z)/\Z_p}^1\xrightarrow{\cong}\omega_{(Z,M_Z)/W}^1$.

Suppose that there is a closed immersion $(Y,M_Y)\hookrightarrow (Z,M_Z)$. 
For $i\geq 1$, let $(Z^{(i)},M_{Z^{(i)}})$ be the $(i+1)$-st self-product of $(Z,M_Z)$ in the category of $p$-adic fine log formal schemes. 
Let 
\[
(D,M_D)\coloneqq (D_{(Y,M_Y)}(Z,M_Z), M_D) 
\]
denote the $p$-adically completed log PD-envelope as in \cite[Def.~5.4]{Kato-log}. Similarly, let $(D^{(i)},M_{D^{(i)}})$ denote the $p$-adically completed log PD-envelope attached to the closed immersion $(Y,M_Y)\hookrightarrow (Z^{(i)},M_{Z^{(i)}})$. As in Construction~\ref{construction: eval of isocrystal at log PD-formal scheme}, set $D_n\coloneqq D\otimes_{\Z_p}\Z_p/p^n$ and $D_n^{(i)}\coloneqq D^{(i)}\otimes_{\Z_p}\Z_p/p^n$. It is straightforward to see $(Y,D_n,M_{D_n}), (Y,D_n^{(i)},M_{D_n^{(i)}})\in (Y,M_Y)_\CRIS$ and $(Y,D_n^{(i)},M_{D_n^{(i)}})$ represents $(i+1)$-st self-product of $(Y,D_n,M_{D_n})$ in $(Y,M_Y)_\CRIS$. We also regard $(Y,D)$ as an ind-object $\varinjlim_n (Y,D_n)$ of $(Y,M_Y)_{\CRIS}$.

\begin{lem}\label{lem:weakly-final}
With the notation as above, for any $(U,T) \in (Y,M_Y)_{\CRIS}^{\mathrm{aff}}$, there is a morphism $(U, T)\rightarrow (Y,D)$ in $(Y,M_Y)_{\CRIS}$. More generally, for any affine $p$-adic log PD-formal scheme $T$ in $(Y,M_Y)_\CRIS$ as in Construction~\ref{construction: eval of isocrystal at log PD-formal scheme}, there is a morphism $(\overline{T},T)\rightarrow (Y,D)$ of ind-objects in $(Y,M_Y)_\CRIS$.
\end{lem}

\begin{proof}
For the first statement, write $T=\Spec A$ and $U=\Spec(A/I)$ with $I$ a PD-ideal of $A$. Note that $p$ is nilpotent in $A$ and thus $I$ is nilpotent. Since the map $Z\rightarrow \Spf W$ is a log smooth morphism of $p$-adic fine log formal schemes, the map $\Spec (A/I)\rightarrow Y\rightarrow Z$ lifts to a map $\Spec A\rightarrow Z$ compatible with log structures by \cite[Prop.~IV.3.1.4.2]{Ogus-log}. Finally, the universal property of log PD-envelope induces a morphism $(U,T)=(\Spec (A/I), \Spec A)\rightarrow (Y,D)$. In general, if $T$ is an affine $p$-adic log formal scheme in $(Y,M_Y)_\CRIS$, the above log smoothness allows one to obtain morphisms $(\overline{T},T_n)\rightarrow (Y,D)$ compatibly when varying $n$, which yields a morphism $(\overline{T},T)\rightarrow (Y,D)$.
\end{proof}

\begin{eg}[Small affine case]\label{eg:BreuilSR}
Consider a $p$-adically completed \'etale map 
\[
\square\colon R^0\coloneqq \mathcal{O}_K \langle T_1, \ldots, T_m, T_{m+1}^{\pm 1}, \ldots, T_d^{\pm 1}\rangle / (T_1\cdots T_m - \pi) \rightarrow R
\]
together with the prelog structure $\mathbf{N}^d \rightarrow R^0\rightarrow R$ sending $e_i \mapsto T_i$ for $1 \leq i \leq d$. Set $(X,M_X)=(\Spf R, \N^d)^a$ and let $(Y,M_Y)\coloneqq (X_1,M_{X_1})$ denote the log scheme associated to the mod $p$ fiber. By \cite[04D1]{stacks-project}, the map $R^0/p\rightarrow R/p$ lifts to a $p$-adically completed \'etale morphism
\[
\widetilde{R}^0\coloneqq W(k)\langle T_1, \ldots, T_m, T_{m+1}^{\pm 1}, \ldots, T_d^{\pm 1}\rangle\rightarrow \widetilde{R}. 
\]
We set $Z\coloneqq \Spf \widetilde{R}$ and equip it with log structure $M_Z$ attached to $\N^d\rightarrow \widetilde{R}: e_i\mapsto T_i$. Then $(Y,M_Y)\hookrightarrow (Z,M_Z)$ satisfies the assumptions of this subsection. 

Let $S$ denote the $p$-adic completion of the (log) PD-envelope of $\widetilde{R}$ with respect to the kernel of the exact surjection $\widetilde{R}\rightarrow R/p$; then $D$ is identified with $\Spf S$. Furthermore, with the notation of Section~\ref{sec:Breuil--Kisin log prism}, the obvious map $\widetilde{R}^0\rightarrow \fkS_{R^0}$ extends to a map $\widetilde{R}\rightarrow \fkS_{R,\square}$, which identifies $S$ with the $p$-adic completion of the log PD-envelope of $\fkS_{R,\square}$ with respect to the kernel of the exact surjection $\fkS_{R,\square} \to R/p$. One can also check that the $W(k)$-algebra endomorphism $\varphi_{\widetilde{R}^0}$ on $\widetilde{R}^0$ given by $T_i\mapsto T_i^p$ uniquely extends to $\varphi_S\colon S\rightarrow S$, giving a lift of Frobenius on $S/p$.

For $D^{(i)}$, observe that the exactification of the closed immersion $(Y,M_Y)\hookrightarrow (Z^{(i)},M_{Z^{(i)}})$ is affine, which we write $\Spf \widetilde{R}^{(i)}$. Set $\widetilde{I}^{(i)}\coloneqq \Ker(\widetilde{R}^{(i)}\twoheadrightarrow R)$ and let $S^{(i),\mathrm{nc}}$ denote the PD-envelope of $\widetilde{R}^{(i)}$ with respect to $(p,\widetilde{I}^{(i)})$. Define $S^{(i)}$ to the $p$-adic completion of $S^{(i),\mathrm{nc}}$. By construction, $D^{(i)}$ is given by $\Spf S^{(i)}$. The endomorphism $\varphi_{\widetilde{R}^0}$ also yields a lift of Frobenius $\varphi_{S^{(i)}}$ on $S^{(i)}$.

Let us also discuss the relation with objects in Appendix~\ref{sec:log convergent site}: let $(X^{(1)},M_{X^{(1)}})$ denote the exactification of the closed immersion $(X,M_X)\rightarrow (X,M_X)\times_{(\Spf \calO_K,M_\mathrm{can})}(X,M_X)$ as a $p$-adic log formal scheme. For each $n \geq 1$, let $(X^{(1),n}, M_{X^{(1),n}})$ denote the $n$-th infinitesimal neighborhood of $(X, M_X)$ in $(X^{(1)}, M_{X^{(1)}})$. In this case, $X^{(1)}$ and $X^{(1),n}$ are all affine. So we write $X^{(1)}=\Spf R^{(1)}$ and $X^{(1),n}=\Spf R^{(1),n}$: more concretely, $R^{(1),n}=R^{(1)}/(I^{(1)})^{n+1}$ for $I^{(1)}\coloneqq \Ker(R^{(1)}\twoheadrightarrow R)$. 
Consider the two maps $S\rightarrow R\rightrightarrows R^{(1)}\rightarrow R^{(1),n}$. Since the induced two log structures on $X^{(1),n}$ from $D$ coincide, we obtain a map $\widetilde{R}^{(1)}\rightarrow R^{(1),n}$, which maps $\tilde{I}^{(1)}$ to $I^{(1)}$. This yields a map $S^{(1),\mathrm{nc}}\rightarrow R^{(1),n}[p^{-1}]$ and further extends to $S^{(1)}\rightarrow R^{(1),n}[p^{-1}]$ because $I^{(1)}$ is nilpotent in $R^{(1),n}[p^{-1}]$. We remark that when $n=1$, the map $S^{(1)}\rightarrow R^{(1),n}[p^{-1}]$ factors as $S^{(1)}\rightarrow R^{(1),1}\subset R^{(1),1}[p^{-1}]$ since $(I^{(1)})^2=0$ in $R^{(1),1}$; for a general $n$, if we take $m\geq 1$ such that $p^m\geq n+1$, then the composite $S^{(1)}\xrightarrow{\varphi^m_{S^{(1)}}}S^{(1)}\rightarrow R^{(1),n}[p^{-1}]$ factors through $R^{(1),n}$.
\end{eg}

Let us go back to the general set-up and keep the notation before Lemma~\ref{lem:weakly-final}.

\begin{defn}[{cf.~\cite[Def.~4.3.1]{Shiho-I}}]
Write $p_i\colon D^{(1)}\rightarrow D$ ($i=1,2$) for the $j$-th projection.
    An \emph{HPD-stratification} on an $\calO_D$-module $\calM$ on $D_\et$ is an isomorphism\footnote{The standard convention in the literature is $p_2^\ast \calM \xrightarrow{\cong}p_1^\ast\calM$, but we use the opposite convention.} 
    \[
    \eta\colon p_1^\ast \calM \xrightarrow{\cong}p_2^\ast \calM 
    \]
    on $D^{(1)}_\et$ that satisfies the cocycle condition on $D^{(2)}_\et$. Note that the cocycle condition implies that the pullback of $\eta$ along the diagonal $D\rightarrow D^{(1)}$ is $\mathrm{id}_{\calM}$; see Remark~\ref{rem:conditions on stratification} below.    
\end{defn}

\begin{rem}\label{rem:conditions on stratification}
Let us spell out the conditions in the definition of HPD-stratifications.
Consider the following projection and diagonal maps of $p$-adic fine log schemes induced by the universality:
\[
p_{12},p_{13}, p_{23}\colon D^{(2)} \rightarrow D^{(1)}, \;
q_1, q_2, q_3\colon D^{(2)}\rightarrow D, \;
\Delta^{(2)}\colon D\rightarrow D^{(2)}, \;\text{and}\;
\Delta\colon D\rightarrow D^{(1)}.
\]
Then the condition on the isomorphism $\eta\colon p_1^\ast \calM \xrightarrow{\cong}p_2^\ast \calM$ is 
\[
p_{13}^\ast\eta=p_{23}^\ast\eta \circ p_{12}^\ast\eta 
\]
under the identifications (that are unique up to a unique isomorphism) $p_{12}^\ast p_1^\ast \calM\cong q_1^\ast\calM \cong p_{13}^\ast p_1^\ast\calM$, and so on. Note that pulling back $p_{13}^\ast\eta=p_{23}^\ast\eta \circ p_{12}^\ast\eta$ along $\Delta^{(2)}$ yields
$\Delta^\ast\eta=\Delta^\ast\eta\circ \Delta^\ast\eta$. Since $\Delta^\ast\eta$ is an isomorphism as $\eta$ is, we conclude $\Delta^\ast\eta=\mathrm{id}_\calM$. 
\end{rem}

Finally, we review the notion of log connections. The map $d\colon \calO_Z\rightarrow \omega^1_{(Z,M_Z)/\Z_p}$ extends to $d\colon \calO_D\rightarrow \calO_D\otimes_{\calO_Z}\omega^1_{(Z,M_Z)/\Z_p}$ (e.g., by a log analogue of \cite[Cor.~6.3, Exer.~6.4]{berthelot-ogus-book} using \cite[Prop.~6.5]{Kato-log}). Recall that a \emph{log connection} on a quasi-coherent $\calO_D$-module $\calM$ is an additive map 
\[
\nabla\colon \calM\rightarrow\calM\otimes_{\calO_Z}\omega^1_{(Z,M_Z)/\Z_p}=\calM\otimes_{\calO_D}(\calO_D\otimes_{\calO_Z}\omega^1_{(Z,M_Z)/\Z_p})
\]
on $D_\et$ such that $\nabla(ae)=a\nabla(e)+e\otimes da$ for $a\in \calO_Z$ and $e\in \calM$. We call $\nabla$ \emph{integrable} if the composite $\calM\xrightarrow{\nabla}\calM\otimes_{\calO_Z}\omega^1_{(Z,M_Z)/\Z_p}\xrightarrow{\nabla}\calM\otimes_{\calO_Z}\omega^2_{(Z,M_Z)/\Z_p}$ is zero. We say that an integrable log connection $\nabla$ on $D_{\et}$ is \emph{quasi-nilpotent} if the pullback of $(\calM,\nabla)$  to $D_{n,\et}$ is quasi-nilpotent in the sense of \cite[Thm.~6.2(iii)]{Kato-log} and \cite[p.~19]{ogus-griffiths} for each $n$.

\begin{prop} \label{prop: crystals and quasi-nilpotent connections}
With the notation as above, the following categories are equivalent:
\begin{enumerate}
    \item the category $\CR((Y,M_Y)_\CRIS)$ of coherent crystals on $(Y,M_Y)_\CRIS$;
    \item the category of finitely generated quasi-coherent $\calO_D$-modules together with HPD-stratifications;
    \item the category of finitely generated quasi-coherent $\calO_D$-modules $\calM$ on $D_{\et}$ together with a quasi-nilpotent integrable log connection
    \[
    \nabla\colon \calM\ra \calM\otimes_{\calO_Z}\omega^1_{(Z,M_Z)/\Z_p}.
    \]
\end{enumerate}
\end{prop}

Recall that coherent crystals mean finitely generated quasi-coherent crystals.

\begin{proof}
This is essentially proved in \cite[Thm.~6.2]{Kato-log}. See also \cite[Prop.~6.8]{du-moon-shimizu-cris-pushforward}.
\end{proof}

By passing to isogeny, the category coherent isocrystals on $(Y,M_Y)_\CRIS$ is equivalent to the isogeny category of (2) or (3) in the above proposition. 

\begin{cor} \label{cor:iso-crystals and quasi-nilpotent connections} 
The equivalence in Proposition~\ref{prop: crystals and quasi-nilpotent connections} gives rise to a fully faithful functor from the category $\CR_\Q((Y,M_Y)_\CRIS)$ of isocrystals to the category of isocoherent sheaves $\calM$ on $D_{\et}$ together with an integrable log connection $\nabla\colon \calM\rightarrow \calM\otimes_{\calO_Z}\omega^1_{(Z,M_Z)/\Z_p}$.
\end{cor}

\subsection{Log convergent site}\label{sec:log convergent site}

When we work on a semistable $p$-adic log formal scheme $(X,M_X)$, we associate to an $F$-isocrystal on $(X_1,M_{X_1})_\CRIS$ a vector bundle with integral connection on the generic fiber $X_\eta$. This is best explained by introducing log convergent sites. The convergent site was introduced by Ogus in \cite{ogus-F-converg-isocryst-de-rham-cohom-II, Ogus-convergenttopos}, and its logarithmic counterpart was developed by Shiho in \cite{Shiho-I, Shiho-II}.

Let $V$ be a CDVR of mixed characteristic $(0,p)$ with perfect residue field $k$. Later we consider $V=W(k)$ or $\calO_K$.
In this subsection, all the $p$-adic log formal $V$-schemes are assumed to be topologically of finite type over $V$, in particular, \emph{noetherian}.
For such a $(T,M_T)$, set $T_0\coloneqq (T\otimes_Vk)_\mathrm{red}$ (resp. $T_1\coloneqq (T\otimes_V V/p)$) and equip it with the pullback log structure $M_{T_0}$ (resp. $M_{T_1}$) from $M_T$.

Let $(Y,M_Y)\rightarrow (\Sigma,M_\Sigma)$ be a morphism of fine $p$-adic log formal $V$-schemes.

\begin{defn}
The \emph{log convergent site} $((Y,M_Y)/(\Sigma,M_\Sigma))_{\mathrm{conv}}$
is defined as follows:
\begin{itemize}
   \item an \textit{enlargement} of $(Y, M_Y)$ over $(\Sigma, M_{\Sigma})$ is a triple $(T, M_T, f)$ where $(T, M_T)$ is a fine $p$-adic log formal $V$-scheme over $(\Sigma, M_{\Sigma})$ and $f\colon (T_0, M_{T_0}) \rightarrow (Y, M_Y)$ is a morphism of log formal $V$-schemes over $(\Sigma, M_{\Sigma})$;
    \item the \textit{log convergent site} $((Y, M_Y)/(\Sigma,M_\Sigma))_{\mathrm{conv}}$ consists of enlargements of 
    $(Y, M_Y)$ over $(\Sigma, M_{\Sigma})$ with the obvious notion of morphisms. We often denote an enlargement $(T,M_T,f)$ by $T$.
    A cover over $(T,M_T,f)$ is a family of morphisms that are jointly surjective and strict \'etale over $(T,M_T)$.
\end{itemize}
The assignment $T \mapsto \Gamma(T, \mathcal{O}_T)$ (resp. $T \mapsto  \Gamma(T, \mathcal{O}_T)\otimes_\Z\Q$) defines a sheaf on $((Y, M_Y)/(\Sigma,M_\Sigma))_{\mathrm{conv}}$, for which we write $\mathcal{O}_{Y/\Sigma}$ (resp. $\mathcal{K}_{Y/\Sigma}$).

When $\Sigma=\Spf V$, we also write $((Y, M_Y)/(V,M_{\Spf V}))_{\mathrm{conv}}$ for $((Y, M_Y)/(\Spf V,M_{\Spf V}))_{\mathrm{conv}}$.
\end{defn}

\begin{defn}
An \emph{isocrystal} on $((Y, M_Y)/(\Sigma,M_\Sigma))_{\mathrm{conv}}$ is a sheaf $\mathcal{E}$ of finite locally free $\mathcal{K}_{Y/\Sigma}$-modules such that for any morphism $g\colon (T', M_{T'}) \rightarrow (T, M_T)$ in $((Y, M_Y)/(\Sigma,M_\Sigma))_{\mathrm{conv}}$, the map $g^* \mathcal{E}_T \rightarrow \mathcal{E}_{T'}$ is an isomorphism.   
\end{defn}

\begin{rem}\label{rem:first remarks on convergent site}
\hfill
\begin{enumerate}
    \item If $(Y',M_Y')\hookrightarrow (Y,M_Y)$ is an exact closed immersion of fine $p$-adic log formal $V$-schemes whose underlying morphism of topological spaces is a homeomorphism, then the canonical functor $((Y', M_{Y'})/(\Sigma,M_\Sigma))_{\mathrm{conv}}\rightarrow((Y, M_Y)/(\Sigma,M_\Sigma))_{\mathrm{conv}}$
    is an equivalence of categories \cite[Rem.~5.1.3]{Shiho-I}.
    \item The log convergent site enjoys the obvious functoriality for $(Y,M_Y)\rightarrow (\Sigma,M_\Sigma)$ and the category of convergent isocrystals satisfies the \'etale descent: see \cite[Rem.~5.1.5, 5,1,7]{Shiho-I}.
    \item The above definition of the log convergent site follows \cite[Def.~5.1.1]{Shiho-I}, and a different definition is given in \cite[Def.~2.1.1]{Shiho-II}. See \cite[Rem.~2.1.2, Prop.~2.1.7]{Shiho-II} for the comparison of the two notions.
\end{enumerate}
\end{rem}

\noindent
\textbf{Semistable formal schemes.}
Assume that $\Sigma=\Spf \calO_K$ is equipped with the canonical log structure $M_\mathrm{can}\coloneqq (\calO_K\smallsetminus\{0\})^a$. Let $(X, M_X)$ be a semistable $p$-adic log formal scheme over $\mathcal{O}_K$.
Write $X_\eta$ for the adic generic fiber of $X$, and identify the category of finite locally free isocoherent sheaves on $X$ with that of vector bundles on $X_\eta$.
Set $\omega^1_{X/\calO_K}\coloneqq\omega^1_{(X,M_X)/(\Spf\calO_K,M_\mathrm{can})}$.

\begin{prop} \label{prop:functor-converg-isocrys-to-vector-bundle}
There is a natural functor from the category of isocrystals on the log convergent site $((X, M_X) / (\mathcal{O}_K, M_\mathrm{can}))_{\mathrm{conv}}$ to the category of vector bundles with integrable connection on $X_{\eta}$.    
\end{prop}

\begin{proof}
This follows from \cite[Thm.~3.2.15, p.~624]{Shiho-I}:
let $(X^{(1)},M_{X^{(1)}})$ denote the exactification of the closed immersion $(X,M_X)\rightarrow (X,M_X)\times_{(\Spf \calO_K,M_\mathrm{can})}(X,M_X)$ as a $p$-adic log formal scheme. For each $n \geq 1$, let $(X^{(1),n}, M_{X^{(1),n}})$ denote the $n$-th infinitesimal neighborhood of $(X, M_X)$ in $(X^{(1)}, M_{X^{(1)}})$, and write $p_{i}^n\colon (X^{(1),n}, M_{X^{(1),n}}) \rightarrow (X, M_X)$ ($i = 1, 2$) for the projections. Then the log scheme $((X^{(1),n})_0,M_{(X^{(1),n})_0})$ is identified with $(X_0,M_{X_0})$, and 
$(X^{(1),n}, M_{X^{(1),n}})$ defines an object of $((X, M_X) / (\mathcal{O}_K, M_\mathrm{can}))_{\mathrm{conv}}$. 

Let $\mathcal{E}$ be an isocrystal on $((X, M_X) / (\mathcal{O}_K, M_\mathrm{can}))_{\mathrm{conv}}$. Evaluation on $(X,M_X)$ gives the isocoherent sheaf $\mathcal{E}_X$ on $X$, which is regarded as a vector bundle on $X_{\eta}$. We then have an infinitesimal log stratification on $\calE_X$ (which is called a formal log stratification in \cite[Def.~3.2.10]{Shiho-I}), namely, 
 the isomorphisms
\[
\eta_{\mathrm{inf}}^n\colon p_{1}^{n,\ast} \mathcal{E}_X \cong \mathcal{E}_{X^{(1),n}} \cong p_{2}^{n,\ast} \mathcal{E}_X
\]
satisfying the compatibility and cocycle conditions. Under the identification $\omega^1_{X/\calO_K}[p^{-1}]=\Omega^1_{X_\eta/K}$, this is equivalent to an integrable connection $\nabla\colon \mathcal{E}_X \rightarrow \mathcal{E}_X\otimes_{X_{\eta}} \Omega^1_{X_{\eta}/K}$ by \cite[Thm.~3.2.15]{Shiho-I}.
\end{proof}

\begin{rem}
With the notation as in the proof of Proposition~\ref{prop:functor-converg-isocrys-to-vector-bundle}, let $\calI$ denote the ideal of definition of $X$ in $X^{(1)}$. Let $B^{(1),n}\rightarrow X^{(1)}$ be the admissible blow-up of $X^{(1)}$ along $p\calO_{X^{(1)}}+\calI^{n+1}$ and define $T^{(1),n}\subset B^{(1),n}$ to be the maximal open $p$-adic formal subscheme on which $(p\calO_{X^{(1)}}+\calI^{n+1})\calO_{T^{(1),n}}=p\calO_{T^{(1),n}}$. We equip $T^{(1),n}$ with the pullback log structure $M_{T^{(1),n}}$ along $T^{(1),n}\rightarrow X^{(1)}$ and write $p_{T,i}^n\colon (T^{(1),n}, M_{T^{(1),n}}) \rightarrow (X, M_X)$ ($i = 1, 2$) for the projections.
Then $(T^{(1),n}, M_{T^{(1),n}})$ defines an object of $((X, M_X) / (\mathcal{O}_K, M_\mathrm{can}))_{\mathrm{conv}}$, and the category of isocrystals on $((X, M_X) / (\mathcal{O}_K, M_\mathrm{can}))_{\mathrm{conv}}$ is equivalent to the category of finite locally free isocoherent sheaves $\calE_X$ on $X$ together with compatible families of isomorphisms $p_{T,1}^{n,\ast} \mathcal{E}_X \cong p_{T,2}^{n,\ast} \mathcal{E}_X$ satisfying the compatibility and cocycle conditions (called convergent log stratifications): see \cite[pp.~626-631]{Shiho-I}.
The projection $p_{i}^n\colon (X^{(1),n}, M_{X^{(1),n}}) \rightarrow (X, M_X)$ factors
through $p_{T,i}^n$ and the infinitesimal log stratification on $\calE_X$ in the above proof is the pullback of the convergent log stratification on $\calE_X$ corresponding to $\calE$. One can express this in terms of a convergence condition of the resulting integrable connection $\nabla$ on $\calE_X$. 
\end{rem}

Using Dwork's trick, we now construct a canonical functor from the category 
$\Vect_\Q^\varphi((X_1,M_{X_1})_\CRIS)$ of $F$-isocrystals on $(X_1,M_{X_1})_\CRIS$ to the category of isocrystals on $((X, M_X) / (\mathcal{O}_K, M_\mathrm{can}))_{\mathrm{conv}}$. 
Let $(T, M_T, f) \in ((X, M_X) / (\mathcal{O}_K, M_\mathrm{can}))_{\mathrm{conv}}$ be an enlargement. Write $i_T\colon T_0 \rightarrow T_1$ for the natural inclusion. For $m_T \gg 0$, there exists a unique map $\rho_T^{(m_T)}\colon T_1 \rightarrow T_0$ such that $i_T \circ \rho_T^{(m_T)} = F_{T_1}^{m_T}$ and $\rho_T^{(m_T)} \circ i_T = F_{T_0}^{m_T}$. Set $f_T^{(m_T)} \coloneqq f \circ \rho_T^{(m_T)} \colon (T_1, M_{T_1}) \rightarrow (X_1, M_{X_1})\subset (X,M_X)$. Note $f_T^{(m_T+1)}=f_T^{(m_T)}\circ F_{T_1}=F_{X_1}\circ f_T^{(m_T)}$. Consider $(T, M_T)$ as an ind-object in $(T_1, M_{T_1})_{\CRIS}$. Then Construction~\ref{construction:strict log structure for functoriality} gives an ind-object $(T, M_{(T, f_T^{(m_T)})})$ in $(X_1,M_{X_1})_\CRIS$.

\begin{prop}\label{prop: F-isocrystals to convergent isocrystals}
There is a canonical functor from $\Vect_\Q^\varphi((X_1,M_{X_1})_\CRIS)$ to the category of isocrystals on $((X, M_X) / (\mathcal{O}_K, M_\mathrm{can}))_{\mathrm{conv}}$.    
\end{prop}

\begin{proof}
Let $(\calE_\Q,\varphi_{\calE_\Q}\colon F^\ast_{X_1,\CRIS}\calE_\Q\xrightarrow{\cong}\calE_\Q)$ be an $F$-isocrystal on $(X_1,M_{X_1})_\CRIS$. For each $(T, M_T, f) \in ((X, M_X) / (\mathcal{O}_K, M_\mathrm{can}))_{\mathrm{conv}}$, consider the evaluation $\mathcal{E}_\Q(T, M_{(T, f_T^{(m_T)})})$. By Proposition~\ref{prop:functoriality of log crystalline topoi}, we have the canonical isomorphisms
\[
\mathcal{E}_\Q(T, M_{(T, f_T^{(m_T)})})=(f_{T,\CRIS}^{(m_T),\ast}\calE_\Q)(T,M_T)\underset{\cong}{\xleftarrow{\varphi_{\calE_\Q}}}(f_{T,\CRIS}^{(m_T),\ast}F^\ast_{X_1,\CRIS}\calE_\Q)(T,M_T)=\mathcal{E}_\Q(T, M_{(T, f_T^{(m_T+1)})}).
\]
We set $\mathcal{E}_\Q(T, M_T, f) \coloneqq \mathcal{E}_\Q(T, M_{(T, f_T^{(m_T)})})$. By the above identification, it is well-defined and independent of the choice of $m_T$. By construction, this assignment is compatible with morphisms of enlargements in $(X_1,M_{X_1})_\CRIS$. Since $\calE$ is finite locally free and $T$ is noetherian, we easily see that $(T,M_T,f)\mapsto \calE(T, M_T,f)$ defines an isocrystal on $((X, M_X) / (\mathcal{O}_K, M_\mathrm{can}))_{\mathrm{conv}}$.
\end{proof}

\begin{rem}\label{rem:integrable connection associated to F-isocrystal}
By Propositions~\ref{prop:functor-converg-isocrys-to-vector-bundle} and \ref{prop: F-isocrystals to convergent isocrystals}, there is a functor from $\Vect_\Q^\varphi((X_1,M_{X_1})_\CRIS)$ to the category of vector bundles with integrable connection on $X_{\eta}$. We give a more explicit description.

With the notation as in the proof of Proposition~\ref{prop:functor-converg-isocrys-to-vector-bundle}, the map $X\rightarrow X^{(1),1}$ is a closed immersion with the same underlying topological space, and
$\omega^1_{X/\calO_K}$ is identified with $\Ker(\calO_{X^{(1),1}}\twoheadrightarrow\calO_X)$ (cf.~\cite[Prop.~3.2.5]{Shiho-I}).
Moreover, the exact closed immersions 
\[
(X_1,M_{X_1})\hookrightarrow (X,M_X) \quad\text{and}\quad (X_1,M_{X_1})\hookrightarrow (X^{(1),1},M_{X^{(1),1}})
\]
admit the canonical PD-structure on their defining ideals (by \cite[Ex.~3.2.4]{berthelot-ogus-book} for the latter) and yield ind-objects of $(X_1,M_{X_1})_\CRIS$, which we still write $(X,M_X)$ and $(X^{(1),1},M_{X^{(1),1}})$ by abuse of notation.

Let $(\mathcal{E}_{\mathbf{Q}},\varphi_{\mathcal{E}_{\mathbf{Q}}})$ be an $F$-isocrystal on $(X_1,M_{X_1})_{\CRIS}$, and write $(\mathbf{E},\nabla_{\mathbf{E}})$ for the associated vector bundle with an integrable connection on $X_{\eta}$. 
Then $\nabla_{\mathbf{E}}\colon \mathbf{E}\rightarrow \mathbf{E}\otimes_{\calO_X}\omega^1_{X/\calO_K}\subset p_{1}^{1,\ast}\mathbf{E}$ is given as 
\begin{equation}\label{eq:connection and stratification}
\nabla(e)=(\eta_{\mathrm{inf}}^1)^{-1}(1\otimes e)-e\otimes 1.
\end{equation}
With the notation as in Construction~\ref{construction: eval of isocrystal at log PD-formal scheme}, the proof of Proposition~\ref{prop: F-isocrystals to convergent isocrystals} provides the canonical identifications
\[
\mathbf{E}\cong\calE_{\Q,X}\quad\text{and}\quad
\eta_{\mathrm{inf}}^1\colon p_{1}^{1,\ast}\mathbf{E}\cong p_{1}^{1,\ast}\calE_{\Q,X}\cong \calE_{\Q,X^{(1),1}}\cong p_{2}^{1,\ast}\calE_{\Q,X}\cong p_{2}^{1,\ast}\mathbf{E}.
\]
More precisely, the first identification is given as $\mathbf{E}\cong F_{X_1,\CRIS}^{m,\ast}\calE_{\Q,X}\xrightarrow{\varphi_{\calE_\Q}^m}\calE_{\Q,X}$ for some $m$ such that $F_{X_1}^m\colon X_1\rightarrow X_1$ factors as $X_1\rightarrow X_0\rightarrow X_1$, and the second identification is given similarly. 

We also explain how to obtain the isomorphisms $\eta_{\mathrm{inf}}^n\colon p_{1}^{n,\ast} \mathcal{E}_X \cong \mathcal{E}_{X^{(1),n}} \cong p_{2}^{n,\ast} \mathcal{E}_X$ for $n\geq 2$. For this, we assume that $X=\Spf R$ is small affine and use the notation in Example~\ref{eg:BreuilSR}; we have the following commutative diagram 
\[
\xymatrix{
R[p^{-1}]\ar[r]^-{p_1^{n}} & R^{(1),n}[p^{-1}] & R[p^{-1}]\ar[l]_-{p_2^{n}}\\
S[p^{-1}]\ar[r]^-{p_{S,1}}\ar[u]^-f & S^{(1)}[p^{-1}]\ar[u]^-{f^{(1),n}} & S[p^{-1}].\ar[l]_-{p_{S,2}}\ar[u]_-f
}
\]
By Proposition~\ref{prop: crystals and quasi-nilpotent connections}, $(\calE_\Q,\varphi_{\calE_\Q})$ defines a finite projective $S[p^{-1}]$-module, an $S[p^{-1}]$-linear isomorphism $\varphi_\calM\colon \varphi_S^\ast\calM\rightarrow \calM$, and $S^{(1)}[p^{-1}]$-linear isomorphism $\eta_\calM\colon p_{S,1}^\ast\calM\xrightarrow{\cong}p_{S,2}^\ast\calM$. Pullback along the diagram gives a finite projective $R[p^{-1}]$-module $f^\ast\calM$ together with an isomorphism $f^{(1),n,\ast}\eta_\calM\colon p_1^{n,\ast}f^\ast\calM\xrightarrow{\cong}p_2^{n,\ast}f^\ast\calM$. The previous paragraph shows that the isocoherent sheaf on $X$ corresponding to $f^\ast\calM$ is identified with $\mathbf{E}$. Unwinding the construction and using $\varphi_\calM$ and the last sentence of Example~\ref{eg:BreuilSR}, one can check that $f^{(1),n,\ast}\eta_\calM$ yields $\eta_{\mathrm{inf}}^n$.
\end{rem}

Finally, let us give a definition of filtered $F$-isocrystals on $(X,M_X)$, which is suitable for our paper:

\begin{defn}\label{def:filtered F-isocrystals}
A \emph{filtered $F$-isocrystal} on a semistable $p$-adic formal log scheme $(X,M_X)$ over $\calO_K$ is a pair $((\calE_\Q,\varphi_{\calE_\Q}),(\mathbf{E},\nabla,\Fil^\bullet\mathbf{E}))$ where
\begin{itemize}
\item $(\calE_\Q,\varphi_{\calE_\Q})$ is an $F$-isocrystal on $(X_1,M_{X_1})_\CRIS$ (Definition~\ref{def:F-isocrystals});
\item $(\mathbf{E},\nabla)$ is the vector bundle with integral connection on $X_\eta$ attached to $(\calE_\Q,\varphi_{\calE_\Q})$ (Remark~\ref{rem:integrable connection associated to F-isocrystal});
\item $\Fil^\bullet\mathbf{E}$ is a $\Z$-indexed decreasing separated and exhaustive filtration of $\mathcal{O}_{X_{\eta}}$-submodules of $\mathbf{E}$ satisfying the Griffiths transversality such that $\mathrm{Fil}^i \mathbf{E} / \mathrm{Fil}^{i+1} \mathbf{E}$ is locally free over $\mathcal{O}_{X_{\eta}}$ for each $i$.
\end{itemize}
\end{defn}

\section{Functoriality of absolute strict log prismatic sites}\label{app: Functoriality}
 
We state the functoriality of the absolute strict log prismatic sites and refer the reader to \cite[\S 3]{du-moon-shimizu-prismatic-crystalline-comparison} for the detail.
Let $f\colon (X,M_X)\rightarrow (Y,M_Y)$ be a morphism of $p$-adic log formal schemes. We will discuss a morphism of topoi 
\[
f_\Prism=(f^{-1}_\Prism, f_{\Prism,\ast})\colon \operatorname{Sh}((X,M_X)_\Prism^{\str})\rightarrow\operatorname{Sh}((Y,M_Y)_\Prism^{\str}).
\]
If $f$ is strict, we have a fully faithful cocontinuous functor $(X,M_X)_\Prism^{\str} \rightarrow (Y,M_Y)_\Prism^{\str}$, which yields $f_\Prism$. Below we will construct $f_\Prism$ when $M_Y$ is integral and quasi-coherent (Proposition~\ref{prop:functoriality of log prismatic site}).

An \emph{$f$-prismatic morphism} (or an $f$-$\Prism$ morphism for short) from $(\Spf B, J, M_{\Spf B})\in (X,M_X)_\Prism^{\str}$ to $(\Spf A, I, M_{\Spf A})\in (Y,M_Y)_\Prism^{\str}$ is a morphism $g\colon (\Spf B, J, M_{\Spf B}) \rightarrow (\Spf A, I, M_{\Spf A})$ of log prisms that is compatible with $f$, namely, makes the following induced diagram commutative:
\[
\xymatrix{
(X,M_X)\ar[d]^-f & (\Spf (B/J), M_{\Spf (B/J)})\ar[l]\ar[d]^-{\overline{g}}\ar@{^{(}->}[r] & (\Spf B, M_{\Spf B})\ar[d]^-{g}\\
(Y,M_Y) & (\Spf (A/I), M_{\Spf (A/I)})\ar[l]\ar@{^{(}->}[r] & (\Spf A, M_{\Spf A}).
}
\]

\begin{lem}\label{lem:pullback of representable prismatic sheaf}
For $(\Spf A, I, M_{\Spf A})\in (Y,M_Y)_\Prism^{\str}$, the presheaf $f^{-1} h_A\coloneqq f^{-1}(\Spf A, I, M_{\Spf A})$ on $(X,M_X)_\Prism^{\str}$ defined by
\[
(X,M_X)_\Prism^{\str}\ni(\Spf B, J, M_{\Spf B})\mapsto \Mor_{f\text{-}\Prism}\bigl((\Spf B, J, M_{\Spf B}),(\Spf A, I, M_{\Spf A})\bigr)
\]
is a sheaf of sets. 
\end{lem} 

When $f=\mathrm{id}$, the lemma claims that the site $(X,M_X)_\Prism^{\str}$ is subcanonical.

\begin{proof}
See \cite[Lem.~3.3]{du-moon-shimizu-prismatic-crystalline-comparison}.
\end{proof}

\begin{construction}\label{construction:initial object for f-prismatic morphism}
Assume that $M_Y$ admits an integral chart $P\rightarrow \Gamma(Y,M_Y)$. For every prism $(\Spf B, J, M_{\Spf B})\in (X,M_X)_\Prism^{\str}$, let $M_{\Spf (B/J)}^P$ and $M_{\Spf B}^P$ denote the log structures on $\Spf(B/J)$ and $\Spf B$ associated to
the monoid homomorphisms
\begin{align*}
&P\rightarrow \Gamma(Y,M_Y) \rightarrow \Gamma(X,M_X)\rightarrow \Gamma(\Spf (B/J), M_{\Spf (B/J)}) \rightarrow B/J \quad\text{and}    \\
&P_B\coloneqq P\times_{\Gamma(\Spf (B/J), M_{\Spf (B/J)})}\Gamma(\Spf B, M_{\Spf B})\rightarrow \Gamma(\Spf B, M_{\Spf B}) \rightarrow B.
\end{align*}
\end{construction}

\begin{lem}\label{lem:f-prismatic morphism induced by chart}
The induced map $(\Spf (B/J), M_{\Spf (B/J)}^P)\rightarrow (\Spf B, M_{\Spf B}^P)$ is an exact closed immersion of integral and quasi-coherent $p$-adic log formal schemes, and the induced map $M_{\Spf B}^P\rightarrow M_{\Spf (B/J)}^P\times_{M_{\Spf (B/J)}}M_{\Spf B}$  in $\operatorname{Sh}((\Spf (B/J))_\et)=\operatorname{Sh}((\Spf B)_\et)$ is an isomorphism.
Moreover, the $\delta_\log$-structure on $\Gamma(\Spf B, M_{\Spf B})$ defines a $\delta$-log structure on $P_B$, and $(\Spf B, J, M_{\Spf B}^P)$ becomes a log prism in $(Y,M_Y)_\Prism^{\str}$ with an $f$-prismatic morphism $(\Spf B, J, M_{\Spf B})\rightarrow (\Spf B, J, M_{\Spf B}^P)$.
\end{lem}

\begin{proof}
See \cite[Lem.~3.8]{du-moon-shimizu-prismatic-crystalline-comparison}.
\end{proof}

\begin{prop}\label{prop:functoriality of log prismatic site}
If $M_Y$ is integral and quasi-coherent, then there exists a morphism of topoi
\[
f_\Prism\coloneqq (f^{-1}_\Prism,f_{\Prism,\ast})\colon \operatorname{Sh}((X,M_X)_\Prism^{\str})\rightarrow\operatorname{Sh}((Y,M_Y)_\Prism^{\str})
\]
satisfying the following properties.
\begin{itemize}
 \item For $\calF\in \operatorname{Sh}((X,M_X)_\Prism^{\str})$ and $(\Spf A, I, M_{\Spf A})\in (Y,M_Y)_\Prism^{\str}$, we have
\[
(f_{\Prism,\ast}\calF)(\Spf A, I, M_{\Spf A})=\Hom_{\operatorname{PSh}((X,M_X)_\Prism^{\str})}(f^{-1}h_A,\calF).
\]
 \item Assume that $M_Y$ admits an integral chart $P\rightarrow \Gamma(Y,M_Y)$. For $\calG\in \operatorname{Sh}((Y,M_Y)_\Prism^{\str})$ and $(\Spf B, J, M_{\Spf B})\in (X,M_X)_\Prism^{\str}$, we have 
\[
(f_\Prism^{-1}\calG)(\Spf B, J, M_{\Spf B})=\calG(\Spf B, J, M_{\Spf B}^P).
\]
\end{itemize}
\end{prop}

\begin{proof}
See \cite[Prop.~3.6]{du-moon-shimizu-prismatic-crystalline-comparison}.
Note that $M_Y$ admits an integral chart \'etale locally.
\end{proof}

\begin{rem}\label{rem:functoriality of perfect log prismatic site}
Assume that $M_X$ is integral and quasi-coherent. It is straightforward to check that the above construction also gives a morphism of topoi
\[
f_\Prism^\perf\coloneqq (f^{\perf,-1}_\Prism,f_{\Prism,\ast}^{\perf})\colon \operatorname{Sh}((X,M_X)_\Prism^{\str, \perf})\rightarrow\operatorname{Sh}((Y,M_Y)_\Prism^{\str, \perf})
\]
such that $(f_{\Prism,\ast}^\perf\calF)(\Spf A, I, M_{\Spf A})=\Hom_{\operatorname{PSh}((X,M_X)_\Prism^{\str, \perf})}(f^{-1}h_A,\calF)$.

When both $X$ and $Y$ are semistable or the log formal spectrum of a CDVR, Remark~\ref{rem:log prismatic site when admitting finite free chart}(2) gives $\operatorname{Sh}((X,M_X)_\Prism^{\str, \perf})\cong \operatorname{Sh}(X_\Prism^\perf)$ and $\operatorname{Sh}((Y,M_Y)_\Prism^\perf)\cong \operatorname{Sh}(Y_\Prism^\perf)$. In this case, $f_\Prism^\perf$ agrees with the morphism of topoi $\operatorname{Sh}(X_\Prism^\perf)\rightarrow \operatorname{Sh}(Y_\Prism^\perf)$ induced from the cocontinuous functor
\[
X_\Prism^\perf\rightarrow Y_\Prism^\perf: (\Spf A, I, \Spf (A/I)\rightarrow X) \mapsto (\Spf A, I, \Spf (A/I)\rightarrow X\rightarrow Y).
\]
This follows from the above description of $f_\Prism^\perf$.
\end{rem}

The pullback $f_\Prism^{-1}$ preserves various notions of crystals.

\begin{prop} \label{prop: prismatic pullback of crystals}
Assume that $M_Y$ is integral and quasi-coherent. 
\begin{enumerate}
    \item We have
\[
f_\Prism^{-1}\calO_\Prism=\calO_\Prism, \quad f_\Prism^{-1}\calO_\Prism[p^{-1}]^\wedge_{\calI_\Prism}=\calO_\Prism[p^{-1}]^\wedge_{\calI_\Prism}, \quad\text{and}\quad
f_\Prism^{-1}\calO_\Prism[\calI_\Prism^{-1}]^\wedge_p=\calO_\Prism[\calI_\Prism^{-1 }]^\wedge_p.
\]
 \item
Via the identification in (1),  $f_\Prism^{-1}$ induces
\[
f_\Prism^{-1}\colon \Vect((Y,M_Y)_\Prism^{\str},\calO_\Prism)
\rightarrow \Vect((X,M_X)_\Prism^{\str},\calO_\Prism).
\]
Moreover, if $M_Y$ admits an integral chart $P\rightarrow \Gamma(Y,M_Y)$, then for $\calG\in\Vect((Y,M_Y)_\Prism^{\str},\calO_\Prism)$  and $(\Spf B,J, M_{\Spf B})\in (X,M_X)_\Prism^{\str}$
\[
(f_\Prism^{-1}\calG)(\Spf B,J, M_{\Spf B})= \calG(\Spf B,J, M_{\Spf B}^{P}) \in \Vect(B).
\]
A similar assertion holds if $\Vect$ is replaced with $\Vect^\varphi$ or  $\calO_\Prism$ is replaced with $\calO_\Prism[p^{-1}]^\wedge_{\calI_\Prism}$ or $\calO_\Prism[\calI_\Prism^{-1}]^\wedge_p$.
\end{enumerate}
\end{prop}

\begin{proof}
This is a direct consequence of Proposition~\ref{rem:functoriality of perfect log prismatic site}.
\end{proof}

\begin{defn} \label{defn:pullback-analytic-prismatic-F-crystals}
Assume that $M_Y$ is integral and quasi-coherent.  Define the pullback functor of analytic prismatic $F$-crystals
\[
f_\Prism^{-1}\colon \Vect^{\an,\varphi}((Y,M_Y)_\Prism^{\str})\rightarrow \Vect^{\an,\varphi}((X,M_X)_\Prism^{\str})
\]
by associating to $(\calE_\Prism,\varphi_{\calE_\Prism})$
\[
\bigl((\calE_\Prism,\varphi_{\calE_\Prism})(\Spf B,J,M_{\Spf B}^P)\bigr)\in \varprojlim_{(\Spf B,J,M_{\Spf B})}\Vect^{\an,\varphi}(B,J)= \Vect^{\an,\varphi}((X,M_X)_\Prism^{\str}),
\]
where the limit is taken over the full subcategory of $(X,M_X)_\Prism^{\str}$ consisting of objects $(\Spf B,J,M_{\Spf B})$ such that $\Spf B/J\rightarrow Y$ factors through an \'etale $Y$-formal scheme $Y'$ admitting an integral chart $P\rightarrow \Gamma(Y',M_Y)$. Obviously, this construction is well-defined and compatible with Proposition~\ref{prop: prismatic pullback of crystals}(2) under the restriction functor $\Vect^{\varphi}((-,M_-)_\Prism^{\str})\rightarrow \Vect^{\mathrm{an},\varphi}((-,M_-)_\Prism^{\str})$ or the scalar extension functor $\Vect^{\an,\varphi}((-,M_-)_\Prism^{\str}) \to \Vect((-,M_-)_\Prism^{\str}, \calO_\Prism[\calI_\Prism^{-1}]^\wedge_p)^{\varphi=1}$.
\end{defn}

\bibliographystyle{amsalpha}
\bibliography{library}
	
\end{document}